\documentclass[11pt]{amsart}
\usepackage{amsmath,amssymb,amsthm}
\usepackage[margin=1.15in]{geometry}
\usepackage{booktabs}
\usepackage{graphicx}
\usepackage{xcolor}
\usepackage[colorlinks=true,linkcolor=blue,citecolor=blue,urlcolor=blue]{hyperref}

\usepackage{dsfont}
\newcommand{\id}{\mathds{1}}

\newtheorem{proposition}{Proposition}[section]

\newtheorem{lettertheorem}{Theorem}%

\newtheorem{lemma}[proposition]{Lemma}
\newtheorem{corollary}[proposition]{Corollary}
\newtheorem{definition}[proposition]{Definition}

\theoremstyle{definition}
\newtheorem{example}[proposition]{Example}
\theoremstyle{remark}
\newtheorem{remark}[proposition]{Remark}

\makeatletter
\NewDocumentEnvironment{restatable}{o m m +b}{%
  \IfNoValueTF{#1}%
    {\expandafter\gdef\csname rst@open@#3\endcsname{\begin{#2}}}%
    {\expandafter\gdef\csname rst@open@#3\endcsname{\begin{#2}[#1]}}%
  \expandafter\long\expandafter\gdef\csname rst@body@#3\endcsname{#4}%
  \def\@upn##1{\hyperref[rst:#3]{\textup{##1}}}%
  \csname rst@open@#3\endcsname
    \expandafter\protected@xdef\csname rst@num@#3\endcsname{\csname the#2\endcsname}%
    #4%
  \end{#2}%
  \expandafter\gdef\csname #3\endcsname{\@ifstar{\rst@restate{#2}{#3}}{\rst@restate{#2}{#3}}}%
}{}
\newcommand*\rst@restate[2]{%
  \begingroup
  \let\rst@label\label
  \let\label\@gobble
  \expandafter\let\csname the#1\expandafter\endcsname\csname rst@num@#2\endcsname
  \expandafter\def\csname theH#1\endcsname{restated.#2}%
  \addtocounter{#1}{-1}%
  \csname rst@open@#2\endcsname\rst@label{rst:#2}\csname rst@body@#2\endcsname\end{#1}%
  \endgroup}
\makeatother

\newcommand{\BF}{\bf\boldmath }
\newcommand{\one}{\mathbf 1}
\newcommand{\R}{\mathbb R}
\newcommand{\Z}{\mathbb Z}
\newcommand{\N}{\mathbb N}
\newcommand{\RP}{\mathbb{RP}}
\newcommand{\cN}{\mathcal N}
\newcommand{\cI}{\mathcal I}
\newcommand{\cF}{\mathcal F}

\newcommand{\cT}{\mathcal T}

\newcommand{\cL}{\mathcal L}
\newcommand{\GN}{\Gamma_{\mathcal N}}
\newcommand{\PGN}{{\mathbb P}\Gamma_{\mathcal N}}
\newcommand{\cG}{\mathcal G}
\newcommand{\GG}{\Gamma_{\mathcal G}}
\newcommand{\Oct}{\mathrm{O}}
\newcommand{\udimB}{\overline{\dim}_{\mathrm B}}
\newcommand{\sN}{s_{\mathrm N}}
\newcommand{\sA}{s_{\mathrm A}}
\newcommand{\dimH}{\dim_{\mathrm H}}
\newcommand{\dimB}{\dim_{\mathrm B}}
\newcommand{\abs}[1]{\lvert #1\rvert}
\newcommand{\norm}[1]{\lVert #1\rVert}
\newcommand{\pr}[1]{\langle #1\rangle}
\newcommand{\cK}{\mathcal K}
\newcommand{\PT}{\mathbb P(T)}

\DeclareMathOperator{\midd}{mid}
\DeclareMathOperator{\diam}{diam}

\begin{document}

\title[The Levitt--Yoccoz gasket]{The Levitt--Yoccoz gasket:\\ dynamics, Hausdorff dimension and norm exponent}
\author{Roberto De Leo}
\address{Department of Mathematics, Howard University, Washington, DC 20059, USA}
\date{\today}

\begin{abstract}
In S.\,P.~Novikov's problem on plane sections of triply periodic surfaces, the author and I.\,A.~Dynnikov showed that, for the $\mu$-cube, the directions that lie in no stability zone form the union of a fractal set $\cN$, the Levitt--Yoccoz gasket (projectively equivalent to the Rauzy gasket), and of its images under the sign changes of the coordinates. The gasket $\cN$ is the limit set of a parabolic projective iterated function system $\PGN$ generated by a semigroup $\GN\subset SL_3(\Z)$, and each zone is labelled by a primitive integer vector, its soul.
We show that $\cN$ is the global attractor of $\PGN$, and that the chain graphs of $\PGN$ and of its Hutchinson map consist of a single node.
Building on the theorem of Jiao--Li--Pan--Xu and Jurga that $\dimH\cN$ equals the affinity exponent of $\GN$, we prove $1.6853<\dimH\cN=\dimB\cN<1.6854$.
For gaskets $K$ cut out of an $n$-simplex by convex zones labelled by the elements of a matrix semigroup $\Gamma$, we bound $\udimB K$ from both sides by the norm exponent $\sN(\Gamma)$.
The norm exponent $\sN$ of $\GN$, the counting exponent of the souls, is Baragar's exponent of the Markoff--Hurwitz equation $x_1^2+x_2^2+x_3^2+x_4^2=ax_1x_2x_3x_4$. We prove $2.4411\le\sN\le2.4465$, improving Baragar's bounds, and $\frac23\sN\le\dimH\cN\le\frac12(1+\sN)$, which proves a conjecture recorded by Gamburd, Magee and Ronan; and we show that the exponent of the equation in six variables exceeds $3$, answering a question of Baragar.
For the three-dimensional Levitt--Yoccoz gasket $\cG\subset\RP^3$ we prove $2.74\le\sN(\GG)\le2.777$ and $\dimH\cG\le2.3885$.
Finally, for the renormalization flow of Avila, Hubert and Skripchenko, we show that the invariant measure they introduce is the Gibbs measure at the parameter $s_N$ and that $s_N$ is the maximal entropy per unit time of the return map, attained only by this measure. 
As a consequence, the diffusion exponent of their chaotic plane sections is less than 0.7766.
\end{abstract}

\maketitle

\section{Introduction}

We met what we call the \textbf{Levitt--Yoccoz} gasket (also called the Rauzy gasket in the literature) in 2006 in the study of the Novikov problem of planar sections of triply periodic surfaces.

\subsection*{The Novikov problem.}
Let $M\subset\R^3$ be a triply periodic surface and $H\in\RP^2$ a direction in $\R^3$.
Novikov's problem asks for the asymptotic behaviour of the open components of the sections of $M$ by the planes orthogonal to $H$. 
By Dynnikov's theory~\cite{Dyn99}, under some genericity assumptions for $M$, for each direction $H$ outside a set of measure zero, either there is no unbounded section or all unbounded sections 
lie in periodic components $C_i$ of $M$ that are warped planes with holes homotopic to zero.
Moreover, in that case there is a primitive integer vector $L$ (determined up to sign), that depends continuously on $H$ (and so is locally constant), such that each of these warped planes $C_i$ lies in a slab of finite width perpendicular to $L$.
In particular, each open trajectory has a strong asymptotic direction.
We call $L$ the soul of the pair $(M,H)$.
Given a primitive vector $L$ which is the soul for some direction $H$, we denote by $Z_L\subset\RP^2$ the set of all directions $H$ with soul $L$.
$Z_L$ is open and we call it the \textbf{stability zone} of $L$.

By~\cite{Dyn99} (see~\cite[Props.~1--3]{DD09}, where these results are summarized), any two mutually disjoint triply periodic surfaces $M_1,M_2$ have stability zones that are compatible: if $L_1$ is the soul of $(M_1,H)$ and $L_2$ is the soul of $(M_2,H)$, then $L_1=L_2$.
Hence, given a triply periodic function $f:\R^3\to\R$, it is possible to associate stability zones to $f$ as well.
Moreover, for a generic $f$, if there is more than a single stability zone, then: (a) the union of the stability zones is dense in $\RP^2$, (b) each zone has a piecewise smooth boundary and (c) the closures of two distinct zones meet at most in countably many points of their boundaries, so that the complement of the union of all stability zones is a cut-out fractal $K_f$~\cite[Prop.~3 and the remark after it]{DD09}.
Novikov conjectured that, in case of functions, almost all directions give rise, at some level surface $M_c=f^{-1}(c)$, to open sections with a strong asymptotic direction and $K_f$ has a Hausdorff dimension strictly between 1 and 2 (personal communication; see also \cite{MN04}).

\medskip
In case of the {\BF$\mu$-cube}, namely the infinite regular skew polyhedron $\{4,6\,|\,4\}$, the set of its stability zones, due to its particularly rich symmetry group, coincides with the set of the stability zones of the function 
\begin{equation}
\label{eq:mucube}
   G(x,y,z)=\midd\big(g(x),g(y),g(z)\big)-\tfrac12,
   \qquad g(t)=|2\{t\}-1|,
\end{equation}
where $\midd(a,b,c)$ is the median of three numbers and $\{t\}=t-\lfloor t\rfloor$ is the fractional part of $t$.
For this case the author, jointly with Ivan Dynnikov, was able to determine an algorithm to generate all stability zones~\cite{DD07,DD09}: three of them are squares, and every other one is, up to sign changes of the coordinates, an open triangle produced by an explicit recursion (Section~\ref{sec:souls}). 
We call Levitt--Yoccoz gasket, and we denote it by $\cN$, the set of directions in the first quadrant of $\RP^2$ that lie in no stability zone. 
The cut-out fractal $K_G$ corresponding to $G$, and so to the $\mu$-cube, is the union of $\cN$ and of its three images under the sign changes of the coordinates, one in each quadrant (Figure~\ref{fig:zonesRP2}).

\begin{figure}[t]
\centering
\includegraphics[width=0.55\textwidth]{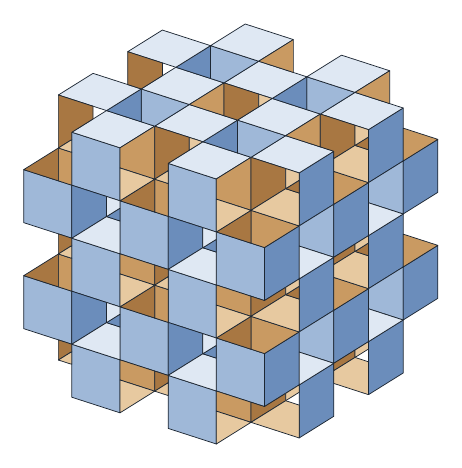}
\caption{A piece of the $\mu$-cube $\{4,6\,|\,4\}$, the surface $G=0$ with $G$ as in~\eqref{eq:mucube}, inside the cube $[-\frac14,\frac94]^3$. Its faces are squares of side $\frac12$, six at each vertex. The surface divides $\R^3$ into two congruent labyrinths: $\{G>0\}$, the union of the square tubes of side $\frac12$ around the lines of $\Z^3$ parallel to the axes, and $\{G<0\}=\{G>0\}+\frac12\one$. The side of each face looking into $\{G>0\}$ is brown, the one looking into $\{G<0\}$ is blue; at the boundary of the box the tubes are left open.}
\label{fig:mucube}
\end{figure}

\subsection*{A set discovered four times}

The set $\cN$ first appeared in literature in 1989 in an article by Meester on percolation theory: up to the boundaries of the stability zones and to permutations of the coordinates, $\cN$ is the set of parameters of Meester's circle percolation on $\Z^2$ for which his algorithm for the critical value never stops~\cite{Mee89}.
The algorithm subtracts the smallest coordinate from the other two until the largest one exceeds their sum, so it runs backwards the recursion that produces the stability zones, until the point lands in one of them (Section~\ref{sec:souls}).
In the same year, in the last section of a paper on the infinite clusters of circle percolation, Meester and Nowicki proved that this set has zero area~\cite[Thm.~5.1]{MN89}.

In 1991, Arnoux and Rauzy described $\cN$, without proof, as a set of measure zero ``du genre tapis de Sierpinski'' in an article on the representation of sequences of complexity $2n+1$ by exchanges of six intervals on the circle~\cite[\S4]{AR91}.
To our knowledge, Rauzy did not publish any other result on $\cN$ while Arnoux dedicated a full article to it, joint work with Starosta, in 2013 under the name of {\em Rauzy gasket}~\cite{AS13}.

In 1993, $\cN$ showed up again in an article of Levitt on pseudogroups of rotations of the circle~\cite{Lev93}: its points give minimal pseudogroups~\cite[Remarque~V.10]{Lev93}, and it is the set where Meester's algorithm never stops in his proof that, for almost every pair of rotation angles, the set of pseudogroups with a dense orbit in a two-parameter family is a Sierpi\'nski curve (a set different from $\cN$, cut out by squares).
In the same work, Yoccoz gave another proof that $\cN$ has zero area~\cite[Lemme~V.9]{Lev93}, and with it gave the first description of $\cN$ as a cut-out set: starting from the triangle with vertices $(0,1,1)$, $(1,0,1)$, $(1,1,0)$, which is our triangle $\cT$ (Section~\ref{sec:souls}), each triangle with vertices $X,Y,Z$ is replaced by the three triangles with vertices $X,X+Y,X+Z$, $Y,Y+Z,Y+X$ and $Z,Z+X,Z+Y$, a construction topologically analogous to that of the Sierpi\'nski triangle; the triangles removed are the zones of $\cN$.
It is for this reason that we believe that the name of the gasket should refer to these two authors, as suggested to the author by Novikov.

In 1995, Kraaikamp and Meester extended the result of Meester and Nowicki to circle percolation on $\Z^d$, whose set of parameters is, in the same way, the $d$-dimensional analogue $\mathcal G^d$ of $\cN$ (Section~\ref{sec:MH})~\cite{KM95}.
Their Lemma~1 is, in our language, the statement that $\mathcal G^d$ has zero volume (Theorem~\ref{thm:KM}).

To our knowledge, since then $\cN$ did not resurface until the author rediscovered it numerically in 2006 in connection with the $\mu$-cube, and then studied its analytic properties jointly with Ivan Dynnikov in~\cite{DD07} (2007) and~\cite{DD09} (2009).
In~\cite{DD07} the zones were described explicitly for the first time as the images of a single zone under the compositions of three projective maps, the maps $\psi_1,\psi_2,\psi_3$ of the present paper, and the points of $\cN$ as the limits of these compositions, so that $\cN$ is the limit set of a self-projective iterated function system (IFS).
In the same work its dimension was also estimated numerically for the first time, as approximately $1.7$, and it was announced that $\cN$ has zero area, which was proved in~\cite{DD09} by an argument different from those of~\cite{MN89} and~\cite{Lev93}.
The present author was immediately fascinated by the non-trivial structure of this set and started its study \emph{per se}, together with that of its analogues $\mathcal G^d$ in every dimension, with the 2008 preprint~\cite{DL08}, whose results were published in 2015 in~\cite{DL15EM,DL15JGA}.
These two articles, together with, and independently of, the works of Barnsley and Vince~\cite{BV12} (2012) and of Vince~\cite{Vin13} (2013), are the first published works devoted to self-projective IFSs.

\begin{figure}[t]
\centering
\begin{minipage}[b]{0.48\textwidth}\centering
\includegraphics[width=\textwidth]{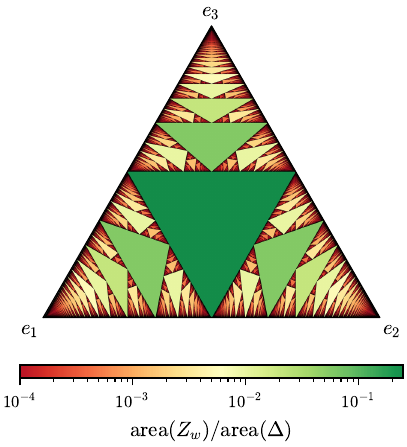}\\[2pt]
{\small (a) The Levitt--Yoccoz gasket $\cN$}
\end{minipage}\hfill
\begin{minipage}[b]{0.48\textwidth}\centering
\includegraphics[width=\textwidth]{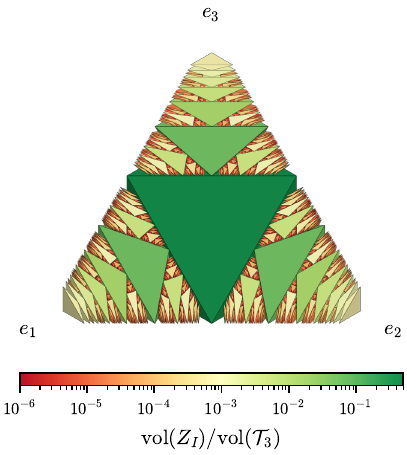}\\[2pt]
{\small (b) The three-dimensional Levitt--Yoccoz gasket $\cG$}
\end{minipage}
\caption{The two gaskets and their zones, coloured from green to red as the size of the zone decreases. (a) The image of $\cN\subset\cT$ under the similarity $x\mapsto\one-2x$ of the plane $\{\one^*x=1\}$, which maps the triangle $\cT$ of \S\ref{sec:souls} onto the triangle $\Delta$ with vertices $e_1,e_2,e_3$: the images of all the zones $Z_I$ larger than the resolution of the picture, coloured by $\mathrm{area}(Z_I)/\mathrm{area}(\cT)$; their boundaries, drawn in black, lie in $\cN$, and the black set is $\cN$ at the resolution of the picture. (b) The image of $\cG\subset\cT_3$ under the similarity $x\mapsto\one-3x$, which maps the tetrahedron $\cT_3$ of \S\ref{sec:3d} onto the $3$-simplex $\Delta$, drawn in $\R^3$ by the projection $x\mapsto(x_1,x_2,x_3)$ of the hyperplane $x_1+x_2+x_3+x_4=1$, as in \cite[Fig.~9]{DL08}, seen from the direction $(1,1,1)$ (the vertex $e_4$ is hidden behind the centre): the images of the octahedral zones $Z_I$ with $\abs I\le6$, drawn as solid bodies and coloured by $\mathrm{vol}(Z_I)/\mathrm{vol}(\cT_3)$. The front face is the face $\{x_4=0\}$ of $\Delta$, the image of the face of $\cT_3$ opposite to $\pr{\one-e_4}$, which meets $\cG$ in a copy of $\cN$; near the vertices and edges, where the zones with $\abs I>6$ are missing, the picture shows the zones lying behind. $\cG$ is the complement of the union of all the zones.}
\label{fig:gaskets}
\end{figure}

After that, several authors studied this set over the years.
Arnoux and Starosta studied it through Arnoux--Rauzy words and multidimensional continued fraction algorithms~\cite{AS13} (2013); Avila, Hubert and Skripchenko proved that $\dimH\cN<2$~\cite{AHS16} and constructed a natural invariant measure on it~\cite{AHSinv} (2016); Fougeron proved that the Levitt-Yoccoz gasket of every dimension has Hausdorff dimension smaller than that of the ambient space~\cite{Fou20} (2020); and Dynnikov, Hubert and Skripchenko gave a dictionary between several of the dynamical systems in which it appears~\cite{DHS23} (2023).
The known bounds on $\dimH\cN$ are recalled in {\em Main result~2} below.
To our knowledge, $\cN$, or the semigroup $\GN$ that generates it, appears in at least the following contexts:
\begin{itemize}
\item percolation theory: circle percolation on $\Z^d$~\cite{Mee89,MN89,KM95} (1989--1995);
\item symbolic dynamics and multidimensional continued fractions: Arnoux--Rauzy sequences \cite{AR91} (1991) and their letter frequencies~\cite{AS13} (2013); the fully subtractive algorithm~\cite{KM95,Sch00,FKN11} (1995--2011) and its dual, the Arnoux--Rauzy algorithm~\cite{AS13}; simplicial systems~\cite{Fou20} (2020);
\item pseudogroups of rotations, systems of isometries and band complexes~\cite{Lev93,Dyn08,Skr12,DS15,AHS16,AHSinv,DHS23} (1993--2023);
\item Novikov's problem on the plane sections of triply periodic surfaces~\cite{DD07,DD09,DS15,AHSinv} (2007--2016);
\item self-projective IFSs~\cite{DD07,DL08,DL15JGA,DL15EM} (2007--2015);
\item digital geometry: up to permutations of the coordinates and the boundaries of the stability zones, the cone over $\cN$ is the set of normal vectors of the arithmetical discrete planes for which the algorithm of Jamet and Toutant for the connecting thickness never stops~\cite{JT09,DJT09} (2009); the $2$-connectedness of these planes is characterized in~\cite{BJJP13} (2013);
\item tiling billiards: the exceptional trajectories of triangle tiling billiards are parametrized by the Levitt-Yoccoz gasket~\cite{HPR19} (2019);
\item Diophantine equations and hyperbolic geometry: the norm exponent of $\GN$ is the exponent of the number of integer solutions of the Markoff--Hurwitz equation in four variables~\cite{Bar98} (1998); up to transposition and permutation matrices, $\GN$ is the semigroup of~\cite{GMR19} (2019), where Gamburd, Magee and Ronan observed that its attracting set is the Rauzy gasket and, in more variables, its higher-dimensional versions; and, by results of Huang and Norbury~\cite{HN17} (2017) and of Magee~\cite{Mag20} (2020), the norm exponent of $\GN$ is the growth exponent of the number of one-sided simple closed geodesics on a thrice-punctured projective plane (Section~\ref{sec:MH}).
\end{itemize}
Figure~\ref{fig:gaskets}(a) shows $\cN$ together with its zones, and Figure~\ref{fig:gaskets}(b) the three-dimensional Levitt--Yoccoz gasket $\cG$ discussed below.

\subsection*{Souls}
In~\cite{DL03} we proved that, for a triply periodic Morse function $f$ with more than one stability zone, every point of $K_f$ is an accumulation point of souls; the argument rests on the fact that there is a constant $c>0$ such that $Z_L$ is contained in a disc of radius $c/\|L\|$ centered at the direction $[L]$ (see~\cite{DL05}).
The function $G$ is not a Morse function, but the same holds for $K_G$: by Corollary~\ref{cor:accumulation}, the gasket $\cN$ is the set of accumulation points of the souls of the zones in the first quadrant, and so, by the symmetry under the sign changes, $K_G$ is the set of accumulation points of all the souls. So every quantity related to the gasket can be extracted from the set of its souls.
In~\cite{DL08} we evaluated numerically the box dimension of $\cN$, finding $\dimH\cN\le\dimB\cN\simeq1.7$ (two box counts gave $1.69$ and $1.7$; see Remark~\ref{rem:172}; for the first estimates, see~\cite{DD07,DD09}), and started the study of the dimension of the three-dimensional gasket (Main result~4).

\medskip
Before going over the structure of the article, we now illustrate the main results of this work. 

\subsection*{Main result 1: dynamics of the gasket}

In~\cite{DD07,DL08} (see also~\cite{DL15JGA,DL15EM}) we showed that $\cN$, as a cut-out fractal, can be seen as a compact invariant subset of an IFS, that we denote by $\PGN$, generated by three real projective maps of a triangle $\cT$ into itself (Section~\ref{sec:souls}), and in~\cite[Prop.~2]{DL08} that the same holds for its analogues in every dimension. This description is similar
to the way the Sierpi\'nski gasket can be defined as the (unique) compact invariant subset of an IFS generated by three contractions defined on a triangle and the Apollonian gasket can be defined as the (unique) compact invariant subset of an IFS generated by three M\"obius maps.
Unlike the Sierpi\'nski case, though, $\PGN$ is not contractive, namely its maps are not contractions (they have a parabolic fixed point, see Section~\ref{sec:souls}), and, unlike the Apollonian gasket, $\PGN$ is not conformal either.

It is natural to ask what the qualitative dynamics of $\PGN$ is.
In~\cite{DLstreams}, generalizing an analogue construction on dynamical systems introduced in~\cite{DLY26DCDS,DeLeoYorke2026JDE} (joint work with Jim Yorke), the author studied some fundamental properties of the qualitative dynamics of a general (i.e.\ not necessarily contractive) IFS $\cF$ on a locally compact metric space $X$ and showed how to encode, using a generalization of Charles Conley's chains, this behavior into a graph. Chain recurrence and the Conley decomposition had been extended before to closed relations on compact spaces by McGehee and Wiandt~\cite{MW06}, and to iterated function systems on compact metric spaces, through Conley attractors, by Barnsley and Vince~\cite{BV13}, building on Kieninger's study of IFSs on compact spaces~\cite{Kie02}.
The essential dynamical condition we require on $\cF$ is that it has a global attractor, namely a compact subset which is invariant under the Hutchinson map of $\cF$ and attracts, under $\cF$, every compact subset of $X$.

In case of dynamical systems, the relevance of the graph is given by the following property: every point in the global attractor is either chain-recurrent, in which case its forward orbit lies entirely inside a node of the graph, or belongs to a trajectory of the system that is defined for all times and whose $\alpha$-limit set lies in one node and whose $\omega$-limit set lies in a different node.
In this case we say that there is an edge from the first node to the second.
In short, each point in the global attractor has either a recurrent dynamics, in which case it belongs to a node, or a gradient dynamics, in which case it defines an edge.
In the case of IFSs, we show~\cite{DLstreams} that there are three options: a point has either recurrent dynamics under each generator of the IFS, or a mix of recurrent and gradient dynamics, or gradient dynamics only.

Our first main result is that $\cN$ is the global attractor of $\PGN$ and, in fact, its unique compact $H$-invariant set, and that both the graph of $\PGN$ and the graph of its Hutchinson map consist of a single node. Below, $\cT$ is the triangle of Section~\ref{sec:souls} that carries the zones of the first quadrant, on which $\PGN$ is defined, and $\mathcal K(\cT)$ is the space of the non-empty compact subsets of $\cT$ with the Hausdorff distance $d_H$; the notions involved are recalled in Section~\ref{sec:souls}.

\begin{restatable}{theorem}{thmsinglenode}
\label{cor:singlenode}
    Let $\mathcal H$ be the Hutchinson map of $\PGN$.
    Then:
    \begin{enumerate}
    \item[(i)] $\cN$ is the only non-empty $H$-invariant compact subset of $\cT$ under $\PGN$;
    \item[(ii)] $\mathcal H^n(K)\to\cN$ in $d_H$ for every $K\in\mathcal K(\cT)$;
    \item[(iii)] the global attractor of $\PGN$ is $\cN$;
    \item[(iv)] the global attractor of $\mathcal H$ is $\{\cN\}$;
    \item[(v)] the chain graph of $\PGN$ consists of the single node $\cN$;
    \item[(vi)] the chain graph of $\mathcal H$ consists of the single node $\{\cN\}$.
\end{enumerate}
\end{restatable}

Hence, $\PGN$ has the same qualitative dynamics as a contractive IFS. Items (i)--(iv) can also be obtained from Edelstein's fixed point theorem~\cite{Ede62}, because the maps $\psi_j$ shrink the distance between any two distinct points of $\cT$ (Proposition~\ref{lem: Dpsi} and Remark~\ref{rem:edelstein}); the new content of Theorem~\ref{cor:singlenode} is (v)--(vi), the structure of the chain graphs. Their proof needs only that the cylinders $\psi_I(\cT)$ with $\abs I=n$ shrink uniformly as $n\to\infty$; this also follows from Proposition~\ref{lem: Dpsi}, by compactness (Remark~\ref{rem:edelstein}), and Lemma~\ref{lem:geom} gives it with the explicit bound $6\sqrt3/(3+2n)$ on their diameters.

\subsection*{Main result 2: Hausdorff dimension of the gasket}

Our second main result is the evaluation of new rigorous bounds for the Hausdorff dimension of $\cN$.

\begin{restatable}{theorem}{thmdimH}
    \label{thm: dimH}
     The following rigorous bounds hold:
     \[
     1.6853<\dimH\cN=\dimB\cN<1.6854.
     \]
     So $\dimH\cN=1.6853\ldots$, correct to four decimal places.
\end{restatable}

The following bounds were known previously: $\dimH\cN>1.19$ \cite{GRM20};
$\dimH\cN\ge3/2$ \cite{JLPX23} (see also \cite{LPX23}); $\dimH\cN\ge1.6196$ \cite{Sew23} (preprint), by an argument also based on the norm exponent, via the box dimension;
$\dimH\cN<2$ \cite{AHS16};
$\dimH\cN<1.825$ \cite{Fou20} (preprint);
$\dimH\cN\le1.7415$ \cite{PS23}.
Theorem~\ref{thm: dimH} improves both the lower and the upper bound, and it reduces the width of the known interval from $0.24$ ($0.12$ with the preprint bound of~\cite{Sew23}) to $10^{-4}$.
This range is consistent with the first numerical estimates of its box dimension in~\cite{DD07,DL08,DD09}, namely $\simeq1.7$ for the set and $\simeq1.69$ for its barycentres, with the conclusion ``probably about 1.7''.

The proof rests on the recent theorem of Jiao, Li, Pan and Xu~\cite{JLPX23,LPX23} and Jurga~\cite{Jur26}, obtained independently: the Hausdorff dimension of $\cN$ equals the affinity exponent $\sA(\GN)$, the critical exponent of the series of the singular value functions of the elements of $\GN$ (Section~\ref{sec:affinity}). We compute rigorous bounds for $\sA(\GN)$. The diameter of the zone $Z_I$ is comparable to $\sigma_2(A_I)/\sigma_1(A_I)$, and so is the norm of the derivative $D_{x_0}\psi_I$ at the centre $x_0=\pr\one$ (Lemmas~\ref{lem:diamZ} and~\ref{lem:DI}). We split every multiindex into runs of equal indices: each run multiplies the size of a soul by a factor that depends only on the soul point before the run, and the same holds for the derivative once its direction is added to the state (Lemma~\ref{lem:Dcocycle}). So $\sA$ is the critical exponent of a positive transfer operator on a three-dimensional space, and it suffices to exhibit positive test functions $h$ with $\tilde\cK h\ge h$, or $\tilde\cK h\le\theta h$ with $\theta<1$ (Proposition~\ref{cor:certA}). Ours are continuous and piecewise linear on a graded mesh of $3.0\times10^7$ cells. They are approximate eigenfunctions, corrected so that the extremes of $(\tilde\cK h)/h$ between the nodes, rather than its values at the nodes, are nearly constant. The two inequalities are verified at every point, in interval arithmetic with directed rounding: the operator is evaluated directly at the vertices of the cells, or of the boxes of an adaptive subdivision (about $1.5\times10^9$ boxes in all), and between the vertices the interpolation error is bounded through enclosures of second derivatives, of increments or of oscillations (Proposition~\ref{prop:P1} and Section~\ref{sec:comp}). The verification took about three days on six machines with $16$ cores each; a coarser version of it, which proves $1.68364<\sA(\GN)<1.68637$, takes a few hours on two cores (Remark~\ref{rem:coarse}). Runs are the right steps because $\GN$ is a fast semigroup in the sense of~\cite{DL15JGA}: norms add within a run, where the generators are parabolic, but multiply up to a bounded factor across runs (Remark~\ref{rem:fast}).

It is instructive to compare with the Apollonian gasket, whose Hausdorff dimension Vytnova and Wormell recently computed to $128$ decimal places~\cite{VW24}.
Their setting is close to ours in several respects: 
the Apollonian gasket is also the limit set of a parabolic system of three maps; they also replace it by an induced, uniformly contracting system with infinitely many branches, the analogue of our runs; their bounds also come from an approximate eigenfunction of a finite-rank approximation of a transfer operator; and they also have to bound infinite sums over the branches, by the Euler--Maclaurin formula. 
The derivative of a conformal map is a similarity, so the size of the image of a small set is a single number, the modulus of the derivative, and by the chain rule the size of a composition is exactly the product of the factors contributed by its maps. For the induced system, which is uniformly contracting, the sum of the $s$-th powers of the sizes of its pieces of generation $n$ then grows or decays exponentially in $n$, at a rate given by the leading eigenvalue of a transfer operator, and the dimension is the exponent $s$ at which it does neither (Bowen's formula~\cite{Bow79}).
Moreover, the maps are holomorphic, so the transfer operator acts compactly on a Hardy space of holomorphic functions; there, Chebyshev--Lagrange approximation converges exponentially fast, with explicit constants, and a matrix of moderate size gives many digits. 

None of this is available for $\cN$. 
The maps $\psi_j$ are projective and not conformal: the two singular values of $D_x\psi_I$ can be of different orders of magnitude, the zones $Z_I$, of area of order $\|A_I\|_1^{-3}$, can have diameter of order $\|A_I\|_1^{-1}$ and so be far from round (Lemma~\ref{lem:geom} and Remark~\ref{rem:geomsharp}), and the size of a zone is not a single number: it has a length and a width, the two singular values. Under composition these satisfy only inequalities, such as $\sigma_1(AB)\le\sigma_1(A)\sigma_1(B)$, because the length of the image depends on how the long direction of one map lines up with the other; so $\dimH\cN$ is the affinity exponent, defined through these inequalities, rather than an exponent given by exact products. Exact products come back if the direction of the derivative is carried along with the point, and this is why our transfer operator acts on functions of three real variables.
Its eigenfunction is not smooth: near two lines of the state space it varies on scales proportional to the distance from them (Remark~\ref{rem:fan}). 
Piecewise-linear approximations of its leading eigenvalue converge only with order about $1.7$ in the mesh size, and, what matters for a certificate, so do the deviations of $(\tilde\cK h)/h$ between the nodes, for piecewise-quadratic functions as well (Remark~\ref{rem:sAnum}). 
This is why we certify the bounds by inequalities for piecewise-linear test functions on graded meshes, verified at every point, and why each further digit is expensive: the fourth decimal place of $\dimH\cN$ needed a computation on a cluster, and certifying the fifth would need meshes with about a hundred times as many cells (Remark~\ref{rem:sAnum}).

The lower bound of Theorem~\ref{thm: dimH} uses the equality $\dimH\cN=\sA$; the upper bound does not, because a cut-out argument with the exact shape of the zones gives $\udimB\cN\le\sA$ (Proposition~\ref{prop:cutout}). 
By~\cite{Sew23} together with~\cite{JLPX23,Jur26}, $\dimB\cN$ exists and equals $\dimH\cN$; since~\cite{Sew23} is currently still unpublished, we provide a proof in Corollary~\ref{cor: dimB}.
Finally, we provide a new numerical estimate: $\dimH\cN\approx1.685342$ (Remark~\ref{rem:sAnum}).

\subsection*{Main results 3: the norm exponent}
The souls of the zones in $\cT$ are the vectors $L_I=A_I\one$, where $A_I$ ranges over the elements of the free semigroup $\GN$ of non-negative matrices of $SL_3(\Z)$ generated by the three unipotent matrices $A_i=\id+(\one-e_i)e_i^*$ (Section~\ref{sec:souls}). Since all entries of the generators are non-negative, so are all entries of each element of $\GN$; on non-negative matrices $\|\cdot\|_1=\one^*(\cdot)\one$ is a \emph{linear} functional, and
\[
  \|L_I\|_1=\|A_I\|_1=\one^* A_I\one .
\]

\medskip
Counting souls by their $\ell^1$-norm is therefore counting the elements of $\GN$ by their norm. The exponent that governs such counts is defined for every finitely generated semigroup of matrices. Let $\Gamma$ be the semigroup generated by the matrices $A_1,\dots,A_m$, let $\cI^m$ be the set of the multiindices in $\{1,\dots,m\}$, let $A_I$ be the product of the generators along $I\in\cI^m$, and fix a norm $\norm\cdot$ on matrices such that $\nu_I=\norm{A_I}\ge1$ for every $I\in\cI^m$ and $N_\Gamma(r)=\#\{I\in\cI^m:\ \nu_I\le r\}<\infty$ for every $r$. The \textbf{zeta function} and the \textbf{norm exponent} of $\Gamma$ \cite{DL12,DL15JGA} are
\[
  \zeta_\Gamma(t)=\sum_I\nu_I^{-t},\qquad
  \sN(\Gamma)=\inf\Bigl\{t:\ \zeta_\Gamma(t)<\infty\Bigr\}.
\]
Since any two norms are comparable up to constant factors, $\sN(\Gamma)$ does not depend on the norm; and, since $\zeta_\Gamma$ is a Dirichlet series with positive coefficients, $\sN(\Gamma)=\limsup_{r\to\infty}\log N_\Gamma(r)/\log r$ \cite[Thm.~3]{DL15JGA}. For $\GN$ we always use the norm $\norm\cdot_1$, and we write simply
\[
  \zeta(t)=\zeta_{\GN}(t)=\sum_I \|A_I\|_1^{-t}=\sum_I \|L_I\|_1^{-t},\qquad \sN=\sN(\GN).
\]

\medskip\noindent\textbf{Souls and Markoff--Hurwitz quadruples.}
The series $\zeta$ was introduced, in a different language, by Baragar~\cite{Bar98}, in his study of the integer solutions of the Markoff--Hurwitz equation
\begin{equation}\label{eq:MH}
  x_1^2+x_2^2+x_3^2+x_4^2=a\,x_1x_2x_3x_4 .
\end{equation}
Baragar reduced the count of these solutions to a sum over a tree whose branching rule is the action of the generators $A_j$ on the souls (Lemma~\ref{lem:baragartree}), and he proved that, if~\eqref{eq:MH} has an integer solution other than $0$ (for instance $(1,1,1,1)$ for $a=4$), then the number of its integer solutions with $\max_i\abs{x_i}\le R$ is $(\log R)^{\sN+o(1)}$ as $R\to\infty$ (Theorem~\ref{thm:baragar}). So $\sN$ is the exponent denoted $\alpha(3)$ in~\cite{Bar98} and $\beta(4)$ in~\cite{GMR19}, and the bounds $2.430<\sN<2.477$ of~\cite[Thm.~0.1]{Bar98} were the best known so far. In the same way, the norm exponent of the semigroup $\GG$ of Main result~4 below is the exponent of the equation in five variables. By results of Huang and Norbury~\cite{HN17} and of Magee~\cite{Mag20}, $\sN$ is also the growth exponent of the number of one-sided simple closed geodesics of length at most $L$ on any hyperbolic thrice-punctured projective plane of finite area (Remark~\ref{rem:geodesics}). So the norm exponent of $\GN$ is at the same time the counting exponent of the souls in Novikov's problem, of the integer solutions of~\eqref{eq:MH}, and of these geodesics.

The norm exponent is a tool to bound the dimension of gaskets, and the first result says when it works. Let $\Delta$ be an $n$-simplex, $n\ge2$, in an affine $n$-plane of $\R^{n+1}$, and let $(Z_I)$ be pairwise disjoint convex subsets of $\Delta$ with non-empty interior in the $n$-plane, the \textbf{zones} (they need not be open: the zones of the three-dimensional gasket $\cG$ below contain some of their faces), indexed by the $I\in\cI^m$, where $m$ is the number of generators of a semigroup $\Gamma$ with a norm as above (that is, by the elements $A_I$ of $\Gamma$, when $\Gamma$ is free). We call $K=\Delta\setminus\bigcup_IZ_I$ a \textbf{cut-out gasket}. Volumes and areas are Euclidean, in the $n$-plane, and $\partial$ denotes the boundary relative to it.

\begin{restatable}{theorem}{thmcutout}\label{thm:cutout}
Let $K$ be a cut-out gasket, and assume $n-1<\sN(\Gamma)<n+1$.
\begin{enumerate}
\item[(a)] If $K$ has measure zero, and $\mathrm{vol}(Z_I)\le C\,\nu_I^{-(n+1)}$ and $\mathrm{area}(\partial\overline{Z_I})\le C\,\nu_I^{-(n-1)}$ for every $I$, then
\[
  \udimB K\ \le\ \tfrac12\bigl(n-1+\sN(\Gamma)\bigr).
\]
\item[(b)] If $\mathrm{vol}(Z_I)\ge c\,\nu_I^{-(n+1)}$ and $\partial\overline{Z_I}\subset K$ for every $I$, then
\[
  \udimB K\ \ge\ \tfrac{n}{n+1}\,\sN(\Gamma).
\]
\end{enumerate}
\end{restatable}

Part (a) is an elementary cut-out argument. Part (b) rests on the isoperimetric inequality: among zones of a given volume the round ones have the smallest boundary, so round zones are the worst case for a lower bound. 
The volume hypotheses hold, with $\nu_I=\norm{A_I}_1$, for the gaskets $\cN$ and $\cG$ discussed below (Lemma~\ref{lem:geom}(iii) and Lemma~\ref{lem:geom3}). 
The bound in (b) is the box-dimension counterpart of the inequality $\sA\ge\min\{n,\frac n{n+1}\sN\}$, which holds for every finitely generated semigroup in $SL_{n+1}(\R)$ (Proposition~\ref{prop:main}). 
Hence our Conjecture~1.4 of~\cite{DL15EM}, namely that $(n+1)\dimH\Lambda\ge n\,\sN$ for the limit set $\Lambda$ of such a semigroup, holds whenever $\dimH\Lambda=\sA<n$, and its box-dimension version holds for every cut-out gasket satisfying (b). Without a restriction of this kind the inequality can fail, because $\sN$ can exceed $n+1$ (Remark~\ref{rem:minneeded}); the condition $\sA<n$ excludes this, since by Proposition~\ref{prop:main} it forces $\sN<n+1$.

For the Levitt--Yoccoz gasket both parts apply.

\begin{restatable}{theorem}{thmsN}
    \label{thm: sN}
    $\displaystyle\tfrac23\, \sN\;\le\;\dimH\cN\;=\;\dimB\cN\;\le\;\frac{1+ \sN}{2}$.
\end{restatable}

This result is the case $n=2$ of Theorem~\ref{thm:cutout}, which gives $\frac23\sN\le\udimB\cN\le\frac12(1+\sN)$, together with Corollary~\ref{cor: dimB}, by which the box dimension of $\cN$ exists and equals its Hausdorff dimension. That corollary rests on the equality $\dimH\cN=\sA(\GN)$ proved by Jiao, Li, Pan and Xu~\cite{JLPX23} and by Jurga~\cite{Jur26} (see Remarks~\ref{rem:conj} and~\ref{rem:conventions} for the conventions). Given this equality, both bounds of Theorem~\ref{thm: sN} also follow from the bounds $\sigma_1(A_I)^{-1/2}\le\sigma_2(A_I)\le3$ on the middle singular value (Lemma~\ref{lem:sv} and Remark~\ref{rem:ceiling}). What Theorem~\ref{thm:cutout} adds is that it bounds the upper box dimension from both sides without using that equality, and in every dimension.

Theorem~\ref{thm: sN} relates the exponent of the Markoff--Hurwitz equation~\eqref{eq:MH} to the dimension of the gasket. Gamburd, Magee and Ronan characterized this exponent, $\beta(4)=\sN$, as the unique parameter for which the action of the accelerated semigroup on the projective plane has a conformal measure, one that transforms under each map by the power $\beta/3$ of its Jacobian~\cite[Thm.~10 and Lemma~44]{GMR19} (in our coordinates, the Jacobian of $\psi_I$ at $x$ is $\norm{A_Ix}_1^{-3}$, so this power is $\norm{A_Ix}_1^{-\sN}$). They noted that this is a partial analogue of the Patterson--Sullivan theory, which identifies the critical exponent of a convex cocompact Fuchsian group with the Hausdorff dimension of its limit set, but that here ``the parameter $\beta$ is not in any obvious way connected to the Hausdorff dimension'' of the gasket~\cite[Remark~11]{GMR19}. Right after this remark they recorded our conjecture~\cite{DL15EM} that the Hausdorff dimension $\delta$ of the Rauzy gasket satisfies $\delta\ge\frac23\beta(4)$, and observed that, with Baragar's bounds, it would give $\delta>1.62$. The lower bound of Theorem~\ref{thm: sN} proves this conjecture (its box-dimension version is due to Sewell~\cite[Thm.~1.8]{Sew23}), and the upper bound gives the inequality $\delta\le\frac12(1+\beta(4))$ in the other direction.
The value $\frac23\beta(4)$ is also what the Patterson--Sullivan analogy predicts: for conformal maps of the plane the Jacobian is the square of the modulus of the derivative, so a measure that transforms by the power $\beta/3$ of the Jacobian transforms by the power $\frac23\beta$ of the modulus of the derivative, and the analogy would give $\delta=\frac23\beta$. By Theorems~\ref{thm: dimH} and~\ref{thm: delta}, $\frac23\beta(4)\le1.6310<1.6853<\delta$: the conformal prediction is a strict lower bound for the dimension, and the gap comes from the non-conformality of the maps $\psi_j$. Accordingly, the lower bound $\frac23\sN$ of Theorem~\ref{thm:cutout}(b) is the one given by round zones, which are the worst case for the isoperimetric argument, and it corresponds to the extreme case $\sigma_2(A_I)=\sigma_1(A_I)^{-1/2}$ of the middle singular value (Remark~\ref{rem:ceiling}).

\medskip
Then, we provide rigorous bounds for $\sN$.

\begin{restatable}{theorem}{thmdelta}
    \label{thm: delta}
    The following rigorous bounds hold:
    $2.4411\le \sN\le 2.4465$.
\end{restatable}

The method works because $\GN$ is a fast semigroup in the sense of~\cite{DL15JGA}: when the last index of $I$ differs from the first index of $J$, the ratio $\|A_{IJ}\|_1/(\|A_I\|_1\|A_J\|_1)$ lies in $[\frac14,\frac12]$, and the constant $\frac14$ is sharp (Lemma~\ref{lem:junction}).
We compute this ratio exactly, as a function of the direction of the soul of $J$, and this turns $\zeta$ into the sum of the iterates of a transfer operator on three triangles, on which the runs of equal generators act by uniform contractions (Lemma~\ref{lem:contr}).
The bounds of Theorem~\ref{thm: delta} are Collatz--Wielandt bounds for its discretizations; in practice their width decreases like the mesh size, while the bounds given by \cite[Thm.~1]{DL15JGA}, which replace the ratio by a constant, improve only logarithmically. The same holds for Baragar's method~\cite{Bar98}, a truncation of his tree refined by Boyd's convexity argument, which converges like $1/\log y$ in the truncation height $y$~\cite[Thm.~2.3]{Bar98} and gave $2.430<\sN<2.477$; Theorem~\ref{thm: delta} narrows this interval by a factor of about $9$ (Remark~\ref{rem:JGA}).

Set
\[
  N(r)=\#\{I:\ \|A_I\|_1\le r\}
  =\#\{\text{zones }Z_L\text{ in $\cT$ with }\|L\|_1\le r\}.
\]
Over all of $\RP^2$ the number of zones with $\|L\|_1\le r$ is $4N(r)+3$ for $r\ge1$ (four sign classes, three squares; Theorem~\ref{thm:DD} and Figure~\ref{fig:zonesRP2}). As recalled above, $\sN=\limsup_{r\to\infty}\log N(r)/\log r$, and by Baragar's theorem the $\limsup$ is a limit. Much more is true. Reading a theorem of Gamburd, Magee and Ronan~\cite[Thm.~13]{GMR19} through the dictionary of Section~\ref{sec:MH} (Proposition~\ref{lem:GMRdict}), we obtain the strong asymptotics of the number of souls.

\begin{restatable}[counting the souls]{theorem}{thmcount}\label{thm:count}
There is $c_N>0$ such that
\[
  N(r)\ \sim\ c_N\,r^{\sN}\qquad(r\to\infty).
\]
Hence the number of stability zones of the $\mu$-cube in $\RP^2$ whose soul has $\ell^1$-norm at most $r$ is asymptotic to $4c_N\,r^{\sN}$. More generally, for every $d\ge2$, let $\Gamma^{(d)}\subset SL_{d+1}(\Z)$ be the semigroup generated by the matrices $\id+(\one-e_j)e_j^*$, $j=1,\dots,d+1$, and $N_d(r)=\#\{I\in\cI^{d+1}:\ \norm{A_I}_1\le r\}$. Then there is $c_d>0$ such that $N_d(r)\sim c_d\,r^{\sN(\Gamma^{(d)})}$ as $r\to\infty$.
\end{restatable}

For $d=2$ and $d=3$ the semigroups $\Gamma^{(d)}$ are $\GN$ and the semigroup $\GG$ of Main result~4, so $N=N_2$ and $c_N=c_2$. The limit set of the projective action of $\Gamma^{(d)}$ on the $d$-simplex $\cT_d$ is the $d$-dimensional Levitt--Yoccoz gasket $\mathcal G^d$ of~\cite{DL08,DL15JGA} (Section~\ref{sec:MH}; see also~\cite{AS13,Fou20}), so $\cN=\mathcal G^2$, and the gasket $\cG$ of Main result~4 is $\mathcal G^3$.

This is the analogue, for the souls, of the count of the circles of an Apollonian packing by Kontorovich and Oh~\cite{KO11}: there the number of circles of curvature at most $r$ is asymptotic to a constant times $r$ raised to the Hausdorff dimension of the residual set. Here the counting exponent is not the dimension of the gasket, but it bounds it from both sides (Theorem~\ref{thm: sN}). By Theorem~\ref{thm:count}, the answer to the question raised in~\cite[Open Question~1]{DL12}, whether the number of elements of norm at most $r$ of a semigroup with finite norm exponent is always strongly asymptotic to $r^{\sN}$, is positive for $\GN$ and the norm $\norm\cdot_1$. In general it is negative: already for the Sierpi\'nski gasket, whose norms are exactly multiplicative, the normalized count oscillates between two values whose ratio is $3$ (\cite[Motivational Example~3]{DL12} and Proposition~\ref{prop:lattice}), and the answer may depend on the norm (Remark~\ref{rem:norms}). The proof of~\cite{GMR19} is the renewal method of Lalley~\cite{Lal89}, applied to a transfer operator whose steps are the runs of equal generators, the same steps that we use for the certificates of Theorems~\ref{thm: dimH} and~\ref{thm: delta}; what distinguishes $\GN$ from the Sierpi\'nski gasket is that its run factors are non-lattice (Remark~\ref{rem:GMRproof}).
Numerically, counting all $1.6\times10^{11}$ multiindices with $\norm{A_I}_1\le 2^{17}$, the
local slopes $\log_2(N(2r)/N(r))$ decrease to $2.4437$
and a fit gives $c_N\approx0.0511$ (Remark~\ref{rem:countnum}).

The identification of the norm exponent with Baragar's exponent holds in every dimension: $\sN(\Gamma^{(d)})$ is the exponent of the Markoff--Hurwitz equation in $d+2$ variables (Lemma~\ref{lem:baragartree} and Theorem~\ref{thm:baragar}). Baragar observed that his bounds $2.963<\sN(\Gamma^{(4)})<3.048$ raise the question whether the exponent for six variables is exactly $3$~\cite[\S4]{Bar98}. With the run operator of Section~\ref{sec:fast} and a polynomial test function, we answer this question.

\begin{restatable}[Baragar's question]{theorem}{thmalphafive}\label{thm:alpha5}
$\sN(\Gamma^{(4)})>3.0002$. Hence, if $a\in\Z$ and the equation $x_1^2+\dots+x_6^2=a\,x_1\cdots x_6$ has an integer solution other than $0$, the number of its integer solutions with $\max_i\abs{x_i}\le R$ is $(\log R)^{\alpha+o(1)}$ with $\alpha>3$.
\end{restatable}

Numerically $\sN(\Gamma^{(4)})\approx3.0051$. The value $3$ is also the threshold below which the hypothesis $n-1<\sN$ of Theorem~\ref{thm:cutout} fails for the four-dimensional gasket $\mathcal G^4$ (Remark~\ref{rem:higher}).

\medskip
Theorems~\ref{thm: sN} and~\ref{thm: delta} give at once the rigorous bounds $1.6274\le\dimH\cN\le1.7233$, which already improve the lower and upper bounds previously available in the literature, recalled above. Norm bounds cannot reach the precision of Theorem~\ref{thm: dimH}: the two bounds of Theorem~\ref{thm: sN} come from the pointwise bounds $\sigma_1(A_I)^{-1/2}\le\sigma_2(A_I)\le3$ on the middle singular value, and each of them is attained along an explicit family of elements (Lemma~\ref{lem:sv} and Remark~\ref{rem:ceiling}). The norm exponent remains the right tool in four respects: it is the counting exponent of the souls, a quantity of Novikov's problem in its own right, and also of the Markoff--Hurwitz quadruples and of the one-sided simple closed geodesics of thrice-punctured projective planes; its computation needs a state space of the same dimension as the gasket, while the affinity exponent needs additional variables, so it scales to higher dimension; it gives a plausible explanation of the value $\approx1.72$ suggested by the numerics of~\cite{DD09} (Remark~\ref{rem:172}); and it is the entropy per unit time of the natural measure of the dynamics of chaotic sections (Main result~5 below).

\medskip
\subsection*{Main result 4: the three-dimensional Levitt--Yoccoz gasket}
This result concerns the analogue of the Levitt--Yoccoz gasket in dimension 3.
Let $\GG$ be the semigroup of non-negative matrices of $SL_4(\Z)$ generated by the four matrices
\[
  A_j=\id+(\one-e_j)e_j^*,\qquad j=1,\dots,4,
\]
(the generators of $\GN$ are the $3\times3$ matrices of the same form), acting projectively, by $\psi_I(x)=\pr{A_Ix}$, on the tetrahedron $\cT_3=\{x\in\R^4:\ \one^*x=1,\ 0\le x_m\le\frac13\}$, whose vertices are $\frac13(\one-e_m)$, and let
\[
  \cG=\bigcap_{n\ge0}\ \bigcup_{\abs I=n}\psi_I(\cT_3)
\]
be its limit set, the \textbf{three-dimensional Levitt--Yoccoz gasket} (Figure~\ref{fig:gaskets}(b)). Up to a projective change of coordinates, $\cG$ is the case $n=3$ of the fractals $F(\mathcal E)$ of~\cite{DL08}, which is also the cubic gasket in $\RP^3$ of~\cite{DL15JGA}: $F(\mathcal E)$ is the image of $\cG$ under the projective map with matrix $\mathcal E(\one\one^*-3\id)$ (Section~\ref{sec:MH}), and for the basis $\mathcal E_T$ of \cite[\S4.3]{DL08} it is, in the chart $x_4=1$, the image of $\cG$ under the similarity $x\mapsto\one-3x$ of the plane $\one^*x=1$, followed by the projection on the first three components.
See also~\cite{AS13,Fou20}; in~\cite{Fou20}, up to the same similarity, $\cG$ is called the Rauzy gasket of the $3$-simplex and denoted $\mathcal G^3$.
Its complement in $\cT_3$ is the disjoint union of the sets $Z_I=\psi_I(\Oct)$, where $\Oct=\{x\in\cT_3:\ x_m>\frac16,\ m=1,\dots,4\}$ is the octahedron with vertices at the midpoints of the edges of $\cT_3$, without its four faces that face the vertices of $\cT_3$ (Lemma~\ref{lem:oct}).
It is known that $\cG$ has zero volume (Theorem~\ref{thm:KM}; see also~\cite[Thm.~1]{DL08}), that $\dimH\cG<3$ \cite{Fou20,PS24} (sketched earlier in~\cite[\S11]{AHS16}), and, explicitly, that $\dimH\cG<2.7$ \cite{Fou20} (preprint); a box count of the barycentres of the zones gave the numerical estimate $2.20$ in~\cite{DL08}, where it was suggested that $\dimH\cG>2$.

\begin{restatable}{theorem}{thmthreed}
    \label{thm:3d}
    The norm exponent of $\GG$ satisfies $2.74\le\sN(\GG)\le2.777$, and
    \[
      1.6853\ <\ \dimH\cG\ \le\ \udimB\cG\ \le\ 1+\tfrac12\,\sN(\GG)\ \le\ 2.3885 .
    \]
\end{restatable}

The upper bound is Theorem~\ref{thm:cutout}(a) with $n=3$: the zone $Z_I$ has volume of order $\|A_I\|_1^{-4}$ and surface area $O(\|A_I\|_1^{-2})$ (Lemma~\ref{lem:geom3}). The lower bound only uses that $\cG$ contains a copy of $\cN$ in each face of $\cT_3$, together with Theorem~\ref{thm: dimH}. Theorem~\ref{thm:cutout}(b) does not apply to $\cG$, because the closures of its zones meet along faces (Remark~\ref{rem:octfaces}).
Numerically $\sN(\GG)\approx2.7578$. The bounds for $\sN(\GG)$ improve Baragar's bounds $2.730<\sN(\GG)<2.798$ for the exponent of the Markoff--Hurwitz equation in five variables~\cite[Thm.~0.1]{Bar98} (Lemma~\ref{lem:baragartree}), and by~\cite[Thm.~13]{GMR19} the number of $I$ with $\norm{A_I}_1\le r$ is asymptotic to $c_3\,r^{\sN(\GG)}$ for some $c_3>0$ (Theorem~\ref{thm:count}). A shell estimate of the affinity exponent gives $\sA(\GG)\approx2.26$, and a box count of $\cG$ gives $\approx2.21$; each value varies by less than $0.01$ across scales, and we do not know the origin of the gap between them (Remark~\ref{rem:3dnum}). By Proposition~\ref{prop:main} with $n=3$ and Theorem~\ref{thm:3d}, $\sA(\GG)\ge\frac34\sN(\GG)\ge2.055$; so an analogue in $\RP^3$ of the equality $\dimH=\sA$ of \cite{JLPX23,Jur26}, which is not available, would give $\dimH\cG>2$.

\medskip
\subsection*{Main result 5: diffusion of chaotic sections}
Avila, Hubert and Skripchenko~\cite{AHSinv} introduced an invariant probability measure $\mu_0$ on $\cN$, which they identify as the measure of maximal entropy of a renormalization flow~\cite[Thm.~25]{AHSinv}, and proved that, for $\mu_0$-almost every chaotic regime, the trajectories of the plane sections diffuse with an exponent $\delta=-\lambda_3/\lambda_1\in(\frac12,1)$, where $\lambda_1>\lambda_2>\lambda_3$ are the Lyapunov exponents of their cocycle (the product of the matrices of the successive steps, each depending only on the current point). Their flow is a suspension over the first return, to the set of the points whose induction begins with a fixed complete path, of an accelerated induction whose steps are, in our notation, the runs of equal generators. We show that, after the projective change of coordinates of Remark~\ref{rem:conj}, it is the suspension flow of the shift of the runs, with the logarithm of our run factors as roof function, seen on the corresponding cross-section, and that $\sN$ is the maximal entropy per unit time of its return map: the maximum, over the invariant probability measures $\nu$ of the return map with $\int r\,d\nu<\infty$, where $r$ is the roof function, of the entropy $h_\nu$ divided by the mean return time $\int r\,d\nu$, attained only at $\nu=\mu_0$ (Proposition~\ref{prop:AHS}). By Abramov's formula~\cite{Abr59b}, $h_\nu/\int r\,d\nu$ is the entropy of the flow for the measure induced by $\nu$ when the flow is built over an invertible map; the return map of~\cite{AHSinv} is not invertible (Remark~\ref{rem:maxentflow}).
Since $\sN<3$, this value lies outside the range in which~\cite{AHSinv} prove that the pressure is finite, and the existence of $\kappa_0$ and $\mu_0$, claimed in~\cite[Lemma~24]{AHSinv}, needs a proof: for the accelerated induction, before passing to the first returns, it follows from results of Gamburd, Magee and Ronan~\cite{GMR19}, who also noted that their conformal measure is essentially $\mu_0$~\cite[Remark~12]{GMR19}, and Proposition~\ref{prop:AHS} completes it for the first-return system of~\cite{AHSinv} (see Section~\ref{rem:AHS}).
Combined with Theorem~\ref{thm: dimH}, this gives a quantitative bound.

\begin{restatable}{theorem}{thmdiffusion}\label{cor:diffusion}
The diffusion exponent $\delta$ of~\cite{AHSinv} satisfies
\[
  \delta\ \le\ \frac{1+\sA-\sN}{2-\sA}\ <\ 0.7766 .
\]
\end{restatable}

A heuristic numerical estimate gives $\delta\approx0.74$ (Remark~\ref{rem:diffnum}).

\subsection*{Organization}
Section~\ref{sec:souls} identifies the zones and souls with the elements of $\GN$, studies the dynamics of $\PGN$ and proves Theorem~\ref{cor:singlenode}; as a corollary, it gives, for $\cN$, a direct proof of the accumulation statement of~\cite{DL03,DL05} (Corollary~\ref{cor:accumulation}).
Section~\ref{sec:affinity} relates $\sA(\GN)$ to the shape of the zones and proves the cut-out bound $\udimB\cN\le\sA$.
Section~\ref{sec:cocycle} introduces the run factors and the needle, and reduces bounds on $\sA$ to positive test functions for a transfer operator.
Section~\ref{sec:disc} reduces the operator by symmetry and gives a pointwise certificate for continuous piecewise-linear test functions, and Section~\ref{sec:comp} describes the computation and proves Theorem~\ref{thm: dimH}.
Section~\ref{sec:general} proves Theorem~\ref{thm:cutout}, the inequality between the affinity and norm exponents, and Theorem~\ref{thm: sN}.
Section~\ref{sec:fast} treats fast semigroups and the growth rate per run of their zeta function, and rewrites the zeta function of $\GN$ by splitting the multiindices into runs: the $\ell^1$-size of a soul is the product of its run factors, which turns $\zeta$ into the sum of the iterates of a transfer operator on the first-level cylinders. It also explains why this works: $\GN$ is a fast semigroup, with sharp junction constant $\frac14$, and the run factors compute exactly the junction factors that \cite[Thm.~1]{DL15JGA} bounds by the fastness constant.
Section~\ref{sec:thmB} proves Theorem~\ref{thm: delta}.
Section~\ref{sec:MH} puts $\GN$ and $\GG$ in the family of the semigroups $\Gamma^{(d)}$ generated by the matrices $\id+(\one-e_j)e_j^*$ in every dimension, whose limit sets are the $d$-dimensional Levitt--Yoccoz gaskets of~\cite{DL08,DL15JGA}, identifies their norm exponents with Baragar's exponents of the Markoff--Hurwitz equations, deduces Theorem~\ref{thm:count} from the theorem of Gamburd, Magee and Ronan, and contrasts it with the lattice case of the Sierpi\'nski gasket.
Section~\ref{sec:3d} treats the three-dimensional Levitt--Yoccoz gasket and proves Theorem~\ref{thm:3d}. Section~\ref{sec:AHS} relates $\sN$ to the measure of~\cite{AHSinv} and proves Theorem~\ref{cor:diffusion}.
The programs used for the certificates, the logs of the runs and the test functions are deposited at Zenodo~\cite{DLdata}.

\section{Zones, souls and the semigroup}\label{sec:souls}

In this section we describe the zones and the souls of $\cN$ through the semigroup $\GN$, prove the geometric estimates on the zones used throughout the paper (Lemma~\ref{lem:geom}), and study the dynamics of $\PGN$ (Theorem~\ref{cor:singlenode}).
The description of the zones and of the souls is in~\cite{DD09,DL08,DL15EM}, and we recall it for the reader's convenience; Lemma~\ref{lem:geom}, Proposition~\ref{lem: Dpsi} and Theorem~\ref{cor:singlenode} are new.

\begin{definition}
    Let $m>0$. 
    We denote by $\cI^m$ the set of multiindices in $\{1,\dots,m\}$.
    Given a metric space $X$, an {\bf Iterated Function System \boldmath $\cF$} on $X$ is a family $\{f_1,\dots,f_m\}$ of continuous maps $f_i:X\to X$ acting on $X$ by the left semigroup action
    \[
    \Phi_\cF:\cI^m\times X\to X
    \]
    defined by 
    \[
    \Phi^\emptyset_\cF=id_X,\quad\Phi_\cF^I=f_I,
    \]
    where $I=i_1\dots i_r\in\cI^m$ and $f_I(x)=f_{i_1}\circ\cdots\circ f_{i_r}(x)$.
    We say that $\cF$ is {\bf generated} by the $f_i$.

    We denote by {\boldmath$\mathcal K(X)$} the space of all compact non-empty subsets of $X$ endowed with the Hausdorff distance {\boldmath$d_H$} and, given an IFS $\cF$, we call the {\bf Hutchinson map} of $\cF$ the map
    \[
    \text{\boldmath$H_\cF$}:\mathcal K(X)\to\mathcal K(X)
    \] 
    defined by $H_\cF(A)=f_1(A)\cup\dots\cup f_m(A)$.
    We say that a set $A$ is {\bf\boldmath $H$-invariant} under $\cF$ if $H_\cF(A)=A$.
\end{definition}
Notice that, when there are no relations among the $f_i$, the action of $\cF$ coincides with the action of the semigroup of continuous self-maps of $X$ freely generated by the $f_i$.

\medskip
We denote by $e_1,e_2,e_3$ the canonical basis of $\R^3$ and by $\|\cdot\|_1$ both the $\ell^1$ norm for vectors and the entrywise $\ell^1$ norm for matrices.
We denote by $\abs{\cdot}$ the Euclidean norm, by $x\times y$ the cross product, so that $\abs{x\times y}=\abs x\abs y\sin\angle(x,y)$, and by $A^*$ the transpose of a matrix or vector $A$.

The largest stability zones are the square with vertices at $[e_1\pm e_2]$, $[e_1\pm e_3]$, centered at $[e_1]$ and with soul $e_1$, and the other similar squares centered respectively at $[e_2]$ and $[e_3]$.
All other zones are obtained in the following way.

\begin{figure}[t]
  \centering
  \includegraphics[width=0.8\textwidth]{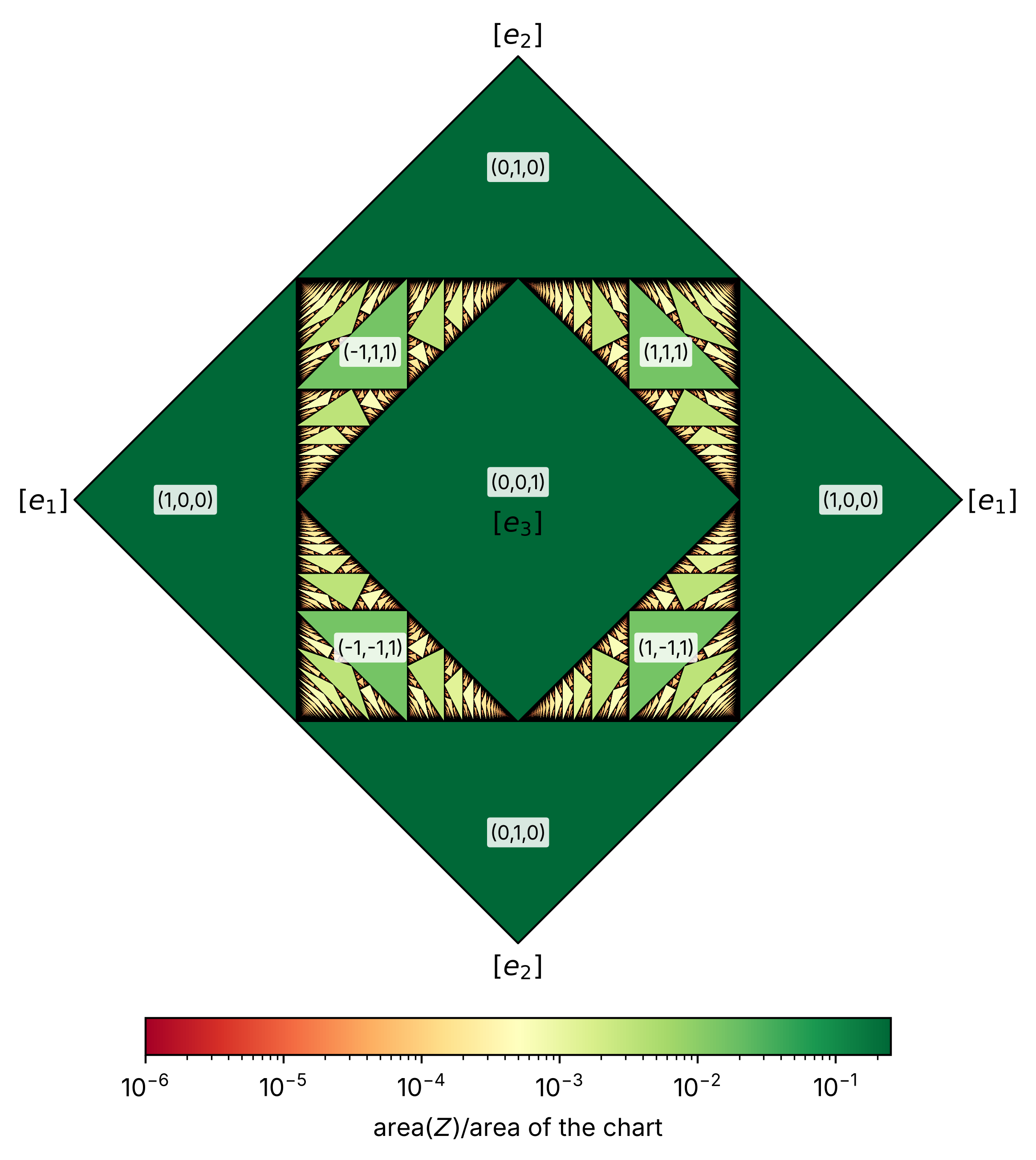}
  \caption{The stability zones of the $\mu$-cube $\{4,6\,|\,4\}$ in $\mathbb{RP}^2$,
    drawn in the octahedral chart $[v]\mapsto (v_1,v_2)/\|v\|_1$, $v_3\ge 0$: the square
    $|x|+|y|\le 1$ with antipodal points of its boundary identified. On each of its four
    quadrants the chart is projective, so the zones are drawn with straight sides.
    The three square zones have souls $e_1$, $e_2$, $e_3$ and are centred at $[e_1]$,
    $[e_2]$, $[e_3]$; the squares about $[e_1]$ and $[e_2]$ are each split into two
    halves, glued by the antipodal identification. Each quadrant contains a copy of the
    triangle $\mathcal T$, the image of $\mathcal T$ under a change of sign of the
    coordinates; its largest zone is labelled by its soul $(\pm1,\pm1,1)$.
    Zones are coloured from green to red as their area, relative to that of the chart,
    decreases. All zones of diameter at least $1.5\cdot10^{-3}$ are drawn; the black
    set is the cut-out fractal $K_G$ at this resolution, the union of the
    Levitt--Yoccoz gasket $\mathcal N\subset\mathcal T$ and of its three images under
    the sign changes.}
  \label{fig:zonesRP2}
\end{figure}

Denote by $\GN$ the IFS generated by the linear maps $A_i=\id+(\one-e_i)e_i^*$, where $\one=e_1+e_2+e_3$.
In coordinates
\[
  A_1=\begin{pmatrix}1&0&0\\1&1&0\\1&0&1\end{pmatrix},\quad
    A_2=\begin{pmatrix}1&1&0\\0&1&0\\0&1&1\end{pmatrix},\quad
  A_3=\begin{pmatrix}1&0&1\\0&1&1\\0&0&1\end{pmatrix}.
\]
Since there are no relations among the $A_i$ (see Corollary~\ref{cor: free}), the action of $\GN$ coincides with the action of the subsemigroup of $SL_3(\Z)$ generated by the $A_i$, which we still denote by $\GN$. 
Finally, denote by $\psi_i:\RP^2\to\RP^2$ the corresponding projective automorphisms, defined by $\psi_i([x])=[A_ix]$.
For $I=i_1\dots i_r\in\cI^3$ we write $A_I=A_{i_1}\cdots A_{i_r}$, so that $\psi_I([x])=[A_Ix]$; $A_\emptyset=\id$.
For a permutation $\pi$ of $\{1,2,3\}$ we denote by the same letter the permutation matrix with $\pi e_j=e_{\pi(j)}$, so that $(\pi x)_{\pi(j)}=x_j$ and $\pi A_j\pi^{-1}=A_{\pi(j)}$; hence $\pi\circ\psi_j=\psi_{\pi(j)}\circ\pi$. We use the same notation for the permutations of $\{1,\dots,4\}$ in Section~\ref{sec:3d}.
Notice that
\[
\psi_1([e_1])=[\one],\;
\psi_1([e_2])=[e_2],\;
\psi_1([e_3])=[e_3],\;
\psi_1([\one])=[1:2:2]
\]
and similarly for the other two.

All zones in the first quadrant are obtained through the action of $\GN$ on the triangle $\cT$ with vertices $[e_1+e_2], [e_2+e_3], [e_3+e_1]$, and the Levitt--Yoccoz gasket $\cN$ is the complement in $\cT$ of their union.
All other zones are the images of these under the sign changes of the coordinates (Figure~\ref{fig:zonesRP2}), so we can focus our attention solely on $\cT$.

Notice that $\cT$ is invariant under $\GN$ and that $\psi_i$ keeps its vertex $[\one-e_i]$ fixed.
We denote by $\PGN$ the IFS (i.e., in this case, semigroup) on $\cT$ generated by the restrictions to $\cT$ of the $\psi_i$, which we still denote by $\psi_i$.
As for the contractions of an equilateral triangle that generate the Sierpi\'nski gasket, each $\psi_i(\cT)$ is the corner of $\cT$ at the vertex fixed by $\psi_i$, cut off by a side of the central triangle described below, and this decomposition is repeated inside every $\psi_I(\cT)$ (Lemma~\ref{lem:zones}); so the two cut-out fractals are combinatorially the same object.
As we will see, their geometry is, on the contrary, quite different.

More in detail, $\psi_i(\mathcal T)$ is the triangle with vertices $[\mathbf 1-e_i],[\mathbf 1+e_j],[\mathbf 1+e_k]$. 
The complement of $\bigcup_i\psi_i(\mathcal T)$ is the open triangle with vertices $[\mathbf 1+e_1],[\mathbf 1+e_2],[\mathbf 1+e_3]$, i.e. $[2:1:1],[1:2:1],[1:1:2]$, which is the largest stability zone inside $\cT$ which we denote by $Z_\emptyset$.
Its soul is the vector $\one$.

More generally, for every multiindex $I$ we call the triple of primitive
vectors $a_{I,k}=A_I(\one-e_k)$, $k=1,2,3$, which span the vertices of $\psi_I(\mathcal T)$, the
\emph{frame} of $I$. 
The stability zone $Z_I=\psi_I(Z_\emptyset)$ is the open triangle with vertices $[a_{I,j}+a_{I,k}]$,
$\{i,j,k\}=\{1,2,3\}$, and its soul is half the sum of the frame,
\[
  L_I=\tfrac12\,(a_{I,1}+a_{I,2}+a_{I,3})=A_I\one .
\]
That these, together with the three squares and the images of the triangles under the sign changes, are all the stability zones is the following theorem of Dynnikov and the author.

\begin{lettertheorem}[De Leo--Dynnikov~{\cite[Thm.~1]{DD09}}; see also~\cite{DD07}]\label{thm:DD}
The stability zones of the $\mu$-cube $\{4,6\,|\,4\}$ are the interiors of the three squares $\{[h]:\ \abs{h_i}\ge\abs{h_j}+\abs{h_k}\}$, $\{i,j,k\}=\{1,2,3\}$, with souls $e_1,e_2,e_3$; the open triangles $\psi_I(Z_\emptyset)$, $I\in\cI^3$, with souls $A_I\one$; and the images of these triangles under the sign changes of the coordinates.
\end{lettertheorem}

By Theorem~\ref{thm:DD}, the stability zones in $\cT$ are the $Z_I$, the orbit of $Z_\emptyset$ under the action of $\PGN$, and their souls are the orbit of $\one$ under $\GN$. We use this identification throughout, in particular in the proof of Lemma~\ref{lem:zones}.

\begin{figure}
  \centering
  \includegraphics[width=\textwidth]{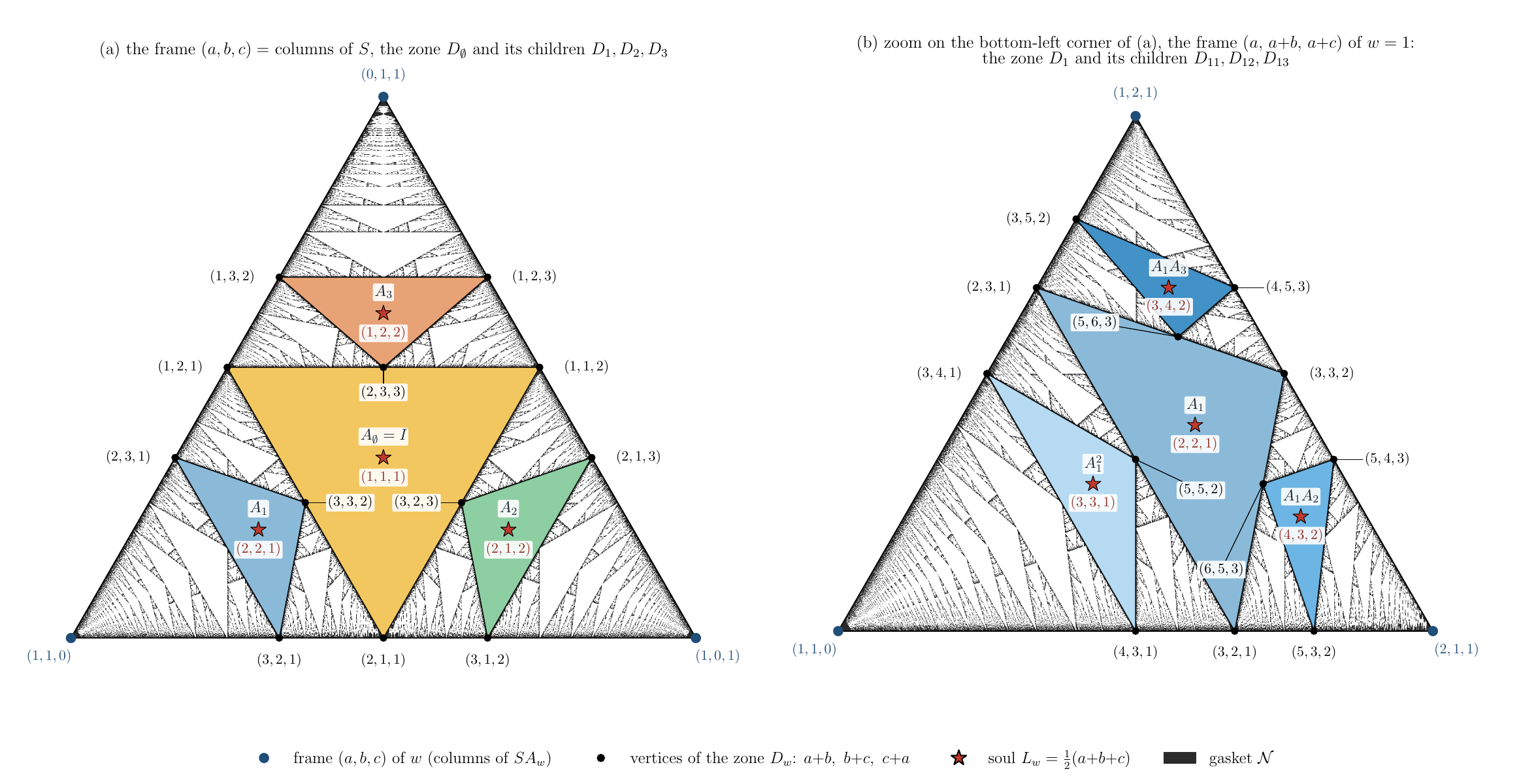}
  \caption{The Levitt--Yoccoz gasket $\mathcal N$ (black) in the triangle $\mathcal T$,
  drawn in the $\ell^1$ chart $\pr x=x/\one^*x$. Each highlighted zone $Z_I$ is labelled
  by its matrix $A_I$. Blue discs: the frame $a_{I,k}=A_I(\one-e_k)$, $k=1,2,3$, of $I$,
  which spans the vertices of the cylinder $\psi_I(\mathcal T)$; black discs: the
  vertices $\pr{a_{I,j}+a_{I,k}}$, $j\neq k$, of $Z_I$; stars: the soul point
  $\pr{L_I}$, where $L_I=\frac12(a_{I,1}+a_{I,2}+a_{I,3})=A_I\one$.
  (a)~The frame $(\one-e_1,\one-e_2,\one-e_3)$ of the empty multiindex (the vertices of
  $\mathcal T$), the central zone $Z_\emptyset$ and its children $Z_1,Z_2,Z_3$.
  (b)~The corner of (a) at the vertex $\pr{\one-e_1}$. By~\eqref{eq:tree}, the frame of
  $I=1$ is $(a,a+b,a+c)$, where $(a,b,c)$ is the frame of the empty multiindex, and the
  same rule produces $Z_1$ and the zones $Z_{11},Z_{12},Z_{13}$. Iterating the rule
  removes all the zones; what remains is $\mathcal N$.}
  \label{fig:gasket}
\end{figure}

\subsection*{The geometry of $\cN$}
From now on we work in the $\ell^1$ chart $\pr x=x/\one^*x$ (so $\pr x=x/\norm x_1$ for $x\ge0$, $x\neq0$): distances and areas are Euclidean, in the plane $\{\one^*x=1\}$. A non-negative invertible matrix $g$ acts by $x\mapsto\pr{gx}$, and it maps the triangle with vertices $\pr p,\pr q,\pr r$, $p,q,r\ge0$, onto the triangle with vertices $\pr{gp},\pr{gq},\pr{gr}$.

Since $A_j(\one-e_j)=\one-e_j$ and $A_j(\one-e_k)=(\one-e_j)+(\one-e_k)$ for $k\neq j$, the frames are produced by the recursion
\begin{equation}\label{eq:tree}
  a_{Ij,j}=a_{I,j},\qquad a_{Ij,k}=a_{I,j}+a_{I,k}\quad(k\neq j),
\end{equation}
rooted at $(\one-e_1,\one-e_2,\one-e_3)$: up to the order, a step replaces a triple $(x,y,z)$ by $(x,x+y,x+z)$. Similarly, since $A_jA_I=A_I+(\one-e_j)\,e_j^*A_I$, prepending $j$ to $I$ adds the $j$-th row of $A_I$ to the other two rows, so the rows of the $A_I$ are produced by the same rule from the rows of $\id$.

The cylinder $\psi_I(\cT)$ is the closed triangle with vertices $\pr{a_{I,1}},\pr{a_{I,2}},\pr{a_{I,3}}$, and $Z_I\subset\psi_I(\cT)$ is the open triangle with vertices $\pr{A_I(\one+e_i)}=\pr{a_{I,j}+a_{I,k}}$, $\{i,j,k\}=\{1,2,3\}$. Since $A_I\one=\frac14\sum_iA_I(\one+e_i)$, the zone $Z_I$ contains the \textbf{soul point}
\[
  \pr{L_I}=\pr{A_I\one}\in Z_I ,
\]
which is the soul direction $[L_I]$ in the chart. Finally, $a_{I,1}+a_{I,2}+a_{I,3}=2A_I\one$, so
\[
  \norm{a_{I,1}}_1+\norm{a_{I,2}}_1+\norm{a_{I,3}}_1=2\norm{A_I}_1 .
\]

\begin{lemma}\label{lem:zones}
    The triangles $Z_I$ are pairwise disjoint and $\cN=\bigcap_n\bigcup_{\abs I=n}\psi_I(\cT)$. More precisely, for every $n\ge0$,
\[
  \cT=\bigcup_{\abs I<n}Z_I\ \cup\ \bigcup_{\abs I=n}\psi_I(\cT),
\]
where the closed triangles $\psi_I(\cT)$, $\abs I=n$, have pairwise disjoint interiors and the $Z_I$, $\abs I<n$, are pairwise disjoint and disjoint from all of them.
\end{lemma}

\begin{proof}
    The triangle $\psi_j(\cT)$, with vertices $\pr{\one-e_j},\pr{\one+e_k},\pr{\one+e_l}$, is the corner of $\cT$ at $\pr{\one-e_j}$ cut off by a side of $Z_\emptyset$.
    Hence $\cT=Z_\emptyset\cup\psi_1(\cT)\cup\psi_2(\cT)\cup\psi_3(\cT)$, the three closed triangles having pairwise disjoint interiors and not meeting $Z_\emptyset$.
    Each $\psi_I$ is a homeomorphism of $\RP^2$ mapping $\cT$ into itself, so applying this decomposition inside every $\psi_I(\cT)$, $\abs I=n$, gives the displayed decomposition by induction. If $y\in\cT\setminus\cN$, then $y$ lies in some zone $Z_J$, which does not meet $\bigcup_{\abs I=\abs J+1}\psi_I(\cT)$. Conversely, if $y\notin\bigcup_{\abs I=n}\psi_I(\cT)$ for some $n$, then $y\in Z_I$ for some $\abs I<n$.
\end{proof}

\begin{lettertheorem}[Meester--Nowicki~{\cite[Thm.~5.1]{MN89}}; see also~{\cite[Lemme~V.9]{Lev93}} and~{\cite[Thm.~2]{DD09}}]\label{thm:yoccoz}
$\mathrm{area}(\cN)=0$.
\end{lettertheorem}

\begin{lemma}
    \label{lem:geom}
    Let $a,b,c$ be the frame of $I$ and $r_1,r_2,r_3$ the rows of $A_I$.
    \begin{enumerate}
    \item[(i)] For every permutation $(x,y,z)$ of $(a,b,c)$, $\|x\|_1\le \|y\|_1+\|z\|_1$; equivalently, $\|x\|_1\le\|A_I\|_1\le\|y\|_1+\|z\|_1$\\ (since $\pr{L_I}=\bigl(\|a\|_1\pr a+\|b\|_1\pr b+\|c\|_1\pr c\bigr)/(2\|A_I\|_1)$, this says that the soul point $\pr{L_I}$ lies in the Euclidean medial triangle of the cylinder $\psi_I(\cT)$);
    \item[(ii)] for every permutation $(x,y,z)$ of $(a,b,c)$, $\|x\times y\|_1+\|x\times z\|_1\le3\|x\|_1$, and for every permutation $(x,y,z)$ of $(r_1,r_2,r_3)$, $\|x\times y\|_1+\|x\times z\|_1\le2\|x\|_1$; in particular $\norm{A_I^{-1}}_1\le\|A_I\|_1$, with equality when $I=j^k$;
    \item[(iii)] $\mathrm{area}(Z_I)=\frac{2\sqrt3}{(\|a\|_1+\|b\|_1)(\|b\|_1+\|c\|_1)(\|c\|_1+\|a\|_1)}$, so
    \[
    \tfrac{\sqrt3}4\,\|A_I\|_1^{-3}\le\mathrm{area}(Z_I)\le2\sqrt3\,\|A_I\|_1^{-3};
    \]
    \item[(iv)] $\abs{\pr x-\pr y}\le6\sqrt3\,\|A_I\|_1^{-1}$ for $x,y$ in the frame, so
  \[\|A_I\|_1^{-3/2}\le\mathrm{diam}\,Z_I\le\mathrm{diam}\,\psi_I(\cT)\le6\sqrt3\,\|A_I\|_1^{-1}\]
  and
  \[3\,\|A_I\|_1^{-3/2}\le\mathrm{per}(Z_I)\le\mathrm{per}(\psi_I(\cT))\le18\sqrt3\,\|A_I\|_1^{-1};\]
  \item[(v)] $\|A_{Ij}\|_1=\|A_I\|_1+\|a_{I,j}\|_1\ge\|A_I\|_1+2$ for $j=1,2,3$; hence $\|A_I\|_1\ge3+2\abs I$.
\end{enumerate}
\end{lemma}

\begin{proof}
    (i) and (ii) are proved by induction along \eqref{eq:tree}, and along the same rule for the rows.
    At the root, the frame vectors have $\|x\|_1=2$ and $\norm{(\one-e_i)\times(\one-e_j)}_1=3$, and the rows have $\|x\|_1=1$ and $\norm{e_i\times e_j}_1=1$. A step replaces the triple by $(x,x+y,x+z)$.
    For (i), $\|x\|_1\le \|x+y\|_1+\|x+z\|_1$ is clear, and
\[
\|x+y\|_1=\|x\|_1+\|y\|_1\le \|x\|_1+(\|x\|_1+\|z\|_1)=\|x\|_1+\|x+z\|_1,
\]
and symmetrically.
For (ii), with $C=3$ for the frame and $C=2$ for the rows,
\[
x\times(x+y)=x\times y,\;
x\times(x+z)=x\times z\text{ and }(x+y)\times(x+z)=x\times z+y\times x+y\times z;
\]
so the left side of (ii) does not change for $x$, while for $x+y$ it is at most $\norm{x\times y}_1+\norm{x\times z}_1+\norm{y\times x}_1+\norm{y\times z}_1\le C\|x\|_1+C\|y\|_1=C\|x+y\|_1$, and the same for $x+z$.
The two forms of (i) are equivalent because $\|a\|_1+\|b\|_1+\|c\|_1=2\|A_I\|_1$, and the description of $\pr{L_I}$ holds because $A_I\one=\frac12(a+b+c)$ and $x=\|x\|_1\pr x$; a point of a triangle lies in its medial triangle if and only if all its barycentric coordinates are at most $\frac12$.
Since $\det A_I=1$, the columns of $A_I^{-1}$ are $r_2\times r_3$, $r_3\times r_1$ and $r_1\times r_2$; summing (ii) for the rows over $x=r_1,r_2,r_3$ counts each of them twice, so $\norm{A_I^{-1}}_1\le\|A_I\|_1$. For $I=j^k$, $A_j^{\pm k}=\id\pm k\,(\one-e_j)e_j^*$, and both norms equal $3+2k$.

(iii) For $p,q,r\ge0$ the triangle with vertices $\pr p,\pr q,\pr r$ has area
\[
\frac{\sqrt3}{2}\frac{\abs{\det(p,q,r)}}{\|p\|_1\|q\|_1\|r\|_1}:
\]
the determinant is six times the volume
of the tetrahedron with apex $0$ over the triangle, whose height is $1/\sqrt3$. Here
$\abs{\det(b+c,c+a,a+b)}=2\abs{\det(a,b,c)}=2\abs{\det(\one-e_1,\one-e_2,\one-e_3)}=4$ and $\|a+b\|_1=\|a\|_1+\|b\|_1$, etc.; the bounds
follow from (i), since $\|A_I\|_1\le \|y\|_1+\|z\|_1\le2\|A_I\|_1$.

    (iv) For $p,q\ge0$ non-zero, $\pr p-\pr q=\one\times(p\times q)/(\|p\|_1\|q\|_1)$, so
\[
\abs{\pr p-\pr q}\le\sqrt3\,\norm{p\times q}_1/(\|p\|_1\|q\|_1).
\]
    For two frame vectors, (ii) gives $\norm{x\times y}_1\le3\min(\|x\|_1,\|y\|_1)$, so that, by (i),
\[
\abs{\pr x-\pr y}\le3\sqrt3/\max(\|x\|_1,\|y\|_1)\le6\sqrt3/(\|x\|_1+\|y\|_1)\le6\sqrt3/\|A_I\|_1.
\]
    Hence every side of $\psi_I(\cT)$ has length at most $6\sqrt3/\|A_I\|_1$, and $Z_I\subset\psi_I(\cT)$ is convex, so its perimeter is at most that of $\psi_I(\cT)$.
For the lower bounds, let $Z$ be a triangle and $PQ$ a longest side. The third vertex lies within $\mathrm{diam}\,Z$ of both $P$ and $Q$, hence at distance at most $\frac{\sqrt3}2\,\mathrm{diam}\,Z$ from the line $PQ$, so $\mathrm{area}(Z)\le\frac{\sqrt3}4(\mathrm{diam}\,Z)^2$. Also $\mathrm{per}(Z)^2\ge12\sqrt3\,\mathrm{area}(Z)$, by Heron's formula and the AM--GM inequality. Both lower bounds now follow from (iii).

    (v) Since $A_j\one=\one+(\one-e_j)$, $\norm{A_{Ij}}_1=\one^*A_I\one+\one^*A_I(\one-e_j)=\norm{A_I}_1+\norm{a_{I,j}}_1$. Each $A_I$ is a product of non-negative matrices whose entries are at least those of $\id$, so $A_I\ge\id$ entrywise and $\norm{a_{I,j}}_1\ge\norm{\one-e_j}_1=2$. As $\norm{\id}_1=3$, induction gives $\norm{A_I}_1\ge3+2\abs I$.
\end{proof}

\begin{remark}[sharpness of (iv)]\label{rem:geomsharp}
Both exponents in Lemma~\ref{lem:geom}(iv) are attained. For the cusp zones, $\mathrm{diam}\,Z_{1^k}=1/(\sqrt2\,(k+2))$, the distance between $\pr{A_1^k(\one+e_2)}=\pr{(1,k+2,k+1)}$ and $\pr{A_1^k(\one+e_3)}=\pr{(1,k+1,k+2)}$, while $\|A_1^k\|_1=2k+3$. For $I=(123)^n$, $\mathrm{diam}\,Z_I\asymp\sigma_2(A_I)/\sigma_1(A_I)\asymp\|A_I\|_1^{-3/2}$ by Lemmas~\ref{lem:sv} and~\ref{lem:diamZ}. So no power of $\|A_I\|_1$ alone controls the size of the zones; the missing information is $\sigma_2(A_I)$ (Remark~\ref{rem:ceiling}). The point called the barycenter of $Z_I$ in~\cite{DL08} is the soul direction $[A_I\one]$, and $\|A_I\one\|_1=\|A_I\|_1$; so Lemma~\ref{lem:geom}(iv) proves Conjecture~1 of~\cite{DL08} for $n=2$, with optimal exponents.
\end{remark}

\subsection*{The dynamics of \texorpdfstring{$\PGN$}{the IFS}}
Here we briefly discuss the qualitative dynamics of $\PGN$, through the graph of~\cite{DLstreams} recalled in the introduction.
The essential condition is that the IFS has a global attractor: a compact set that is invariant under the Hutchinson map and attracts every compact subset.

We recall the relevant notions of~\cite{DLstreams}. 
An {\bf\boldmath$\varepsilon$-chain} from $x$ to $y$ is a finite sequence $x=x_0,x_1,\dots,x_k=y$ together with a sequence of indices $i_0,\dots,i_{k-1}$ such that, for every $r<k$, $d(\psi_{i_r}(x_r),x_{r+1})<\varepsilon$. 
A point $y$ is \emph{\bf\boldmath chain-downstream} from $x$ if, for every $\varepsilon>0$, there is an $\varepsilon$-chain from $x$ to $y$;  $x$ and $y$ are \emph{\bf\boldmath chain-equivalent} if each of them is chain-downstream from the other. 
Notice that every point on the orbit of $x$ is chain-downstream from $x$.
A point is \emph{\bf\boldmath chain-recurrent} if there is a $\psi_i$ such that $x$ is chain-downstream from $\psi_i(x)$ (and so $x$ is chain-equivalent to $\psi_i(x)$).

\begin{figure}
  \centering
  \includegraphics[width=\textwidth]{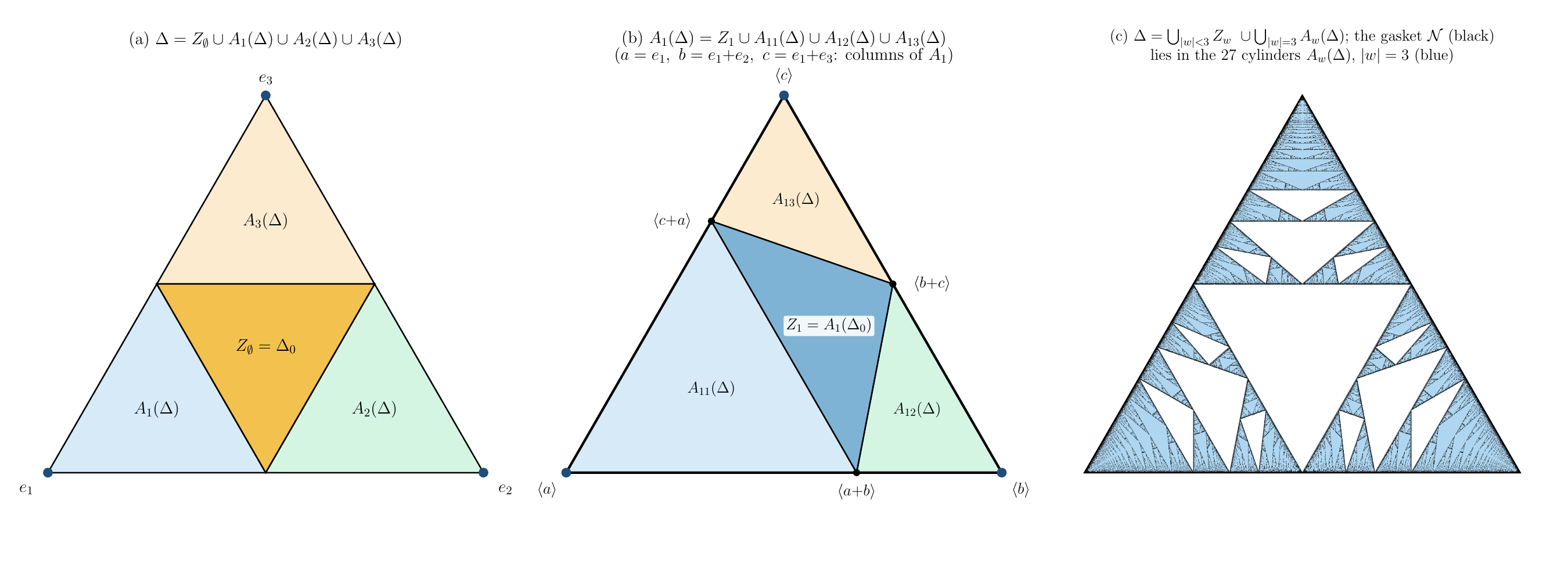}
    \caption{Cylinders $\psi_I(\mathcal T)$ and zones $Z_I=\psi_I(Z_\emptyset)$.
  (a)~$\mathcal T$ splits into $Z_\emptyset$ and the three cylinders $\psi_j(\mathcal T)$.
  (b)~The same splitting inside $\psi_1(\mathcal T)$, enlarged. (c)~After three steps:
  the zones $Z_I$, $|I|<3$ (white), and the cylinders $\psi_I(\mathcal T)$, $|I|=3$
  (blue), which contain $\mathcal N$ (black).}
  \label{fig:cylinders}
\end{figure}

Through chains we encode the qualitative behavior of an IFS into a graph in the following way.
The {\bf chain nodes} of the graph are the classes of chain-equivalent chain-recurrent points; the \textbf{chain graph} has the chain nodes as vertices, and there is an {\bf edge} from a chain node to a different chain node if the points of the second are chain-downstream from those of the first. 
Nodes generalize periodic orbits, and edges the trajectories connecting them. 
The same definitions, with the single map $\mathcal H$ in place of the $\psi_j$, 
give the chain graph of the discrete semiflow of $\mathcal H$ on $\mathcal K(\cT)$.

Most of the IFS literature is about contractive IFSs.
None of the $\psi_j$ is a contraction, though: $\psi_j$ fixes the vertex $\pr{\one-e_j}$, its derivative there is the identity, and so that vertex is a parabolic fixed point.
Hence, $\PGN$ is not covered by Hutchinson's theorem~\cite{Hut81}. It is, however, weakly contractive (Proposition~\ref{lem: Dpsi}), and weakly contractive IFSs on compact spaces are well understood~\cite{Kie02}; so items (i)--(iv) of Theorem~\ref{cor:singlenode} also follow from classical results (Remark~\ref{rem:edelstein}).
Our main new result is the structure of the chain graphs of $\PGN$ and of its Hutchinson map, items (v)--(vi) of Theorem~\ref{cor:singlenode}.

\medskip
In the $\ell^1$ chart, the tangent plane of $\{\one^*x=1\}$ is $T=\{v\in\R^3:\one^*v=0\}$, with the Euclidean norm. In this chart $\cT=\{x:\ \one^*x=1,\ 0\le x_i\le\frac12,\ i=1,2,3\}$ is the equilateral triangle with vertices $\pr{\one-e_i}=\frac12(\one-e_i)$, $i=1,2,3$. Since $A_jx=x+x_j(\one-e_j)$ and $\one^*A_jx=1+2x_j$ for $\one^*x=1$, in the chart the generators act on $\cT$ by
\begin{equation}\label{eq:fj}
  \psi_j(x)=\pr{A_jx}=\frac{x+x_j(\one-e_j)}{1+2x_j}=\pr{\one-e_j}+\frac{x-\pr{\one-e_j}}{1+2x_j},\qquad j=1,2,3 :
\end{equation}
$\psi_j$ moves every point of $\cT$ radially towards the vertex $\pr{\one-e_j}$, by the factor $1/(1+2x_j)$. In $\RP^2$, the set of points fixed by $\psi_j$ is the line $\{x_j=0\}$; in the chart, it is the line through $\pr{\one-e_j}$ parallel to the opposite side $\{x_j=\frac12\}$ of $\cT$, which meets $\cT$ only at $\pr{\one-e_j}$.

\medskip
Although the $\psi_i$ are not contractions, the next result shows that they shrink the distance between any two points of $\cT$.

\begin{proposition}
    \label{lem: Dpsi}
For $j=1,2,3$ and $x\in\cT$, let $D_x\psi_j:T\to T$ be the derivative of $\psi_j$ at $x$, $D_x\psi_j(v)$ its value on $v\in T$, and $\norm{D_x\psi_j}$ its Euclidean operator norm. Then $D_p\psi_j=\id$ at $p=\pr{\one-e_j}$, and $\norm{D_x\psi_j}<1$ for $x\neq p$. Consequently,
\[
  \abs{\psi_j(x)-\psi_j(y)}<\abs{x-y}\qquad\text{for all }x\neq y\text{ in }\cT .
\]
\end{proposition}

\begin{proof}
Let $p=\pr{\one-e_j}$, $w=x-p\in T$, $u=2x_j$ and $d=1+u$. Since $p_j=0$ we have $w_j=x_j$; if $\{j,l,m\}=\{1,2,3\}$, then $w_l,w_m\le0$ (because $x\in\cT$) and $w_l+w_m=-x_j$, so $\abs w^2\le2x_j^2=\frac12u^2$. Let $g=e_j-\frac13\one\in T$, so that $v_j=g^*v$ for $v\in T$ and $\abs g^2=\frac23$. Differentiating~\eqref{eq:fj},
\[
  D_x\psi_j=\frac1d\,M,\qquad M=\id-\frac2d\,w\,g^*\quad\text{on }T;
\]
in particular $D_p\psi_j=\id$. On the plane $T$, $\det M=1-\frac2d\,g^*w=1-\frac ud=\frac1d$ and
\[
  \norm M_F^2=2-\frac{2u}d+\frac4{d^2}\abs w^2\abs g^2\ \le\ \frac2d+\frac{4u^2}{3d^2}.
\]
Let $x\neq p$, so $u>0$, and let $s_1\ge s_2$ be the singular values of $M$. Then $s_2^2\le s_1s_2=\frac1d<d^2$, so $s_1<d$ if and only if $(d^2-s_1^2)(d^2-s_2^2)=d^4-d^2\norm M_F^2+d^{-2}>0$. Multiplying by $d^2$,
\[
  d^6-d^4\norm M_F^2+1\ \ge\ (d^3-1)^2-\tfrac43u^2d^2=u^2\Bigl((3+3u+u^2)^2-\tfrac43(1+u)^2\Bigr)>0 .
\]
Hence $\norm{D_x\psi_j}=s_1/d<1$. Finally, for $x\neq y$ in $\cT$ the segment from $x$ to $y$ meets the vertex $p$ at most at an endpoint, so $\abs{\psi_j(y)-\psi_j(x)}\le\abs{y-x}\int_0^1\norm{D_{x+t(y-x)}\psi_j}\,dt<\abs{y-x}$.
\end{proof}

\begin{remark}[Edelstein's theorem]\label{rem:edelstein}
Proposition~\ref{lem: Dpsi} already gives items (i)--(iv) of Theorem~\ref{cor:singlenode}. Let $A\ne B$ in $\mathcal K(\cT)$. For $a\in A$ and each $j$, choose $b\in B$ with $\abs{a-b}=\operatorname{dist}(a,B)$; then $\operatorname{dist}(\psi_j(a),\psi_j(B))\le\abs{\psi_j(a)-\psi_j(b)}<d_H(A,B)$ if $a\ne b$, and the left side vanishes if $a=b$. Since $A$ is compact and these functions of $a$ are continuous, the supremum over $a\in A$ is attained, and, exchanging the roles of $A$ and $B$, $d_H(\mathcal H(A),\mathcal H(B))<d_H(A,B)$. So $\mathcal H$ is a contractive map, in the sense of Edelstein, of the compact metric space $(\mathcal K(\cT),d_H)$, and by Edelstein's theorem~\cite{Ede62} it has a unique fixed point, which attracts every $K\in\mathcal K(\cT)$.

Items (v)--(vi) need more: that the cylinders shrink uniformly, $\max_{\abs I=n}\diam\psi_I(\cT)\to0$ as $n\to\infty$, which makes the convergence uniform in $K$ (see~\eqref{eq:rate}). This too follows from Proposition~\ref{lem: Dpsi}, by compactness. For an infinite sequence $\mathbf i=i_1i_2\cdots$ of indices, the cylinders $\psi_{i_1\cdots i_n}(\cT)$ form a decreasing sequence of compact sets; let $D(\mathbf i)$ be the diameter of their intersection, which is the limit of their diameters. Since $\psi_{i_1}$ is injective, this intersection is the image under $\psi_{i_1}$ of the intersection of the cylinders $\psi_{i_2\cdots i_n}(\cT)$ of the sequence $\mathbf i'=i_2i_3\cdots$. On $\{1,2,3\}^{\N}$, $D$ is the limit of the decreasing sequence of the continuous functions $\mathbf i\mapsto\diam\psi_{i_1\cdots i_n}(\cT)$, so it is upper semicontinuous and attains its maximum, at some $\mathbf i$. If this maximum were positive, then, choosing two points $p\ne q$ of the intersection for $\mathbf i'$ whose images under $\psi_{i_1}$ are at distance $D(\mathbf i)$, Proposition~\ref{lem: Dpsi} would give $D(\mathbf i)<\abs{p-q}\le D(\mathbf i')$, against maximality. So $D=0$, and, by Dini's theorem, $\max_{\abs I=n}\diam\psi_I(\cT)\to0$. The proof below does not use this argument: Lemma~\ref{lem:geom}(iv),(v) gives the explicit bound $\delta_n$ in~\eqref{eq:rate}.
\end{remark}

\medskip
We can now prove our first main result, Theorem~\ref{cor:singlenode}: the qualitative dynamics of $\PGN$ is that of a contractive iterated function system, although none of its maps is a contraction.

\thmsinglenode*

\begin{proof}
Let $\delta_n=6\sqrt3/(3+2n)$. Since $\mathcal H^n(\cT)=\bigcup_{\abs I=n}\psi_I(\cT)$,
Lemma~\ref{lem:zones} gives $\cN=\bigcap_n\mathcal H^n(\cT)$. The sets
$\mathcal H^n(\cT)$ are compact and decreasing, so they converge to $\cN$ in $d_H$
(for every $\varepsilon>0$, the compact sets of the points of $\mathcal H^n(\cT)$ at distance at least $\varepsilon$ from $\cN$ decrease and have empty intersection, so they are empty for large $n$),
and, since $\mathcal H$ is continuous, $\mathcal H(\cN)=\cN$. For $\abs I=n$
the cylinder $\psi_I(\cT)$ has diameter at most $\delta_n$
(Lemma~\ref{lem:geom}(iv),(v)) and contains the non-empty set
$\psi_I(\cN)\subset\cN$. Since $\mathcal H^n(K)=\bigcup_{\abs I=n}\psi_I(K)$
and $\cN=\mathcal H^n(\cN)=\bigcup_{\abs I=n}\psi_I(\cN)$, every point of
either set lies in a cylinder that meets the other set, so
\begin{equation}\label{eq:rate}
  d_H\bigl(\mathcal H^n(K),\cN\bigr)\le\delta_n
  \qquad\text{for all }K\in\mathcal K(\cT)\text{ and }n\ge0 .
\end{equation}

(i) The set $\cN$ is compact, non-empty and satisfies
$\mathcal H(\cN)=\cN$. If $K\in\mathcal K(\cT)$ and $\mathcal H(K)=K$, then
$d_H(K,\cN)=d_H(\mathcal H^n(K),\cN)\le\delta_n$ for every $n$, so
$K=\cN$.

(ii) This is~\eqref{eq:rate}, since $\delta_n\to0$; the convergence is uniform
in $K$.

(iii) By (i) and (ii), $\cN$ is a compact set with
$\mathcal H(\cN)=\cN$ that attracts every compact subset of $\cT$. So it is
the global attractor of $\PGN$.

(iv) By~\eqref{eq:rate}, $\mathcal H^n$ maps all of $\mathcal K(\cT)$ into the
closed $\delta_n$-ball about $\cN$. So the fixed point $\cN$ of $\mathcal H$
attracts every subset of $\mathcal K(\cT)$, uniformly, and $\{\cN\}$ is the
global attractor of $\mathcal H$.

(v) That every chain-recurrent point lies in the global attractor was proven recently by the author in ~\cite[Thm.~3]{DLstreams}; since~\cite{DLstreams} is not yet published, we give the argument for $\PGN$ directly. 

Let $x\in\mathcal T$, $y\in \cN$ and $\varepsilon>0$,
and choose $n\geq 1$ with $\delta_n<\varepsilon$.
Since $\cN\subset\mathcal H^n(\mathcal T)$, we have
$y\in\psi_I(\mathcal T)$ for some $|I|=n$, and therefore
$|\psi_I(x)-y|<\varepsilon$.
The orbit segment from $x$ to $\psi_I(x)$, with its
last point replaced by $y$, is an $\varepsilon$-chain
from $x$ to $y$.
Hence any two points of $\cN$ are chain-equivalent
and, since $\psi_j(\cN)\subset \cN$, they are
chain-recurrent.

Conversely, let $x\in\mathcal T\setminus \cN$ and put
$\eta=\operatorname{dist}(x,\cN)>0$.
Choose $n\geq 1$ with $\delta_n<\eta/3$.
By the uniform continuity of the $\psi_j$ on
$\mathcal T$, choose $\varepsilon>0$ such that every
$\varepsilon$-chain $z_0,\ldots,z_n$ with indices
$i_0,\ldots,i_{n-1}$ satisfies $\bigl|z_n-\psi_{i_{n-1}}\circ\cdots\circ\psi_{i_0}(z_0)\bigr|<\eta/3$.
Consider any $\varepsilon$-chain
$z_0=x,z_1,\ldots,z_r$ of length $r\geq n$.
Applying the preceding estimate to its last $n$ steps, we find a multiindex $I$ of length $n$ such that $|z_r-\psi_I(z_{r-n})|<\eta/3$.
By~\eqref{eq:rate}, applied to the singleton $\{z_{r-n}\}$, $\operatorname{dist}(\psi_I(z_{r-n}),\cN)\leq\delta_n$.
Consequently, $\operatorname{dist}(z_r,\cN)<\delta_n+\eta/3<\eta$.
In particular, no $\varepsilon$-chain of length $r\geq n$ can return from $x$ to itself.
Nor can there be a return chain of positive length $r<n$, since repeating it sufficiently many times would produce a return chain of length at least $n$.
Thus, $x$ is not chain-recurrent.

Therefore the chain-recurrent set of $\PGN$
is $\cN$, and its points are pairwise
chain-equivalent: the chain graph has the single
node $\cN$.

(vi) This follows at once from (iv) and the fact that, for a semiflow  with {\em compact dynamics}, i.e. with a global attractor, acting on a locally compact space, its chain-recurrent set is contained in the global attractor~\cite[Thm.~5.24]{DLY26DCDS}.
\end{proof}

Notice that the fact that $\mathcal H$ has a single node does not follow from $\PGN$ having a single node. 
Indeed there are IFSs (in fact, even standard dynamical systems, namely IFSs with $m=1$) whose graph consists of a single node and whose Hutchinson map has more than one node (in fact, even infinitely many).
An elementary example is the dyadic odometer, namely the map $x\mapsto x+1$ on the ring $\Z_2=\{0,1\}^\N$ of $2$-adic integers ($x\leftrightarrow\sum_kx_k2^k$),
with the $2$-adic metric $d(x,y)=2^{-\min\{k:\,x_k\neq y_k\}}$.
This map is minimal, because the orbit $x+\N$ of every $x$ is dense ($\N$ meets every residue class modulo $2^k$), and so its chain graph has the single node $\Z_2$. Its Hutchinson map is $K\mapsto K+1$, and its chain graph has infinitely many nodes. Indeed, $d(x,y)\le2^{-k}$ if and only if $x\equiv y\pmod{2^k}$, so $d_H(K,K')\le2^{-k}$ if and only if $K$ and $K'$ have the same reduction modulo $2^k$. Hence, if $\varepsilon\le2^{-k}$, the reductions modulo $2^k$ of the successive terms of an $\varepsilon$-chain of the Hutchinson map are successive translates by $1$, and they all have the same number of elements. For $n\ge0$, the set $K_n=\{0,2^n\}$ is chain-recurrent, because the continuous map $x\mapsto\{x,x+2^n\}$ conjugates the odometer to the Hutchinson map, so that the orbit of $K_n+1$ accumulates at $K_n$; and, for $n<n'$, $K_n$ and $K_{n'}$ are not chain-equivalent, because their reductions modulo $2^{n'}$ have two elements and one element. So the $K_n$ lie in pairwise distinct chain nodes of the Hutchinson map.

\medskip
As a corollary of Theorem~\ref{cor:singlenode}, we get a new proof of the fact that $\cN$ is the set of limit points of the souls, namely of the orbit of $\pr\one$ under $\PGN$.

\begin{corollary}
\label{cor:accumulation}
    The soul points $\pr{L_I}=\pr{A_I\one}$, that is, the soul directions $[L_I]$ in the chart, are pairwise distinct, and $\cN$ is exactly the set of their accumulation points.
\end{corollary}

\begin{proof}
Since $\one=\frac14\sum_i(\one+e_i)$, the point $\pr\one$ lies in $Z_\emptyset$, so $\pr{L_I}=\psi_I(\pr\one)\in\psi_I(Z_\emptyset)=Z_I$; the zones $Z_I$ are pairwise disjoint by Lemma~\ref{lem:zones}, so the $\pr{L_I}$ are pairwise distinct and none of them lies in $\cN$.
The set $\mathcal H^n(\{\pr\one\})$ is the set of the soul points $\pr{L_I}$ with $\abs I=n$, and, by~\eqref{eq:rate}, $d_H\bigl(\mathcal H^n(\{\pr\one\}),\cN\bigr)\le\delta_n\to0$.
Hence, every $y\in\cN$ is within $\delta_n$ of some $\pr{L_{I_n}}$ with $\abs{I_n}=n$, and these points are pairwise distinct, so $y$ is an accumulation point.
Conversely, there are finitely many soul points of each generation, so an accumulation point $y$ is a limit of soul points $\pr{L_{I_m}}$ with $\abs{I_m}\to\infty$; since $\operatorname{dist}(\pr{L_{I_m}},\cN)\le\delta_{\abs{I_m}}\to0$ and $\cN$ is closed, $y\in\cN$.
\end{proof}

Finally, we get the following.

\begin{corollary}
    \label{cor: free}
    $\GN$ is free (and so $\PGN$ is free as well).
\end{corollary}

\begin{proof}
    If $A_I=A_J$ with $I\neq J$, then $Z_I=\psi_I(Z_\emptyset)=\psi_J(Z_\emptyset)=Z_J$, against Lemma~\ref{lem:zones}. If $\psi_I=\psi_J$ on $\cT$, then, since $\cT$ has non-empty interior, the projective maps $[x]\mapsto[A_Ix]$ and $[x]\mapsto[A_Jx]$ coincide on $\RP^2$, so $A_I=cA_J$ for some $c\neq0$ with $c^3=\det A_I/\det A_J=1$, so $A_I=A_J$; hence $\PGN$ is free as well.
\end{proof}

\begin{remark}[the action on $\RP^2$]\label{rem:RP2}
Notice that, since each $\psi_j$ is a projective automorphism of $\RP^2$, $\RP^2$ is invariant under each $\psi_j$ and so the global attractor of the IFS generated by the $\psi_j$ on $\RP^2$ is $\RP^2$ itself: $\bigcup_j\psi_j(\RP^2)=\RP^2$.
In case of the first quadrant, namely the triangle $\Delta$ with vertices $[e_1],[e_2],[e_3]$, this set is only forward-invariant, namely $\psi_j(\Delta)\subset\Delta$, but nevertheless $\bigcup_j\psi_j(\Delta)=\Delta$, so $\Delta$ is the global attractor of the IFS generated by the $\psi_j$ on $\Delta$. 
This is why we define $\PGN$ on the triangle $\cT$, its natural phase space. In particular, $\cN$ is not an attractor of the IFS generated by the $\psi_j$ on $\RP^2$ in the sense of Barnsley and Vince, who characterize the attractors of projective IFSs~\cite{BV12}: $\psi_j$ fixes every point of the line $\{x_j=0\}$, which passes through the point $\pr{\one-e_j}$ of $\cN$, so no neighbourhood of $\cN$ in $\RP^2$ is attracted to $\cN$.
\end{remark}

\begin{remark}[$\cN$ and the Rauzy gasket of~\cite{JLPX23,Jur26}]
\label{rem:conj}
In~\cite{AS13,JLPX23,Jur26} the Rauzy gasket $\mathcal R$ is the limit set, in the triangle $\Delta$ with vertices $[e_1],[e_2],[e_3]$, of the IFS generated by the maps $[x]\mapsto[A_i^*x]$. 
Since $(\one\one^*-\id)A_i^*=A_i(\one\one^*-\id)$, the projective map $[x]\mapsto[(\one\one^*-\id)x]$, which in the $\ell^1$ charts is the similarity $x\mapsto\frac12(\one-x)$ of ratio $\frac12$ from $\Delta$ onto $\cT$, conjugates $[x]\mapsto[A_i^*x]$ to $\psi_i$ for each $i$. Hence it maps the Hutchinson-invariant non-empty compact subsets of $\Delta$ for the maps $[x]\mapsto[A_i^*x]$ to those of $\cT$ for the $\psi_i$, and, by Theorem~\ref{cor:singlenode}(i), it maps $\mathcal R$ onto $\cN$. So $\cN$ and $\mathcal R$ are similar, and in particular $\dimH\cN=\dimH\mathcal R$.
Moreover, $A^*_{i_1}\cdots A^*_{i_n}=(A_{i_n}\cdots A_{i_1})^*$ has the same singular
values as $A_{i_n}\cdots A_{i_1}$, so the two semigroups have the same affinity exponent.
Hence the theorem of~\cite{JLPX23,LPX23} and~\cite{Jur26} reads $\dimH\cN=\sA(\GN)$ (see Theorem~\ref{thm:jurga} below).
\end{remark}

As in the introduction, all entries of every element of $\GN$ are non-negative, so $\|\cdot\|_1=\one^*(\cdot)\one$ is linear, hence additive, on $\GN$, and the soul $L_I=A_I\one$ of $Z_I$ has $\norm{L_I}_1=\norm{A_I}_1$. Counting souls by their $\ell^1$-norm is therefore counting the elements of $\GN$ by their $\ell^1$ entrywise norm; the corresponding growth exponent, the norm exponent $\sN$ of $\GN$, is the subject of Sections~\ref{sec:general}--\ref{sec:3d}.

\section{The affinity exponent and the shape of the zones}\label{sec:affinity}

In this section $\sigma_1(A)\ge\sigma_2(A)\ge\sigma_3(A)$ are the singular values of $A\in SL_3(\R)$ for the Euclidean norm. 
Inside proofs, when $A$ is fixed, we sometimes write simply $\sigma_k$ for $\sigma_k(A)$. 
The singular value function of the projective action of $A$ on $\RP^2$ is~\cite[Eq. 1.2]{Jur26} (see also~\cite[Eq. 1.1]{JLPX23})
\[
  \varphi^s(A)=
  \begin{cases}
    \Bigl(\dfrac{\sigma_2(A)}{\sigma_1(A)}\Bigr)^s, & 0\le s\le1,\\[2ex]
    \dfrac{\sigma_2(A)}{\sigma_1(A)}\Bigl(\dfrac{\sigma_3(A)}{\sigma_1(A)}\Bigr)^{s-1}=\sigma_2(A)^{2-s}\sigma_1(A)^{1-2s}, & 1\le s\le2,
  \end{cases}
\]
where the last equality is due to the fact that $\det A=1$, and the affinity exponent of $\GN$ is defined as 
\[
  \sA=\sA(\GN)=\inf\Bigl\{s\in[0,2]:\ \sum_I\varphi^s(A_I)<\infty\Bigr\}
\]
(with $\sA=2$ if no such $s$ exists). 
Notice that this is a generalization to projective maps of the singular value function defined by Falconer for affine actions~\cite{Fal88}.

Since $\sigma_3\le\sigma_2\le\sigma_1$, $\varphi^s(A)$ is non-increasing in $s$, so the series converges for $s>\sA$ and diverges for $s<\sA$. 
Recall that $T=\{v\in\R^3:\ \one^*v=0\}$ is the tangent plane of the $\ell^1$ chart (Section~\ref{sec:souls}); we denote by $\PT$ its projective line, the set of lines through $0$ in $T$, that is, of unoriented directions in $T$, and by $x_0=\pr\one$ the center of $\cT$.

\begin{remark}[conventions]\label{rem:conventions}
The affinity exponent in \cite{JLPX23,Jur26} is defined for the projective action of a semigroup that may differ from $\GN$ by a transposition or by a fixed change of coordinates.
This does not change $\sA$: the singular values of $A^*$ are those of $A$, and the elements of the semigroup generated by the $A_j^*$ are the $A_I^*$; and for a fixed $C\in GL_3(\R)$, $\sigma_k(CAC^{-1})$ lies between $\sigma_k(A)/\mu$ and $\mu\,\sigma_k(A)$, with $\mu=\norm C\norm{C^{-1}}$, so that $\varphi^s(CA_IC^{-1})\asymp\varphi^s(A_I)$ uniformly in $I$ and $s\in[0,2]$.
\end{remark}

The equality between the Hausdorff dimension of $\cN$ and the affinity exponent of $\GN$, on which Theorem~\ref{thm: dimH} rests, is the following theorem of Jiao, Li, Pan and Xu and of Jurga, obtained independently (see Remarks~\ref{rem:conj} and~\ref{rem:conventions}, and~\cite[Remark~1.6]{LPX23}).

\begin{lettertheorem}[Jiao--Li--Pan--Xu~{\cite[\S5]{JLPX23}, \cite[Thm.~1.4]{LPX23}}; Jurga~{\cite[Thm.~1.1]{Jur26}}]\label{thm:jurga}
$\dimH\cN=\sA(\GN)$.
\end{lettertheorem}

\subsection*{Singular values and the shape of the zones}

\begin{lemma}\label{lem:sv}
For every multiindex $I$,
\[
  \tfrac13\norm{A_I}_1\le\sigma_1(A_I)\le\norm{A_I}_1,\qquad
  \sigma_1(A_I)^{-1/2}\le\sigma_2(A_I)\le3 .
\]
Both bounds on $\sigma_2$ are sharp up to multiplicative constants, along families
with $\sigma_1(A_I)\to\infty$:
\begin{enumerate}
\item[(a)] $\sigma_2(A_j^k)=1$ for all $j$ and $k\ge1$;
\item[(b)] $\sigma_2(A_I)\asymp\sigma_1(A_I)^{-1/2}$ for $I=(123)^n$, $n\ge1$.
\end{enumerate}
So neither bound can be replaced by a better power of $\sigma_1(A_I)$.
\end{lemma}

\begin{proof}
For a non-negative matrix $A$, $\norm A_1=\one^*A\one\le\abs{\one}\,\abs{A\one}\le3\sigma_1(A)$
and $\sigma_1(A)\le\norm A_F\le\norm A_1$, where $\norm{\cdot}_F$ is the Frobenius norm.
Since $\det A_I=1$, $\sigma_2^2\ge\sigma_2\sigma_3=\sigma_1^{-1}$. For the upper bound,
by Lemma~\ref{lem:geom}(ii),
\[
  \sigma_1\sigma_2=\sigma_3^{-1}=\norm{A_I^{-1}}\le\norm{A_I^{-1}}_F
  \le\norm{A_I^{-1}}_1\le\norm{A_I}_1\le3\sigma_1 .
\]

(a) Let $\{j,l,m\}=\{1,2,3\}$ and $v=e_l-e_m$. Since $e_j^*v=0$ and
$(\one-e_j)^*v=0$, both $A_j^k=\id+k(\one-e_j)e_j^*$ and its transpose fix $v$.
Hence $(A_j^k)^*A_j^kv=v$, and $1$ is a singular value of $A_j^k$. The other two
singular values have product $1$, because $\det A_j^k=1$, and they are not both
equal to $1$, because $\sigma_1(A_j^k)\ge\frac13\norm{A_j^k}_1=\frac13(3+2k)>1$. So they
are $t$ and $1/t$ for some $t>1$, and $1$ is the middle one: $\sigma_2(A_j^k)=1$.

(b) Let $B=A_1A_2A_3$, so that $A_I=B^n$. Recall that a matrix $M$ that is
diagonalizable over $\mathbb C$ satisfies $\norm{M^n}\asymp r(M)^n$, where $r(M)$ is
its spectral radius. The characteristic polynomial of $B$ is
$\lambda^3-7\lambda^2+5\lambda-1$, whose discriminant $-44$ is negative. So $B$ has
one real eigenvalue $\rho$ and two complex-conjugate ones $\mu,\bar\mu$. They are
distinct, so $B$ is diagonalizable over $\mathbb C$. Here $\rho=\tau^3\approx6.2223$,
where $\tau\approx1.8393$ is the tribonacci constant, and $\abs\mu^2=\rho^{-1}$ because
$\det B=1$. Hence $\sigma_1(B^n)=\norm{B^n}\asymp\rho^n$. The eigenvalues of
$B^{-1}$ have moduli $\rho^{-1}$, $\rho^{1/2}$, $\rho^{1/2}$, so
$\sigma_1\sigma_2(B^n)=\sigma_3(B^n)^{-1}=\norm{B^{-n}}\asymp\rho^{n/2}$. Dividing,
$\sigma_2(B^n)\asymp\rho^{-n/2}\asymp\sigma_1(B^n)^{-1/2}$.
\end{proof}

\begin{remark}[norm bounds and $\sigma_2$]\label{rem:ceiling}
For $1\le s\le2$, the bound $\sigma_2\le3$ and $\sigma_1\ge\frac13\norm{A_I}_1$ give
$\varphi^s(A_I)\le3^{s+1}\norm{A_I}_1^{1-2s}$, hence $\sA\le\frac12(1+\sN)$; the bound $\sigma_2\ge\sigma_1^{-1/2}$ gives
$\varphi^s(A_I)\ge\sigma_1(A_I)^{-3s/2}\ge\norm{A_I}_1^{-3s/2}$, hence $\sA\ge\frac23\sN$.
These are the two bounds of Theorem~\ref{thm: sN} for $\dimH\cN$.
By Lemma~\ref{lem:sv}, each of the two pointwise bounds on $\sigma_2$ is attained, up to constants,
along an explicit family. So neither can be improved term by term, and a precise value of $\sA$ requires $\sigma_2$ itself. (Compare the bounds for the areas and perimeters of the zones in~\cite[Prop.~9]{DD09}.)
\end{remark}

\begin{lemma}\label{lem:diamZ}
For every multiindex $I$,
\[
  \frac1{18}\,\frac{\sigma_2(A_I)}{\sigma_1(A_I)}\ \le\ \diam Z_I\ \le\
  3\sqrt{33}\,\frac{\sigma_2(A_I)}{\sigma_1(A_I)} .
\]
\end{lemma}

\begin{proof}
For $x,y\in\cT$ let $\varrho(x,y)=\abs{x\times y}/(\abs x\abs y)$, the sine of the angle
between $x$ and $y$. We first show that
\begin{equation}\label{eq:round}
  \tfrac1{\sqrt3}\abs{x-y}\ \le\ \varrho(x,y)\ \le\ \sqrt3\,\abs{x-y}\qquad(x,y\in\cT).
\end{equation}
Indeed $\abs{x\times y}=\abs{x\times(y-x)}$ and $y-x\in T$. 
The angle between $x\ge0$ and any non-zero $v\in T$ differs from a right angle by at most $\angle(x,\one)$, because $\one^*v=0$; so its sine is at least $\cos\angle(x,\one)=\one^*x/(\sqrt3\abs x)\ge1/\sqrt3$.
So $\abs{x\times y}\ge\abs x\abs{y-x}/\sqrt3$, and also $\abs{x\times y}\le\abs x\abs{y-x}$.
Now use $1/\sqrt3\le\abs y\le1$ on $\cT$.

Let $A=A_I$ and $p_m=\one+e_m$, $m=1,2,3$. The vertices of $Z_I$ are the points $\pr{Ap_m}$, and $\diam Z_I$ is the
largest of the three distances between them. Since $\det A=1$,
$(Ap)\times(Aq)=A^{-*}(p\times q)$, where $A^{-*}=(A^{-1})^*$ is the cofactor matrix. The vectors $n_1=p_2\times p_3=(3,-1,-1)$,
$n_2=p_3\times p_1=(-1,3,-1)$ and $n_3=p_1\times p_2=(-1,-1,3)$ have Gram matrix $16\,\id-5J$, with
$J$ the all-ones matrix, whose smallest eigenvalue is $1$. So every unit vector is
$\sum_m c_mn_m$ with $\abs c\le1$, and therefore
$\max_m\abs{A^{-*}n_m}\ge\norm{A^{-*}}/\sqrt3=\sigma_1\sigma_2/\sqrt3$.

Moreover $\abs{Ap_m}\le\sqrt6\,\sigma_1$, and, by
Lemma~\ref{lem:geom}(i) applied to the two frame vectors of $I$ whose sum is $Ap_m$,
\[
  \abs{Ap_m}\ge\norm{Ap_m}_1/\sqrt3\ge\norm A_1/\sqrt3\ge\sigma_1/\sqrt3.
\]
By \eqref{eq:round},
\[
  \abs{\pr{Ap_m}-\pr{Ap_{m'}}}\le\sqrt3\,\frac{\abs{A^{-*}n}}{\abs{Ap_m}\abs{Ap_{m'}}}
  \le\sqrt3\cdot\frac{\sigma_1\sigma_2\sqrt{11}}{\sigma_1^2/3}=3\sqrt{33}\,\frac{\sigma_2}{\sigma_1},
\]
where $n=\pm n_{m''}$ and $\abs n=\sqrt{11}$. For the pair $(m,m')$ that maximizes
$\abs{A^{-*}n}$,
\[
  \abs{\pr{Ap_m}-\pr{Ap_{m'}}}\ge\frac{\sigma_1\sigma_2/\sqrt3}{\sqrt3\cdot6\sigma_1^2}
  =\frac1{18}\,\frac{\sigma_2}{\sigma_1}. \qedhere
\]
\end{proof}

\begin{lemma}\label{lem:DI}
Let $v_0,v_1$ be an orthonormal basis of $T$. There is $c_0>0$ such that, for every multiindex $I$, the derivative $D_{x_0}\psi_I:T\to T$ satisfies
\[
  c_0\,\frac{\sigma_2(A_I)}{\sigma_1(A_I)}\ \le\ \max\bigl(\abs{D_{x_0}\psi_I(v_0)},\abs{D_{x_0}\psi_I(v_1)}\bigr)
  \ \le\ \norm{D_{x_0}\psi_I}\ \le\ 27\,\frac{\sigma_2(A_I)}{\sigma_1(A_I)} .
\]
Consequently $\diam Z_I\asymp\max_m\abs{D_{x_0}\psi_I(v_m)}\asymp\norm{D_{x_0}\psi_I}$, uniformly in $I$.
\end{lemma}

\begin{proof}
    Let $g$ be a non-negative invertible matrix. 
    For $x\in\cT$ and $v\in T$, differentiating $t\mapsto g(x+tv)/\one^*g(x+tv)$ at $t=0$ shows that the derivative of $x\mapsto\pr{gx}$ on the plane $\{\one^*x=1\}$ is
\begin{equation}\label{eq:Dproj}
  v\ \longmapsto\ \frac{gv-(\one^*gv)\,\pr{gx}}{\one^*gx} .
\end{equation}
    Let $A=A_I$. 
    By~\eqref{eq:Dproj} with $g=A$ and $x=x_0$,
\[
  D_{x_0}\psi_I(v)=\frac{P(Av)}{\one^*Ax_0},\qquad P(w)=w-(\one^*w)\,\pr{Ax_0} .
\]
    The vector $P(w)$ lies in $T$ and differs from $w$ by a multiple of $Ax_0$. 
    Hence, the distance $\delta$ from $w$ to the line $\R Ax_0$ satisfies $\delta\le\abs{P(w)}$ and, by the angle bound used for~\eqref{eq:round}, $\abs{P(w)}\le\sqrt3\,\delta$. 
    Since $\mathrm{dist}(Av,\R Ax_0)=\abs{Av\times Ax_0}/\abs{Ax_0}$ and $\abs{Ax_0}\le\one^*Ax_0\le\sqrt3\abs{Ax_0}$,
\begin{equation}\label{eq:Dv}
  \frac{\abs{Av\times Ax_0}}{\sqrt3\,\abs{Ax_0}^2}\ \le\ \abs{D_{x_0}\psi_I(v)}\ \le\
  \frac{\sqrt3\,\abs{Av\times Ax_0}}{\abs{Ax_0}^2}.
\end{equation}
    Also $\sigma_1/(3\sqrt3)\le\norm A_1/(3\sqrt3)=\norm{Ax_0}_1/\sqrt3\le\abs{Ax_0}\le\sigma_1\abs{x_0}=\sigma_1/\sqrt3$.

    \emph{Upper bound.} For a unit vector $v\in T$, $v\perp x_0$ and $\abs{Av\times Ax_0}\le\sigma_1\sigma_2\abs{x_0}=\sigma_1\sigma_2/\sqrt3$. So, by~\eqref{eq:Dv}, $\abs{D_{x_0}\psi_I(v)}\le\sqrt3\cdot\frac{\sigma_1\sigma_2}{\sqrt3}\cdot\frac{27}{\sigma_1^2}=27\sigma_2/\sigma_1$.

    \emph{Lower bound.} First, $\abs{Av\times Ax_0}=\abs{A^{-*}(v\times x_0)}$. On $T$, the map $v\mapsto v\times x_0=\frac13v\times\one$ is $1/\sqrt3$ times an isometry of $T$ onto itself. 
    So $\max_m\abs{A^{-*}(v_m\times x_0)}\ge\norm{A^{-*}|_T}/\sqrt6$.
    Next, let $q$ be a unit right singular vector of $A$ for $\sigma_3$. 
    Then, $q$ is a right singular vector of $A^{-*}$ for its largest singular value $\sigma_1\sigma_2$. 
    Testing $A^{-*}$ on the unit vector of $T$ in the direction of the orthogonal projection $q_T$ of $q$ on $T$ gives $\norm{A^{-*}|_T}\ge\sigma_1\sigma_2\abs{q_T}$.
    Let $\theta$ be the angle between $q$ and $\one$, so that $\abs{q_T}=\sin\theta$. 
    Writing $\one/\sqrt3=\cos\theta\,q+\sin\theta\,q'$ with $q'\perp q$ a unit vector,
\[
  \frac{\sigma_1}3\le\frac{\norm A_1}3=\frac{\norm{A\one}_1}{3}\le\frac{\abs{A\one}}{\sqrt3}
  \le\sigma_3+\sigma_1\sin\theta .
\]
    Since $\sigma_3\le\sigma_1^{-1/2}$, this gives $\sin\theta\ge\frac13-\sigma_1^{-3/2}\ge\frac5{24}$ when $\sigma_1\ge4$. 
    In that case, by \eqref{eq:Dv},
\[
  \max_m\abs{D_{x_0}\psi_I(v_m)}\ \ge\ \frac{3}{\sqrt3\,\sigma_1^2}\cdot\frac{\sigma_1\sigma_2}{\sqrt6}\cdot\frac5{24}
  \ =\ \frac5{24\sqrt2}\,\frac{\sigma_2}{\sigma_1}.
\]
    If $\sigma_1<4$ then $\norm A_1<12$, so $\abs I\le4$ (Lemma~\ref{lem:geom}(v)) and there are finitely many such $I$.
    For each of them $D_{x_0}\psi_I$ is invertible, and $c_0$ is the minimum of $5/(24\sqrt2)$ and of these finitely many ratios. 
    The last statement follows from Lemma~\ref{lem:diamZ}.
\end{proof}

\begin{proposition}\label{prop:phiD}
For $0\le\kappa\le1$ and every multiindex $I$, with constants independent of $I$ and $\kappa$,
\[
  \varphi^{2-\kappa}(A_I)\ \asymp\ \norm{A_I}_1^{-3(1-\kappa)}\,
  \max\bigl(\abs{D_{x_0}\psi_I(v_0)},\abs{D_{x_0}\psi_I(v_1)}\bigr)^{\kappa} .
\]
Consequently $\sum_I\varphi^{2-\kappa}(A_I)<\infty$ if and only if
$\sum_I\norm{A_I}_1^{-3(1-\kappa)}\abs{D_{x_0}\psi_I(v)}^\kappa<\infty$ for every unit vector $v\in T$,
and if and only if this holds for $v=v_0$ and $v=v_1$.
\end{proposition}

\begin{proof}
Since $\sigma_1\sigma_2\sigma_3=1$,
$\varphi^{2-\kappa}=\frac{\sigma_2}{\sigma_1}\bigl(\frac{\sigma_3}{\sigma_1}\bigr)^{1-\kappa}
=\bigl(\frac{\sigma_2}{\sigma_1}\bigr)^\kappa\sigma_1^{-3(1-\kappa)}$. The first statement
then follows from Lemmas~\ref{lem:sv} and~\ref{lem:DI}. For the second, use
$\abs{D_{x_0}\psi_I(v)}\le\norm{D_{x_0}\psi_I}$ and 
\[
\max_m\abs{D_{x_0}\psi_I(v_m)}^\kappa\le\abs{D_{x_0}\psi_I(v_0)}^\kappa+\abs{D_{x_0}\psi_I(v_1)}^\kappa.
\]
\end{proof}

\begin{corollary}\label{cor:cutoutA}
For $0\le\kappa\le1$, uniformly in $I$ and $\kappa$,
$\mathrm{area}(Z_I)^{1-\kappa}\mathrm{per}(Z_I)^\kappa\asymp\varphi^{2-\kappa}(A_I)$.
Hence the series $\sum_I\mathrm{area}(Z_I)^{1-\kappa}\mathrm{per}(Z_I)^\kappa$
converges if $\kappa<2-\sA$ and diverges if $\kappa>2-\sA$.
\end{corollary}

\begin{proof}
    By Lemmas~\ref{lem:geom}(iii) and~\ref{lem:sv}, $\mathrm{area}(Z_I)\asymp\norm{A_I}_1^{-3}\asymp\sigma_1(A_I)^{-3}$. 
    Moreover, we have that $2\,\diam Z_I\le\mathrm{per}(Z_I)\le3\,\diam Z_I$, and, by Lemma~\ref{lem:diamZ}, $\diam Z_I\asymp\sigma_2(A_I)/\sigma_1(A_I)$. 
    Since $\varphi^{2-\kappa}=(\sigma_2/\sigma_1)^\kappa\sigma_1^{-3(1-\kappa)}$, this proves the comparison, and the second statement follows from the definition of $\sA$.
\end{proof}

\begin{proposition}[the cut-out bound]\label{prop:cutout}
$\udimB\cN\le\sA$.
\end{proposition}

\begin{proof}
    Recall that $\mathrm{area}(\cN)=0$ (Theorem~\ref{thm:yoccoz}).
    Let $s\in(\sA,2]$ and $\kappa=2-s\in[0,1)$. 
    For $0<\varepsilon<1$, let $F_\varepsilon$ be the set of points of the plane $\{\one^* x=1\}$ at distance less than $\varepsilon$ from $\cN$.
    Then:

    \noindent
    {\em(a)} $F_\varepsilon\setminus\cT$ lies in the $\varepsilon$-neighbourhood of $\cT$   minus $\cT$, whose area is $\frac3{\sqrt2}\,\varepsilon+\pi\varepsilon^2<6\varepsilon$.

    \noindent
    {\em(b)} Let $\xi\in(F_\varepsilon\cap\cT)\setminus\cN$. By Lemma~\ref{lem:zones}, $\xi\in Z_I$ for some $I$. 
    Let $q\in\cN$ with $\abs{\xi-q}<\varepsilon$; since $q\notin Z_I$, the segment from $\xi$ to $q$ meets $\partial Z_I$. 
    So $\xi$ lies in the collar $C_\varepsilon(Z_I)=\{\xi\in Z_I:\mathrm{dist}(\xi,\partial Z_I)<\varepsilon\}$.

    \noindent
    {\em(c)} For a triangle $Z$, $\mathrm{area}(C_\varepsilon(Z))\le\min(\mathrm{area}(Z),\varepsilon\,\mathrm{per}(Z))$:
    if $\xi\in Z$ and $d=\mathrm{dist}(\xi,\partial Z)<\varepsilon$, the closed disc of radius $d$ about $\xi$ lies in $Z$ and touches $\partial Z$ at a point $p$ of a side $\Sigma$; it lies on one side of the line of $\Sigma$, so $p$ is the foot of the perpendicular from $\xi$, and $\xi$ lies in the rectangle of width $\varepsilon$   erected inwards on $\Sigma$.

Since $\cN$ has measure zero and $\min(\alpha,\beta)\le\alpha^{1-\kappa}\beta^\kappa$,
Corollary~\ref{cor:cutoutA} gives
\[
  \mathrm{area}(F_\varepsilon)\le6\varepsilon+\varepsilon^\kappa\sum_I\mathrm{area}(Z_I)^{1-\kappa}\mathrm{per}(Z_I)^\kappa
  \le6\varepsilon+C\varepsilon^\kappa\sum_I\varphi^s(A_I)=O(\varepsilon^\kappa),
\]
the last series being convergent because $s>\sA$. For a bounded set $F$ in the plane
\cite[Prop.~3.2]{FalBook},
\[
  \udimB F=2-\liminf_{\varepsilon\to0}\frac{\log\mathrm{area}(F_\varepsilon)}{\log\varepsilon}.
\]
Hence $\udimB\cN\le2-\kappa=s$, and letting $s\downarrow\sA$ gives the claim.
\end{proof}

The following result follows from~\cite{Sew23} together with~\cite{JLPX23,Jur26}; since~\cite{Sew23} is still unpublished, we provide here a proof.

\begin{corollary}
    \label{cor: dimB}
    The box dimension of $\cN$ exists and $\dim_B\cN=\dim_H\cN=s_A$.
\end{corollary}

\begin{proof}
    By Theorem~\ref{thm:jurga} and Proposition~\ref{prop:cutout},
    $s_A=\dimH\cN\le\underline{\dim}_{\mathrm B}\cN\le\udimB\cN\le\sA$.
\end{proof}

\begin{remark}[a direct estimate]\label{rem:cutout}
Lemma~\ref{lem:geom}(iv) bounds every perimeter by a multiple of $\norm{A_I}_1^{-1}$, the order of the perimeters of the cusp zones, while the perimeters of other zones can be as small as the order $\norm{A_I}_1^{-3/2}$ (Remark~\ref{rem:geomsharp}).
We compared half-octave shells $2^{j/2}<\|A_I\|_1\le2^{(j+1)/2}$ of the series of Corollary~\ref{cor:cutoutA} over the $3.4\times10^7$ multiindices with $\|A_I\|_1\le4096$: the last four shells give critical values of $\kappa$ between $0.311$ and $0.315$, that is, $2-\kappa\approx1.686$, a direct numerical estimate of $\sA$. The perimeter, unlike $\|\cdot\|_1$, is not a product of run factors (Section~\ref{sec:cocycle}), but it is comparable to the norm of the derivative $D_{x_0}\psi_I$, which becomes a product of one factor per run once its direction is added to the state (Lemma~\ref{lem:Dcocycle}). This is the route to Theorem~\ref{thm: dimH}.
\end{remark}

\section{Runs and the needle}\label{sec:cocycle}

In this section we split every multiindex into maximal runs of equal indices. The $\ell^1$-size of a soul is the product of one factor per run, the run factors, each depending only on the soul point before the run; and so is the length of the image of a tangent vector, once its direction is carried along, the needle (Lemma~\ref{lem:Dcocycle}). This turns the series that defines $\sA$ into the sum of the iterates of a positive transfer operator, and Proposition~\ref{cor:certA} reduces bounds on $\sA$ to positive test functions for it.

We call {\bf run} a multiindex $J=j\dots j=j^k$ where only one type of index appears, so that $A_J=A^k_j$.
Since $e_j^*(\one-e_j)=0$, we have $A_j^k=\id+k\,(\one-e_j)e_j^*$. 
Therefore, given an $x\in\cT$, we have that
\begin{equation}\label{eq:run}
  A_j^kx=x+kx_j(\one-e_j),\qquad
  \norm{A_j^kx}_1=1+2kx_j,\qquad
  \psi_j^k(x)=\pr{A_j^kx}=(1-\beta)\,x+\beta\,\pr{\one-e_j},
\end{equation}
with $\beta=2kx_j/(1+2kx_j)$: a run of $j$ pushes $x$ toward the parabolic
point $\pr{\one-e_j}$. Throughout the paper, the \textbf{first index} of a multiindex $I=i_1\cdots i_n$ is $i_1$, the leftmost, and its \textbf{last index} is $i_n$, the rightmost, which is the first to act on a vector: $A_Ix=A_{i_1}\bigl(\cdots(A_{i_n}x)\bigr)$. Every multiindex is uniquely a sequence of runs; we list them from the
right, $I=j_n^{k_n}\cdots j_1^{k_1}$ with $k_r\ge1$ and $j_{r+1}\neq j_r$, so that
$j_1^{k_1}$ is the first run to act on a vector, and $j_1$ is the last index of $I$. Since $\norm{A_I}_1=\norm{A_I\one}_1
=3\norm{A_Ix_0}_1$, iterating \eqref{eq:run} with $x^{(0)}=x_0$ and
$x^{(r)}=\psi_{j_r}^{k_r}(x^{(r-1)})$ gives
\begin{equation}\label{eq:cocycle}
  \norm{L_I}_1=\norm{A_I}_1=3\prod_{r=1}^n\Bigl(1+2k_r\,x^{(r-1)}_{j_r}\Bigr),
  \qquad x^{(n)}=\pr{L_I},
\end{equation}
exactly. Geometrically, $Z_{j^kJ}=\psi_j^k(Z_J)$ and $\pr{L_{j^kJ}}=\psi_j^k(\pr{L_J})$:
prepending a run to $J$ moves the zone of $J$ and its soul point by the map
$\psi_j^k$ of the iterated function system, and multiplies the size of the soul by
$1+2k\,(\pr{L_J})_j$. The states $x^{(r)}$ are the soul points of the
suffixes of $I$, and the last one is the soul point of $I$ itself. We call the factors $1+2k_r\,x^{(r-1)}_{j_r}$ in~\eqref{eq:cocycle} the \textbf{run factors} of $I$; each depends only on the run $j_r^{k_r}$ and on the state $x^{(r-1)}$ before it.

\begin{example}\label{ex:run}
Take $I=1^2\,2\,3^2=11233$. Read from the right, its runs are $3^2$, $2$, $1^2$:
\[
\renewcommand{\arraystretch}{1.2}
\begin{array}{c|c|c|c|c}
  r & j_r^{k_r} & x^{(r-1)} & 1+2k_r\,x^{(r-1)}_{j_r} & x^{(r)}\\\hline
  1 & 3^2 & (\frac13,\frac13,\frac13) & \frac73 & (\frac37,\frac37,\frac17)\\
  2 & 2 & (\frac37,\frac37,\frac17) & \frac{13}7 & (\frac6{13},\frac3{13},\frac4{13})\\
  3 & 1^2 & (\frac6{13},\frac3{13},\frac4{13}) & \frac{37}{13} & (\frac6{37},\frac{15}{37},\frac{16}{37})
\end{array}
\]
and \eqref{eq:cocycle} reads $37=3\cdot\frac73\cdot\frac{13}7\cdot\frac{37}{13}$. The
final state is the soul point $\pr{L_I}=\pr{A_{11233}\one}=\pr{(6,15,16)}=\frac1{37}(6,15,16)$.
\end{example}

\begin{figure}[t]
\centering
\includegraphics[width=\textwidth]{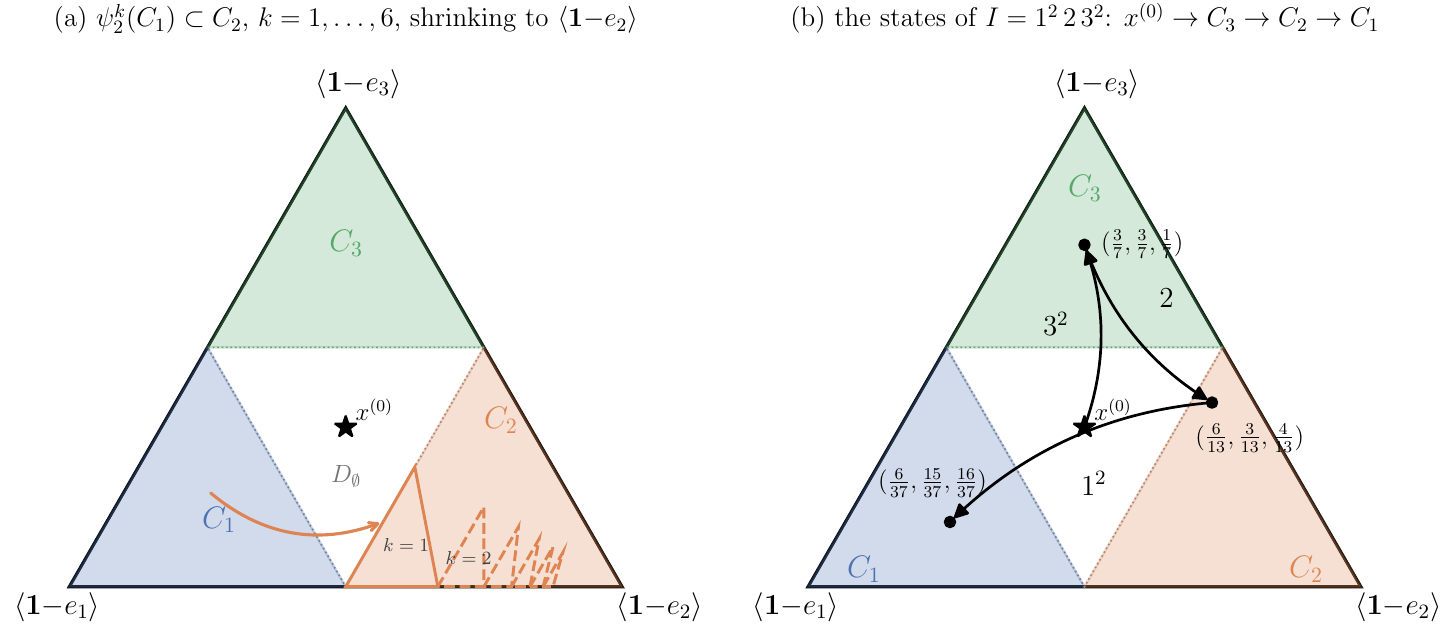}
\caption{The cylinders $C_i=\psi_i(\cT)$ of Lemma~\ref{lem:trap}. Left: the images $\psi_2^k(C_1)=\psi_2^k\psi_1(\cT)\subset C_2$, $k=1,\dots,6$, have pairwise disjoint interiors (consecutive ones share a vertex) and shrink to the parabolic point $\pr{\one-e_2}$. Right: the states of Example~\ref{ex:run}, $I=1^2\,2\,3^2$: starting from $x^{(0)}=\pr{\one}$, the runs $3^2$, $2$, $1^2$ move the state to $x^{(1)}\in C_3$, $x^{(2)}\in C_2$ and $x^{(3)}=\pr{L_I}\in C_1$.}
\label{fig:runs}
\end{figure}

\begin{lemma}\label{lem:trap}
Let $C_i=\psi_i(\cT)=\{x\in\cT:\ x_i\le\frac14\}$, the first-level cylinder of
$i$. Then $\psi_j^k(\cT)\subset C_j$ for every $j$ and $k\ge1$, and
$2x_j\ge\frac12$ for all $x\in C_i$ and $j\neq i$.
\end{lemma}

\begin{proof}
$\psi_j^k(\cT)\subset\psi_j(\cT)$. If $x\in C_i$, $j\ne i$ and $l$ is the third index, then
$x_j=1-x_i-x_l\ge1-\frac14-\frac12=\frac14$.
\end{proof}

\begin{remark}[runs and fastness]\label{rem:fast}
Each generator $A_j$ is parabolic: the map $\psi_j$ fixes the vertex $\pr{\one-e_j}$, and within a run norms add rather than multiply, $\norm{A_j^{k+k'}}_1=3+2(k+k')$, while $\norm{A_j^k}_1\norm{A_j^{k'}}_1\asymp kk'$. A transfer operator that adds one index at a time therefore has weights arbitrarily close to $1$ near the parabolic points, and its iterates decay only polynomially. Across runs, instead, norms multiply up to a bounded factor. If $J$ begins with $i\neq j$, then $\pr{L_J}\in\psi_J(\cT)\subset C_i$, so $2(\pr{L_J})_j\ge\frac12$ by Lemma~\ref{lem:trap}, and by \eqref{eq:run}, since $A_J\one=\norm{A_J}_1\pr{L_J}$,
\[
  \norm{A_j^kA_J}_1=\norm{A_J}_1\bigl(1+2k(\pr{L_J})_j\bigr)\ \ge\ \norm{A_J}_1\bigl(1+\tfrac k2\bigr)\ \ge\ \tfrac14\,\norm{A_j^k}_1\,\norm{A_J}_1 .
\]
Together with Lemma~\ref{lem:junction}, this makes $\GN$ fast in the sense of \cite[Def.~3]{DL15JGA}, recalled in Section~\ref{sec:fast}.
We therefore take each step of the transfer operator to be a maximal run. Grouping the steps of a parabolic system into runs is the inducing scheme of parabolic iterated function systems~\cite{MU00}, and, for the Rauzy gasket, the acceleration of~\cite{AS13,AHS16,GMR19}. Each step then starts where a run of another generator ended: the state lies in $C_i$, away from the parabolic point $\pr{\one-e_j}$ of the run $A_j^k$, and the weights $\bigl(1+2kx_j\bigr)^{-t}\le(1+k/2)^{-t}$ are uniformly summable in $k$ for $t>1$. The same method gives the norm exponent (Section~\ref{sec:fast}); here we apply it to the derivative.
\end{remark}

The derivative of a run can be written explicitly as follows.

\begin{lemma}\label{lem:contraction}
    Fix a $j\in\{1,2,3\}$ and let $k\ge1$ and $x\in\cT$ with $x_j>0$. 
    Set $u=2x_j$, $d=1+ku$ and $c=\bigl(\pr{\one-e_j}-x\bigr)/x_j$. 
    Then $c\in T$ and the derivative of $\psi_j^k$ at $x$ is, for $v\in T$,
\begin{equation}\label{eq:DR}
  D_x\psi_j^k(v)=\frac{v-2k\,v_j\,\bigl(\psi_j^k(x)-\pr{\one-e_j}\bigr)}{d}
  =\frac1d\Bigl(v+\frac{ku}{d}\,v_j\,c\Bigr).
\end{equation}
    If $x\in C_i$ for some $i\neq j$, or $x=x_0$, then $\norm{D_x\psi_j^k}\le(1+2/\sqrt3)/d$.
\end{lemma}

\begin{proof}
    For $v\in T$ we have $A_j^kv=v+k\,v_j(\one-e_j)$ and $\one^*A_j^kv=2k\,v_j$, and $\one^*A_j^kx=1+2kx_j=d$ by~\eqref{eq:run}. 
    So, as in the proof of Lemma~\ref{lem:DI}, formula~\eqref{eq:Dproj} with $g=A_j^k$ gives 
    \[
    D_x\psi_j^k(v)=\frac{A_j^kv-2k\,v_j\,\psi_j^k(x)}{d}. 
    \]
    This is the first expression in~\eqref{eq:DR}, since $\one-e_j=2\pr{\one-e_j}$. 
    Moreover, by~\eqref{eq:run}, we have that 
    \[
    \psi_j^k(x)-\pr{\one-e_j}=\bigl(x-\pr{\one-e_j}\bigr)/d=-x_jc/d,
    \]
    which gives the second.

    Now, let $\beta=ku/d\in(0,1)$, and let $g=e_j-\frac13\one$ be the orthogonal projection of $e_j$ on $T$.
    Then $v_j=g\cdot v$ for $v\in T$, $\abs g^2=\frac23$ and, since the $j$-th coordinate of $\pr{\one-e_j}$ is $0$, $g\cdot c=c_j=-1$, so that, on $T$,
\begin{equation}\label{eq:DRB}
  D_x\psi_j^k=d^{-1}B,\qquad B=\id+\beta\,c\,g^* .
\end{equation}
    The function $x\mapsto\abs{x-\pr{\one-e_j}}/x_j$ is quasi-convex on $C_i$, because its sublevel sets $\{\abs{x-\pr{\one-e_j}}\le s\,x_j\}$ are convex, so its maximum on the triangle $C_i$ is attained at a vertex. 
    The values at the vertices $\pr{\one-e_i}$, $\pr{\one+e_l}$ and $\pr{\one+e_j}$ of $C_i$, where $l$ is the third index, are $\sqrt2$, $\sqrt2$ and $\sqrt{3/2}$; hence $\abs c^2\le2$ on $C_i$, and $\abs c^2=\frac32$ at $x_0$. So $\norm B\le1+\beta\abs c\abs g\le1+\sqrt2\sqrt{2/3}=1+2/\sqrt3$.
\end{proof}

\begin{lemma}[the needle along runs]\label{lem:Dcocycle}
Let $I=j_n^{k_n}\cdots j_1^{k_1}$ be written in runs from the right, and let $x^{(r)}$ be
as in \eqref{eq:cocycle}. For a unit vector $v\in T$ let $w^{(0)}=v$ and
\[
  w^{(r)}=\frac{D_{x^{(r-1)}}\psi_{j_r}^{k_r}(w^{(r-1)})}{\abs{D_{x^{(r-1)}}\psi_{j_r}^{k_r}(w^{(r-1)})}}
  \qquad(r=1,\dots,n).
\]
Then
\[
  \abs{D_{x_0}\psi_I(v)}=\prod_{r=1}^n\abs{D_{x^{(r-1)}}\psi_{j_r}^{k_r}(w^{(r-1)})} ,
\]
where each factor is given explicitly by \eqref{eq:DR}.
\end{lemma}

\begin{proof}
$\psi_I=\psi_{j_n}^{k_n}\circ\cdots\circ\psi_{j_1}^{k_1}$ and $x^{(r)}$ is the image of $x_0$
after $r$ runs, so this is the chain rule.
\end{proof}

We call $(x^{(r)},[w^{(r)}])\in\cT\times\PT$ the \textbf{needle} after $r$ runs: each run moves its base point, turns its direction and multiplies its length by the factor $\abs{D_{x^{(r-1)}}\psi_{j_r}^{k_r}(w^{(r-1)})}$, which depends only on the run and on the needle before it. So the length of the image of a unit needle at $x_0$ is the product of these factors, as $\norm{A_I}_1$ is the product of the run factors in~\eqref{eq:cocycle}.

Bounds on $\sA$ will be decided by positive test functions, through the following elementary lemma.

\begin{lemma}\label{lem:test}
Let $X$ be a set, $\mathcal B(X)$ the space of bounded real functions on $X$, and
$\mathcal L:\mathcal B(X)\to\mathcal B(X)$ a positive linear operator, i.e.\
$\mathcal Lh\ge0$ whenever $h\ge0$. Let $h\in\mathcal B(X)$ with
$0<\inf h\le\sup h<\infty$.
\begin{enumerate}
\item[(i)] If $\mathcal Lh\ge h$, then $\mathcal L^m\mathds 1\ge\inf h/\sup h$ for all $m\ge0$.
\item[(ii)] If $\mathcal Lh\le\theta h$ for some $\theta<1$, then
  $\mathcal L^m\mathds 1\le\theta^m\sup h/\inf h$ for all $m\ge0$.
\end{enumerate}
\end{lemma}

\begin{proof}
    By positivity and linearity, $\mathcal Lh\ge h$ gives $\mathcal L^mh\ge h$ and
    $\mathcal Lh\le\theta h$ gives $\mathcal L^mh\le\theta^mh$. 
    Since $h/\sup h\le\mathds 1\le h/\inf h$, positivity also gives
    $\mathcal L^mh/\sup h\le\mathcal L^m\mathds 1\le\mathcal L^mh/\inf h$.
\end{proof}

For a non-negative matrix $W$ with spectral radius $r(W)$ and a positive vector
$h$, the lemma reduces to the Collatz--Wielandt bounds (see \cite{Col42,Wie50} and
\cite[Thm.~8.1.26]{HJ85})
\[
  \min_k\frac{(Wh)_k}{h_k}\ \le\ r(W)\ \le\ \max_k\frac{(Wh)_k}{h_k}.
\]
The same principle underlies rigorous computations of Hausdorff dimension; see, for
instance, \cite[Lemma~2.2]{FN18}.

\begin{proposition}[certificates for the affinity exponent]
\label{cor:certA}
For $0<\kappa<1$ let $\cK_\kappa$ act on bounded functions on $\mathcal C\times\PT$, where $\mathcal C=\bigsqcup_iC_i\times\{i\}$, by
\[
  (\cK_\kappa h)(x,[w],i)=\sum_{j\neq i}\sum_{k\ge1}
  \bigl(1+2kx_j\bigr)^{-3(1-\kappa)}\,\abs{D_x\psi_j^k(\hat w)}^{\kappa}\,
  h\bigl(\psi_j^k(x),[D_x\psi_j^k(w)],j\bigr),
\]
where $\hat w=w/\abs w$. Then $\cK_\kappa$ is a bounded positive operator. 
Let $h$ be bounded with $\inf h>0$.
\begin{enumerate}
\item[(i)] If $\cK_\kappa h\ge h$, then $\sA(\GN)\ge2-\kappa$.
\item[(ii)] If $\cK_\kappa h\le\theta h$ for some $\theta<1$, then $\sA(\GN)\le2-\kappa$.
\end{enumerate}
\end{proposition}

\begin{proof}
By Lemma~\ref{lem:contraction} and $2x_j\ge\frac12$ on $C_i$ ($j\neq i$), the weights are
at most $(1+2/\sqrt3)^\kappa(1+k/2)^{-(3-2\kappa)}$. These are summable in $k$ because
$\kappa<1$; so $\cK_\kappa$ is bounded.

For $(x,[w],i)\in\mathcal C\times\PT$, induction on $m$ shows that $(\cK_\kappa^m\mathds 1)(x,[w],i)$ is the sum, over the multiindices $I'$ with exactly $m$ runs whose last index is not $i$, of the product over the runs of $I'$ of the weights $\bigl(1+2k_rx^{(r-1)}_{j_r}\bigr)^{-3(1-\kappa)}\abs{D_{x^{(r-1)}}\psi_{j_r}^{k_r}(\hat w^{(r-1)})}^\kappa$, computed as in \eqref{eq:cocycle} and Lemma~\ref{lem:Dcocycle} but starting from $x^{(0)}=x$ and $w^{(0)}=w$. Let $v\in T$ be a unit vector. Every non-empty multiindex is $I=I'j^k$, where $j^k$ is its first run to act and $I'$ has last index different from $j$; at $x_0$ the first factor of \eqref{eq:cocycle} is $1+\frac{2k}3$. Hence, by \eqref{eq:cocycle} and Lemma~\ref{lem:Dcocycle},
\begin{multline*}
  \sum_I\norm{A_I}_1^{-3(1-\kappa)}\abs{D_{x_0}\psi_I(v)}^\kappa
  =3^{-3(1-\kappa)}\Bigl(1+\sum_{j}\sum_{k\ge1}\bigl(1+\tfrac{2k}3\bigr)^{-3(1-\kappa)}
  \abs{D_{x_0}\psi_j^k(v)}^\kappa\\
  \cdot\sum_{m\ge0}(\cK_\kappa^m\mathds 1)\bigl(\psi_j^k(x_0),[D_{x_0}\psi_j^k(v)],j\bigr)\Bigr),
\end{multline*}
an identity in $[0,\infty]$. The bound of Lemma~\ref{lem:contraction} also holds at $x_0$, so the sum over $k$ of the outer weights
is finite.
If $\cK_\kappa h\ge h$, Lemma~\ref{lem:test}(i) shows that the series diverges, for every $v$. By
Proposition~\ref{prop:phiD}, $\sum_I\varphi^{2-\kappa}(A_I)=\infty$, so
$\sA\ge2-\kappa$. If $\cK_\kappa h\le\theta h$, Lemma~\ref{lem:test}(ii) shows that the series
converges for every unit vector $v$, and Proposition~\ref{prop:phiD} gives $\sA\le2-\kappa$.
\end{proof}

\section{The algorithm: from infinitely many points to a finite check}
\label{sec:disc}

By Proposition~\ref{cor:certA}, to prove $\sA(\GN)\ge2-\kappa$ it suffices to find a bounded function $h$ on $\mathcal C\times\PT$ with $\inf h>0$ and $\cK_\kappa h\ge h$, and to prove $\sA(\GN)\le2-\kappa$ it suffices to find one with $\cK_\kappa h\le\theta h$ for some $\theta<1$. Each condition is an inequality at every point of $\mathcal C\times\PT$, a space of three real dimensions. In this section we reduce it to finitely many inequalities between numbers that can be enclosed in interval arithmetic. This takes three steps.

First, we use the fact that the permutations of the coordinates permute the generators $A_1,A_2,A_3$ of $\GN$ by conjugation, so $\cK_\kappa$ maps functions invariant under them to functions invariant under them. Such a function $h$ is determined by its restriction $h_1$ to $C_1\times\PT$, which is invariant under the transposition $\tau_{23}$ of $x_2$ and $x_3$; conversely, every $\tau_{23}$-invariant function on $C_1\times\PT$ is such a restriction. We look for $h$ among the invariant functions. On them, $\cK_\kappa$ acts as an explicit operator $\tilde\cK_\kappa$ on functions on $C_1\times\PT$, and $h$ satisfies either inequality if and only if $h_1$ satisfies it with $\tilde\cK_\kappa$ in place of $\cK_\kappa$ (Lemma~\ref{lem:sym}).

Second, we take $h_1$ positive and piecewise linear on a mesh of $C_1\times\PT$ whose cells are prisms over triangles. 
To check an inequality on a set $Q$ that covers part of a cell, we evaluate both sides at the vertices of $Q$, and require the inequality to hold there with a margin. The margin bounds how far $\tilde\cK_\kappa h_1$ can deviate, inside $Q$, from its interpolation between the vertices: we split $\tilde\cK_\kappa h_1$ on $Q$ into three parts, bound the deviation of each through an enclosure of its Hessian, of its increments, or of its oscillation on $Q$ (Lemma~\ref{lem:interp}), and take the sum of the three bounds. If the inequality holds at the vertices of $Q$ with this margin, then it holds at every point of $Q$ in the cell (Proposition~\ref{prop:P1}). The margin is smaller on smaller sets: we first take $Q$ to be the cell itself and, where this test fails, products of a rectangle and an arc of directions, obtained by repeatedly bisecting the smallest such product that contains the cell. By the symmetry $\tau_{23}$, it suffices to do this for the cells on one side of the diagonal $x_2=x_3$.

Third, we explain how to obtain the three parts and their enclosures run by run. The term of a run is bounded through its Hessian when the run maps $Q$ into a single cell, where $h_1$ is a polynomial; otherwise only the part of the term due to the kinks of $h_1$ needs to be bounded through its increments. Groups of runs, in practice those with large $k$, are bounded together through their increments, and any term may instead be bounded through its oscillation.

\subsection{Symmetry reduction}
For $i=1,2,3$ let $\tau_i$ be the transposition of $1$ and $i$ ($\tau_1=\id$), and let $\tau_{23}$ be the transposition of $2$ and $3$. Every permutation $\pi$ (Section~\ref{sec:souls}) preserves the plane $\{\one^*x=1\}$, the triangle $\cT$, the plane $T$ and its Euclidean norm, and maps $C_i$ onto $C_{\pi(i)}$, because $(\pi x)_{\pi(i)}=x_i$; in particular $\tau_i(C_1)=C_i$. So the permutations act on $\mathcal C\times\PT$ by
\[
  \pi\cdot(y,[v],i)=(\pi y,[\pi v],\pi(i)).
\]
A function $g$ on $\mathcal C\times\PT$ is \emph{invariant} if $g(\pi\cdot z)=g(z)$ for every permutation $\pi$ and every $z\in\mathcal C\times\PT$.

\begin{lemma}\label{lem:sym}
\begin{enumerate}
\item[(i)] If $g$ is a bounded invariant function on $\mathcal C\times\PT$, then so is $\cK_\kappa g$.
\item[(ii)] The restriction $h\mapsto h_1=h(\cdot,\cdot,1)$ is a bijection from the bounded invariant functions on $\mathcal C\times\PT$ onto the bounded $\tau_{23}$-invariant functions on $C_1\times\PT$, with inverse $h(y,[v],i)=h_1(\tau_iy,[\tau_iv])$. In particular $\inf h=\inf h_1$ and $\sup h=\sup h_1$.
\item[(iii)] If $h$ is bounded and invariant, then $(\cK_\kappa h)(x,[w],1)=(\tilde\cK_\kappa h_1)(x,[w])$ for $(x,[w])\in C_1\times\PT$, where
\[
(\tilde\cK_\kappa h_1)(x,[w])=\sum_{j=2}^3\sum_{k\ge1}\bigl(1+2kx_j\bigr)^{-3(1-\kappa)}\abs{D_x\psi_j^k(\hat w)}^\kappa\,h_1\bigl(\tau_j\psi_j^k(x),[\tau_jD_x\psi_j^k(w)]\bigr).
\]
\end{enumerate}
Consequently, $\tilde\cK_\kappa$ maps $\tau_{23}$-invariant functions to $\tau_{23}$-invariant functions, and $\cK_\kappa h\ge h$ (resp.\ $\cK_\kappa h\le\theta h$) on $\mathcal C\times\PT$ if and only if $\tilde\cK_\kappa h_1\ge h_1$ (resp.\ $\tilde\cK_\kappa h_1\le\theta h_1$) on $C_1\times\PT$.
\end{lemma}

\begin{proof}
Let $\pi$ be a permutation. Since $\pi\circ\psi_j=\psi_{\pi(j)}\circ\pi$, also $\pi\circ\psi_j^k=\psi_{\pi(j)}^k\circ\pi$, and, differentiating, $\pi\,D_x\psi_j^k(w)=D_{\pi x}\psi_{\pi(j)}^k(\pi w)$ for $w\in T$. Moreover $(\pi x)_{\pi(j)}=x_j$, and $\pi$ is an isometry of $T$.

(i) Let $z=(x,[w],i)$. In the sum defining $(\cK_\kappa g)(\pi\cdot z)$ write $j=\pi(j')$; then $j\neq\pi(i)$ if and only if $j'\neq i$, and by the identities above the term of $(j,k)$ is
\[
\bigl(1+2kx_{j'}\bigr)^{-3(1-\kappa)}\abs{D_x\psi_{j'}^k(\hat w)}^\kappa\,g\bigl(\pi\cdot(\psi_{j'}^k(x),[D_x\psi_{j'}^k(w)],j')\bigr).
\]
Since $g$ is invariant, this is the term of $(j',k)$ in $(\cK_\kappa g)(z)$. So $(\cK_\kappa g)(\pi\cdot z)=(\cK_\kappa g)(z)$.

(ii) Since $\tau_i$ is an involution and $\tau_i(1)=i$, every $z=(y,[v],i)$ equals $\tau_i\cdot(\tau_iy,[\tau_iv],1)$, with $\tau_iy\in C_1$. So an invariant $h$ satisfies $h(y,[v],i)=h_1(\tau_iy,[\tau_iv])$ and is determined by $h_1$, and $h_1$ is $\tau_{23}$-invariant because $\tau_{23}$ fixes the index $1$. Conversely, let $h_1$ be $\tau_{23}$-invariant and define $h$ by this formula. For a permutation $\pi$ and $z=(y,[v],i)$, put $x=\tau_iy$, $w=\tau_iv$ and $\pi'=\tau_{\pi(i)}\pi\tau_i$; then $h(\pi\cdot z)=h_1(\pi'x,[\pi'w])$. Since $\pi'(1)=1$, $\pi'$ is $\id$ or $\tau_{23}$, so $h(\pi\cdot z)=h_1(x,[w])=h(z)$. In both directions $h$ and $h_1$ take the same values.

(iii) For $(x,[w])\in C_1\times\PT$ the sum defining $(\cK_\kappa h)(x,[w],1)$ runs over $j\in\{2,3\}$, and its terms contain $h\bigl(\psi_j^k(x),[D_x\psi_j^k(w)],j\bigr)$, where $\psi_j^k(x)\in C_j$ by Lemma~\ref{lem:trap}. By (ii) this equals $h_1\bigl(\tau_j\psi_j^k(x),[\tau_jD_x\psi_j^k(w)]\bigr)$.

Finally, by (i) and (iii), $\tilde\cK_\kappa h_1$ is the restriction of the invariant function $\cK_\kappa h$, so it is $\tau_{23}$-invariant by (ii). By (ii) again, $(\cK_\kappa h-\theta h)(y,[v],i)=(\tilde\cK_\kappa h_1-\theta h_1)(\tau_iy,[\tau_iv])$ for every $\theta$, and $(y,[v],i)\mapsto(\tau_iy,[\tau_iv])$ maps $\mathcal C\times\PT$ onto $C_1\times\PT$; this gives the last statement.
\end{proof}

\subsection{Piecewise-linear test functions}\label{ssec:P1}
We use the coordinates $\lambda=(\lambda_1,\lambda_2)=(2-4x_2,2-4x_3)$, in which $C_1$ is the triangle $\Lambda=\{\lambda\ge0,\ \lambda_1+\lambda_2\le1\}$ and $\tau_{23}$ exchanges $\lambda_1$ and $\lambda_2$.
We parametrize $\PT$ by $\vartheta\in\R/4\Z$ as follows. For $v\in T\setminus\{0\}$ let $\xi=v\cdot g_1=3v_1$ and $\eta=v\cdot g_2=v_2-v_3$, where $g_1=(2,-1,-1)$ and $g_2=(0,1,-1)$ form an orthogonal basis of $T$, and set
\[
  \vartheta=1+\eta/\xi\ \text{ if }\abs\eta\le\abs\xi,\qquad \vartheta=3-\xi/\eta\ \text{ otherwise}.
\]
This is a homeomorphism of $\PT$ onto $\R/4\Z$, monotone with respect to the angle, and $\tau_{23}$ acts by $\vartheta\mapsto2-\vartheta$. Across the lines $\abs\eta=\abs\xi$ the two expressions have the same first derivatives, so $\vartheta$ is $C^{1,1}$ as a function of $(\xi,\eta)\ne0$. The vectors $w(\vartheta)=(2,-1+3s,-1-3s)$, $s=\vartheta-1$, for $\vartheta\in[0,2]$, and $w(\vartheta)=(2r,3-r,-3-r)$, $r=3-\vartheta$, for $\vartheta\in[2,4]$, represent the directions with parameter $\vartheta$, and they are affine in $\vartheta$ on each of the two intervals.

The formula defining $(\tilde\cK_\kappa h_1)(x,[w])$ makes sense for every $x\in\cT$ with $x_2,x_3\ge\frac14$, that is, for $\lambda$ in the square $[0,1]^2$: by Lemma~\ref{lem:trap}, $\tau_j\psi_j^k(x)\in C_1$, and $2x_j\ge\frac12$, so the series converges as in the proof of Proposition~\ref{cor:certA}. 
We use this extension of $\tilde\cK_\kappa h_1$ to $[0,1]^2\times\PT$ below.

\emph{Meshes.} A \emph{mesh} is given by numbers $0=\ell_0<\ell_1<\dots<\ell_M=1$ with $\ell_m+\ell_{M-m}=1$, and numbers $0=\vartheta_0<\vartheta_1<\dots<\vartheta_{N-1}<4$ such that $2$ is one of them and $\vartheta\mapsto2-\vartheta\pmod 4$ permutes them; we put $\vartheta_N=4$. For $m+n\le M-1$ let $Q^+_{mn}$ be the triangle with vertices $(\ell_m,\ell_n)$, $(\ell_{m+1},\ell_n)$, $(\ell_m,\ell_{n+1})$, and for $m+n\le M-2$ let $Q^-_{mn}$ be the triangle with vertices $(\ell_{m+1},\ell_n)$, $(\ell_m,\ell_{n+1})$, $(\ell_{m+1},\ell_{n+1})$. Since $\ell_m+\ell_{M-m}=1$, the nodes $(\ell_m,\ell_{M-m})$ lie on the side $\lambda_1+\lambda_2=1$, and these triangles tile $\Lambda$. A \emph{cell} is a product $\gamma=Q\times\Theta_t$ of one of these triangles and of an arc $\Theta_t=[\vartheta_t,\vartheta_{t+1}]$, $0\le t<N$; we denote by $R_\gamma$ the rectangle $[\ell_m,\ell_{m+1}]\times[\ell_n,\ell_{n+1}]$ that contains $Q=Q^\pm_{mn}$. The map $\tau_{23}$ sends the cell $Q^\pm_{mn}\times\Theta_t$ to the cell $Q^\pm_{nm}\times(2-\Theta_t)$.

A \emph{piecewise-linear function} on the mesh is a continuous function $h_1$ on $C_1\times\PT$ such that on every cell $\gamma=Q\times\Theta_t$
\[
  h_1(\lambda,\vartheta)=h_\gamma(\lambda,\vartheta):=(1-s)\,a_t(\lambda)+s\,a_{t+1}(\lambda),\qquad s=\frac{\vartheta-\vartheta_t}{\vartheta_{t+1}-\vartheta_t},
\]
where $a_t$ and $a_{t+1}$ are the affine functions on $\R^2$ that interpolate the values of $h_1$ at the points $(v,\vartheta_t)$ and $(v,\vartheta_{t+1})$, $v$ a vertex of $Q$. It is determined by its values at the nodes $(\ell_m,\ell_n,\vartheta_t)$, $m+n\le M$, and it is $\tau_{23}$-invariant if and only if these values are. The polynomial $h_\gamma$, defined on all of $\R^3$, is the \emph{extension} of $h_1|_\gamma$; it is affine in $\lambda$ and in $\vartheta$ separately, so its Hessian is constant, with zero diagonal.

\emph{Interpolation.} Let $Q$ be either a box $\prod_{i=1}^3[p_i,p_i+\Delta_i]$ in the coordinates $(\lambda_1,\lambda_2,\vartheta)$, or a prism $P\times[p_3,p_3+\Delta_3]$ over a right triangle $P$ whose legs, of lengths $\Delta_1$ and $\Delta_2$, are parallel to the coordinate axes, as are the triangles $Q^\pm_{mn}$. Let $\mathcal V$ be the set of its $8$ or $6$ vertices, and let $If=\sum_{v\in\mathcal V}\phi_vf(v)$ be the interpolation of $f$ at $\mathcal V$ which is multilinear on boxes, and on prisms is affine on $P$ and linear in the third variable. The weights $\phi_v$ are non-negative, $\sum_v\phi_v=1$ and $\sum_v\phi_v(y)\,v=y$ for $y\in Q$. Moreover $Ih_\gamma=h_\gamma$ for every cell $\gamma$ and every such $Q$.

\begin{lemma}\label{lem:interp}
Let $Q$ and $I$ be as above, and let $y\in Q$.
\begin{enumerate}
\item[(i)] If $f(q)-f(p)\in\sum_{i=1}^3[G_i^-,G_i^+]\,(q_i-p_i)$ for all $p,q\in Q$, then
$\abs{f(y)-(If)(y)}\le\frac12\sum_i(G_i^+-G_i^-)\Delta_i$.
\item[(ii)] If $f\in C^1(Q)$ has a Lipschitz gradient, and its Hessian, which exists almost everywhere, lies in an interval matrix $\bigl([H_{ik}^-,H_{ik}^+]\bigr)$, then $-A^+\le f(y)-(If)(y)\le A^-$, where
\[
  A^\pm=\tfrac18\Bigl(\sum_i(H_{ii}^\pm)_\pm\,\Delta_i^2+2\sum_{i<k}\abs{H_{ik}}\Delta_i\Delta_k\Bigr),
\]
$t_+=\max(t,0)$, $t_-=\max(-t,0)$ and $\abs{H_{ik}}=\max\bigl(\abs{H_{ik}^-},\abs{H_{ik}^+}\bigr)$.
\end{enumerate}
\end{lemma}

\begin{proof}
(i) There is a vertex $v$ with $\abs{y_i-v_i}\le\Delta_i/2$ for $i=1,2,3$. On a box, take the nearer endpoint in each coordinate. On a prism over the triangle with vertices $0$, $\Delta_1e_1$, $\Delta_2e_2$ (the other position is symmetric), take $0$ if $y_1\le\Delta_1/2$ and $y_2\le\Delta_2/2$, take $\Delta_1e_1$ if $y_1>\Delta_1/2$, in which case $y_2<\Delta_2/2$, and take $\Delta_2e_2$ otherwise; in the third coordinate take the nearer end. For each $i$ the differences $u_i-v_i$, $u\in\mathcal V$, have the same sign, and $\sum_u\phi_u(y)(u_i-v_i)=y_i-v_i$. Hence
\[
  (If)(y)-f(v)=\sum_u\phi_u(y)\bigl(f(u)-f(v)\bigr)\in\sum_i[G_i^-,G_i^+]\,(y_i-v_i),
\]
and also $f(y)-f(v)\in\sum_i[G_i^-,G_i^+](y_i-v_i)$; subtracting gives the claim.

(ii) For $u\in\mathcal V$ let $d=u-y$. Since $\nabla f$ is Lipschitz, $t\mapsto f(y+td)$ has an absolutely continuous derivative $t\mapsto\nabla f(y+td)\cdot d$, whose derivative is $d^*D^2f(y+td)\,d$ at every $t$ at which $D^2f(y+td)$ exists. The Hessian exists almost everywhere, but the segment from $y$ to $u$ is a null set; however, for each vertex $u$, by Fubini's theorem in polar coordinates centred at $u$, for almost every $y\in Q$ the Hessian exists at almost every point of the segment from $y$ to $u$. So, for almost every $y\in Q$ and every $u\in\mathcal V$,
$f(u)=f(y)+\nabla f(y)\cdot d+\int_0^1(1-t)\,d^*D^2f(y+td)\,d\,dt$; we prove the bounds for these $y$, and they extend to all $y\in Q$ because both sides are continuous in $y$. Multiplying by $\phi_u(y)$ and summing over $u$, the linear terms cancel, and
\[
  (If)(y)-f(y)=\sum_u\phi_u(y)\int_0^1(1-t)\,d^*D^2f(y+td)\,d\,dt .
\]
For every $H$ in the interval matrix, $d^*Hd\le\sum_i(H^+_{ii})_+d_i^2+2\sum_{i<k}\abs{H_{ik}}\abs{d_i}\abs{d_k}$; moreover $\int_0^1(1-t)\,dt=\frac12$, and
\[
  \sum_u\phi_u(y)\,d_i^2\le\frac{\Delta_i^2}4,\qquad \sum_u\phi_u(y)\,\abs{d_i}\abs{d_k}\le\frac{\Delta_i\Delta_k}4\quad(i\ne k).
\]
On a box the weights are products of one-dimensional ones, and these sums are $(y_i-p_i)(p_i+\Delta_i-y_i)$ and $\prod_{l\in\{i,k\}}2(y_l-p_l)(p_l+\Delta_l-y_l)/\Delta_l$. On the triangle with vertices $0$, $\Delta_1e_1$, $\Delta_2e_2$, writing $y=(\alpha\Delta_1,\beta\Delta_2)$ with $\alpha,\beta\ge0$ and $\alpha+\beta\le1$, they are $\alpha(1-\alpha)\Delta_1^2$, $\beta(1-\beta)\Delta_2^2$ and $\alpha\beta(3-2\alpha-2\beta)\Delta_1\Delta_2\le\frac14\Delta_1\Delta_2$, and $\sum_u\phi_u\abs{d_1}=2\alpha(1-\alpha)\Delta_1\le\frac12\Delta_1$; on a prism the mixed sums with the third coordinate factor. This proves the lower bound, and the upper bound follows by applying it to $-f$.
\end{proof}

\begin{proposition}\label{prop:P1}
Let $0<\kappa<1$ and $\theta>0$, and let $h_1$ be a positive, $\tau_{23}$-invariant piecewise-linear function on a mesh. Let $\gamma=Q^\pm_{mn}\times\Theta_t$ be a cell, and let $Q$ be either the prism $\gamma$ itself or a box contained in $R_\gamma\times\Theta_t$. Suppose that on $Q$
\[
  \tilde\cK_\kappa h_1=F_S+F_K+F_R ,
\]
where $F_S$ satisfies the hypothesis of Lemma~\ref{lem:interp}(ii) with an interval matrix $[H]$, $F_K$ satisfies that of Lemma~\ref{lem:interp}(i) with intervals $[G_i^-,G_i^+]$, and $\sup_QF_R-\inf_QF_R\le o$. Let $A^\pm$ be as in Lemma~\ref{lem:interp}(ii), let $E=\frac12\sum_i(G_i^+-G_i^-)\Delta_i$, and $L^\pm=A^\pm+E+o$.
\begin{enumerate}
\item[(i)] If $(\tilde\cK_\kappa h_1)(v)-h_\gamma(v)\ge L^+$ at every vertex $v$ of $Q$, then $\tilde\cK_\kappa h_1\ge h_1$ on $Q\cap\gamma$.
\item[(ii)] If $(\tilde\cK_\kappa h_1)(v)-\theta h_\gamma(v)\le-L^-$ at every vertex $v$ of $Q$, then $\tilde\cK_\kappa h_1\le\theta h_1$ on $Q\cap\gamma$.
\end{enumerate}
If, for every cell $\gamma=Q^\pm_{mn}\times\Theta_t$ with $m\le n$, the hypothesis of~(i) holds on each set $Q$ of a finite cover of $\gamma$, then $\sA(\GN)\ge2-\kappa$. If the hypothesis of~(ii) holds in the same way with $\theta<1$, then $\sA(\GN)\le2-\kappa$.
\end{proposition}

\begin{proof}
Let $\mu\in\{1,\theta\}$ and $F=\tilde\cK_\kappa h_1-\mu h_\gamma=(F_S-\mu h_\gamma)+F_K+F_R$ on $Q$. Since $Ih_\gamma=h_\gamma$, Lemma~\ref{lem:interp} gives
$(F_S-\mu h_\gamma)+F_K\ge I\bigl(F_S-\mu h_\gamma+F_K\bigr)-A^+-E$ on $Q$, and the right side is at least $\min_v(F-F_R)(v)-A^+-E$, because the weights of $I$ are non-negative with sum $1$. Let $v_*$ be a vertex at which $F-F_R$ is minimal. Since $F_R(y)\ge F_R(v_*)-o$,
\[
  F(y)\ \ge\ (F-F_R)(v_*)+F_R(v_*)-A^+-E-o\ =\ F(v_*)-L^+\qquad(y\in Q).
\]
So under~(i) $F\ge0$ on $Q$, and $h_\gamma=h_1$ on $Q\cap\gamma$. The proof of~(ii) is the same, with maxima.
For the last statement, $\tilde\cK_\kappa h_1-\mu h_1$ is $\tau_{23}$-invariant by Lemma~\ref{lem:sym}, and $\tau_{23}$ maps the cells with $m\le n$ onto those with $m\ge n$. So $\tilde\cK_\kappa h_1\ge h_1$, respectively $\tilde\cK_\kappa h_1\le\theta h_1$, on all of $C_1\times\PT$. Since $h_1$ is continuous and positive on a compact set, $0<\inf h_1\le\sup h_1<\infty$, and Lemma~\ref{lem:sym} and Proposition~\ref{cor:certA} conclude.
\end{proof}

The proposition leaves the decomposition free; it may differ from one set $Q$ to another, and it is chosen term by term. A term is a single run, $t_{j,k}=W_{j,k}\cdot(h_1\circ\Phi_{j,k})$, where
\begin{align*}
  W_{j,k}(x,\vartheta)&=\bigl(1+2kx_j\bigr)^{-3(1-\kappa)}\abs{D_x\psi_j^k(\hat w(\vartheta))}^\kappa,\\
  \Phi_{j,k}(x,\vartheta)&=\bigl(\tau_j\psi_j^k(x),[\tau_jD_x\psi_j^k(w(\vartheta))]\bigr)
\end{align*}
in the coordinates $(\lambda,\vartheta)$, or a group of runs. The functions $W_{j,k}$ and $\Phi_{j,k}$ are real-analytic on $Q$, except that the $\vartheta$-component of $\Phi_{j,k}$ is only $C^{1,1}$ across the lines where the chart changes. Three cases occur.
\begin{enumerate}
\item[(a)] If $\Phi_{j,k}(Q)$ lies in a single cell $\gamma'$, then $t_{j,k}=W_{j,k}\cdot(h_{\gamma'}\circ\Phi_{j,k})$ on $Q$, a $C^{1,1}$ function, and it is put in $F_S$.
\item[(b)] Otherwise, let $B$ be a box containing $\Phi_{j,k}(Q)$. On $B\cap(\Lambda\times\R)$ the function $h_1$ is continuous and piecewise bilinear, so its increments lie in $[\nabla h_1]\cdot(q-p)$, where $[\nabla h_1]$ is the hull of the gradients of the bilinear functions $h_{\gamma''}$ on $B\cap\gamma''$, over the cells $\gamma''$ meeting $B$. By the mean value theorem for $W_{j,k}$ and $\Phi_{j,k}$, the increments of $t_{j,k}$ lie in $\bigl([\nabla W_{j,k}][h_1]+[W_{j,k}][\nabla h_1][D\Phi_{j,k}]\bigr)(q-p)$, where $[\,\cdot\,]$ denotes enclosures over $Q$, respectively over $B$, and $t_{j,k}$ may be put in $F_K$. Better, for a cell $\gamma'$ write $t_{j,k}=W_{j,k}\cdot(h_{\gamma'}\circ\Phi_{j,k})+W_{j,k}\cdot\bigl((h_1-h_{\gamma'})\circ\Phi_{j,k}\bigr)$: the first term goes in $F_S$, and the second in $F_K$, with the same enclosure in which $h_1$ is replaced by $h_1-h_{\gamma'}$. The enclosures of $h_1-h_{\gamma'}$ and of its gradient are small when $B$ is small: they are controlled by the differences between the gradients of $h_1$ on neighbouring cells, that is, by the second differences of its nodal values.
\item[(c)] The sum over a group of runs $k\in E$ goes in $F_K$, with the increment enclosure
\[
  [W_E]\,\bigl([\nabla\log W]\,[h_1]+[\nabla h_1][D\Phi]\bigr),
\]
where $[W_E]$ is an interval containing $\sum_{k\in E}W_{j,k}(y_k)$ for all choices of points $y_k\in Q$, $[\nabla\log W]$ and $[D\Phi]$ enclose the gradients of $\log W_{j,k}$ and the derivatives of $\Phi_{j,k}$ over $Q$ uniformly in $k\in E$, and $[h_1]$, $[\nabla h_1]$ are hulls over the cells meeting a box that contains all the images. Indeed, the increment of $W_{j,k}\cdot(h_1\circ\Phi_{j,k})$ between $p$ and $q$ is $\bigl(W_{j,k}(q)-W_{j,k}(p)\bigr)h_1(\Phi_{j,k}(q))+W_{j,k}(p)\bigl(h_1(\Phi_{j,k}(q))-h_1(\Phi_{j,k}(p))\bigr)$, and $W_{j,k}(q)-W_{j,k}(p)=W_{j,k}(y)\nabla\log W_{j,k}(y)\cdot(q-p)$ for some $y$ on the segment.
\end{enumerate}
Any term may instead be put in $F_R$, with its oscillation over $Q$.

\section{The computation}
\label{sec:comp}

In this section we describe how we use the finite check of Section~\ref{sec:disc} to prove Theorem~\ref{thm: dimH}, namely that $1.6853<\dimH\cN<1.6854$. 
Since $\dimH\cN=\sA(\GN)$ (Corollary~\ref{cor: dimB}), it suffices to prove $\sA(\GN)\ge2-\kappa_-$ and $\sA(\GN)\le2-\kappa_+$ for two binary64 numbers $\kappa_\pm$ with $2-\kappa_->1.6853$ and $2-\kappa_+<1.6854$. For each of them we compute in floating point a positive, $\tau_{23}$-invariant, piecewise-linear function $h^\pm$, an approximate eigenfunction of $\tilde\cK_{\kappa_\pm}$; then a program, the \emph{verifier}, proves in interval arithmetic that $\tilde\cK_{\kappa_-}h^-\ge h^-$ and $\tilde\cK_{\kappa_+}h^+\le\theta h^+$ with $\theta<1$, by checking the hypotheses of Proposition~\ref{prop:P1} on every cell with $m\le n$. How the test functions were obtained plays no role in the proof; only the checks made by the verifier do.

We describe in turn: the mesh, which is graded towards two lines near which the eigenfunction varies on small scales (Remark~\ref{rem:fan}); the test functions; how the verifier encloses $\tilde\cK_\kappa h$ at a point, summing the runs $k\le1024$ one by one and grouping the longer ones into blocks; how, on a set $Q$, it splits $\tilde\cK_\kappa h$ into the three parts of Proposition~\ref{prop:P1} and encloses each of them; when it subdivides a cell into boxes; how rounding errors are controlled; how the programs were tested; and the result. The remarks at the end describe a coarser certificate that runs in a few hours on two cores (Remark~\ref{rem:coarse}), a cruder but independent certificate with piecewise-constant test functions (Remark~\ref{rem:P0}), and the numerical value $\sA(\GN)\approx1.6853416$, together with what certifying a fifth decimal place would cost (Remark~\ref{rem:sAnum}).

\medskip\noindent\textbf{The mesh.} The nodes $\ell_m\le\frac12$ are $\ell_0=0$, the numbers
\[
  \ell_{m+1}=\ell_m+\min\bigl(5\cdot10^{-5}+\tfrac18\,\ell_m,\ 10^{-2}\bigr),\qquad 0\le m\le85,
\]
and $\ell_{87}=\frac12$. So the spacing starts at $5\cdot10^{-5}$ and grows by the factor $\frac98$ until it reaches $10^{-2}$, at $\ell\approx0.08$; the last two cells have widths $0.0100$ and $0.0102$. These nodes are computed first, and then $\ell_{174-m}=1-\ell_m$ exactly, so $M=174$. Let $\vartheta_*$ be the binary64 number nearest to $\frac43$. The nodes $\vartheta_t$ in $[1,2]$ are $1$, $\vartheta_*$ and $2$, the $96$ numbers $\vartheta_*-d_i$, $1\le i\le96$, where
\[
  d_0=0,\qquad d_{i+1}=d_i+\min\bigl(5\cdot10^{-6}+\tfrac18\,d_i,\ \tfrac1{128}\bigr),
\]
and the $266$ numbers $\vartheta_*+d'_i$, $1\le i\le266$, where $d'_0=0$ and $d'_{i+1}=d'_i+s_i$, with $s_i=\min(5\cdot10^{-6}+\frac1{20}d'_i,\ 2.5\cdot10^{-3})$ while $d'_i<\frac14$, and $s_i=\min(1.075\,s_{i-1},\ \frac1{128})$ afterwards. So on both sides of $\vartheta_*$ the spacing starts at $5\cdot10^{-6}$: towards $1$ it grows by the factor $\frac98$ up to $\frac1{128}$; towards $2$ it grows by the factor $\frac{21}{20}$ up to $2.5\cdot10^{-3}$, keeps this value up to $\vartheta\approx\vartheta_*+\frac14$, and then grows by the factor $1.075$ up to $\frac1{128}$. The last cells before $1$ and $2$ have widths $0.0088$ and $0.0091$. The nodes in $[3,4)$ are the multiples of $\frac1{128}$. The nodes in $[1,2]\cup[3,4)$ are computed in binary64, and the others are their images under $\vartheta\mapsto2-\vartheta\pmod4$, which are computed exactly; so $N=984$. The mesh has $29\,791\,584$ cells and $15\,153\,600$ nodes. It is stored in the files of the test functions, from which the verifier reads it.

\begin{remark}[the mesh]\label{rem:fan}
    The refinement is dictated by how the runs turn the needle. 
    At a point $x$ of the side $\{x_3=\frac12\}$ of $C_1$ and for $j=2$, the vector $c$ of Lemma~\ref{lem:contraction} is parallel to $e_1-e_2$, the direction with $\vartheta=\frac23$. 
    The runs of $2$, followed by $\tau_2$, map this side into itself and fix the direction $[e_1-e_2]$. 
    This direction spans the kernel of $\id+c\,g^*$ in \eqref{eq:DRB}, and it is a repelling fixed point of the projective action of $D_x\psi_2^k$, with expansion factor $d$. 
    Accordingly, the eigenfunction of $\tilde\cK_\kappa$ varies on scales proportional to the distance from the line $\{\lambda_2=0,\ \vartheta=\frac23\}$. 
    The same happens near the line $\{\lambda_1=0,\ \vartheta=\frac43\}$, where the kernel direction for $j=3$ is $[e_1-e_3]$. 
    The mesh is graded geometrically towards these two lines. 
    In $\vartheta$ it is finer on the side $\vartheta>\frac43$ of the second line, and on the side $\vartheta<\frac23$ of the first, where, roughly along $\vartheta-\frac43\approx\frac12\lambda_1$ and its image under $\tau_{23}$, the eigenfunction passes from its behaviour near the line to its behaviour away from it. 
    With uniform meshes, the oscillation of $(\tilde\cK_\kappa h)/h$ between the nodes stays above $10^{-2}$ up to $M=128$; on the graded mesh of Remark~\ref{rem:coarse} it is about $2.5\times10^{-3}$; on the present mesh, after the correction of the test functions described below, the sampled extremes of the ratio lie within $4.1\times10^{-5}$ of the nodal eigenvalue (Remark~\ref{rem:sAnum}).
\end{remark}

\medskip\noindent\textbf{Test functions.} We use two piecewise-linear functions $h^-$ and $h^+$ on this mesh, computed in floating point at the values $\kappa_-$ and $\kappa_+$ of $\kappa$ given below. Each is obtained in two stages.
First, the power iteration $h\leftarrow\tilde\cK_\kappa h/\max\tilde\cK_\kappa h$ on the nodal values (a collocation at the nodes, in which the runs $k>64$ are grouped into blocks of ratio $1.05$, each evaluated at its middle run), symmetrized under $\tau_{23}$ after each step. It makes $(\tilde\cK_\kappa h)/h$ nearly constant at the nodes, but between the nodes, near the lines of Remark~\ref{rem:fan}, the ratio still departs from its nodal value by more than the margins available here (Remark~\ref{rem:sAnum}).
Second, a correction that equalizes instead the extremes of the ratio around each node: for $h^-$,
\[
  h(v)\ \leftarrow\ h(v)\,\min\Bigl\{\frac{\tilde\cK_\kappa h}{h}(p):\ p=v\ \text{or}\ p\in P_v\Bigr\},
\]
followed by the $\tau_{23}$-symmetrization and a normalization, where $P_v$ is the set of the sample points whose interpolation stencil contains the node $v$: the vertices, the midpoints of the edges and the centroids of the triangles $Q^\pm_{mn}$, at the values of $\vartheta$ at $\frac16$, $\frac12$ and $\frac56$ of each arc $\Theta_t$. For $h^+$ the minimum is replaced by the maximum. The correction reduces to a step of the power iteration when the ratio has no extremes between the nodes, and it converges in a few steps.
We used twenty power iterations followed by twelve corrections, starting from the test functions $h^\pm_0$ of Remark~\ref{rem:coarse}, interpolated to this mesh (program \texttt{p1solve2.c}, an extension of \texttt{p1solve.c}, with the scripts that call it).
The nodal values of $h^\pm$ lie in $[0.5187,1]$. The verifier checks that they are finite, positive and exactly $\tau_{23}$-invariant, and that the mesh has the symmetries required in \S\ref{ssec:P1}; the identities $\ell_m+\ell_{M-m}=1$ and $\vartheta_t+\vartheta_{t^*}\in\{2,6\}$ are checked with error-free transformations of the sums.
The certificate does not depend on how the test functions were obtained.

\medskip\noindent\textbf{Values at the vertices.} At a point $(x,\vartheta)$, the runs $k\le1024$ are summed one by one. Let $u=2x_j$ and $d_k=1+ku$. By \eqref{eq:DR}, and since $ku\,c=(d_k-1)c$, the image direction is that of
\[
  V_k(x)w=d_k^2\,D_x\psi_j^k(w)=d_kw+ku\,w_jc=d_k\,p_w-w_jc,\qquad p_w=w+w_jc,
\]
the weight is $W_{j,k}=(d_k\abs{V_k(x)\hat w})^\kappa d_k^{-3}$, and the image point is $\tau_j\psi_j^k(x)$, with $\psi_j^k(x)_j=x_j/d_k$ and $\psi_j^k(x)_m=\frac12-(\frac12-x_m)/d_k$ for $m\ne j$. The value of $h_1$ at the image is enclosed by the hull of the values of the $h_{\gamma'}$ over the cells $\gamma'$ meeting the enclosure of the image.

The runs $k\ge1025$ are grouped into blocks $E=\{k_1,\dots,k_2\}$ with $k_2=\lfloor1.25\,k_1\rfloor$, up to $2^{28}$, followed by the tail. For a block, \eqref{eq:DR} gives $D_x\psi_j^k(\hat w)=(p_{\hat w}-\alpha_{\hat w}/d_k)/d_k$ with $\alpha_w=w_jc$, so the sum of the weights is
\[
  \sum_{k\in E}d_k^{-(3-2\kappa)}\,\bigl|p_{\hat w}-\alpha_{\hat w}/d_k\bigr|^\kappa .
\]
We enclose $\abs{p_{\hat w}-s\,\alpha_{\hat w}}$ for $s\in[1/d_{k_2},1/d_{k_1}]$, or $s\in[0,1/d_{k_1}]$ for the tail. The function $f(k)=(1+ku)^{-(3-2\kappa)}$ is convex and decreasing, so $\sum_{k\in E}f(k)$ lies between the trapezoid and the midpoint bounds
\[
  \int_{k_1}^{k_2}f+\tfrac12\bigl(f(k_1)+f(k_2)\bigr)\ \le\ \sum_{k=k_1}^{k_2}f(k)\ \le\ \int_{k_1-1/2}^{k_2+1/2}f ;
\]
for the tail, $k_2=\infty$ and $f(\infty)=0$. Since $\psi_j^k(x)$ is coordinatewise monotone in $k$, the image points of a block lie between those of $k_1$ and $k_2$, or of $k_1$ and the limit point $\pr{\one-e_j}$ for the tail. As for the directions, $V_k(x)w$ is affine in $k$. So every $V_k(x)w$ with $k\in E$ is a positive combination of $V_{k_1}(x)w$ and $V_{k_2}(x)w$, and for the tail of $V_{k_1}(x)w$ and $U(x)w=u\,p_w$. We enclose these vectors, for all the $x$ and $w$ under consideration, in boxes of the $(\xi,\eta)$-plane. If the corners of the boxes lie in an open half-plane, we take the arc of directions that they span; otherwise we take all of $\PT$. The relative width of the resulting enclosures of $(\tilde\cK_\kappa h)(v)$ is at most $2.4\times10^{-6}$ at the nodes, about one eighteenth of the smaller margin of the final test functions (Remark~\ref{rem:sAnum}); with the runs $k\le256$ summed one by one, as in Remark~\ref{rem:coarse}, it would be ten times larger.

\medskip\noindent\textbf{Enclosures on a prism or a box.} On a set $Q$ as in Proposition~\ref{prop:P1}, the runs $k\le64$ are treated one by one as follows.
\begin{itemize}
\item \emph{Runs $k\le16$.} Interval arithmetic replaces each number by an interval that contains it, and each of the operations $+$, $-$, $\times$, $/$, $\sqrt{\ }$ by an operation on intervals whose result contains every possible result of the operation on numbers taken from the operands~\cite{MKC09}. Evaluating a formula in interval arithmetic, with the coordinates $(\lambda_1,\lambda_2,\vartheta)$ replaced by their ranges over the smallest box containing $Q$, therefore gives an interval that contains the values of the formula at all points of $Q$. To enclose derivatives as well, we carry along with every intermediate quantity $g$ its gradient $\nabla g$ and its Hessian $D^2g$, and we update all three at each operation by the rules of calculus; for a product, for instance,
\[
  \nabla(g_1g_2)=g_1\nabla g_2+g_2\nabla g_1,
  \]
  \[
  D^2(g_1g_2)=g_1D^2g_2+g_2D^2g_1+\nabla g_1(\nabla g_2)^*+\nabla g_2(\nabla g_1)^* .
\]
This is forward-mode automatic differentiation~\cite{GW08}: the derivatives are computed exactly by the chain rule, not approximated by difference quotients, and in interval arithmetic each resulting interval contains the corresponding derivative at every point of $Q$. In this way we obtain enclosures of the values, gradients and Hessians of $W_{j,k}$ and $\Phi_{j,k}$ over $Q$.

Interval evaluation overestimates ranges, because each occurrence of a variable is treated independently: on $[0,1]$, for instance, the expression $x-x$ gives $[-1,1]$. To reduce this, every intermediate quantity $g$ (the vector $V_k(x)w$, its coordinates $(\xi,\eta)$, the ratio $\abs{V_k(x)w}^2/\abs w^2$ and $d_k\abs{V_k(x)\hat w}$) is also evaluated at the centre $y_0$ of $Q$. By the mean value theorem, $g(y)\in g(y_0)+[\nabla g](y-y_0)$ and $\nabla g(y)\in\nabla g(y_0)+[D^2g](y-y_0)$ for every $y\in Q$, where $[\nabla g]$ and $[D^2g]$ are the enclosures over $Q$. So the enclosures of $g$ and $\nabla g$ may be intersected with $g(y_0)+[\nabla g](Q-y_0)$ and $\nabla g(y_0)+[D^2g](Q-y_0)$, and the results are still enclosures, usually much narrower ones.

Across a change of chart, the enclosures of the derivatives of $\vartheta$ are the hulls of those of the two expressions, which is legitimate because $\vartheta$ is $C^{1,1}$. The term is then treated as in case~(a) if the enclosure of $\Phi_{j,k}(Q)$ meets a single cell. Otherwise it is treated as in case~(b), with $\gamma'$ the cell containing the centre of the enclosure, choosing among the two forms of case~(b) and $F_R$ the one with the smallest contribution to $L^\pm$.
\item \emph{Runs $17\le k\le64$.} We use first-order automatic differentiation and the first form of case~(b).
\item \emph{Runs $k\ge65$.} They are grouped into blocks of ratio $1.25$ and a tail, as above, and treated as in case~(c). The enclosure of
\[
  \nabla\log W_{j,k}=-(3-2\kappa)\,\frac k{d_k}\,\nabla u+\kappa\,\frac{\nabla\abs{p_{\hat w}-s\,\alpha_{\hat w}}}{\abs{p_{\hat w}-s\,\alpha_{\hat w}}},\qquad s=\frac1{d_k},
\]
uses $k/d_k\in[k_1/d_{k_1},k_2/d_{k_2}]$ and $\nabla s=-s\,(k/d_k)\nabla u$; those of $\nabla p_{\hat w}$ and $\nabla\alpha_{\hat w}$ are obtained by first-order automatic differentiation.
\end{itemize}
The three treatments trade accuracy for cost. The runs $k\le16$ carry roughly $95\%$ of the weight, and only for them is it worth enclosing Hessians, which makes the second-order bounds of cases~(a) and~(b) available; the runs $17\le k\le64$ carry a few per cent, and a first-order bound suffices; the infinitely many runs $k\ge65$ must be grouped, and they carry about $1\%$. Any other choice of the thresholds gives a valid certificate; only the margins and the running time depend on it.

A term or a group is moved to $F_R$ whenever its oscillation over $Q$ is smaller than its contribution to $L^\pm$. Finally, the program adds to $L^\pm$ a bound for the interpolation error of $\mu h_\gamma$ through its Hessian. This term is redundant, because $Ih_\gamma=h_\gamma$, and it only strengthens the inequalities that are checked.

\medskip\noindent\textbf{Subdivision.} First $\tilde\cK_\kappa h$ is evaluated at the nodes that are vertices of cells with $m\le n$, and each such cell is tested as a prism, by Proposition~\ref{prop:P1} with $Q=\gamma$. If the test fails, the box $R_\gamma\times\Theta_t$ is tested, and then bisected recursively, up to $18$ times, along the coordinate with the largest contribution to $L^\pm$. Its eight corners are nodes or, near the side $\lambda_1+\lambda_2=1$, points of $[0,1]^2\times\PT$ outside $C_1\times\PT$, where $\tilde\cK_\kappa h$ is evaluated as well. Boxes that do not meet the triangle of $\gamma$ are discarded, and the values of $\tilde\cK_\kappa h$ at new corners are computed once and passed on to the children. A box is not bisected further if the inequality to be proved already fails at one of its vertices. A cell for which the test fails in this way, or still fails after the last bisection, is a failed cell, and the bound is certified only if there is no failed cell; in the two runs that prove Theorem~\ref{thm: dimH} there was none. Since the children of a box cover its intersection with $\gamma$, the bound kept for a box is the better of its own and the worst of its children's.

\medskip\noindent\textbf{Rigour.} All enclosures are computed in IEEE~754 binary64 arithmetic \cite{IEEE754} with rounding towards $+\infty$. Write $\mathrm{fl}^{\uparrow}$ and $\mathrm{fl}^{\downarrow}$ for the results rounded towards $+\infty$ and towards $-\infty$. Since $\mathrm{fl}^{\downarrow}(x)=-\mathrm{fl}^{\uparrow}(-x)$, the lower endpoints are computed as
\[
\begin{aligned}
  \mathrm{fl}^{\downarrow}(a+b)&=-\mathrm{fl}^{\uparrow}\bigl((-a)-b\bigr), &\qquad \mathrm{fl}^{\downarrow}(a-b)&=-\mathrm{fl}^{\uparrow}(b-a),\\
  \mathrm{fl}^{\downarrow}(a\times b)&=-\mathrm{fl}^{\uparrow}\bigl((-a)\times b\bigr), &\qquad \mathrm{fl}^{\downarrow}(a/b)&=-\mathrm{fl}^{\uparrow}\bigl((-a)/b\bigr),
\end{aligned}
\]
which are the formulas coded in \texttt{p1verify.c}; the lower endpoint of $\sqrt a$ is obtained by stepping down from $\mathrm{fl}^{\uparrow}(\sqrt a)$ to the preceding binary64 number. Since these operations are correctly rounded, every operation encloses its exact result \cite{Tuc11,MKC09}. The program checks at start-up that the directed rounding is honoured, and it is compiled without floating-point contraction. No function of the mathematical library enters any bound.

The following quantities are computed with Arb~\cite{Arb} (FLINT~3.6.0 \cite{FLINT}) and rounded outwards to binary64:
\begin{itemize}
\item the powers $g^\kappa$, for $g$ on a geometric grid of $2^{11}$ binary64 numbers per octave in $[2^{-40},2^{24}]$;
\item the integrals and the values of $f$ for the blocks and the tail, at the grid values $u=1-i/2^{11}$, $0\le i\le2^{10}$.
\end{itemize}
For a number $m$ between consecutive grid values $g<g'$, the concavity of $t\mapsto t^\kappa$ gives
\[
  g^\kappa+(m-g)\,\frac{g'^\kappa-g^\kappa}{g'-g}\ \le\ m^\kappa\ \le\ g^\kappa\Bigl(1+\kappa\,\frac{m-g}{g}\Bigr),
\]
which we evaluate with the outward-rounded table values. A power of an enclosed number is bounded by monotonicity, and the derivative $\kappa m^{\kappa-1}$ as $\kappa m^\kappa/m$. The sums over the blocks decrease in $u$, so an enclosure $[u^-,u^+]$ of $u$ is handled with the grid values below $u^-$ and above $u^+$. The values of $\kappa$ are the binary64 number $\kappa_-$ nearest to $0.3147$ and the binary64 number $\kappa_+$ that follows $0.3146$, namely $\kappa_-=0.314699999999999979\ldots$ and $\kappa_+=0.314600000000000046\ldots$, and they are used exactly. (The binary64 number nearest to $0.3146$ is smaller than $0.3146$, and it would only give $\sA(\GN)\le1.685400000000000009\ldots$.) When selecting the cells that meet an enclosure, every cell within relative distance $10^{-9}$ of it is included, so the selection is conservative despite rounding.

\medskip\noindent\textbf{Checks.} An independent floating-point routine, in which the formulas of Section~\ref{sec:cocycle} are coded directly, verified the following for the test functions of Remark~\ref{rem:coarse}, with the same verifier.
\begin{itemize}
\item The values and the gradients of $3000$ random terms $t_{j,k}$, $k\le64$, each at $10$ random points of a random cell ($30\,000$ values and $90\,000$ gradient components), lie in the computed enclosures.
\item At $300$ random points, the values of $\tilde\cK_\kappa h$ computed by direct summation up to $k=2\times10^5$, plus an asymptotic estimate of the remainder, lie in the computed enclosures, up to a relative error of $10^{-7}$.
\item At $19\,627$ random points of random boxes contained in the cells, the values of $(\tilde\cK_\kappa h)/h$ lie within the certified bounds for the boxes.
\item On those boxes, at $16\,100$ random points, the deviation of $\tilde\cK_\kappa h-h$ from its interpolant never exceeds $27\%$ of the certified bound.
\end{itemize}
No violation was found. Earlier versions of the code were tested in the same way on coarser meshes. For the final test functions, a floating-point survey of $(\tilde\cK_\kappa h)/h$, at several points of every cell and more densely on the $6000$ cells with the smallest margins, gives values in $[1.000044,1.000387]$ for $h^-$ at $\kappa_-$ and in $[0.999558,0.999916]$ for $h^+$ at $\kappa_+$; and the verifier, run before the full computation on every $4000$th cell and on each of the $20$ cells with the smallest sampled margin, proves both inequalities there.

\medskip\noindent\textbf{The result.} For each test function, the cells with $m\le n$ were split into $1499$ ranges of $10\,000$ consecutive cells (the last one has $1\,400$), each checked by a separate run of the verifier. The union of the ranges was checked to be the set of all the $14\,981\,400$ cells with $m\le n$, and every log was checked to come from the right test function, through its SHA-256 checksum, and from the same arguments of the verifier, which every log records; in particular, $\kappa$ is passed as the decimal strings \texttt{0.3147} and \texttt{0.31460000000000005}, whose binary64 values are $\kappa_-$ and $\kappa_+$. On six machines with $16$ cores and two hardware threads per core, the two computations took $8\,125$ and $3\,936$ thread-hours, about three days in all. Besides the $14\,981\,400$ prisms, the subdivision tested the following numbers of boxes:
\begin{center}
\begin{tabular}{ccccc}
\toprule
$\kappa$ & test function & inequality verified & boxes & conclusion\\
\midrule
$\kappa_-$ & $h^-$ & $\tilde\cK_\kappa h^-\ge h^-$ & $1\,025\,808\,647$ & $\sA(\GN)\ge2-\kappa_->1.6853$\\
$\kappa_+$ & $h^+$ & $\tilde\cK_\kappa h^+\le0.99999999\,h^+$ & $454\,613\,766$ & $\sA(\GN)\le2-\kappa_+<1.6854$\\
\bottomrule
\end{tabular}
\end{center}
The programs \texttt{p1verify.c}, \texttt{p1solve.c} and \texttt{p1solve2.c}, the scripts, the commands and the logs of the runs, the test functions $h^\pm$ and those of Remark~\ref{rem:coarse}, and the program \texttt{p2solve.c} of the collocations of Remark~\ref{rem:sAnum} are deposited at Zenodo~\cite{DLdata}. The files of $h^-$ and $h^+$ have the SHA-256 checksums
\begin{center}\small
$h^-$:\ \texttt{8c587e7cd1461153b6f39916d0caa2b5913664966945a3fc8f93cff38721928b}\\
$h^+$:\ \texttt{da0e2e05008f9dbeaebdcae48a930ba48948c22913b2eb02b1ad40c0c284afcc}
\end{center}
and every log begins with the checksum of the file it checked. In the programs a point $x\in\cT$ is stored through the vector $y=\one-2x$, which lies in the triangle with vertices $e_1,e_2,e_3$; this affine map is a similarity of ratio $2$ of the plane $\{\one^*x=1\}$, it commutes with the $\tau$'s, its linear part on $T$ is $-2\,\id$, so that the derivatives and the parameter $\vartheta$ of the directions are unchanged, and in these coordinates $2x_j=1-y_j$, $\pr{\one-e_j}$ becomes $e_j$, and $\lambda=2(y_2,y_3)$.

\thmdimH*

\begin{proof}
By Proposition~\ref{prop:P1} and the table, $1.6853<2-\kappa_-\le\sA(\GN)\le2-\kappa_+<1.6854$; indeed, $2-\kappa_-=1.68530000000000002\ldots$ and $2-\kappa_+=1.68539999999999995\ldots$.
By Corollary~\ref{cor: dimB}, the same bounds hold for $\dimH\cN=\dimB\cN$.
\end{proof}

\begin{remark}[a certificate on two cores]\label{rem:coarse}
A coarser version of the same computation proves $1.68364<\sA(\GN)<1.68637$ in a few hours on two cores, and it is a convenient test of the programs. Its mesh has the nodes $\ell_m\le\frac12$ equal to $0$, to the sixteen numbers $5\cdot10^{-4}\cdot1.4^m$, $0\le m\le15$, to the thirteen numbers $\ell_{16}+i\,(\frac12-\ell_{16})/14$, $1\le i\le13$, where $\ell_{16}=5\cdot10^{-4}\cdot1.4^{15}\approx0.0778$, and to $\frac12$, so that $M=60$. Its nodes $\vartheta_t$ form two clusters, around $\vartheta_c=\vartheta_*$ and $\vartheta_c=2-\vartheta_*$: the numbers $\vartheta_c$ and $\vartheta_c\pm10^{-4}\cdot1.25^m$, $0\le m\le28$; outside the clusters they are the multiples of $\frac1{64}$ at distance at least $0.0595$ from both values of $\vartheta_c$. The nodes in $[1,2]\cup[3,4)$ are computed in binary64, and the others are their images under $\vartheta\mapsto2-\vartheta\pmod4$; so $N=358$, and the mesh has $1\,288\,800$ cells and $676\,978$ nodes. Its test functions $h^\pm_0$ are obtained by the power iteration alone, at $\kappa=0.31635$ and $\kappa=0.31364$, starting from the constant function on a coarser mesh, with eight iterations on the final mesh (program \texttt{p1solve.c}, script \texttt{make\_testfunctions.sh}, which reproduce $h^\pm_0$ bit for bit); their nodal values lie in $[0.518,1]$. At the vertices, the runs $k\le256$ are summed one by one, and the relative width of the enclosures at the nodes is at most $2.2\times10^{-5}$. Let $\kappa^0_-=0.316350000000000020\ldots$ and $\kappa^0_+=0.313639999999999974\ldots$ be the binary64 numbers nearest to $0.31635$ and $0.31364$. The verifier proves $\tilde\cK_{\kappa^0_-}h^-_0\ge h^-_0$ and $\tilde\cK_{\kappa^0_+}h^+_0\le0.99999999\,h^+_0$ on all the $655\,140$ cells with $m\le n$, after testing $2\,173\,814$ and $3\,970\,640$ boxes; on two cores the two computations took $3.4$ and $4.5$ hours. So $1.68364<2-\kappa^0_-\le\sA(\GN)\le2-\kappa^0_+<1.68637$.
\end{remark}

\begin{remark}[piecewise-constant test functions]\label{rem:P0}
Let $H=(H_\gamma)$ be a positive, $\tau_{23}$-invariant vector indexed by the cells of a mesh, and let $h_1(p)=\min\{H_\gamma:\ p\in\gamma\}$. On each cell $\gamma$, $\tilde\cK_\kappa h_1$ is bounded below and above by
\[
  \sum_{j,E}\omega^-(\gamma,j,E)\min_{\gamma'\in S(\gamma,j,E)}H_{\gamma'}\qquad\text{and}\qquad\sum_{j,E}\omega^+(\gamma,j,E)\max_{\gamma'\in S(\gamma,j,E)}H_{\gamma'} ,
\]
where:
\begin{itemize}
\item the entries $E$ are single runs or groups of runs;
\item $\omega^\pm(\gamma,j,E)$ enclose the sums of the weights over $\gamma$;
\item $S(\gamma,j,E)$ contains the cells that meet the images of $\gamma$.
\end{itemize}
These are Collatz--Wielandt bounds for two non-negative matrices. Take a uniform mesh with $M=128$ and $256$ arcs, that is, $4.2\times10^6$ cells, and $\kappa=0.3209$ and $0.3086$. The bounds then give $1.6790\le\sA(\GN)\le1.6915$, and the verification takes about two minutes on two cores (program \texttt{rsa.c}, also in~\cite{DLdata}). The error of this discretization is of first order in the mesh size: the width of the certified interval halves when $M$ and the number of arcs are doubled, at eight times the cost. For piecewise-linear test functions the interpolation error is of second order, away from the lines of Remark~\ref{rem:fan}, and this is what makes the interval of Theorem~\ref{thm: dimH} affordable. The two programs are independent implementations apart from the interval primitives, so the piecewise-constant computation is an independent check of the piecewise-linear one, but only at the level of its own width, about $10^{-2}$.
\end{remark}

\begin{remark}[the numerical value]\label{rem:sAnum}
Let $F(\sigma)$ be the mesh defined as above with every spacing multiplied by $2\sigma$: for $\lambda$ the spacing is $\min\bigl(2\sigma(5\cdot10^{-5}+\frac18\ell),\ 2\sigma\cdot10^{-2}\bigr)$, the spacings in $\vartheta$ and their growth rates are rescaled in the same way, the threshold $\frac14$ is kept, and $[3,4)$ is divided into $\mathrm{round}(64/\sigma)$ equal arcs. So the mesh of Theorem~\ref{thm: dimH} is $F(\frac12)$, and $F(1)$ and $F(1/\sqrt2)$ have $(M,N)=(90,498)$ and $(124,702)$. On each mesh, the power iteration of the collocation at two values of $\kappa$ near $0.3146$ gives the value $\kappa_{\rm eig}$ at which the nodal eigenvalue of the collocation equals $1$; the last iteration brackets this eigenvalue to about $10^{-8}$, by the Collatz--Wielandt bounds for the collocation matrix, and the eigenvalue increases by $1.877$ per unit of $\kappa$. The values of $s_{\rm eig}=2-\kappa_{\rm eig}$ are:
\begin{center}
\begin{tabular}{lcccc}
\toprule
mesh & Remark~\ref{rem:coarse} & $F(1)$ & $F(1/\sqrt2)$ & $F(\frac12)$\\
\midrule
$s_{\rm eig}$ & $1.6853907$ & $1.6853613$ & $1.6853514$ & $1.6853460$\\
\bottomrule
\end{tabular}
\end{center}
Over the family, $\kappa_{\rm eig}$ converges with order about $1.7$ in $\sigma$. Fitting $\kappa_{\rm eig}=\kappa_0+C\sigma^p$ on the three meshes $F(\sigma)$ gives $2-\kappa_0=1.6853394$ with $p=1.73$, and $1.6853410$ with $p=2$. In this collocation, as in the computation of the test functions, the blocks of runs are evaluated at their middle run. With a finer rule (the runs $k\le1024$ one by one, then blocks of ratio $1.003$), $s_{\rm eig}$ increases by $8.8\times10^{-7}$, and a coarser rule (runs $k\le256$, ratio $1.01$) already gives all but $10^{-8}$ of this; with this correction the two extrapolations become $1.6853403$ and $1.6853419$. A collocation with continuous functions that are quadratic on each triangle $Q^\pm_{mn}$ and in $\vartheta$ on each arc converges much faster: on $F(2\sqrt2)$, $F(2)$ and $F(\sqrt2)$, with $0.9$, $2.4$ and $6.0$ million nodes, it gives $s_{\rm eig}=1.6853416$ within $3\times10^{-8}$, with the finer quadrature. So numerically $\sA(\GN)=\dimH\cN\approx1.6853416$, at distance about $4\times10^{-5}$ from $1.6853$ and $6\times10^{-5}$ from $1.6854$. This agrees with three other estimates:
\begin{itemize}
\item the Richardson extrapolation of piecewise-constant computations on three meshes, which gives $1.6855$;
\item the critical values of $\kappa$ of the cut-out series in Remark~\ref{rem:cutout};
\item a box count of $\cN$ at scales $2^{-k}$, $k\le11$, which gives $1.686$.
\end{itemize}

The ends of the interval of Theorem~\ref{thm: dimH} were fixed in advance, so the mesh had to be fine enough for $(\tilde\cK_\kappa h)/h$ to stay on the right side of $1$ between the nodes as well. On $F(\frac12)$, at $\kappa_-$ and $\kappa_+$ the nodal eigenvalue is about $1+8.6\times10^{-5}$ and $1-1.0\times10^{-4}$. Between the nodes the ratio departs from its nodal value, with the extremes in the bands of Remark~\ref{rem:fan}: for the eigenfunctions of the collocation on the mesh of Remark~\ref{rem:coarse}, it dips by up to about $2\times10^{-3}$ below its nodal value and rises by up to about $5\times10^{-4}$ above it. After the correction, the sampled minimum of the ratio for $h^-$ lies $4.1\times10^{-5}$ below the nodal eigenvalue, and the sampled maximum for $h^+$ lies $1.7\times10^{-5}$ above it, against $4.6\times10^{-4}$ and $1.4\times10^{-4}$ for the corrected functions on the mesh of Remark~\ref{rem:coarse}; along the family $F(\sigma)$ these two losses decrease roughly like $\sigma^{1.8}$ and $\sigma^{1.4}$. What is left are the margins $4.4\times10^{-5}$ and $8.4\times10^{-5}$ of the survey above, about eighteen and thirty-five times the width of the enclosures at the nodes.

By the estimates above, $\sA(\GN)$ lies about $1.6\times10^{-6}$ above $1.68534$. Certifying the fifth decimal place, that is, $\sA(\GN)>1.68534$, would require a lower certificate at $\kappa=2-1.68534$ with a margin of about $3\times10^{-6}$ in the eigenvalue, fifteen times smaller than the margin of $h^-$; for piecewise-linear test functions this means meshes with about a hundred times as many cells as $F(\frac12)$. Piecewise-quadratic functions do not change this much: for the eigenfunctions of the collocations, before any correction, the ratio $(\tilde\cK_\kappa h)/h$ between the nodes deviates from the nodal eigenvalue by $4.5\times10^{-4}$ on $F(2)$ and $2.6\times10^{-4}$ on $F(\sqrt2)$, against $8.2\times10^{-4}$ for piecewise-linear functions on $F(1)$, with a similar number of nodes; and these deviations decrease at about the same rate, roughly $\sigma^{1.7}$. So the behaviour of the eigenfunction in the bands of Remark~\ref{rem:fan}, rather than the degree of the test functions, limits the certificate.
\end{remark}

\section{Cut-out gaskets and the norm exponent}\label{sec:general}

In this section we prove Theorem~\ref{thm:cutout}, which bounds the upper box dimension of a cut-out gasket from both sides through the norm exponent; the inequality $\sA\ge\min\{n,\frac n{n+1}\sN\}$ between the affinity and the norm exponents of a finitely generated semigroup of $SL_{n+1}(\R)$ (Proposition~\ref{prop:main}); and Theorem~\ref{thm: sN}.

Throughout this section $K=\Delta\setminus\bigcup_IZ_I$ is a cut-out gasket as in the introduction, with zones indexed by the $I\in\cI^m$, where $m$ is the number of generators of $\Gamma$; $\nu_I$, $N_\Gamma$ and $\zeta_\Gamma$ are as there, and we write $\sN$ for $\sN(\Gamma)$. For $\varepsilon>0$, $K_\varepsilon$ is the set of points of the $n$-plane at distance less than $\varepsilon$ from $K$, and for a bounded set $F$ in the $n$-plane \cite[Prop.~3.2]{FalBook}
\begin{equation}\label{eq:minkowski}
  \udimB F=n-\liminf_{\varepsilon\to0}\frac{\log\mathrm{vol}(F_\varepsilon)}{\log\varepsilon}.
\end{equation}
For a convex body $Z$ and $t\ge0$ let $Z_{-t}=\{z\in Z:\ \mathrm{dist}(z,\partial Z)\ge t\}$, and let $C_\varepsilon(Z)=Z\setminus Z_{-\varepsilon}$ be its inner collar of width $\varepsilon$.

\begin{lemma}\label{lem:collar}
Let $Z$ be a convex body in the $n$-plane, with volume $V$ and surface area $A$. Then
\[
  \min\Bigl(V,\ \tfrac1n\,\varepsilon A\Bigr)\ \le\ \mathrm{vol}(C_\varepsilon(Z))\ \le\ \min(V,\ \varepsilon A),
  \qquad A\ \ge\ n\,\omega_n^{1/n}\,V^{(n-1)/n},
\]
where $\omega_n$ is the volume of the unit ball of $\R^n$.
\end{lemma}

\begin{proof}
The last inequality is the isoperimetric inequality \cite[(7.20)]{Sch14}. The sets $Z_{-t}$ are convex, and by the coarea formula \cite[3.2.11]{Fed69} for the distance to $\partial Z$, $\mathrm{vol}(C_\varepsilon(Z))=\int_0^\varepsilon\mathrm{area}(\partial Z_{-t})\,dt$; since the surface area of convex bodies is monotone under inclusion (by Cauchy's surface area formula~\cite[(5.73)]{Sch14}), this is at most $\varepsilon A$, which gives the upper bound. For the lower bound, if $a\in Z_{-s}$ and $b\in Z_{-t}$, the balls of radii $s$ and $t$ about $a$ and $b$ lie in $Z$, and so does their convex combination; hence $(1-\lambda)Z_{-s}+\lambda Z_{-t}\subset Z_{-((1-\lambda)s+\lambda t)}$, and by the Brunn--Minkowski inequality \cite[Thm.~7.1.1]{Sch14} the function $f(t)=\mathrm{vol}(Z_{-t})^{1/n}$ is concave where it is positive. Its right derivative at $0$ is $-A/(nV^{(n-1)/n})$, because $\frac d{dt}\mathrm{vol}(Z_{-t})=-\mathrm{area}(\partial Z_{-t})$ by the coarea formula, and $\mathrm{area}(\partial Z_{-t})\to A$ as $t\downarrow0$: the convex bodies $Z_{-t}$ converge to $Z$ in the Hausdorff metric, and the surface area of a convex body, the total mass of its surface area measure, depends continuously on the body in the Hausdorff metric, because the surface area measure depends weakly continuously on it~\cite[\S4.2]{Sch14}. By concavity, $f(\varepsilon)\le V^{1/n}(1-\varepsilon A/(nV))$ whenever the right side is positive, so
\[
  \mathrm{vol}(C_\varepsilon(Z))=V-f(\varepsilon)^n\ \ge\ V\Bigl(1-\bigl(1-y\bigr)_+^n\Bigr)\ \ge\ V\min(1,y),\qquad y=\frac{\varepsilon A}{nV},
\]
since $(1-y)^n\le1-y$ for $0\le y\le1$, while for $y>1$ the right side of the concavity bound is negative, so $f(\varepsilon)=0$ and $\mathrm{vol}(C_\varepsilon(Z))=V$.
\end{proof}

\thmcutout*

\begin{proof}
(a) Fix $t$ with $\sN<t<n+1$, so that $\zeta_\Gamma(t)<\infty$, and let $0<\varepsilon<1$. The set $K_\varepsilon\setminus\Delta$ lies in the $\varepsilon$-neighbourhood of $\Delta$ minus $\Delta$, whose volume is at most $C_\Delta\,\varepsilon$ by Steiner's formula~\cite[(4.8)]{Sch14}. Let $\xi\in(K_\varepsilon\cap\Delta)\setminus K$; then $\xi\in Z_I$ for some $I$. Let $q\in K$ with $\abs{\xi-q}<\varepsilon$; the segment from $\xi$ to $q$ is connected and joins a point of $Z_I$ to a point outside $Z_I$, so it meets the boundary of $Z_I$, which is $\partial\overline{Z_I}$ because $Z_I$ is convex with non-empty interior. Hence $\mathrm{dist}(\xi,\partial\overline{Z_I})<\varepsilon$, and $\xi\in C_\varepsilon(\overline{Z_I})$. Since $K$ has measure zero, Lemma~\ref{lem:collar} gives
\[
  \mathrm{vol}(K_\varepsilon)\le C_\Delta\,\varepsilon+\sum_I\min\bigl(C\nu_I^{-(n+1)},\ C\varepsilon\,\nu_I^{-(n-1)}\bigr).
\]
Split the sum at $X=\varepsilon^{-1/2}$:
\[
  \varepsilon\sum_{\nu_I\le X}\nu_I^{-(n-1)}\le\varepsilon X^{t-(n-1)}\zeta_\Gamma(t)=\varepsilon^{\frac{n+1-t}2}\zeta_\Gamma(t),\qquad
  \sum_{\nu_I>X}\nu_I^{-(n+1)}\le X^{t-(n+1)}\zeta_\Gamma(t)=\varepsilon^{\frac{n+1-t}2}\zeta_\Gamma(t).
\] Since $(n+1-t)/2<1$, $\mathrm{vol}(K_\varepsilon)\le C'\varepsilon^{(n+1-t)/2}$, and by \eqref{eq:minkowski} $\udimB K\le n-\frac{n+1-t}2=\frac{n-1+t}2$. Let $t\downarrow\sN$.

(b) Since $\partial\overline{Z_I}\subset K$, the collars $C_\varepsilon(\overline{Z_I})$ lie in $K_\varepsilon$. They overlap only in a set of measure zero, because two distinct closures $\overline{Z_I}$ and $\overline{Z_J}$ meet only in the boundary of $\overline{Z_I}$, which has measure zero: an interior point of $\overline{Z_I}$ is an interior point of $Z_I$, since $Z_I$ is convex with non-empty interior, so it has a neighbourhood contained in $Z_I$, which does not meet $Z_J$. By Lemma~\ref{lem:collar} and the hypothesis on the volumes, and since $\min(V,\varepsilon c'V^{(n-1)/n})$ is increasing in $V$,
\[
  \mathrm{vol}(K_\varepsilon)\ \ge\ \sum_I\min\bigl(c\,\nu_I^{-(n+1)},\ c''\varepsilon\,\nu_I^{-(n+1)(n-1)/n}\bigr).
\]
Let $X=\varepsilon^{-n/(n+1)}$. For $X<\nu_I\le2X$ the two terms of the minimum differ by a bounded factor, so each such $I$ contributes at least $c'''X^{-(n+1)}$, and
\begin{equation}\label{eq:shell}
  \mathrm{vol}(K_\varepsilon)\ \ge\ c'''\,\bigl(N_\Gamma(2X)-N_\Gamma(X)\bigr)\,X^{-(n+1)} .
\end{equation}
Let $s<s'<\sN$. Then $N_\Gamma(r_m)\ge r_m^{s'}$ along a sequence $r_m\to\infty$. If every dyadic shell $(X,2X]\subset[r_m^{1/2},r_m]$ had $N_\Gamma(2X)-N_\Gamma(X)<X^{s}$, then $N_\Gamma(r_m)\le N_\Gamma(r_m^{1/2})+Cr_m^{s}$, which is impossible for large $m$: since $\sum_I\nu_I^{-t}<\infty$, $N_\Gamma(r)\le C_tr^{t}$ for every $t>\sN$, and choosing $t\in(\sN,2s')$, which is possible if $s'>\sN/2$, gives $N_\Gamma(r_m^{1/2})=o(r_m^{s'})$. So there are $X_m\to\infty$ with $N_\Gamma(2X_m)-N_\Gamma(X_m)\ge X_m^{s}$, and by \eqref{eq:shell}, with $\varepsilon_m=X_m^{-(n+1)/n}$,
\[
  \mathrm{vol}(K_{\varepsilon_m})\ \ge\ c'''X_m^{s-(n+1)}=c'''\,\varepsilon_m^{\,n-\frac{n}{n+1}s}.
\]
By \eqref{eq:minkowski}, $\udimB K\ge\frac n{n+1}s$. Let $s\uparrow\sN$.
\end{proof}

\begin{remark}
Part (b) uses neither fastness nor the measure of $K$, and only a lower bound on the volumes of the zones. Its hypothesis $\partial\overline{Z_I}\subset K$ says that the zones are separated by the gasket; for the Levitt--Yoccoz gasket it holds (see the proof of Theorem~\ref{thm: sN}), and for $\cN$ the bound $\udimB\cN\ge\frac23\sN$ was proved by Sewell~\cite[Thm.~1.8]{Sew23} for self-projective sets with simplicial holes. Part (a) bounds the surface area of every zone by the order $\nu_I^{-(n-1)}$ of the zones of cusp type, so we expect it to be sharp only when the zones of cusp type dominate. For $\cN$ the perimeters of the zones range between the orders $\nu_I^{-3/2}$ and $\nu_I^{-1}$, and both are attained (Lemma~\ref{lem:geom}(iv) and Remark~\ref{rem:geomsharp}); Remark~\ref{rem:172} discusses which of them prevail.
\end{remark}

\begin{remark}[the cut-out argument and the value $1.72$]\label{rem:172}
Theorem~\ref{thm:cutout}(a) treats every zone as a cusp: for $\cN$ it uses only the bound $\mathrm{per}(Z_I)=O(\|A_I\|_1^{-1})$ of Lemma~\ref{lem:geom}(iv), the size of the cusp triangles. The perimeters range between the orders $\|A_I\|_1^{-3/2}$ and $\|A_I\|_1^{-1}$ (Lemma~\ref{lem:geom}(iv) and Remark~\ref{rem:geomsharp}). The same cut-out argument with the actual perimeters gives $\udimB\cN\le2-\kappa$ whenever $\sum_I\mathrm{area}(Z_I)^{1-\kappa}\mathrm{per}(Z_I)^\kappa$ converges, and the critical value of $\kappa$ is exactly $2-\sA(\GN)$ (Corollary~\ref{cor:cutoutA} and Proposition~\ref{prop:cutout}); numerically, half-octave shells of this series give $2-\kappa\approx1.686$ (Remark~\ref{rem:cutout}), and computing $\sA$ along the runs gives $1.6853<\dimH\cN<1.6854$ (Theorem~\ref{thm: dimH}).
The value $\frac12(1+\sN)=1.7219$ is what the cut-out argument gives when every zone is treated as a cusp, with perimeter of order $\nu_I^{-1}$. If, at the scales reachable numerically, the volume of the $\varepsilon$-neighbourhoods is still dominated by zones that look cusp-like, then estimates of this volume at moderate $\varepsilon$ read the cusp exponent $\frac12(1+\sN)$ rather than $\dimH\cN$. This is a plausible explanation of the value $\approx1.72$ suggested by the numerics of~\cite{DD09} and quoted in~\cite{PS23}, which is excluded by Theorem~\ref{thm: dimH}.
The same discrepancy is already visible in~\cite{DL08}, where the estimate from the area of $\varepsilon$-neighbourhoods ($1.7$--$1.8$) lies above the two box counts of the same data ($1.69$ for the barycentres of the zones and $1.7$ for the set).
\end{remark}

\subsection*{Affinity exponent and norm exponent.}
Let $\Gamma\subset SL_{n+1}(\R)$ be a semigroup generated by finitely many matrices, and let $\sigma_1(A)\ge\dots\ge\sigma_{n+1}(A)$ be the singular values of $A\in SL_{n+1}(\R)$. For $0\le s\le n$ and $k=\lfloor s\rfloor$, the singular value function of the projective action of $A$ on $\RP^n$ is
\[
  \varphi^s(A)=\prod_{i=1}^k\frac{\sigma_{i+1}(A)}{\sigma_1(A)}\cdot
  \Bigl(\frac{\sigma_{k+2}(A)}{\sigma_1(A)}\Bigr)^{s-k}
\]
(for $s=n$ the last factor is $1$), that is, Falconer's singular value function~\cite{Fal88} of a linear map of $\R^n$ with singular values $\sigma_2(A)/\sigma_1(A)\ge\dots\ge\sigma_{n+1}(A)/\sigma_1(A)$. For $n=2$ it is the function $\varphi^s$ of Section~\ref{sec:affinity}, and the affinity exponent $\sA=\sA(\Gamma)$ is defined as there, with $[0,2]$ replaced by $[0,n]$. As in Remark~\ref{rem:conventions}, it does not depend on the choice of the inner product.

\medskip
The following is an immediate consequence of Falconer's inequality $\varphi^s(g)\ge\abs{\det g}^{s/n}$ for linear maps $g$ of $\R^n$ and $0\le s\le n$ \cite{Fal88}, applied to the projective singular values; a box-dimension analogue, under geometric assumptions on the limit set, is~\cite[Thm.~1.8]{Sew23}.

\begin{proposition}\label{prop:main}
For every finitely generated semigroup $\Gamma\subset SL_{n+1}(\R)$,
\[
  \sA\ \ge\ \min\Bigl\{n,\ \frac{n}{n+1}\,\sN\Bigr\}.
\]
\end{proposition}

\begin{proof}
Let $g\in SL_{n+1}(\R)$ and $r_i=\sigma_{i+1}(g)/\sigma_1(g)$, $i=1,\dots,n$. The
$r_i$ are non-increasing and $\prod_ir_i=\abs{\det g}/\sigma_1^{n+1}=\sigma_1^{-(n+1)}$.
Now $\log\varphi^s(g)=\sum_{i\le k}\log r_i+(s-k)\log r_{k+1}$ is the integral over
$[0,s]$ of the non-increasing step function $x\mapsto\log r_{\lceil x\rceil}$, which
is at least $s/n$ times its integral over $[0,n]$. Hence
\[
  \varphi^s(g)\ \ge\ \Bigl(\prod_{i=1}^nr_i\Bigr)^{s/n}=\sigma_1(g)^{-\frac{n+1}{n}s}
  \qquad(0\le s\le n).
\]
If $\sA<n$, take $s\in(\sA,n]$. Since $\varphi^s$ is non-increasing in $s$,
$\sum_I\varphi^s(A_I)<\infty$, hence $\sum_I\sigma_1(A_I)^{-\frac{n+1}{n}s}<\infty$
and $\frac{n+1}{n}s\ge\sN$. Letting $s\downarrow\sA$ gives $\sA\ge\frac{n}{n+1}\sN$.
\end{proof}

For $n=2$ the inequality reads $\varphi^s(g)\ge\sigma_1(g)^{-3s/2}$; directly,
$\varphi^s=\sigma_2^{2-s}\sigma_1^{1-2s}$ for $1\le s\le2$ and $\sigma_2\ge\sigma_1^{-1/2}$
because $\sigma_2^2\ge\sigma_2\sigma_3=1/\sigma_1$. For $\GN$ this is Lemma~\ref{lem:sv} and Remark~\ref{rem:ceiling}, where the bound $\sigma_2(A_I)\le3$ gives also $\sA(\GN)\le\frac12(1+\sN)$. For $n=1$ one has
$\varphi^s=\sigma_1^{-2s}$ exactly, so $\sA=\frac12\sN$ (when $\sN\le2$).

\begin{remark}\label{rem:minneeded}
The minimum is needed: since $\sA\le n$ by definition, the bound
$\sA\ge\frac{n}{n+1}\sN$ fails whenever $\sN>n+1$, and this happens. For instance,
let $n=2$ and take five distinct contracting similarities $f_i(x)=\frac12x+t_i$ of
$\R^2$, written as elements of $SL_3(\R)$ by
\[
  A_i=4^{1/3}\begin{pmatrix}\frac12\id&t_i\\0&1\end{pmatrix}.
\]
Every $I\in\cI^5$ of length $k$ gives
$A_I=4^{k/3}\begin{pmatrix}2^{-k}\id&t_I\\0&1\end{pmatrix}$ with $t_I$ bounded, so
$\norm{A_I}\asymp4^{k/3}$, while the singular values of $A_I$ are
$\asymp4^{k/3}(1,2^{-k},2^{-k})$. Summing over the $5^k$ multiindices of length $k$,
$\zeta_\Gamma(t)<\infty$ iff $t>3\log_45\approx3.48$, so $\sN=3\log_45>3$
and $\frac23\sN=\log_25\approx2.32>2$; on the other hand
$\varphi^s(A_I)\asymp2^{-ks}$, and $\sum_I\varphi^s(A_I)$ diverges for all
$s\le2$, so $\sA=2$.
\end{remark}

\begin{corollary}\label{cor:main}
Let $\Gamma\subset SL_{n+1}(\R)$ be a finitely generated semigroup preserving a properly convex cone (this assumption only serves to define its limit set $\Lambda_\Gamma$). If $\dimH\Lambda_\Gamma=\sA<n$, then
$(n+1)\dimH\Lambda_\Gamma\ge n\,\sN$.
\end{corollary}

This proves Conjecture~1.4 of \cite{DL15EM} (stated there for $SL_n$ as $n\dimH R_A\ge(n-1)s_A$, where $s_A$ denotes the norm exponent of the semigroup $A$, our $\sN$) for every gasket whose Hausdorff dimension equals its affinity exponent and is less than $n$; for the Rauzy gasket, by Theorem~\ref{thm:jurga} and Lemma~\ref{lem:baragartree}, it is the inequality $\delta\ge\frac23\beta(4)$ recorded in~\cite{GMR19}, right after Remark~11 there (see the introduction). For the standard Sierpi\'nski gasket the inequality is an equality; for the affine Sierpi\'nski gaskets a box-dimension version was proved in \cite[Prop.~3.1]{DL15EM} for a range of parameters, and for a class of self-projective sets in \cite[Thm.~1.8]{Sew23}.

\subsection*{Elementary bounds for $\sN$.}
For $I\in\cI^2$, let $L_I=A_I\one=(a,b,c)$.
Prepending $1$, resp.\ $2$, maps it to $(a,a+b,a+c)$, resp.\ $(a+b,b,b+c)$, because $A_jL=L+L_j(\one-e_j)$; hence, by induction from $L_\emptyset=\one$, $c=a+b-1$ and $\|A_I\|_1=\norm{L_I}_1=2(a+b)-1$, and $(a,b)$ runs through the Euclidean algorithm backwards from $(1,1)$.
Thus, $I\mapsto(a,b)$ is a bijection from $\cI^2$ onto the coprime pairs of positive integers: the two steps $(a,b)\mapsto(a,a+b)$ and $(a,b)\mapsto(a+b,b)$ preserve coprimality, and a coprime pair $(a,b)\ne(1,1)$ has $a\ne b$ and is obtained from exactly one pair, $(a,b-a)$ if $a<b$ and $(a-b,b)$ if $a>b$, so that induction on $a+b$ gives the bijection; and there are $\varphi_{\mathrm E}(n)$ such pairs with $a+b=n$, where $\varphi_{\mathrm E}$ is Euler's totient function.
Therefore $\zeta(t)\ge\sum_{n\ge2}\varphi_{\mathrm E}(n)(2n-1)^{-t}$, which diverges for
$t\le2$ since $\sum_{n\le X}\varphi_{\mathrm E}(n)\asymp X^2$ \cite[Thm.~330]{HW08}.
This shows that $\sN\geq2$. This argument is the embedding of the case $d=1$ in the case $d=2$ of the monotonicity of Remark~\ref{rem:monotone}: the pairs $(a,b)$ are the souls of the semigroup $\Gamma^{(1)}$ of Section~\ref{sec:MH}, and the bound $\sN\ge2$ corresponds to the exponent $2$ of the Markoff triples~\cite[Thm.~0.1]{Bar98} (see also~\cite[Remark~14]{GMR19}).

Since $\|A_I\|_1\ge3>1$, each term $\|A_I\|_1^{-t}$ is decreasing in $t$, so the set of $t$ with $\zeta(t)<\infty$ is a half-line whose left endpoint is $\sN$.
This half-line is not empty: $I\mapsto L_I$ is injective, and $L_I$ is a vector of positive integers with $\|L_I\|_1=\|A_I\|_1$ (Corollary~\ref{cor:accumulation} and the proof of Lemma~\ref{lem:geom}(v)), so at most $\binom{n-1}{2}<n^2/2$ multiindices have $\|A_I\|_1=n$, and
$\zeta(t)\le\frac12\sum_{n\ge3}n^{2-t}<\infty$ for $t>3$.
Hence, $2\leq \sN\le3$.

\thmsN*

\begin{proof}
The zones $Z_I$ of $\cN$ are open triangles with $\frac{\sqrt3}4\,\nu_I^{-3}\le\mathrm{area}(Z_I)\le2\sqrt3\,\nu_I^{-3}$ and $\mathrm{per}(Z_I)\le18\sqrt3\,\nu_I^{-1}$ (Lemma~\ref{lem:geom}(iii),(iv)), and $\mathrm{area}(\cN)=0$ (Theorem~\ref{thm:yoccoz}). Their sides lie in $\cN$: a point $p$ of $\partial Z_I$ lies in $\cT$, and it lies in no zone, because the zones are open, so a zone containing $p$ would meet $Z_I$; by Lemma~\ref{lem:zones}, $p\in\cN$. By Lemmas~\ref{lem:geom}(iii) and~\ref{lem:zones}, $\zeta(3)\le\frac4{\sqrt3}\sum_I\mathrm{area}(Z_I)\le\frac4{\sqrt3}\,\mathrm{area}(\cT)<\infty$, while $\zeta(\sN)=\infty$ by Proposition~\ref{prop:runpressure} and Remark~\ref{rem:GNfast} below; since $\sN\ge2$ (see above), we get $2\le\sN<3$, and Theorem~\ref{thm:cutout} with $n=2$, applied to the triangle $\cT$ in place of $\Delta$, gives $\frac23\sN\le\udimB\cN\le\frac12(1+\sN)$. By Corollary~\ref{cor: dimB}, $\dimB\cN$ exists and equals $\dimH\cN$. (Alternatively, by Theorem~\ref{thm:jurga}, $\dimH\cN=\sA(\GN)$, and $\frac23\sN\le\sA(\GN)\le\frac12(1+\sN)$ by Proposition~\ref{prop:main} with $n=2$ and by Remark~\ref{rem:ceiling}.) With Theorem~\ref{thm: delta}, $\frac23\cdot2.4411=1.6274$ and $\frac12(1+2.4465)=1.72325$.
\end{proof}

\section{Fast semigroups and the run factors}\label{sec:fast}

In this section we rewrite the zeta function of $\GN$ as the sum of the iterates of a positive transfer operator $\cL_t$ on the first-level cylinders, by splitting the multiindices into runs (Corollary~\ref{cor:test}), and explain why this works: for a fast semigroup the weight of the multiindices with $r$ runs grows or decays exponentially in $r$, and the norm exponent is the value of $t$ at which it does neither, the zero of the growth rate per run (Proposition~\ref{prop:runpressure}); and $\GN$ is fast, with sharp junction constant $\frac14$ (Lemma~\ref{lem:junction}). 
We also prove the uniform contraction of the runs on the cylinders (Lemma~\ref{lem:contr}), used in Sections~\ref{sec:thmB} and~\ref{sec:AHS}.

\subsection*{The zeta function seen from a vector.}
For a non-negative vector $x\neq0$ set
\[
  \zeta_x(t)=\sum_I\norm{A_Ix}_1^{-t},
\]
where the empty multiindex contributes $\norm{x}_1^{-t}$. Clearly $\zeta_{\one}=\zeta$, and $\zeta_{\lambda x}=\lambda^{-t}\zeta_x$ for $\lambda>0$, so $\zeta_x=\norm{x}_1^{-t}\zeta_{\pr x}$.
Since $\norm{\cdot}_1$ is additive on non-negative vectors, $\norm{A_IB}_1=\one^*A_I(B\one)=\norm{A_I(B\one)}_1$ for every non-negative matrix $B$. For $B=A_J$ we have $A_J\one=L_J=\norm{A_J}_1\,\pr{L_J}$, hence
\begin{equation}\label{eq:seen}
  \norm{A_{IJ}}_1=\norm{A_J}_1\,\norm{A_I\pr{L_J}}_1,\qquad
  \sum_I\norm{A_{IJ}}_1^{-t}=\norm{A_J}_1^{-t}\,\zeta_{\pr{L_J}}(t):
\end{equation}
the growth from $J$ to $IJ$ depends on $J$ only through its soul point. This is the exact form of the inequalities~\cite[(3)]{DL15JGA}, in which $\zeta_{\pr{L_J}}(t)$ is replaced by bounds proportional to $\zeta(t)$.
The convergence of such series will be decided with positive test functions, through Lemma~\ref{lem:test}.

\medskip
Splitting off the first index to act on $x$ gives, on $\cT$, the renewal $\zeta_x(t)=1+\sum_j(1+2x_j)^{-t}\zeta_{\psi_j(x)}(t)$, with the operator $(\mathcal M_tf)(x)=\sum_j(1+2x_j)^{-t}f(\psi_j(x))$, for which $(\mathcal M_t^m\mathds 1)(x)=\sum_{\abs I=m}\norm{A_Ix}_1^{-t}$. This operator gives no upper bounds for $\sN$: each $\psi_j$ is parabolic, with fixed point $\pr{\one-e_j}$, where the weight equals $1$, and $(\mathcal M_t^m\mathds 1)(x_0)\ge\norm{A_j^mx_0}_1^{-t}=(1+2m/3)^{-t}$, with $x_0=\pr{\one}$, decays only polynomially, so Lemma~\ref{lem:test}(ii) never applies to $\mathcal M_t$. We therefore take whole runs as the steps of the renewal.

\subsection*{The run factors.}
Recall from Section~\ref{sec:cocycle} that $\norm{A_j^kx}_1=1+2kx_j$ for $x\in\cT$, by~\eqref{eq:run}, and that, by~\eqref{eq:cocycle}, the $\ell^1$-size of a soul is exactly the product of its run factors, and the states $x^{(r)}$ are the soul points of the suffixes of $I$ (Example~\ref{ex:run} and Figure~\ref{fig:runs}). The first-level cylinders $C_i$ trap the states (Lemma~\ref{lem:trap}), and the derivative of a run is given by~\eqref{eq:DR} (Lemma~\ref{lem:contraction}); the latter is used in Lemma~\ref{lem:contr}.

Hence, after the first run, the state $(x^{(r)},j_r)$ lies in
$\mathcal C=\bigsqcup_iC_i\times\{i\}$; the next run uses an index $j\neq j_r$
with $2x_j\ge\frac12$, and the weights $(1+2kx_j)^{-t}\le(1+k/2)^{-t}$
are summable in $k$ for $t>1$. For $t>1$ define on
$\mathcal B=\mathcal B(\mathcal C)$ the positive operator
\[
  (\cL_th)(x,i)=\sum_{j\neq i}\sum_{k\ge1}\bigl(1+2kx_j\bigr)^{-t}\,
  h\bigl(\psi_j^k(x),j\bigr),\qquad h\in\mathcal B;
\]
it maps $\mathcal B$ to itself because
$\sum_{j\neq i}\sum_{k\ge1}\bigl(1+2kx_j\bigr)^{-t}\le2\sum_{k\ge1}(1+k/2)^{-t}<\infty$.
For $i\in\{1,2,3\}$ let $\zeta_{x,i}(t)$ be the sum defining $\zeta_x(t)$ restricted
to the multiindices $I$ whose last index is not $i$. Every non-empty such $I$ ends with a
run $j^k$, $j\neq i$, preceded by a multiindex whose last index is not $j$; hence
\begin{equation}\label{eq:renewal}
  \zeta_{x,i}(t)=\norm{x}_1^{-t}+\sum_{j\neq i}\sum_{k\ge1}\zeta_{A_j^kx,\,j}(t),
\end{equation}
where again $\norm{x}_1^{-t}$ is the contribution of the empty multiindex; similarly
\[
\zeta_x(t)=\norm{x}_1^{-t}+\sum_{j}\sum_{k\ge1}\zeta_{A_j^kx,\,j}(t).
\]
This is the splitting into runs of \cite{DL15JGA}, and it involves no
normalization.
On $\mathcal C$, homogeneity and \eqref{eq:run} give, by induction on $m$,
\begin{equation}\label{eq:Lm}
  (\cL_t^m\mathds 1)(x,i)=\sum_{I}\norm{A_Ix}_1^{-t},
\end{equation}
the sum running over the multiindices $I$ with exactly $m$ runs whose last index is not $i$. Summing over $m$,
\[
\zeta_{x,i}(t)=\sum_{m\ge0}(\cL_t^m\mathds 1)(x,i)\qquad\bigl((x,i)\in\mathcal C\bigr),
\]
an identity in $[0,\infty]$; equivalently, $(x,i)\mapsto\zeta_{x,i}(t)$ is a fixed point of $h\mapsto\mathds 1+\cL_th$ whenever it is finite.
Applying the unrestricted renewal at $x=\one$, where $\norm{A_j^k\one}_1=3+2k$, and homogeneity,
\begin{equation}\label{eq:zeta}
  \zeta(t)=3^{-t}\Bigl(1+\sum_{j=1}^3\sum_{k\ge1}\bigl(1+\tfrac{2k}3\bigr)^{-t}
  \sum_{m\ge0}(\cL_t^m\mathds 1)\bigl(\psi_j^k(x_0),j\bigr)\Bigr).
\end{equation}

\begin{corollary}\label{cor:test}
Let $t>1$ and $h\in\mathcal B$ with $0<\inf h\le\sup h<\infty$. If $\cL_th\ge h$,
then $t\le\sN$; if $\cL_th\le\theta h$ for some $\theta<1$, then $t\ge\sN$.
\end{corollary}

\begin{proof}
If $\cL_th\ge h$, then by Lemma~\ref{lem:test}(i) every term
$(\cL_t^m\mathds 1)(\psi_j^k(x_0),j)$ in \eqref{eq:zeta} is at least $\inf h/\sup h$, so
\eqref{eq:zeta} diverges and $t\le\sN$. If $\cL_th\le\theta h$, then by
Lemma~\ref{lem:test}(ii) the inner sum in \eqref{eq:zeta} is at most
$\sup h/((1-\theta)\inf h)$, and $\sum_k(1+\frac{2k}3)^{-t}<\infty$ for $t>1$, so
\eqref{eq:zeta} converges and $t\ge\sN$.
\end{proof}

\subsection*{Fast semigroups and the growth rate per run.}
The rewriting of $\zeta(t)$ by runs works for a structural reason, which we now
isolate. In this subsection $\Gamma$ is the semigroup generated by invertible matrices $A_1,\dots,A_m$,
$m\ge2$, sums are over $\cI^m$, $\norm{\cdot}$ is a submultiplicative norm with
$\norm{A_I}\ge1$ for every $I\in\cI^m$ (for instance the entrywise $\ell^1$ norm on a semigroup of
invertible non-negative integer matrices), and only finitely many $I\in\cI^m$ have $\norm{A_I}$ below any given bound; $\zeta_\Gamma$ and $\sN=\sN(\Gamma)$ are as in the introduction. Following \cite[Def.~3]{DL15JGA}, $\Gamma$ is \emph{fast} if there is $c>0$
such that
\[
  \norm{A_{IJ}}\ \ge\ c\,\norm{A_I}\,\norm{A_J}
\]
for every multiindex $I$ and every multiindex $J$ whose first two indices are distinct.
As in \eqref{eq:cocycle}, every non-empty multiindex is uniquely $I=j_r^{k_r}\cdots j_1^{k_1}$ with $k_s\ge1$ and
$j_{s+1}\neq j_s$; we call $r=r(I)$ its number of runs. Set
\[
  \zeta^{(r)}(t)=\sum_{r(I)=r}\norm{A_I}^{-t}\in(0,\infty],\qquad
  \zeta_\Gamma(t)=\sum_{r\ge0}\zeta^{(r)}(t),\qquad
  t_0=\inf\{t\ge0:\ \zeta^{(1)}(t)<\infty\}.
\]
We call $\zeta^{(r)}(t)$ the \textbf{weight of the $r$-th generation} of runs: it is the total weight $\norm{A_I}^{-t}$ of the multiindices made of exactly $r$ runs, and $\zeta_\Gamma(t)$ is the sum over all generations. Since $\norm{A_I}\ge1$, each $\zeta^{(r)}$ is non-increasing in $t$; since
$\zeta^{(1)}\le\zeta_\Gamma$, we have $t_0\le\sN$. The number $t_0$ is the exponent of the first generation, that is, of the cusp elements $A_j^k$ alone: $\zeta^{(1)}(t)=\sum_j\sum_{k\ge1}\norm{A_j^k}^{-t}$.

Proposition~\ref{prop:runpressure} below shows that, for a fast semigroup and $t>t_0$, the weight of the $r$-th generation grows or decays exponentially in $r$, at a rate
\[
  p(t)=\lim_{r\to\infty}\frac1r\log\zeta^{(r)}(t)
\]
that we call the \textbf{growth rate per run} of $\zeta_\Gamma$ at $t$, and that $\sN$ is the zero of $p$: below $\sN$ the generations grow, above $\sN$ they decay, and at $t=\sN$ their weight stays bounded away from $0$. In the language of thermodynamic formalism, $p(t)$ is the pressure of the potential $-t\log(\text{run factor})$ on the shift of the runs; this is made precise in Lemma~\ref{lem:pressure}.

\begin{lemma}\label{lem:fastruns}
Let $\Gamma$ be fast and $c_*=c/\max_i\bigl(\norm{A_i}\,\norm{A_i^{-1}}\bigr)$.
Then
\[
  \norm{A_{IJ}}\ \ge\ c_*\,\norm{A_I}\,\norm{A_J}
\]
for all non-empty multiindices $I,J$ such that the last index of $I$ differs from the
first index of $J$. Conversely, if this holds for some $c_*>0$, then $\Gamma$ is
fast.
\end{lemma}

\begin{proof}
Write $I=I'i$ and let $J$ begin with $j\neq i$. The multiindex $iJ$ begins with two
distinct indices, so $\norm{A_{IJ}}=\norm{A_{I'}A_{iJ}}\ge c\norm{A_{I'}}\norm{A_{iJ}}$;
moreover $\norm{A_I}\le\norm{A_{I'}}\norm{A_i}$ and
$\norm{A_J}\le\norm{A_i^{-1}}\norm{A_{iJ}}$. Conversely, if $J=jj'K$ with $j\ne j'$,
the hypothesis applied to $Ij$ and $j'K$ gives
$\norm{A_{IJ}}\ge c_*\norm{A_{Ij}}\norm{A_{j'K}}
\ge c_*\norm{A_I}\norm{A_J}/(\norm{A_j^{-1}}\norm{A_j})$.
\end{proof}

Thus fastness says exactly that concatenating runs costs at most a bounded factor;
only the merging of two runs of the same index is exempt. Note that
$c_*\le c\le1$, by submultiplicativity and $\norm{A_i}\norm{A_i^{-1}}\ge\norm{\id}\ge1$.

\begin{proposition}\label{prop:runpressure}
Let $\Gamma$ be fast and $t>t_0$.
\begin{enumerate}
\item[(i)] $\zeta^{(r)}(t)<\infty$ for every $r$; the limit
  $p(t)=\lim_{r\to\infty}\frac1r\log\zeta^{(r)}(t)$ exists and is finite; and
  $\zeta^{(r)}(t)\ge c_*^{\,t}\,e^{r\,p(t)}$ for every $r\ge1$.
\item[(ii)] $\zeta_\Gamma(t)<\infty$ if and only if $p(t)<0$.
\item[(iii)] $p$ is convex on $(t_0,\infty)$, and there is $b\ge1$ such that
  $p(t')\le p(t)-(t'-t)\,b^{-1}\log2$ for $t_0<t<t'$.
\end{enumerate}
If moreover $\sN>t_0$, then $p(\sN)=0$, $p>0$ on $(t_0,\sN)$, $p<0$ on
$(\sN,\infty)$, and $\zeta_\Gamma(\sN)=\infty$.
\end{proposition}

\begin{proof}
(i) If $I=R_1\cdots R_r$ is the decomposition of $I$ into runs, Lemma~\ref{lem:fastruns}
gives $\norm{A_I}\ge c_*^{\,r-1}\prod_s\norm{A_{R_s}}$, so
$\zeta^{(r)}(t)\le c_*^{-(r-1)t}\zeta^{(1)}(t)^r<\infty$. Cutting a multiindex with $r+r'$ runs
at the run boundary that leaves $r$ runs on the left writes it uniquely as $IJ$ with $r(I)=r$, $r(J)=r'$ and the
last index of $I$ different from the first index of $J$; hence
$\zeta^{(r+r')}(t)\le c_*^{-t}\zeta^{(r)}(t)\,\zeta^{(r')}(t)$, that is,
$a_r=c_*^{-t}\zeta^{(r)}(t)$ is submultiplicative. By Fekete's lemma~\cite{Fek23} (see also \cite[Part~I, Problem~98]{PS98}),
$p(t)=\lim_r\frac1r\log a_r=\inf_r\frac1r\log a_r$ exists in $[-\infty,\infty)$,
and $a_r\ge e^{r\,p(t)}$. Finally $p(t)\ge-t\log\mu$, with $\mu=\max_i\norm{A_i}$,
because the multiindex $1212\cdots$ of length $r$ has $r$ runs and norm at most $\mu^r$.

(ii) If $p(t)\ge0$, then $\zeta^{(r)}(t)\ge c_*^{\,t}$ for every $r$, and $\zeta_\Gamma(t)=\infty$.
If $p(t)<0$, choose $\varepsilon>0$ with $p(t)+\varepsilon<0$; then
$\zeta^{(r)}(t)\le e^{r(p(t)+\varepsilon)}$ for all large $r$, and $\zeta_\Gamma(t)<\infty$.

(iii) Each $t\mapsto\log\zeta^{(r)}(t)$ is convex, by H\"older's inequality, hence so is
$p$. Since only finitely many multiindices $J$ have $\norm{A_J}\le2/c_*$, there is $b\ge1$ such that
$\norm{A_J}>2/c_*$ for every multiindex $J$ of length at least $b$, in particular for
every $J$ with at least $b$ runs. Let $r(I)=r=qb+s$ with $0\le s<b$, and cut $I$
at run boundaries into $q$ multiindices with $b$ runs each, followed by a multiindex with $s$
runs. Lemma~\ref{lem:fastruns} at the (at most $q$) cuts gives
$\norm{A_I}\ge c_*^{\,q}(2/c_*)^q=2^q\ge2^{r/b-1}$, which also holds for $r<b$.
Hence $\zeta^{(r)}(t')\le2^{-(t'-t)(r/b-1)}\zeta^{(r)}(t)$ for $t<t'$, and the inequality
for $p$ follows.

If $\sN>t_0$, then $p$ is finite and convex, hence continuous, on $(t_0,\infty)$.
By (ii), $p<0$ on $(\sN,\infty)$ and $p\ge0$ on $(t_0,\sN)$; so $p(\sN)=0$, and
$p>0$ on $(t_0,\sN)$ because $p$ is strictly decreasing. Finally $\zeta_\Gamma(\sN)=\infty$
by (ii).
\end{proof}

For the grading by the length of the multiindices the situation is different. If $\Gamma$ is parabolic
in the sense of \cite[Def.~1]{DL15JGA}, i.e.\ $\min_{\abs I=k}\norm{A_I}$ grows
subexponentially in $k$, then
\[
  \Bigl(\sum_{\abs I=k}\norm{A_I}^{-t}\Bigr)^{1/k}\ge\Bigl(\min_{\abs I=k}\norm{A_I}\Bigr)^{-t/k}\to1,
\]
while for $t>\sN$ these sums tend to $0$; so their exponential growth rate equals
$1$ for every $t>\sN$ (compare \cite[Thm.~2]{DL15JGA}). This phase transition is
what makes Lemma~\ref{lem:test}(ii) inapplicable index by index, and grading by
runs removes it.
\subsection*{Fastness of $\GN$.}
The rewriting \eqref{eq:zeta} is the splitting into runs used in the proof of \cite[Thm.~1]{DL15JGA}. What is new is that the factor by which the norm grows when two runs are concatenated is computed, not bounded. The property of $\GN$ that makes this effective is the following sharp form of fastness.

\begin{lemma}\label{lem:junction}
Let $I\neq\emptyset$ have last index $j$, let $i\neq j$ and $x\in C_i$. Then
\[
  \tfrac14\norm{A_I}_1\ \le\ \norm{A_Ix}_1\ \le\ \tfrac12\norm{A_I}_1 .
\]
Both constants are optimal.
\end{lemma}

\begin{proof}
Let $y=\one^*A_I$ be the row of the column sums of $A_I$, so that $\norm{A_Ix}_1=yx$ and $\norm{A_I}_1=y\one$. Write $I=I'j$ and
$y'=\one^*A_{I'}$. Since $A_j=\id+(\one-e_j)e_j^*$, appending $j$ replaces the $j$-th column sum by the sum $y'\one$ of all of them and keeps the other two: $y_j=y'\one$ and $y_m=y'_m<y_j$ for $m\neq j$, because $y'_j\ge1$. So $y_j$ is the strict maximum of $y$.

The function $x\mapsto yx$ is linear on the triangle $C_i$, whose vertices are $\pr{\one-e_i}=\frac12(\one-e_i)$,
$\pr{\one+e_j}=\frac14(\one+e_j)$ and $\pr{\one+e_l}=\frac14(\one+e_l)$, where $l$ is the third index. Its values there
are $\frac12(y\one-y_i)$, $\frac14(y\one+y_j)$ and $\frac14(y\one+y_l)$. All of them are at most $\frac12y\one$. The last two are at least $\frac14y\one$, and so is the first, because $y_i\le y_j\le y_j+y_l$ gives $y_i\le\frac12y\one$.

For optimality take $I=j^k$. Then $y=\one^*+2k\,e_j^*$ and $y\one=3+2k$, so the ratio $yx/y\one$ equals $(1+k)/(3+2k)\to\frac12$ at $x=\pr{\one-e_i}$, and $(2+k)/\bigl(2(3+2k)\bigr)\to\frac14$ at $x=\pr{\one+e_l}$.
\end{proof}

\begin{remark}[fastness removes the parabolicity]\label{rem:fastN}
The one-step operator $\mathcal M_t$ fails because each $A_j$ is parabolic: within a run norms add rather than multiply. Across runs they multiply up to the factor $\frac14$ (Remark~\ref{rem:fast}), and Lemma~\ref{lem:junction} extends this, with the same constant, to every $I$ ending in $j$ in place of $A_j^k$: Lemma~\ref{lem:fastruns} holds for $\GN$ with $c_*=\frac14$. Taking whole runs as the steps of the renewal, as in Remark~\ref{rem:fast}, puts every step at a junction, and no state has a weight equal to $1$.
\end{remark}

\begin{remark}[the case of $\GN$]\label{rem:GNfast}
Since $\zeta^{(1)}(t)=3\sum_{k\ge1}(3+2k)^{-t}$, we have $t_0=1<2\le\sN$, and, by Remark~\ref{rem:fastN},
Proposition~\ref{prop:runpressure} applies: $\sN$ is the unique zero of the growth
rate per run, and $\zeta(\sN)=\infty$.

The operator $\cL_t$ of \eqref{eq:zeta} is the exact form of this argument:
by \eqref{eq:Lm}, $(\cL_t^m\mathds 1)(x,i)$ is the sum of $\norm{A_Ix}_1^{-t}$ over the multiindices $I$
with $m$ runs whose last index differs from $i$. By Lemma~\ref{lem:junction},
$\sup_{\mathcal C}\cL_t^m\mathds 1\le4^t\zeta^{(m)}(t)$. Conversely, given $m\ge1$ and a multiindex
$I$ with $m$ runs, choose an index $l$ different from $i$ and from the last index
of $I$; then $Il$ has $m+1$ runs, its last index is not $i$, and, by Lemma~\ref{lem:junction} and Lemma~\ref{lem:geom}(i),(v),
\[
  \norm{A_{Il}x}_1\le\tfrac12\norm{A_{Il}}_1=\tfrac12\bigl(\norm{A_I}_1+\norm{a_{I,l}}_1\bigr)\le\norm{A_I}_1 .
\]
Hence
$\inf_{\mathcal C}\cL_t^{m+1}\mathds 1\ge\zeta^{(m)}(t)\ge4^{-t}e^{m\,p(t)}$ by
Proposition~\ref{prop:runpressure}(i). Consequently, for every $t>1$ with
$t\neq\sN$ there is a test function as in Corollary~\ref{cor:test}:
\begin{itemize}
\item if $t>\sN$, take $e^{p(t)}<\theta<1$ and $h=\sum_{m\ge0}\theta^{-m}\cL_t^m\mathds 1$,
  which is bounded, satisfies $h\ge1$, and $\cL_th=\theta(h-\mathds 1)\le\theta h$;
\item if $1<t<\sN$, then $p(t)>0$, so $\inf_{\mathcal C}\cL_t^N\mathds 1\ge1$ for some
  $N$, and $h=\sum_{m<N}\cL_t^m\mathds 1$ satisfies $h\ge1$ and
  $\cL_th=h-\mathds 1+\cL_t^N\mathds 1\ge h$.
\end{itemize}
So exact test functions exist for every $t>1$ with $t\neq\sN$, and Theorem~\ref{thm: delta} is a matter of approximating them by functions that are constant on the cells of a grid. Remark~\ref{rem:contr} explains why this approximation improves linearly with the mesh; the certificates themselves do not depend on it.
\end{remark}
\begin{remark}[comparison with {\cite[Thm.~1]{DL15JGA}}]\label{rem:JGA}
The proof of \cite[Thm.~1]{DL15JGA} replaces the junction factor by the fastness constant and works with the single unknown $\zeta(t)$. This costs a factor of order $c_*^{-t}$ (Lemma~\ref{lem:fastruns}) at every renewal, so the renewals must run over blocks of norm larger than some $k$, and the bounds converge only like $1/\log k$. For $\GN$, with the norm $\norm{g}_{\infty\to\infty}=\max_i\sum_j\abs{g_{ij}}$, an exact enumeration of the $3.8\times10^{11}$ multiindices of the prefix code $\{I:\ \norm{A_I}_{\infty\to\infty}>65536\ge\norm{A_{I''}}_{\infty\to\infty}\}$, where $I''$ is $I$ without its first index, gives only $\sN\ge2.3444$. Here the junction factor is computed {\em exactly}, as a function of the soul point $\pr{L_J}$, so the renewals are single runs and the unknown is the function $(x,i)\mapsto\zeta_{x,i}(t)$. The only error left is the discretization of $\mathcal C$: the scalar conditions of \cite{DL15JGA} become Collatz--Wielandt conditions for the discretized operators of Section~\ref{sec:thmB}, and the single unknown corresponds to a single cell, $\mathcal C$ itself.
Baragar's bounds $2.430<\sN<2.477$~\cite[Thm.~0.1]{Bar98} come from a related scheme: his tree is the tree of the souls (Lemma~\ref{lem:baragartree}), truncated at height $y$, with the subtrees below the truncation bounded by Hurwitz zeta sums and by Boyd's convexity argument; he proves that the width of his bounds is $O(1/\log y)$~\cite[Thm.~2.3]{Bar98}, and notes that in practice they converge at this slow rate.
\end{remark}

\begin{lemma}[uniform contraction]\label{lem:contr}
Let $\theta_*=\sqrt{160/243}<0.812$, let $i\neq j$ and $k\ge1$. Then $\psi_j^k$ is $\theta_*$-Lipschitz on $C_i$. Moreover, if $I$ is a non-empty multiindex whose last index is not $i$, then $x\mapsto\log\norm{A_Ix}_1$ is $4$-Lipschitz on $C_i$.
\end{lemma}

\begin{proof}
Let $x\in C_i$, and let $u$, $d$ and $c$ be as in Lemma~\ref{lem:contraction}. Then $ku\ge\frac12$, so $\beta=ku/d\in[\frac13,1)$, and $d=1/(1-\beta)$. With $g=e_j-\frac13\one$, the orthogonal projection of $e_j$ on $T=\{v:\one^*v=0\}$, \eqref{eq:DR} reads $D_x\psi_j^k=d^{-1}B$ on $T$, $B=\id+\beta\,c\,g^*$, where $g\cdot c=-1$, $\abs g^2=\frac23$ and $\abs c^2\le2$ (proof of Lemma~\ref{lem:contraction}). So
\[
  \norm{D_x\psi_j^k}^2\le d^{-2}\,\mathrm{tr}(B^*B)=(1-\beta)^2\bigl(2-2\beta+\beta^2\abs c^2\abs g^2\bigr)\le(1-\beta)^2\bigl(2-2\beta+\tfrac43\beta^2\bigr)\le\tfrac{160}{243},
\]
because the last function has derivative $(1-\beta)\bigl(-6+\frac{26}3\beta-\frac{16}3\beta^2\bigr)<0$. Since $C_i$ is convex, $\psi_j^k$ is $\theta_*$-Lipschitz on $C_i$. For the second statement let $y=\one^*A_I$, so that $\log\norm{A_Ix}_1=\log(yx)$. For $x\in C_i$ and a unit vector $v\in T$, Lemma~\ref{lem:junction} gives $yx\ge\frac14y\one$, so $\abs{\partial_v\log(yx)}=\abs{yv}/(yx)\le4$, since $\abs{yv}\le\abs y\le y\one$.
\end{proof}

\begin{remark}\label{rem:contr}
At $k=1$ and $x=\pr{\one+e_l}$, where $l$ is the third index, we have $\beta=\frac13$ and $\abs c^2=2$, so the trace bound in the proof of Lemma~\ref{lem:contr} is an equality, and the largest singular value of $D_x\psi_j$ is $\sqrt{(80+8\sqrt{19})/243}=0.68754\ldots$; numerically, maximizing the operator norm of $D_x\psi_j^k$ over $x\in C_i$ and $k\le100$ gives the same value, which is therefore, numerically, the sharp constant in Lemma~\ref{lem:contr}. 
By the second statement of the lemma, the partial sums of \eqref{eq:Lm} are log-Lipschitz on each $C_i$, with constant $4t$. 
Together, these explain why test functions that are constant on the cells of a grid work well: in practice, the width of the bounds they give decreases linearly with the mesh.
Lemma~\ref{lem:contr} is the analogue, for the soul points, of the uniform contraction of the accelerated maps on the side of the lengths~\cite[Lemma~2]{AS13}, \cite[Lemma~13]{AHS16}, \cite[Prop.~45]{GMR19}, on which the counting theorem of~\cite{GMR19} rests.
We use it for the coding of Section~\ref{sec:AHS}.
\end{remark}

\begin{remark}[the role of the $\ell^1$ norm]\label{rem:l1}
The additivity of $\norm{\cdot}_1$ on non-negative matrices is what makes \eqref{eq:seen} exact and the state space two-dimensional. For the Frobenius norm, $\norm{A_{IJ}}_F^2=\operatorname{tr}\bigl(A_I(A_JA_J^*)A_I^*\bigr)$ depends on $J$ through the Gram matrix $A_JA_J^*$, a point of a $5$-dimensional space, and a grid at the resolution of Section~\ref{sec:thmB} would have about $10^{12}$ cells instead of $2\times10^5$. For the affinity exponent one more variable, the direction of the derivative, is needed (Lemma~\ref{lem:Dcocycle} and Proposition~\ref{cor:certA}).
\end{remark}

\section{Proof of Theorem \ref{thm: delta} }\label{sec:thmB}

In this section we prove Theorem~\ref{thm: delta} with test functions for the operator $\cL_t$ of Section~\ref{sec:fast} that are constant on the cells of a triangulation of each cylinder $C_i$, checked through Corollary~\ref{cor:test}.

\thmdelta*

In the coordinates $\lambda=(2-4x_l,2-4x_m)$ of $C_i$, where $l<m$ are the two indices different from $i$, the cylinder is the triangle $C_i=\{\lambda\ge0,\ \lambda_1+\lambda_2\le1\}$ (as in \S\ref{ssec:P1}). Subdivide it into the $M^2$ closed triangles $Q$ of the standard triangulation of mesh $1/M$ (the squares of side $1/M$ cut by their anti-diagonals), and let $h$ be
constant, $h=H_Q$, on each, with $h(x,i)=\min\{H_Q:\ x\in Q\}$ at points that lie in several cells. For a cell $Q\subset C_i$ and $j\neq i$, the image
$\psi_j^k(Q)$ is the triangle spanned by the images of the vertices of $Q$
($\psi_j^k$ is projective), and, by \eqref{eq:run}, $\bigcup_{k>K}\psi_j^k(Q)$ lies in the convex hull of
$\psi_j^{K+1}(Q)$ and $\pr{\one-e_j}$. Let $\underline H_{j,k}(Q)$ and $\overline H_{j,k}(Q)$ be the
minimum and maximum of $H$ over the cells of $C_j$ meeting the bounding box of that
image, in the coordinates $\lambda$ of $C_j$; any set of cells containing all those that
meet the image may be used (Figure~\ref{fig:thmB}). Then for $x\in Q$
\[
  (\cL_th)(x,i)\ \ge\ \sum_{j\ne i}\Bigl[\sum_{k=1}^K\bigl(1+k\,\bar u_j\bigr)^{-t}
  \underline H_{j,k}(Q)+\frac{\bigl(1+(K+1)\bar u_j\bigr)^{1-t}}{\bar u_j\,(t-1)}\,
  \underline H_{j,\infty}(Q)\Bigr],
\]
with $\bar u_j=\max_Q2x_j$, using $\sum_{k>K}f(k)\ge\int_{K+1}^\infty f$ for
decreasing $f$; similarly
\[
  (\cL_th)(x,i)\ \le\ \sum_{j\ne i}\Bigl[\sum_{k=1}^K\bigl(1+k\,\underline u_j\bigr)^{-t}
  \overline H_{j,k}(Q)+\frac{\bigl(1+K\underline u_j\bigr)^{1-t}}{\underline u_j\,(t-1)}\,
  \overline H_{j,\infty}(Q)\Bigr],
\]
with $\underline u_j=\min_Q2x_j$, using $\sum_{k>K}f(k)\le\int_K^\infty f$.
If the first right-hand side is $\ge H_Q$ for every cell, then $\cL_th\ge h$, and Corollary~\ref{cor:test} gives $t\le\sN$; if the second is $\le\theta H_Q$ with $\theta<1$ for every cell, then $\cL_th\le\theta h$, and it gives $t\ge\sN$. (At a point $x$ lying in several cells, the inequalities hold for each of them, hence for the smallest value $h(x,i)$; and the value of $h$ at a point of an image is the value of one of the cells containing that point, all of which are among the cells used.) We write $\cL_t^-H$ and $\cL_t^+H$ for the vectors of the two right-hand sides.

The program works in these coordinates; it stores a point $x\in\cT$ as the row vector $x^*$, so that $\psi_j^k$ reads $x^*\mapsto x^*(A_j^*)^k/\bigl(x^*(A_j^*)^k\one\bigr)$, and the weights are $1+2kx_j$.

The certificates do not depend on how $H$ was obtained; $H$ enters only as exact data. We use two vectors, $H^-$ for the lower and $H^+$ for the upper bound, computed by the nonlinear power iteration $H\leftarrow\cL_t^\mp H/\max(\cL_t^\mp H)$. So $H^\mp$ approximate the fixed vectors of the monotone, positively homogeneous maps $\cL^\mp_t$, and the margins obtained are close to the best that the discretization allows. (The Perron vector of the matrix obtained by evaluating $\cL_t$ at the centres of the cells does not certify the bounds below.) We bisect in $t$:

\begin{figure}[t]
\centering
\includegraphics[width=\textwidth]{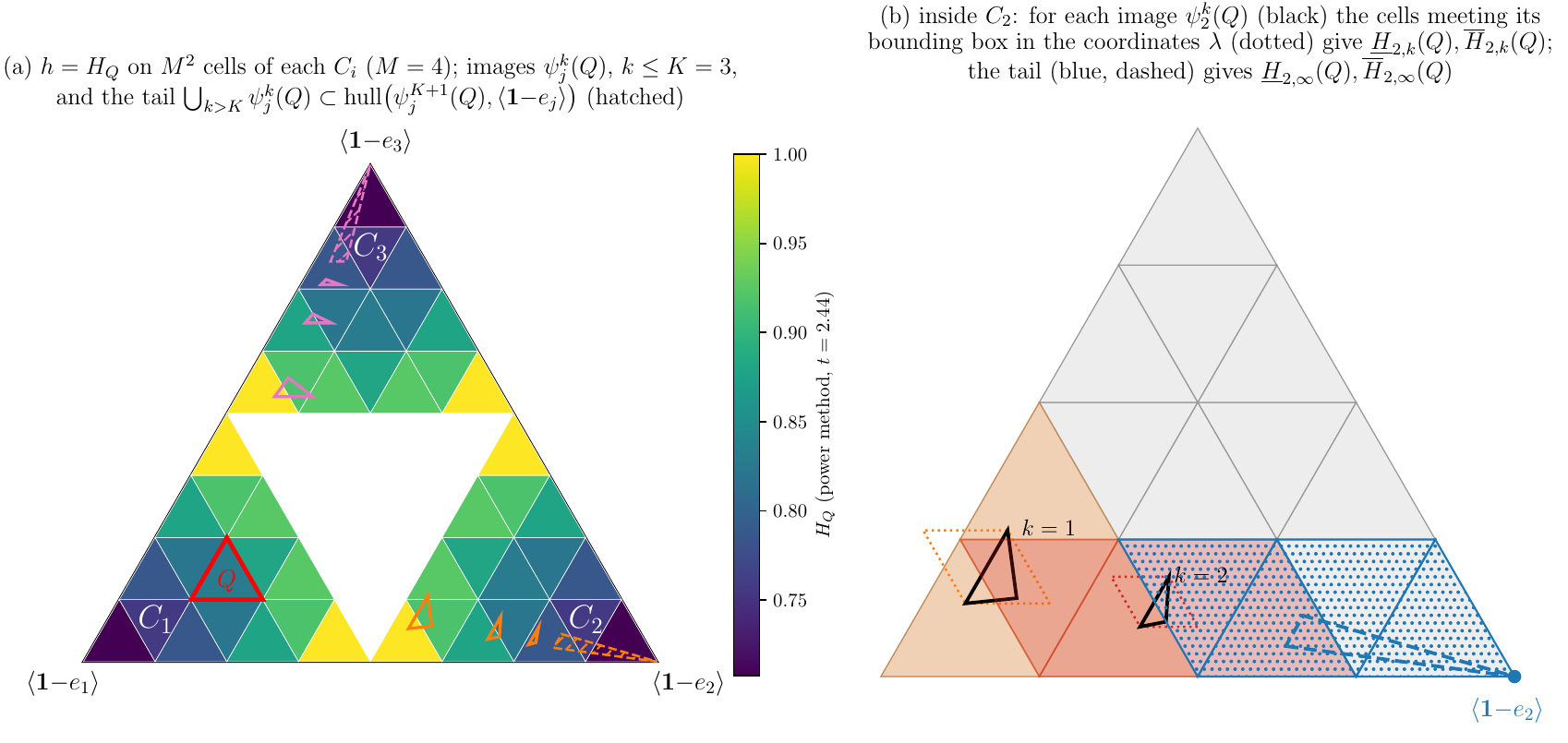}
\caption{The certificate of Section~\ref{sec:thmB}, drawn for $M=4$ and $K=3$, in the coordinates $\lambda$ of the cylinders. (a) The test function $h=H_Q$ is constant on the $M^2$ cells of each cylinder $C_i$ (colour: the values $H_Q$ of a test vector at $t=2.44$). For one cell $Q\subset C_1$ (red) the images $\psi_j^k(Q)$, $k\le K$, are drawn in $C_2$ and $C_3$; the tail $\bigcup_{k>K}\psi_j^k(Q)$ lies in the hatched convex hull of $\psi_j^{K+1}(Q)$ and of the parabolic point $\pr{\one-e_j}$. (b) Inside $C_2$: for each image $\psi_2^k(Q)$ (black) the cells meeting its bounding box (dotted) give $\underline H_{2,k}(Q)$ and $\overline H_{2,k}(Q)$; the cells meeting the bounding box of the tail region (blue, dashed) give $\underline H_{2,\infty}(Q)$ and $\overline H_{2,\infty}(Q)$.}
\label{fig:thmB}
\end{figure}

\begin{center}
\begin{tabular}{rrcc}
\toprule
$M$ & $K$ & certified $\sN\ge$ & certified $\sN\le$\\
\midrule
64 & 150 & 2.4336 & 2.4541\\
128 & 200 & 2.4386 & 2.4489\\
256 & 250 & 2.4411 & 2.4465\\
\bottomrule
\end{tabular}
\end{center}
The verification is done in interval arithmetic with Arb
\cite{Arb}: the images of the cell vertices are enclosed in balls, so that the
selected cells contain all cells meeting the exact images; the weights
and the tail integrals are computed as balls for a ball enclosing
$t$; and each $(\cL^\mp_tH^\mp)_Q$ is evaluated in ball arithmetic. At $M=256$ ($65\,536$ cells per cylinder),
$K=250$, this certifies $\min_Q(\cL^-_tH^-)_Q/H^-_Q\ge1.00014$ at $t=2.4411$ and
$\max_Q(\cL^+_tH^+)_Q/H^+_Q\le0.99976$ at $t=2.4465$.
The code is deposited at Zenodo~\cite{DLdata}.

\begin{proof}[Proof of Theorem~\ref{thm: delta}]
Let $h^-$ and $h^+$ be the functions on $\mathcal C$ with values $H^-$ and $H^+$ on the cells, as above. By the two inequalities above, the certificates at $M=256$ give $\cL_th^-\ge h^-$ at $t=2.4411$ and $\cL_th^+\le0.99976\,h^+$ at $t=2.4465$, and Corollary~\ref{cor:test} gives $2.4411\le\sN\le2.4465$.
\end{proof}

The widths of the certified intervals decrease like $1/M$, and the interval at $M=256$ is consistent with the numerical value $\sN\approx2.4437$. Already at $M=64$ the certified interval lies inside Baragar's interval $(2.430,2.477)$, and at $M=256$ it is about $9$ times narrower (Remark~\ref{rem:JGA}).

\section{Souls, Markoff--Hurwitz equations and higher-dimensional gaskets}\label{sec:MH}

In this section $\cN$ and $\cG$ become the first two members of the family of the Levitt--Yoccoz gaskets $\mathcal G^d$, one in each dimension $d$, and the norm exponents of the corresponding semigroups become the exponents of the Markoff--Hurwitz equations, studied by Baragar~\cite{Bar98} and by Gamburd, Magee and Ronan~\cite{GMR19}. We express their results in our language; in particular, we derive from~\cite[Thm.~13]{GMR19} our Theorem~\ref{thm:count}, the asymptotics of the number of souls.

\subsection*{The semigroups \texorpdfstring{$\Gamma^{(d)}$}{Gamma(d)} and the gaskets \texorpdfstring{$\mathcal G^d$}{G\^{}d}}
For $d\ge1$ let $\Gamma^{(d)}\subset SL_{d+1}(\Z)$ be the semigroup generated by the $d+1$ matrices
\[
  A_j=\id+(\one-e_j)e_j^*,\qquad j=1,\dots,d+1,
\]
so that $\Gamma^{(2)}=\GN$ and $\Gamma^{(3)}=\GG$. As for $\GN$, we call $L_I=A_I\one$ the soul of the multiindex $I$, so that $\norm{A_I}_1=\one^*A_I\one=\norm{L_I}_1$, and we write $N_d(r)=\#\{I\in\cI^{d+1}:\ \norm{A_I}_1\le r\}$.
The semigroup $\Gamma^{(d)}$ acts projectively, by $\psi_I(x)=\pr{A_Ix}$, on the $d$-simplex
\[
  \cT_d=\left\{x\in\R^{d+1}:\ \one^*x=1,\ 0\le x_m\le\frac1d\right\},
\]
whose vertices are the points $\pr{\one-e_m}$, and its limit set
\[
  \mathcal G^d=\bigcap_{n\ge0}\ \bigcup_{\abs I=n}\psi_I(\cT_d)
\]
is the \textbf{$d$-dimensional Levitt--Yoccoz gasket}. Up to a projective change of coordinates, it is the case $n=d$ of the fractals $F(\mathcal E)$ of~\cite{DL08}, which is also the cubic gasket in $\RP^d$ of~\cite{DL15JGA}: $F(\mathcal E)$ is the image of $\mathcal G^d$ under the projective map with matrix $\mathcal E(\one\one^*-d\,\id)$. Indeed, write a basis of $\R^{d+1}$ as the matrix whose columns are its vectors. The tree of bases of~\cite[\S2]{DL08} rooted at $\mathcal E$ consists of the matrices $\mathcal E_I=\mathcal EA_{i_1}^*\cdots A_{i_k}^*$, $I=i_1\cdots i_k$, because replacing $e_m$ by $e_m+e_j$ for $m\ne j$ is the right multiplication by $A_j^*$ (equivalently, the generators of the cubic gasket of~\cite{DL15JGA} over the frame $\mathcal E$ are the matrices $\mathcal EA_j^*\mathcal E^{-1}$), and $F(\mathcal E)=\bigcap_{k\ge0}\bigcup_{\abs I=k}S(\mathcal E_I)$, where $S(\mathcal E_I)$ is the projective simplex of the points $\pr{\mathcal E_Iy}$, $y\ge0$, $y\ne0$. Since $(\one\one^*-\id)A_j^*=\one\one^*-\id+(\one-e_j)(\one-e_j)^*=A_j(\one\one^*-\id)$, as in Remark~\ref{rem:conj}, and $(\one\one^*-d\,\id)(\one\one^*-\id)=d\,\id$, we have $\mathcal E(\one\one^*-d\,\id)A_I(\one\one^*-\id)=d\,\mathcal E_I$. The vertices of $\psi_I(\cT_d)$ are the points $\pr{A_I(\one\one^*-\id)e_m}=\pr{A_I(\one-e_m)}$, so the projective map with matrix $\mathcal E(\one\one^*-d\,\id)$ maps $\psi_I(\cT_d)$ onto $S(\mathcal E_I)$, and $\mathcal G^d$ onto $F(\mathcal E)$; it maps the soul $\pr{L_I}$ to the barycentre $\pr{\mathcal E_I\one}$ of~\cite{DL08}, since $(\one\one^*-\id)\one=d\,\one$. See also~\cite{AS13,Fou20}; in~\cite{Fou20} the image of $\mathcal G^d$ under the similarity $x\mapsto\one-d\,x$ is called the Rauzy gasket of the $d$-simplex.
In every dimension, $\mathcal G^d$ is also, up to the boundaries of its zones and to permutations of the coordinates, the set of parameters of Meester's circle percolation on $\Z^d$ for which his algorithm for the critical value never stops~\cite{Mee89}. This algorithm subtracts the smallest coordinate of a vector of $\R^{d+1}$ from the others, which is $\psi_i^{-1}$ when the $i$-th coordinate is the smallest, as it is on $\psi_i(\cT_d)$; and it stops when the largest coordinate exceeds the sum of the others divided by $d-1$, that is, when the point leaves $\cT_d$, on which $x_m\le\frac1d$, i.e.\ $x_m\le\frac1{d-1}\sum_{l\ne m}x_l$, for every $m$. So $\cN=\mathcal G^2$, with $\cT_2=\cT$, and $\cG=\mathcal G^3$.
For $d=1$, the semigroup $\Gamma^{(1)}$ is generated by $\begin{pmatrix}1&0\\1&1\end{pmatrix}$ and $\begin{pmatrix}1&1\\0&1\end{pmatrix}$; its souls $A_I\one$ are the pairs of coprime positive integers, each obtained once, and $\sN(\Gamma^{(1)})=2$ (Section~\ref{sec:general}).

Kraaikamp and Meester proved that the set of parameters above has zero volume, settling a conjecture of~\cite{Mee89}. Since the boundaries of the zones have zero volume, in our notation their result reads as follows.

\begin{lettertheorem}[Kraaikamp--Meester~{\cite[Lemma~1]{KM95}}; for $d=2$, Meester--Nowicki~{\cite[Thm.~5.1]{MN89}}]\label{thm:KM}
$\mathrm{vol}(\mathcal G^d)=0$ for every $d\ge2$.
\end{lettertheorem}

\begin{remark}[monotonicity in $d$]\label{rem:monotone}
For $d\ge2$, $\sN(\Gamma^{(d)})\ge\sN(\Gamma^{(d-1)})$; this is~\cite[Remark~14]{GMR19} in our language. Directly: let $I\in\cI^d$. For $j\le d$, the matrix $A_j$ adds $x_j$ to the other coordinates of $x\in\R^{d+1}$ and never uses $x_{d+1}$, so the first $d$ coordinates of $L_I=A_I\one$ form the soul of $I$ in $\Gamma^{(d-1)}$. Moreover, each coordinate of a soul of $\Gamma^{(d)}$ is at most the sum of the others. This holds for $\one$, and if it holds for $x\ge0$, it holds for $A_jx$: for $k\ne j$, $(A_jx)_k=x_k+x_j\le x_j+\sum_{m\ne k,j}x_m+(d-1)x_j=\sum_{m\ne k}(A_jx)_m$, because $x_k\le\sum_{m\ne k}x_m$ and $d\ge2$; and $(A_jx)_j=x_j\le\sum_{m\ne j}(x_m+x_j)$. Hence the norm of $A_I$ in $\Gamma^{(d)}$ is at most twice its norm in $\Gamma^{(d-1)}$, so $N_d(2r)\ge N_{d-1}(r)$, and the claim follows. For $d=2$, this is the argument of Section~\ref{sec:general} that gives $\sN\ge2$.
\end{remark}

\subsection*{Baragar's trees}
For $d\ge1$ consider the Markoff--Hurwitz equation in $d+2$ variables
\begin{equation}\label{eq:MHd}
  x_0^2+x_1^2+\dots+x_{d+1}^2=a\,x_0x_1\cdots x_{d+1},
\end{equation}
which for $d=2$ is~\eqref{eq:MH}. It is quadratic in each variable, so it is preserved by the involutions that replace $x_i$ by $a\prod_{m\ne i}x_m-x_i$. For $a\in\mathbb Q$ and a rational solution $P\neq0$, let $N_P(t)$ be the number of points $Q$ of the orbit of $P$ under the group generated by these involutions whose logarithmic height $\max_ih(q_i)$ is less than $t$. Baragar reduces the growth of $N_P(t)$ to that of the following trees~\cite{Bar98}. For $r\in\R^{d+1}$ with positive entries, the tree $E_r$ has root $r$, and every node $b$ has the $d+1$ children
\[
  T_i(b)=\bigl(b_i,\ b_1+b_i,\ \dots,\ \widehat{b_i+b_i},\ \dots,\ b_{d+1}+b_i\bigr),\qquad i=1,\dots,d+1,
\]
where the hat means that the entry is omitted; the nodes are counted with multiplicity, that is, as the words in $T_1,\dots,T_{d+1}$~\cite[\S4]{Bar98}, and the height of a node $b$ is $b_1+\dots+b_{d+1}$. The exponent of $E_r$ is the abscissa of convergence of the series $f_r(s)$ of the heights of its nodes raised to the power $-s$.

\begin{lemma}\label{lem:baragartree}
Let $r\in\R^{d+1}$ have positive entries. For every $m\ge0$, the heights of the nodes of $E_r$ at depth $m$ are the numbers $\one^*A_Ir$, where $I$ runs over the elements of $\cI^{d+1}$ of length $m$, each counted once. In particular $f_\one=\zeta_{\Gamma^{(d)}}$, and the exponent of every tree $E_r$ is $\sN(\Gamma^{(d)})$.
\end{lemma}

\begin{proof}
Since $(A_ib)_m=b_m+b_i$ for $m\ne i$ and $(A_ib)_i=b_i$, we have $T_i=P_iA_i$, where $P_i$ is the permutation matrix that moves the $i$-th coordinate to the first place and keeps the order of the others. For a permutation matrix $P$ with $Pe_k=e_{\pi(k)}$ we have $PA_kP^{-1}=A_{\pi(k)}$, that is, $A_kP=PA_{\pi^{-1}(k)}$. So, by induction on $m$, $T_{i_m}\cdots T_{i_1}=PA_{k_m}\cdots A_{k_1}$ for a permutation matrix $P$, where $k_m=\pi_Q^{-1}(i_m)$ if $T_{i_{m-1}}\cdots T_{i_1}=QA_{k_{m-1}}\cdots A_{k_1}$; conversely $i_m=\pi_Q(k_m)$, so this is a bijection between the words of length $m$ in the $T_i$ and the $K=k_m\cdots k_1\in\cI^{d+1}$ of length $m$. Since $\one^*P=\one^*$, the height of the node $T_{i_m}\cdots T_{i_1}r$ is $\one^*A_Kr$. For $r=\one$ it is $\norm{A_K}_1$, so $f_\one=\zeta_{\Gamma^{(d)}}$; in general $(\min_mr_m)\norm{A_K}_1\le\one^*A_Kr\le(\max_mr_m)\norm{A_K}_1$, because $A_K\ge0$.
\end{proof}

In~\cite{Bar98} the equation~\eqref{eq:MHd} is written with $n+1=d+2$ variables, and its exponent is called $\alpha(n)=\alpha(d+1)$; in~\cite{GMR19} the number of variables is $n=d+2$, and the exponent is called $\beta(n)=\beta(d+2)$. By Lemma~\ref{lem:baragartree}, both are $\sN(\Gamma^{(d)})$, and Baragar's theorem reads as follows.

\begin{lettertheorem}[Baragar~{\cite[Thm.~0.1]{Bar98}}; see also~{\cite[Thm.~2]{GMR19}}]\label{thm:baragar}
Let $d\ge1$, $a\in\mathbb Q$, and let $P\neq0$ be a rational solution of~\eqref{eq:MHd}. Then $\log N_P(t)/\log t\to\sN(\Gamma^{(d)})$ as $t\to\infty$. In particular, if $a\in\Z$ and~\eqref{eq:MHd} has an integer solution other than $0$, then the number of its integer solutions with $\max_i\abs{x_i}\le R$ is $(\log R)^{\sN(\Gamma^{(d)})+o(1)}$ as $R\to\infty$.
\end{lettertheorem}

Baragar also proved the bounds $2.430<\sN(\Gamma^{(2)})<2.477$, $2.730<\sN(\Gamma^{(3)})<2.798$ and $2.963<\sN(\Gamma^{(4)})<3.048$, and $\log(d+1)/\log2<\sN(\Gamma^{(d)})<\log(d+1)/\log2+o(d^{-0.58})$ as $d\to\infty$~\cite[Thm.~0.1 and Table~1]{Bar98}. Theorems~\ref{thm: delta} and~\ref{thm:3d} improve the first two of these bounds. Silverman asked whether these exponents are irrational (see~\cite[p.~753]{GMR19}).

\subsection*{The theorem of Gamburd, Magee and Ronan}
Let $d\ge2$. 
Gamburd, Magee and Ronan use the $(d+1)\times(d+1)$ matrices $\gamma_1,\dots,\gamma_{d+1}$ given by
\[
  \gamma_jw=\bigl(w_1,\dots,\widehat{w_j},\dots,w_{d+1},\ \one^*w\bigr)
\]
\cite[(1.4) and Example~9]{GMR19}: $\gamma_j$ drops the $j$-th coordinate and appends the sum of all of them. 
Since $A_j^*$ replaces $w_j$ by $\one^*w$ and keeps the other coordinates, $\gamma_j=Q_jA_j^*$, where $Q_j$ is the permutation matrix that moves the $j$-th coordinate to the last place and keeps the order of the others; in particular $\gamma_{d+1}=A_{d+1}^*$. 
Gamburd, Magee and Ronan observed that, up to multiplying the generators by permutation matrices, the projective action of this semigroup on the simplex with vertices $\pr{e_1},\dots,\pr{e_{d+1}}$ has the Rauzy gasket as attracting set for $d=2$, and that for $d\ge3$ it is the action that defines the higher-dimensional versions of the Rauzy gasket of~\cite{DL08,AS13}~\cite[\S1, after Example~9]{GMR19}; by Remark~\ref{rem:conj} and the beginning of this section, its attracting set is the image of $\mathcal G^d$ under the similarity $x\mapsto\one-d\,x$. 
They accelerate the semigroup generated by the $\gamma_j$ to the semigroup $\Gamma$ generated by
\[
  \mathcal A=\bigl\{\gamma_{d+1}^{\,b}\gamma_j:\ b\ge0,\ 1\le j\le d\bigr\},
\]
and both semigroups are free on these generators~\cite[p.~756]{GMR19}. On the cone $\mathcal W=\{w\in\R^{d+1}:\ 0\le w_1\le\dots\le w_{d+1}\}\setminus\{0\}$ they count
\[
  N^{\mathrm{GMR}}(w,r)=\#\bigl\{\gamma\in\Gamma\cup\{\id\}:\ \one^*\gamma w\le e^r\,\one^*w\bigr\}
\]
\cite[(1.7)]{GMR19} (in their notation, $\one^*w$ and $\one^*\gamma w$ are the last coordinates $y_n$ and $(\gamma.y)_n$ of the points $y$ and $\gamma.y$ of their hyperplane $H$, written in the basis $e_j+e_n$, $j<n$).

The next result is fundamental for us.

\begin{lettertheorem}[Gamburd--Magee--Ronan~{\cite[Thm.~13]{GMR19}}]
\label{thm:GMR}
    Let $d\ge2$. 
    There are $\beta>1$ and a positive, bounded, $C^1$ function $h$ on $\mathcal W$, constant on the rays, such that, as $r\to\infty$, uniformly in $w\in\mathcal W$,
\[
N^{\mathrm{GMR}}(w,r)=h(w)\,e^{\beta r}\bigl(1+o(1)\bigr).
\]
    Moreover, $\beta$ is the exponent of Theorem~\ref{thm:baragar} for the equation~\eqref{eq:MHd}.
\end{lettertheorem}

\begin{proposition}\label{lem:GMRdict}
Let $d\ge2$, and let $\beta$ and $h$ be as in Theorem~\ref{thm:GMR}. Then, as $r\to\infty$, uniformly in $w\in\mathcal W$,
\[
  \#\bigl\{I\in\cI^{d+1}:\ I=\emptyset\text{ or } i_1\neq d+1,\text{ and }\pr{w,L_I}\le e^r\pr{w,\one}\bigr\}=h(w)\,e^{\beta r}\bigl(1+o(1)\bigr),
\]
where $L_I=A_I\one$ and $\pr{w,v}=w^*v$.
\end{proposition}

\begin{proof}
    It is enough to prove that $N^{\mathrm{GMR}}(w,r)$ equals the left hand side of the equality in the claim.

    Now, for a permutation matrix $P$ with $Pe_k=e_{\pi(k)}$ we have $PA_k^*P^{-1}=A_{\pi(k)}^*$, that is, $A_k^*P=PA^*_{\pi^{-1}(k)}$.
    As in the proof of Lemma~\ref{lem:baragartree}, moving the permutation matrices to the left gives $\gamma_{j_1}\cdots\gamma_{j_m}=PA_{i_m}^*\cdots A_{i_1}^*$ for a permutation matrix $P$, with $i_1=j_m$, and $(j_1,\dots,j_m)\mapsto(i_1,\dots,i_m)$ is a bijection of the words of length $m$.
    Since $\one^*P=\one^*$,
\[
  \one^*\gamma_{j_1}\cdots\gamma_{j_m}w=\one^*A_{i_m}^*\cdots A_{i_1}^*w=\pr{w,A_{i_1}\cdots A_{i_m}\one}=\pr{w,L_I},\qquad I=i_1\cdots i_m,
\]
    and the first index of $I$ is $i_1=j_m$, the index of the first generator to act on $w$.
    The products of elements of $\mathcal A$ are exactly the words in the $\gamma_j$ whose last letter $\gamma_{j_m}$ has $j_m\ne d+1$: such a word splits in a unique way after each of its letters $\gamma_j$ with $j\ne d+1$.
    Since $\Gamma$ is free on $\mathcal A$, the elements of $\Gamma$ correspond bijectively to these words, hence to the $I\in\cI^{d+1}$ with $i_1\neq d+1$; and $\id$ corresponds to $I=\emptyset$, with $\one^*w=\pr{w,L_\emptyset}$.
\end{proof}

\thmcount*

\begin{proof}
Let $d\ge2$, and let $\beta$ and $h$ be as in Theorem~\ref{thm:GMR}. Take $w=\one$, which lies in $\mathcal W$; then $\pr{w,L_K}=\norm{A_K}_1$ and $\pr{w,\one}=d+1=\norm{A_\emptyset}_1$. For a permutation matrix $P$ with $Pe_k=e_{\pi(k)}$ we have $PA_KP^{-1}=A_{\pi(K)}$, where $\pi(K)$ is obtained by applying $\pi$ to each index of $K$, and $\norm{PA_KP^{-1}}_1=\norm{A_K}_1$. So, for $R\ge d+1$, the number of non-empty $K\in\cI^{d+1}$ with first index $i$ and $\norm{A_K}_1\le R$ does not depend on $i$, and it is $(N_d(R)-1)/(d+1)$. Hence the left side of Theorem~\ref{thm:GMR} at $w=\one$ is $1+\frac d{d+1}\bigl(N_d((d+1)e^r)-1\bigr)$, and Theorem~\ref{thm:GMR}, with $R=(d+1)e^r$, gives
\[
  N_d(R)=1+\frac{d+1}d\Bigl(h(\one)\Bigl(\frac R{d+1}\Bigr)^{\beta}\bigl(1+o(1)\bigr)-1\Bigr)\ \sim\ \frac{(d+1)^{1-\beta}}d\,h(\one)\,R^\beta .
\]
Since $\sN(\Gamma^{(d)})=\limsup_{R\to\infty}\log N_d(R)/\log R$ (introduction), $\beta=\sN(\Gamma^{(d)})$, and $N_d(R)\sim c_d\,R^{\sN(\Gamma^{(d)})}$ with $c_d=(d+1)^{1-\sN(\Gamma^{(d)})}h(\one)/d>0$. For $d=2$, $N=N_2$, and the number of zones in $\RP^2$ with $\norm L_1\le r$ is $4N(r)+3$ (introduction).
\end{proof}

\begin{remark}[about the proof of Theorem~\ref{thm:GMR}]\label{rem:GMRproof}
The fact that Theorem~\ref{thm:GMR}, coming from a number theory environment, applies to $\cN$ is not a coincidence.
The proof of Theorem~\ref{thm:GMR} in~\cite{GMR19} works with the same steps as ours, the runs. Indeed, each generator in $\mathcal A$ is a run: since $Q_je_j=e_{d+1}$, we have $\gamma_{d+1}^{\,b}\gamma_j=(A_{d+1}^*)^bQ_jA_j^*=Q_j(A_j^*)^{b+1}$, the transpose of the run $A_j^{b+1}$, of length $b+1$ and index $j$, followed by the relabelling that moves $j$ to the last place, and the condition $j\le d$ says that consecutive runs have different indices. The proof applies the renewal method of Lalley~\cite{Lal89} to the transfer operator whose steps are the elements of $\mathcal A$, with weights $(\one^*\gamma w/\one^*w)^{-s}$, $\gamma\in\mathcal A$, acting on the $C^1$ functions on the projectivization of $\mathcal W$: for real $s>1$ it has a simple leading eigenvalue $\lambda_s$ and a spectral gap~\cite[Thm.~39]{GMR19}, $\lambda_s$ is real-analytic and strictly decreasing, and $\lambda_\beta=1$~\cite[Prop.~43 and \S4.3]{GMR19}, and the weights satisfy the non-lattice condition that the renewal method needs~\cite[Prop.~42]{GMR19}. The main input is the uniform contraction of the maps $\gamma\in\mathcal A$~\cite[Prop.~45]{GMR19}, proved for $d=2$ in~\cite[Lemma~2]{AS13} and~\cite[Lemma~13]{AHS16}. For $\gamma=Q_j(A_j^*)^k$, the weight is $\one^*(A_j^*)^kw/\one^*w=\pr{w,A_j^k\one}/\pr{w,\one}$, the factor by which a run of length $k$ multiplies $\pr{w,L}$ at $L=\one$, on the side of the transposed matrices; its analogue on the side of the souls is the operator $\cL_t$ of Section~\ref{sec:fast}, whose weights are the run factors of Section~\ref{sec:cocycle}, and the analogue of~\cite[Prop.~45]{GMR19} is Lemma~\ref{lem:contr}.
\end{remark}

\begin{remark}[numerics]\label{rem:countnum}
All the column sums of every $A_I\in\GN$ are odd (appending $j$ to $I$ replaces the $j$-th column sum by the sum of all three, as in the proof of Lemma~\ref{lem:junction}), so every $\norm{A_I}_1$ is odd and $N(r)$ jumps only at odd integers. An exact enumeration of the $1.6\times10^{11}$ multiindices with $\norm{A_I}_1\le2^{17}$, by a depth-first search of the tree of the column sums (program \texttt{colsum.c} in~\cite{DLdata}, which also gives the counts of Section~\ref{sec:3d} and Remark~\ref{rem:higher}), gives $N(2^{16})=30\,093\,377\,833$ and $N(2^{17})=163\,722\,940\,573$. The local slopes $\log_2\bigl(N(2r)/N(r)\bigr)$ at $r=2^{10},\dots,2^{16}$ are $2.44435$, $2.44417$, $2.44390$, $2.44386$, $2.44378$, $2.44375$, $2.44374$, and least-squares fits of $\log N(r)=\log c+s\log r$ at the odd $r$ in $[2^k,2^{17}]$, $k=14,15,16$, give $s=2.44368$, $2.44369$, $2.44370$ and $c=0.05112$, $0.05111$, $0.05110$, consistent with Theorem~\ref{thm: delta}; so $c_N\approx0.0511$. This value is fitted jointly with $s$, and it is sensitive to it: $N(2^{17})/2^{17s}=0.05171$, $0.05110$, $0.05050$ for $s=2.4427$, $2.4437$, $2.4447$. With $s=2.4437$, the relative oscillation $(\max-\min)/\mathrm{mean}$ of $N(r)/r^s$ over the odd $r$ in $[2^k,2^{k+1}]$ is $5.0\times10^{-3}$, $2.3\times10^{-3}$, $1.2\times10^{-3}$, $5.5\times10^{-4}$, $2.6\times10^{-4}$, $1.3\times10^{-4}$, $6.2\times10^{-5}$ for $k=10,\dots,16$: it halves at each octave and is of the order $2\sN/r$ of the relative size of a single jump of $N$, with no log-periodic component. For the Sierpi\'nski semigroup of Proposition~\ref{prop:lattice}, the ratio of the largest to the smallest value of $N_{\mathcal S}(r)/r^{s_{\mathcal S}}$ in an octave tends to $3$.

Notice that the parity of $\|L_I\|_1=\|A_I\|_1$ has a topological explanation that does not use the semigroup: the $\mu$-cube is self-complementary, and for every self-complementary triply periodic surface the souls satisfy a parity constraint, which for the $\mu$-cube reads “the sum of the components of each soul is odd” \cite{DL26parity}. Here the same fact is recovered from the recursion on the column sums of $A_I$, and it holds for the souls of all four quadrants, since a sign change of the coordinates does not alter the parity of $a+b+c$.
\end{remark}

\begin{remark}[one-sided geodesics]\label{rem:geodesics}
For $a=4$, Huang and Norbury showed that the positive real solutions of~\eqref{eq:MH} parametrize the Teichm\"uller space of the thrice-punctured projective plane, with coordinates that are functions of the lengths of one-sided simple closed geodesics~\cite{HN17}, and they deduced from Theorem~\ref{thm:baragar} that, for every finite-area hyperbolic structure $J$ on it, the number $n_J(L)$ of one-sided simple closed geodesics of length at most $L$ satisfies $\log n_J(L)/\log L\to\sN$ as $L\to\infty$ (see~\cite[\S1.1]{GMR19}); Magee proved that $n_J(L)\sim c(J)\,L^{\sN}$ for some $c(J)>0$~\cite{Mag20}. So Theorem~\ref{thm: delta} is also a bound for the growth exponent of these geodesics, which, unlike the exponents of Mirzakhani's counts of simple closed geodesics on orientable surfaces, is not an integer.
\end{remark}

\begin{remark}[the gaskets $\mathcal G^d$ for $d\ge4$]\label{rem:higher}
For $x\in\cT_d$ we have $\norm{A_j^kx}_1=1+dkx_j$, so the run factors of $\Gamma^{(d)}$ are as in Sections~\ref{sec:fast} and~\ref{sec:3d}, with state space the first-level cylinders, of dimension $d$; so $\sN(\Gamma^{(d)})$ can be certified as in Sections~\ref{sec:thmB} and~\ref{sec:3d}. Since $\mathcal G^d$ has zero volume (Theorem~\ref{thm:KM}), if the volumes and the surface areas of the zones of $\mathcal G^d$ are bounded as in Lemma~\ref{lem:geom3}, Theorem~\ref{thm:cutout}(a) gives $\udimB\mathcal G^d\le\frac12(d-1+\sN(\Gamma^{(d)}))$, provided that $d-1<\sN(\Gamma^{(d)})<d+1$. A computation of the affinity exponent would need, in addition, the Grassmannians of $k$-planes in the tangent space, $1\le k\le d-1$, so the norm exponent is the quantity that remains computable as $d$ grows.

The condition $\sN(\Gamma^{(d)})>d-1$ fails for every $d\ge5$, by~\cite[Table~1]{Bar98} for $d\le9$ and by the inequality in the proof of~\cite[Thm.~5.1]{Bar98} in general; then $\sum_I\norm{A_I}_1^{-(d-1)}<\infty$, and, under the same bounds on the zones, the proof of Theorem~\ref{thm:cutout}(a) gives only $\udimB\mathcal G^d\le d-1$. The borderline case $d=4$ is the question of Baragar whether $\alpha(5)=3$~\cite[\S4]{Bar98}, that is, whether $\sN(\Gamma^{(4)})=3=d-1$; by Theorem~\ref{thm:alpha5}, $\sN(\Gamma^{(4)})>3$.
Numerically, a spectral collocation of the operator $\tilde\cL_t$ of the proof of Theorem~\ref{thm:alpha5} on symmetric polynomials of degree $4$ and $6$ gives the same leading eigenvalue to within $10^{-5}$, and $\sN(\Gamma^{(4)})\approx3.0051$. This is consistent with the exact count: the local slopes $\log_2\bigl(N_4(2r)/N_4(r)\bigr)$ at $r=2^9,\dots,2^{13}$ are $3.0143$, $3.0115$, $3.0076$, $3.0064$, $3.0058$, still decreasing.
\end{remark}

\subsection*{Baragar's question}
We prove Theorem~\ref{thm:alpha5} with the method of Sections~\ref{sec:fast}, \ref{sec:thmB} and~\ref{sec:3d}, now for the semigroup $\Gamma^{(4)}$ of the $5\times5$ matrices $A_j=\id+(\one-e_j)e_j^*$, acting on $\cT_4=\{x\in\R^5:\ \one^*x=1,\ 0\le x_m\le\frac14\}$. Since the eigenfunction of the run operator is smooth, we use a polynomial test function, and we verify the required inequality at every point of the state space by a mean value argument.

\thmalphafive*

\begin{proof}
For $x\in\cT_4$ we have $A_j^kx=x+kx_j(\one-e_j)$, $\norm{A_j^kx}_1=1+4kx_j$ and $\psi_j^k(x)=(1-\beta)x+\beta\pr{\one-e_j}$ with $\beta=4kx_j/(1+4kx_j)$. The first-level cylinders are $C_i=\psi_i(\cT_4)=\{x\in\cT_4:\ x_i\le\frac18\}$, and, as in Lemma~\ref{lem:trap3}, $\psi_j^k(\cT_4)\subset C_j$ for every $k\ge1$ and $x_j\ge1-\frac18-\frac34=\frac18$ for $x\in C_i$ and $j\ne i$. The renewal~\eqref{eq:renewal} and the expression of the zeta function of Section~\ref{sec:3d} hold verbatim, with $\norm{A_j^k\one}_1=5+4k$, and so does Corollary~\ref{cor:test}, with $\sN(\Gamma^{(4)})$ in place of $\sN$; and, by the argument of Section~\ref{sec:3d}, it suffices to find a positive continuous function $h_1$ on $C_1$, invariant under the permutations of $x_2,\dots,x_5$, such that $\tilde\cL_th_1\ge h_1$ on $C_1$ for some $t>3.0002$, where
\[
  (\tilde\cL_th_1)(x)=\sum_{j=2}^5\sum_{k\ge1}\bigl(1+4kx_j\bigr)^{-t}\,h_1\bigl(\tau_j\psi_j^k(x)\bigr),
\]
$\tau_j$ being the transposition of the coordinates $1$ and $j$. In the coordinates $\lambda=(2-8x_2,\dots,2-8x_5)$, the cylinder $C_1$ is the simplex $\{\lambda\ge0,\ \lambda_1+\dots+\lambda_4\le1\}$, and $\tilde\cL_t$ commutes with the permutations of the $\lambda_i$.

Let $\kappa$ be the binary64 number nearest to $\frac15$, and let $h_1=\sum_\mu c_\mu m_\mu(\lambda-\kappa\one)$, where $m_\mu$ is the monomial symmetric polynomial of the partition $\mu$ in four variables, and the coefficients $c_\mu$ are the binary64 numbers of Table~\ref{tab:h5}, obtained by a spectral collocation of $\tilde\cL_3$. Interval evaluation on a cover of $C_1$ by small cubes gives $0.8611\le h_1\le1.1945$ on $C_1$.

\begin{table}[t]
\centering
\begin{tabular}{cr@{\qquad}cr}
\toprule
$\mu$ & $c_\mu$ & $\mu$ & $c_\mu$\\
\midrule
$\emptyset$ & $0.9361940872141115$ & $(1,1,1)$ & $0.041278142914558065$\\
$(1)$ & $0.13533887045421458$ & $(4)$ & $0.11872677888126022$\\
$(2)$ & $0.20764723855176836$ & $(3,1)$ & $0.012222273707418643$\\
$(1,1)$ & $0.06775138223770444$ & $(2,2)$ & $0.060961678511604646$\\
$(3)$ & $0.04189374357975737$ & $(2,1,1)$ & $-0.008581600298358461$\\
$(2,1)$ & $-0.005183312461752605$ & $(1,1,1,1)$ & $-0.18972389429328756$\\
\bottomrule
\end{tabular}
\caption{The coefficients of the test function $h_1$ in the proof of Theorem~\ref{thm:alpha5}.}
\label{tab:h5}
\end{table}

(a) \emph{$\tilde\cL_3h_1\ge(1+\varepsilon)h_1$ on $C_1$, with $\varepsilon=0.000999$.} By symmetry it suffices to prove this where $\lambda_1\ge\dots\ge\lambda_4$, a set covered by the $38\,097$ cubes $Q=\frac1{64}(c+[0,1]^4)$ with $c\in\Z^4$, $c_1\ge\dots\ge c_4\ge0$ and $c_1+\dots+c_4\le63$; on them the formulas above are defined, since $\frac18\le x_m\le\frac14$ for $m\ge2$. On a cube $Q$ we split $\tilde\cL_3h_1$ into the runs with $k\le40$ and the tail. Since $\psi_j^k(x)$ moves towards $\pr{\one-e_j}$ as $k$ grows, the images $\tau_j\psi_j^k(x)$, $k>40$, lie on the segment from $\tau_j\psi_j^{41}(x)$ to $\tau_j\pr{\one-e_j}=\pr{\one-e_1}$, the point $\lambda=0$; and, the summand being decreasing in $k$,
\[
  \sum_{k>40}\bigl(1+4kx_j\bigr)^{-3}\ \ge\ \int_{41}^\infty\bigl(1+4x_js\bigr)^{-3}ds=\frac1{8x_j\,(1+164\,x_j)^2}.
\]
So the tail is at least $\sum_jT_jm_j$, where $T_j$ is the right side at the largest value of $x_j$ on $Q$, and $m_j$ is a lower bound for $h_1$ on the box hull of $0$ and of an enclosure of the images of $Q$ under $\tau_j\psi_j^{41}$. The remaining part
\[
  F=\sum_{j=2}^5\sum_{k=1}^{40}\bigl(1+4kx_j\bigr)^{-3}\,h_1\circ\tau_j\psi_j^k-(1+\varepsilon)\,h_1
\]
satisfies $F\ge F(\lambda_Q)-\frac1{128}\sum_{i=1}^4\sup_Q\abs{\partial F/\partial\lambda_i}$ on $Q$, where $\lambda_Q$ is the centre of $Q$, by the mean value theorem. The value $F(\lambda_Q)$, enclosures of the four partial derivatives of $F$ over $Q$, computed by forward-mode automatic differentiation~\cite{GW08} in interval arithmetic, and the tail bounds are computed in binary64 arithmetic, every operation being followed by a step of one unit in the last place in the outward direction, so that every computed interval contains the exact result. On every cube the resulting lower bound for $\tilde\cL_3h_1-(1+\varepsilon)h_1$ is at least $1.9\times10^{-3}$; the computation takes four minutes on one core, and the program, with its input and its output, is deposited at Zenodo~\cite{DLdata}.

(b) Let $\delta>0$ and $t=3+\delta$. For $x\in C_1$ and $j\ge2$ we have $\frac12k\le4kx_j\le k$, so $0\le(1+4kx_j)^{-3}-(1+4kx_j)^{-t}\le\delta\,(1+\frac k2)^{-3}\log(1+k)$, and
\[
  \tilde\cL_th_1\ \ge\ \tilde\cL_3h_1-4\delta S\sup h_1,\qquad S=\sum_{k\ge1}\Bigl(1+\frac k2\Bigr)^{-3}\log(1+k)<0.7046 .
\]
By (a), $\tilde\cL_th_1\ge h_1$ on $C_1$ whenever $4\delta S\sup h_1\le\varepsilon\inf h_1$, which holds for $\delta=2.5\times10^{-4}$, since $4\cdot0.00025\cdot0.7046\cdot1.1945<0.000999\cdot0.8611$. Hence $\sN(\Gamma^{(4)})\ge3.00025$. The last statement follows from Theorem~\ref{thm:baragar}.
\end{proof}

\subsection*{The lattice case}

For comparison, let
\[
  B_j=2^{-1/3}\bigl(\id+\pr{\one-e_j}\one^*\bigr),\qquad j=1,2,3 .
\]
Since $\one^*\pr{\one-e_j}=1$, the matrix $\id+\pr{\one-e_j}\one^*$ has eigenvalues $1,1,2$, so $\det B_j=1$. Moreover $B_j\ge0$, and on the plane $\{\one^*x=1\}$ it acts by $\pr{B_jx}=\frac12\bigl(x+\pr{\one-e_j}\bigr)$, the homothety of ratio $\frac12$ centred at the vertex $\pr{\one-e_j}$ of $\cT$, which is also the fixed point of $\psi_j$. So the semigroup $\Gamma_{\mathcal S}$ generated by $B_1,B_2,B_3$ acts on $\cT$ by the iterated function system of the Sierpi\'nski gasket. Let $N_{\mathcal S}(r)=\#\{I:\ \norm{B_I}_1\le r\}$, and let $s_{\mathcal S}$ be the norm exponent of $\Gamma_{\mathcal S}$.

\begin{proposition}\label{prop:lattice}
    $\norm{B_I}_1=3\cdot4^{\abs I/3}$ for every $I\in\cI^3$. 
    Hence, $s_{\mathcal S}=\frac32\log_23$ and
\[
  \liminf_{r\to\infty}\frac{N_{\mathcal S}(r)}{r^{s_{\mathcal S}}}=\frac{3^{-s_{\mathcal S}}}2\ <\ \limsup_{r\to\infty}\frac{N_{\mathcal S}(r)}{r^{s_{\mathcal S}}}=\frac{3^{1-s_{\mathcal S}}}2 .
\]
\end{proposition}

\begin{proof}
Since $\one^*\pr{\one-e_j}=1$, we have $\one^*B_j=2^{2/3}\one^*$, so $\one^*B_I=4^{\abs I/3}\one^*$ and $\norm{B_I}_1=\one^*B_I\one=3\cdot4^{\abs I/3}$. There are $3^n$ multiindices of length $n$, so $\sum_I\norm{B_I}_1^{-t}=3^{-t}\sum_{n\ge0}(3\cdot4^{-t/3})^n$, which is finite if and only if $4^{t/3}>3$, that is, $t>\frac32\log_23$. For $r\ge3$ let $y=\frac32\log_2(r/3)$, so that $r=3\cdot4^{y/3}$ and $r^{s_{\mathcal S}}=3^{s_{\mathcal S}}\,3^y$, because $4^{s_{\mathcal S}/3}=3$. Then $N_{\mathcal S}(r)=\sum_{0\le n\le y}3^n=\frac12\bigl(3^{\lfloor y\rfloor+1}-1\bigr)$, and
\[
  \frac{N_{\mathcal S}(r)}{r^{s_{\mathcal S}}}=\frac{3^{1-s_{\mathcal S}}}2\,3^{-(y-\lfloor y\rfloor)}-\frac{3^{-s_{\mathcal S}-y}}2 .
\]
As $r\to\infty$ the last term tends to $0$, while $y-\lfloor y\rfloor$ takes every value in $[0,1)$ for arbitrarily large $r$.
\end{proof}

The same count, for the real Sierpi\'nski gasket and the maximal row-sum norm, is in~\cite[Motivational Example~3]{DL12}. Here the norms are exactly multiplicative, $\norm{B_{IJ}}_1=\frac13\norm{B_I}_1\norm{B_J}_1$, and each generator multiplies the norm by the constant $4^{1/3}$, so the logarithms of the analogues of the run factors lie in the lattice $\frac23\log2\cdot\Z$; the lattice case of the renewal theorem~\cite[Thm.~3.1(ii)]{KK16} gives a log-periodic function times $r^{s_{\mathcal S}}$, explicit here. Note that $\frac23s_{\mathcal S}=\log_23$ is the dimension of the Sierpi\'nski gasket, so the lower bound of Theorem~\ref{thm:cutout}(b) is an equality for it: its hypotheses hold, because $1<s_{\mathcal S}<3$, the zones are the triangles $B_I(Z_\emptyset)$, with the same central triangle $Z_\emptyset$ as for $\cN$ and area $\frac{\sqrt3}{32}4^{-\abs I}=\frac{27\sqrt3}{32}\norm{B_I}_1^{-3}$, and their sides lie in the gasket. The two gaskets live in the same triangle, and their generators fix the same vertices; what makes the difference is that $\psi_j$ is parabolic at its fixed vertex, so that the norms add within the runs, and the run factors are non-lattice~\cite[Prop.~42]{GMR19}. For instance, the spectral radii of $A_1A_2$ and $A_1A_2A_3$, which are the products of the run factors along the corresponding periodic sequences of runs, are $\frac{3+\sqrt5}2$ and $6.22226\ldots$, the largest root of $x^3-7x^2+5x-1$; this cubic is irreducible over $\mathbb Q$, and its roots are units, so no positive power of $6.22226\ldots$ lies in $\mathbb Q(\sqrt5)$, and the two logarithms are rationally independent.

\begin{remark}[dependence on the norm]\label{rem:norms}
The norm exponent does not depend on the norm, but norm equivalence only gives $N(r/C)\le N'(r)\le N(Cr)$ for some $C\ge1$, hence $N'(r)\asymp N(r)$ when $N$ grows polynomially, and the existence of $\lim N(r)/r^{s}$ may depend on the norm. For instance, for $M_0=\begin{pmatrix}2&0\\0&1\end{pmatrix}$ and $M_1=\begin{pmatrix}2&1\\0&1\end{pmatrix}$ (the maps $x\mapsto2x$ and $x\mapsto2x+1$ in homogeneous coordinates), the product along a multiindex of length $k$ is $\begin{pmatrix}2^k&n\\0&1\end{pmatrix}$, where $n$ runs once over $0,\dots,2^k-1$. For the entrywise $\ell^1$ norm the values $2^k+n+1$ are the integers $\ge2$, each attained once, so $N(r)=\lfloor r\rfloor-1$ for $r\ge2$; for the max-entry norm the norm is $2^k$, and $N(r)=2^{\lfloor\log_2r\rfloor+1}-1$, so that $N(r)/r$ oscillates between $1$ and $2$. Theorem~\ref{thm:count} is therefore a statement about $\norm\cdot_1$, or, more generally, about the linear functionals $L\mapsto\pr{w,L}$ of Theorem~\ref{thm:GMR}. The Euclidean norm of the souls is not of this form, and Theorem~\ref{thm:GMR} does not apply to it as it stands; we do not pursue this.
\end{remark}

\section{The three-dimensional Levitt--Yoccoz gasket}\label{sec:3d}

In this section we prove Theorem~\ref{thm:3d}. We describe the zones of $\cG$ and bound their volumes and surface areas (Lemmas~\ref{lem:oct} and~\ref{lem:geom3}), write the norm of $A_I$ as a product of run factors, certify the bounds for $\sN(\GG)$ as in Section~\ref{sec:thmB}, and apply Theorem~\ref{thm:cutout}(a).

Here $\Delta=\{x\in\R^4:\ x\ge0,\ \one^*x=1\}$ is the $3$-simplex, $\pr x=x/\one^*x$ for $x\ge0$, $x\neq0$, and $A_1,\dots,A_4$, $\GG$, $\cT_3$ and $\cG$ are as in the introduction.
Distances, areas and volumes are Euclidean, in the hyperplane $\{\one^*x=1\}$; in particular $\mathrm{vol}(\Delta)=\frac13$ and $\mathrm{vol}(\cT_3)=\frac1{27}\mathrm{vol}(\Delta)$.
A non-negative invertible matrix $g$ maps the tetrahedron with vertices $\pr{p_1},\dots,\pr{p_4}$, $p_m\ge0$, onto the tetrahedron with vertices $\pr{gp_1},\dots,\pr{gp_4}$, and maps convex subsets of it onto convex sets.
Since $A_jx=x+x_j(\one-e_j)$ and $\one^*A_jx=1+3x_j$ for $\one^*x=1$, the generators act on $\cT_3$ by
\[
  \psi_j(x)=\pr{A_jx}=\frac{x+x_j(\one-e_j)}{1+3x_j}=\pr{\one-e_j}+\frac{x-\pr{\one-e_j}}{1+3x_j},
\]
moving every point radially towards the vertex $\pr{\one-e_j}=\frac13(\one-e_j)$; we write $\psi_I(x)=\pr{A_Ix}$.
As in Corollary~\ref{cor: free}, Lemma~\ref{lem:oct} below shows that $\GG$ is free.
For $I\in\cI^4$ we write $a_m=A_I(\one-e_m)$, $m=1,\dots,4$, for the \textbf{frame} of $I$, so that $\psi_I(\cT_3)$ is the tetrahedron with vertices $\pr{a_1},\dots,\pr{a_4}$; we also write $\nu=\|A_I\|_1$ and $c_m=\one^*A_Ie_m$ for the column sums of $A_I$, so that $\norm{a_m}_1=\nu-c_m$ and $\sum_m\norm{a_m}_1=3\nu$.

\subsection*{Zones}
Let
\[
  \Oct=\bigl\{x\in \cT_3:\ x_m>\tfrac16,\ m=1,\dots,4\bigr\},
\]
a relatively open subset of $\cT_3$ whose closure is the octahedron with vertices $\pr{2\one-e_i-e_k}$, $i<k$, the midpoints of the edges of $\cT_3$; four of its faces lie on $\partial \cT_3$. Let $Z_I=A_I(\Oct)=\psi_I(\Oct)$; its closure is the octahedron with vertices $\pr{a_i+a_k}$, $i<k$.

\begin{lemma}\label{lem:oct}
The sets $Z_I$ are pairwise disjoint, and $\cT_3\setminus\cG=\bigcup_IZ_I$.
\end{lemma}

\begin{proof}
The $j$-th coordinate of $\psi_j(x)$ is $x_j/(1+3x_j)\le\frac16$, and $\psi_j$ maps the vertices $\pr{\one-e_m}$ of $\cT_3$ to $\pr{\one-e_j}$ for $m=j$ and to $\pr{2\one-e_j-e_m}$ for $m\ne j$. So $\psi_j(\cT_3)=\{x\in \cT_3:\ x_j\le\frac16\}$, the corner of $\cT_3$ at $\pr{\one-e_j}$ cut off by the plane through the midpoints of the three edges at $\pr{\one-e_j}$. Two corners meet at most in a vertex of $\Oct$, and $\cT_3$ is the union of $\Oct$ and of the four corners, which do not meet $\Oct$. Each $\psi_I$ is injective on $\cT_3$, so, as in Lemma~\ref{lem:zones}, by induction
\[
  \cT_3=\bigsqcup_{\abs I<n}Z_I\ \sqcup\ \bigcup_{\abs I=n}\psi_I(\cT_3)\qquad(n\ge0),
\]
and the claim follows because $\cG=\bigcap_n\bigcup_{\abs I=n}\psi_I(\cT_3)$.
\end{proof}

\begin{lemma}\label{lem:geom3}
For every $I\in\cI^4$:
\begin{enumerate}
\item[(i)] $\nu\le\norm{a_i+a_k}_1\le2\nu$ for all $i\neq k$;
\item[(ii)] $\norm{\pr{a_i}-\pr{a_k}}_1\le 2/\max_m\norm{a_m}_1$ for all $i,k$; hence
$\operatorname{diam}\psi_I(\cT_3)\le\sqrt2/\max_m\norm{a_m}_1\le\frac{4\sqrt2}3\,\norm{A_I}_1^{-1}$;
\item[(iii)]
\[
  \frac{\mathrm{vol}(Z_I)}{\mathrm{vol}(\Delta)}=\frac2\nu\Bigl(\sum_{\{i,k,l\}}\frac1{(2\nu-c_i-c_k)(2\nu-c_i-c_l)(2\nu-c_k-c_l)}+\sum_{j=1}^4\ \prod_{m\neq j}\frac1{2\nu-c_j-c_m}\Bigr),
\]
  the first sum running over the four $3$-element subsets of $\{1,2,3,4\}$; hence
  \[
    2\,\nu^{-4}\le\mathrm{vol}(Z_I)/\mathrm{vol}(\Delta)\le16\,\nu^{-4};
  \]
\item[(iv)] $\mathrm{area}(\partial\overline{Z_I})\le\mathrm{area}(\partial\psi_I(\cT_3))\le\frac{32\sqrt3}9\,\nu^{-2}$.
\end{enumerate}
\end{lemma}

\begin{proof}
(i) $\norm{a_i+a_k}_1=\one^*A_I(2\,\one-e_i-e_k)=2\nu-c_i-c_k=\nu+c_l+c_m$, where $\{i,k,l,m\}=\{1,2,3,4\}$, and $0\le c_l+c_m\le\nu$.

(ii) Since $A_j(\one-e_j)=\one-e_j$ and $A_j(\one-e_m)=(\one-e_j)+(\one-e_m)$ for $m\ne j$, appending $j$ to $I$ keeps $a_j$ and replaces $a_m$ by $a_m+a_j$ for $m\ne j$, as in~\eqref{eq:tree}. Fix a row vector $f\in[0,1]^4$ and let $t_m=fa_m/\norm{a_m}_1\in[0,1]$, the value of $f$ at the vertex $\pr{a_m}$ of $\psi_I(\cT_3)$. We show that $\max_mt_m-\min_mt_m\le1/\max_m\norm{a_m}_1$. At $I=\emptyset$, $a_m=\one-e_m$, $\norm{a_m}_1=3$ and $t_m=\frac13(f\one-f_m)$, whose spread is at most $\frac13$. Appending $j$, $t_j$ is unchanged and
\[
  t'_m=t_j+c'_m(t_m-t_j),\qquad c'_m=\frac{\norm{a_m}_1}{\norm{a_m}_1+\norm{a_j}_1}\le c'=\frac{b}{b+\norm{a_j}_1},\quad b=\max_{m\neq j}\norm{a_m}_1 .
\]
Since $\min_mt_m\le t_j\le\max_mt_m$, all the new values lie in $t_j+c'\,[\min_mt_m-t_j,\ \max_mt_m-t_j]$, so their spread is at most $c'\,(\max_mt_m-\min_mt_m)\le b/\bigl((b+\norm{a_j}_1)\max_m\norm{a_m}_1\bigr)\le1/(b+\norm{a_j}_1)$, and $b+\norm{a_j}_1$ is the largest new value of $\norm{a_m}_1$. For $u,v\in\Delta$ we have $\norm{u-v}_1=2\max_{f\in\{0,1\}^4}f(u-v)$, which gives the first claim, and, since $\one^*(u-v)=0$, $\abs{u-v}\le\norm{u-v}_1/\sqrt2$; the diameter of the tetrahedron $\psi_I(\cT_3)$ is the largest distance between two of its vertices. Finally $\max_m\norm{a_m}_1\ge\frac14\sum_m\norm{a_m}_1=\frac34\nu$.

(iii) For non-negative non-zero $p_1,\dots,p_4\in\R^4$, the tetrahedron with vertices $\pr{p_1},\dots,\pr{p_4}$ has volume
\[
  \mathrm{vol}(\Delta)\,\frac{\abs{\det(p_1,p_2,p_3,p_4)}}{\one^*p_1\,\one^*p_2\,\one^*p_3\,\one^*p_4}
\]
(cf.~\cite[Lemma~4.1]{PS24}). Indeed, the $4$-simplex with vertices $0$ and $\pr{p_1},\dots,\pr{p_4}$ has $4$-volume $\frac1{4!}\abs{\det(\pr{p_1},\dots,\pr{p_4})}$, and it is a cone of height $\frac12$ (the distance from $0$ to $\{\one^*x=1\}$) over the tetrahedron; moreover $\mathrm{vol}(\Delta)=\frac13$. The closure of $\Oct$ is the union of the eight tetrahedra with apex its centre $\pr{\one}$ over its faces: the four faces on $\partial \cT_3$, with vertices $\pr{2\one-e_i-e_k},\pr{2\one-e_i-e_l},\pr{2\one-e_k-e_l}$, and the four faces facing the corners, with vertices $\pr{2\one-e_j-e_m}$, $m\neq j$. The map $\psi_I$ carries this decomposition to a decomposition of $\overline{Z_I}$ into eight tetrahedra with apex $\pr{A_I\one}$. Now $\abs{\det(2\one-e_i-e_k,2\one-e_i-e_l,2\one-e_k-e_l,\one)}=\abs{\det(2\one-e_j-e_i,2\one-e_j-e_k,2\one-e_j-e_l,\one)}=2$, $\det A_I=1$, $\one^*A_I(2\one-e_i-e_k)=2\nu-c_i-c_k$ and $\one^*A_I\one=\nu$, which gives the formula. By (i), every factor $2\nu-c_i-c_k$ lies in $[\nu,2\nu]$, so each of the eight terms lies between $\frac14\nu^{-4}$ and $2\nu^{-4}$.

(iv) $\overline{Z_I}$ is a convex subset of the tetrahedron $\psi_I(\cT_3)$, and the surface area of convex bodies is monotone under inclusion (by Cauchy's surface area formula~\cite[(5.73)]{Sch14}). Each face of $\psi_I(\cT_3)$ is a triangle with sides at most $D=\mathrm{diam}\,\psi_I(\cT_3)$, so its area is at most $\frac{\sqrt3}4D^2$ (the altitude onto a longest side $c$ is at most $\frac{\sqrt3}2c$); by (ii), $\mathrm{area}(\partial\psi_I(\cT_3))\le\sqrt3\,D^2\le\frac{32\sqrt3}9\,\nu^{-2}$.
\end{proof}

\begin{remark}
The proof of (ii) only uses the rule $a_m\mapsto a_m+a_j$, $m\neq j$, and holds for the semigroup generated by the matrices $\id+(\one-e_j)e_j^*$ in any dimension.
For $\GN$ it gives $\operatorname{diam}\psi_I(\cT)\le\sqrt2/\max_k\norm{a_{I,k}}_1
\le\frac3{\sqrt2}\,\norm{A_I}_1^{-1}$, because the three frame vectors have total $\ell^1$-norm $2\norm{A_I}_1$ (Section~\ref{sec:souls}), which improves the constant $6\sqrt3$ of
Lemma~\ref{lem:geom}(iv); the bound $\sqrt2/\max_k\norm{a_{I,k}}_1$ is attained at $I=\emptyset$.
\end{remark}

\subsection*{The run factors}
As in Section~\ref{sec:fast}, let $x_0=\pr{\one}$. Since $A_j^k=\id+k\,(\one-e_j)e_j^*$, for $x\in \cT_3$
\[
  A_j^kx=x+kx_j\,(\one-e_j),\qquad\norm{A_j^kx}_1=1+3kx_j,\qquad
  \psi_j^k(x)=(1-\beta)\,x+\beta\,\pr{\one-e_j},
\]
with $\beta=3kx_j/(1+3kx_j)$, exactly as in \eqref{eq:run}. Hence, for $I=j_n^{k_n}\cdots j_1^{k_1}$ written in runs from the right,
\begin{equation}\label{eq:cocycle3}
  \|A_I\|_1=4\prod_{r=1}^n\Bigl(1+3k_r\,x^{(r-1)}_{j_r}\Bigr),\qquad
  x^{(0)}=x_0,\quad x^{(r)}=\psi_{j_r}^{k_r}\bigl(x^{(r-1)}\bigr).
\end{equation}

\begin{lemma}\label{lem:trap3}
Let $C_i=\psi_i(\cT_3)=\{x\in \cT_3:\ x_i\le\frac16\}$, the tetrahedron with vertices $\pr{\one-e_i}$ and $\pr{2\one-e_i-e_m}$, $m\neq i$. Then $\psi_j^k(\cT_3)\subset C_j$ for every $j$ and $k\ge1$, and $3x_j\ge\frac12$ for all $x\in C_i$ and $j\neq i$.
\end{lemma}

\begin{proof}
$\psi_j^k(\cT_3)\subset\psi_j(\cT_3)$. If $x\in C_i$ and $j\neq i$, then
$x_j=1-x_i-\sum_{m\ne i,j}x_m\ge1-\frac16-\frac23=\frac16$.
\end{proof}

Consequently, with $\mathcal C=\bigsqcup_{i=1}^4C_i\times\{i\}$ and, for $t>1$ and bounded $h$ on $\mathcal C$,
\[
  (\mathcal L_th)(x,i)=\sum_{j\neq i}\sum_{k\ge1}\bigl(1+3kx_j\bigr)^{-t}\,h\bigl(\psi_j^k(x),j\bigr),
\]
the renewal \eqref{eq:renewal}, which holds verbatim, and $\norm{A_j^k\one}_1=4+3k$ give
\[
  \sum_I\|A_I\|_1^{-t}=4^{-t}\Bigl(1+\sum_{j=1}^4\sum_{k\ge1}\bigl(1+\tfrac{3k}4\bigr)^{-t}
  \sum_{m\ge0}(\mathcal L_t^m\mathds{1})\bigl(\psi_j^k(x_0),j\bigr)\Bigr),
\]
and Corollary~\ref{cor:test} holds verbatim, with $\sN(\GG)$ in place of $\sN$.

\subsection*{The bounds for \texorpdfstring{$\sN(\GG)$}{the norm exponent}}\label{sec:3dcert}
By symmetry the test function can be taken on $C_1$ only. Let $\tau_i$ be the transposition of the coordinates $1$ and $i$ ($\tau_1=\id$), so that $\tau_i(C_1)=C_i$; recall that $\pi\circ\psi_j^k=\psi_{\pi(j)}^k\circ\pi$ for every permutation $\pi$ of the coordinates (Section~\ref{sec:souls}). If $h_1$ is a function on $C_1$ invariant under the permutations of $x_2,x_3,x_4$ and $h(x,i)=h_1(\tau_ix)$, then $(\mathcal L_th)(\tau_ix,i)=(\tilde{\mathcal L}_th_1)(x)$ for $x\in C_1$, where
\[
  (\tilde{\mathcal L}_th_1)(x)=\sum_{j=2}^4\sum_{k\ge1}\bigl(1+3kx_j\bigr)^{-t}\,h_1\bigl(\tau_j\psi_j^k(x)\bigr).
\]
Indeed, writing $j=\tau_i(j')$ in the sum defining $(\mathcal L_th)(\tau_ix,i)$, the term of $(j',k)$ is
\[
  \bigl(1+3kx_{j'}\bigr)^{-t}\,h_1\bigl(\tau_{\tau_i(j')}\,\tau_i\,\psi_{j'}^k(x)\bigr),
\]
and the permutation $\tau_{\tau_i(j')}\tau_i\tau_{j'}$ fixes the coordinate~$1$. So it suffices to find such an $h_1$ with $\tilde{\mathcal L}_th_1\ge h_1$, or with $\tilde{\mathcal L}_th_1\le\theta h_1$ for some $\theta<1$, on $C_1$.

We use the coordinates $\lambda=(2-6x_2,2-6x_3,2-6x_4)$ of $C_1$, in which $C_1=\{\lambda\ge0,\ \lambda_1+\lambda_2+\lambda_3\le1\}$, $x_1=\frac16(\lambda_1+\lambda_2+\lambda_3)$ and $x_m=\frac16(2-\lambda_{m-1})$ for $m\ge2$, and the cubes $Q_c=\frac1M\bigl(c+[0,1]^3\bigr)$, $c\in\Z^3_{\ge0}$, $c_1+c_2+c_3\le M-1$, which cover $C_1$ and are contained in $\cT_3$. We take $h_1=H_c$ on $Q_c$ (the smallest value at points lying in several cubes), with $H$ invariant under the permutations of $c$. The two inequalities of Section~\ref{sec:thmB} hold verbatim, the sums running over $j=2,3,4$, with
\[
  \bar u_j=\max_{Q_c}3x_j=1-\frac{c_{j-1}}{2M},\qquad \underline u_j=\min_{Q_c}3x_j=1-\frac{c_{j-1}+1}{2M}\qquad(j=2,3,4),
\]
and with $\underline H_{j,k}(Q_c)$, $\overline H_{j,k}(Q_c)$ the minimum and the maximum of $H$ over the cubes meeting the closed bounding box, in the coordinates $\lambda$, of the images under $\tau_j\psi_j^k$ of the eight vertices of $Q_c$; for the tail, of the images under $\tau_j\psi_j^{K+1}$ of the vertices and of $\tau_j\pr{\one-e_j}=\pr{\one-e_1}$. These cubes include every cube that meets $\tau_j\psi_j^k(Q_c)$, because $\psi_j^k$ is projective with positive denominator $1+3kx_j$ on $Q_c$. At a vertex $\lambda=c/M$ the vector $6M\,x$ has integer entries, so the coordinates $\lambda$ of $\tau_j\psi_j^k(x)$ are quotients of explicit integers, and the cubes are selected exactly, in integer arithmetic. As in Section~\ref{sec:thmB}, the program stores the points of $\cT_3$ as row vectors.

With $M=48$ ($19\,600$ cubes) and $K=100$, and with $H^-$, $H^+$ the fixed vectors of the lower and the upper discretized operators, computed by nonlinear power iteration as in Section~\ref{sec:thmB} and made exactly symmetric, the lower bound for $\tilde{\mathcal L}_th_1$ on $Q_c$ is at least $1.0015\,H_c$ for every $c$ at $t=2.74$, and the upper bound is at most $0.9982\,H_c$ for every $c$ at $t=2.777$. The weights, the tail integrals and the sums are evaluated in ball arithmetic with Arb~\cite{Arb}, for the exact rational values of $t$. Hence
\[
  2.74\le\sN(\GG)\le2.777 .
\]
The code is deposited at Zenodo~\cite{DLdata}. These bounds improve Baragar's bounds $2.730<\sN(\GG)<2.798$~\cite[Thm.~0.1]{Bar98} (Section~\ref{sec:MH}). They are consistent with the value $\sN(\GG)=2.75784$ given by a spectral collocation discretization of $\tilde{\mathcal L}_t$, and with the exact count $N_\cG(r)=\#\{I:\ \|A_I\|_1\le r\}=N_3(r)$, which is asymptotic to $c_3\,r^{\sN(\GG)}$ by Theorem~\ref{thm:count}: $N_\cG(2^{15})=36\,890\,132\,081$, and the local slopes $\log_2\bigl(N_\cG(2r)/N_\cG(r)\bigr)$ at $r=2^{13}$ and $2^{14}$ are $2.7574$ and $2.7576$.

\thmthreed*

\begin{proof}
The bounds for $\sN(\GG)$ were proved above.

\textbf{Lower bound.} The face $F=\{x\in \cT_3:\ x_4=\frac13\}$ of $\cT_3$, opposite to $\pr{\one-e_4}$, is mapped into itself by $\psi_1,\psi_2,\psi_3$: for $j\le3$ the fourth coordinate of $\psi_j(x)$ is $(x_4+x_j)/(1+3x_j)=\frac13$. In the coordinates $u=\frac32(x_1,x_2,x_3)$, which identify $F$ with the triangle $\cT$ of Section~\ref{sec:souls} through a similarity of ratio $\frac32$, these maps are $u\mapsto(u+u_j(\one-e_j))/(1+2u_j)$, the generators of $\GN$. The limit set of $\{\psi_1,\psi_2,\psi_3\}$ on $F$ is contained in $\cG$, so $\cG$ contains a copy of $\cN$, similar to it, and $\dimH\cG\ge\dimH\cN>1.6853$ by Theorem~\ref{thm: dimH}.

\textbf{Upper bound.} By Theorem~\ref{thm:KM}, $\cG=\mathcal G^3$ has measure zero, and, by Lemma~\ref{lem:oct}, $\cG=\cT_3\setminus\bigcup_IZ_I$, where the $Z_I$ are pairwise disjoint and convex, with non-empty interior; they need not be open (by Remark~\ref{rem:octfaces}, $Z_1$ contains points of the closure of $Z_\emptyset$), and Theorem~\ref{thm:cutout} does not require them to be. By Lemma~\ref{lem:geom3}(iii),(iv), $\mathrm{vol}(Z_I)\le\frac{16}3\,\nu_I^{-4}$ and $\mathrm{area}(\partial\overline{Z_I})\le\frac{32\sqrt3}9\,\nu_I^{-2}$. Theorem~\ref{thm:cutout}(a) with $n=3$, applied to the $3$-simplex $\cT_3$, gives $\udimB\cG\le1+\frac12\sN(\GG)\le1+\frac12\cdot2.777=2.3885$, and $\dimH\cG\le\udimB\cG$.
\end{proof}

\begin{remark}[Theorem~\ref{thm:cutout}(b) and $\cG$]\label{rem:octfaces}
The zones of $\cG$ are not separated by $\cG$. For instance, the closures of $Z_\emptyset=\Oct$ and $Z_1=\psi_1(\Oct)$ meet along a triangle in the plane $\{x_1=\frac16\}$: the point $(\frac13,\frac29,\frac29,\frac29)$ lies in $\Oct$, and $\psi_1$ maps it to $(\frac16,\frac5{18},\frac5{18},\frac5{18})$, a point of $Z_1$ in the closure of $\Oct$. So $\partial\overline{Z_I}\not\subset\cG$, and the collars of the zones do not lie in $\cG_\varepsilon$: near most of its boundary a zone is adjacent to other zones, and $\cG$ meets the faces of the octahedra $\overline{Z_I}$ only in copies of $\cN$. A lower bound $\udimB\cG\ge\frac34\sN(\GG)\ge2.055$, or the corresponding bound for $\dimH\cG$, would need a different argument, or an analogue in $\RP^3$ of the equality $\dimH=\sA$ of \cite{JLPX23,Jur26}.
\end{remark}

\begin{remark}[numerical values]\label{rem:3dnum}
Summing $\varphi^s(A_I)$ over all $I$ with $\|A_I\|_1\le4096$ in half-octave shells of $\sigma_1(A_I)$, and solving for the $s$ at which consecutive shell sums are equal, gives $s=2.263$, $2.256$, $2.263$, $2.261$, $2.260$ on the shells from $\sigma_1\approx128$ to $\sigma_1\approx724$ (these shells are complete, since $\sigma_1(A_I)\ge\|A_I\|_1/4$). A box count of $\cG$ at scales $2^{-k}$ gives the local slopes $2.206$, $2.211$, $2.211$ for $k=9,10,11$; the same code gives $1.686$ for $\cN$ and $2.000$ for the Sierpi\'nski tetrahedron. So the shells give $\sA(\GG)\approx2.26$ and the box count gives $\approx2.21$. Both values lie above the lower bound $\sA(\GG)\ge\frac34\sN(\GG)\ge2.055$ for the affinity exponent and below the upper bound $2.3885$ of Theorem~\ref{thm:3d}, but they differ by about $0.05$, several times more than the variation of each of them across scales. In $\RP^2$ the same box-counting code agrees with $\sA(\GN)$ to within $10^{-3}$; in $\RP^3$ neither the equality $\dimH\cG=\sA(\GG)$ nor the existence of $\dimB\cG$ is known, and the shells above reach only $\sigma_1\approx724$. We do not know whether the gap reflects a difference between $\sA(\GG)$ and the dimensions of $\cG$, or a slow convergence of one of the two estimates. As in Remark~\ref{rem:172}, the cut-out argument loses in the bound $O(\nu^{-2})$ on the surface areas, which is attained by the cusp zones: for the elements that dominate $\sum_I\varphi^s(A_I)$, $-\log(\sigma_m/\sigma_1)/\log\sigma_1\approx1.05$, $1.20$, $1.74$ for $m=2,3,4$, so the corresponding zones are thin plates, whose surface area is much smaller than $\nu_I^{-2}$.
\end{remark}

\section{The diffusion of chaotic sections}\label{sec:AHS}

In~\cite{AHSinv}, Avila, Hubert and Skripchenko construct an invariant probability measure
$\mu_0$ on the Rauzy gasket and show that, for $\mu_0$-almost every chaotic regime, the
trajectories of the plane sections diffuse with exponent
\[
  \delta=\limsup_{t\to\infty}\frac{\log d(x,x_t)}{\log t}=-\frac{\lambda_3}{\lambda_1}
  \in\Bigl(\frac12,1\Bigr)
\]
\cite[Thm.~1 and Thm.~30]{AHSinv} (in the statement of Thm.~30 there, the numerator should read $\log d(x,x_t)$, as in Thm.~1).
In their notations, $\lambda_1>0>\lambda_2>\lambda_3$ are the Lyapunov exponents of their cocycle and $\mu_0$ is the Gibbs measure of the potential $-\kappa_0r$, where $r$ is the roof function of their renormalization flow and $\kappa_0$ is the zero of the pressure of $-\kappa r$.
By~\cite[Thm.~25]{AHSinv}, the measure induced by $\mu_0$ on the flow is its unique measure of maximal entropy (see Remark~\ref{rem:maxentflow}).
Notice that their parameter $\kappa$ plays the role of $t$ in Section~\ref{sec:fast} and has nothing to do with our parameter $\kappa$ in Sections~\ref{sec:cocycle}--\ref{sec:comp}.

\subsection*{Filling a gap in Lemma~24 of~\cite{AHSinv}}
\label{rem:AHS}
Our Theorem~\ref{cor:diffusion} rests on~\cite[Thm.~30]{AHSinv}, which holds for $\mu_0$-almost every point.
The existence of $\kappa_0$ is claimed in~\cite[Lemma~24]{AHSinv}, whose argument needs the pressure to be finite and continuous on an interval on which it changes sign.
The proof of~\cite[Thm.~30]{AHSinv} uses the Oseledets theorem \cite[Thm.~27]{AHSinv}, and so the log-integrability of the cocycle for $\mu_0$ \cite[Lemma~26]{AHSinv}, whose proof derives $\int r\,d\mu_0<\infty$ from $\norm{L_{-\kappa_0r}1}_\infty<\infty$ through \cite[Thm.~20]{AHSinv}.
However, $\norm{L_{-\kappa r}1}_\infty<\infty$ is proved only for $\kappa\ge3$ \cite[Thm.~15 and Thm.~19]{AHSinv}, where the pressure is negative, while $\kappa_0<3$ by Proposition~\ref{prop:AHS}, which shows that $\kappa_0=\sN$, and the proof of Theorem~\ref{thm: sN}, which shows that $\sN<3$.

In order to fix this gap we first notice that, for the accelerated induction, before the passage to the first returns, the corresponding facts were established later by Gamburd, Magee and Ronan~\cite[Thm.~10, Thm.~39 and Prop.~43]{GMR19}, in their study of the Markoff--Hurwitz equations: for $t>1$ the transfer operator of the accelerated induction, with weights the run factors raised to the power $-t$ (the operator of Remark~\ref{rem:GMRproof} for $d=2$), is bounded on $C^1$ functions and has a simple leading eigenvalue $\lambda_t$, which is real-analytic and strictly decreasing in $t$ and equals $1$ at $t=\sN$ (Theorem~\ref{thm:count}), where there is a conformal measure. They also remark that their conformal measure is essentially the measure $\mu_0$ of~\cite{AHSinv}~\cite[Remark~12]{GMR19}. What remains to be done is the passage to the first-return system of~\cite{AHSinv}, on which $\kappa_0$ and $\mu_0$ are defined, together with the integrability that~\cite[Lemma~26]{AHSinv} needs.

We do this in three steps. First, we introduce the shift of the runs and its Gibbs state $\mu_*$, which plays the role of the measure of~\cite{GMR19}: on the shift of the runs, the facts above are Lemma~\ref{lem:pressure}, with $P(-t\xi)=p(t)$ finite for $t>1$ and vanishing only at $\sN$, and Theorem~\ref{thm:RPF}. Then we bound the Lyapunov exponents of the shift of the runs for $\mu_*$, without using~\cite{AHSinv} (Proposition~\ref{prop:lyap}). Finally, we show that $\mu_0$ is the measure induced by $\mu_*$ on the cross-section given by a complete path (Proposition~\ref{prop:AHS}): (i) gives $\kappa_0=\sN$, by comparing the pressure of $-\kappa\rho$ with that of $-\sN\rho$, which vanishes by the Gibbs property of $\mu_*^Y$; (ii) gives $\norm{L_{-\sN r}1}_\infty<\infty$ and the Gibbs measure $\mu_0$, and, granted the two facts on suspension flows of Remark~\ref{rem:maxentflow}, recovers \cite[Thm.~25]{AHSinv}; and (iii) gives the conclusion of \cite[Lemma~26]{AHSinv}.
The Gibbs property of $\mu_0$ in~(ii) is also the input of the bounded distortion property~\cite[Lemma~36]{AHSinv}, on which the simplicity of the Lyapunov spectrum, and so $\delta>\frac12$, rests; Theorem~\ref{cor:diffusion} does not need it.

\subsection*{The symbolic model}
Unlike in~\eqref{eq:cocycle}, where the runs are numbered from the right, in the symbolic model they are numbered from the left: $\omega_1$ is the leftmost run, the last to act.
\emph{Shift spaces.}
We follow~\cite[\S2.1]{MU03} and~\cite[\S2.1]{KK16}, with a different notation. Let $\mathcal E$ be a countable set, the \emph{alphabet}, and let $\mathcal R\subset\mathcal E\times\mathcal E$. The \emph{shift space} of $(\mathcal E,\mathcal R)$ is
\[
  \Sigma=\bigl\{\omega=\omega_1\omega_2\cdots\in\mathcal E^{\N}:\ (\omega_r,\omega_{r+1})\in\mathcal R\text{ for all }r\ge1\bigr\},
\]
with the product of the discrete topologies, its Borel $\sigma$-algebra, and the shift $\operatorname{sh}(\omega_1\omega_2\cdots)=\omega_2\omega_3\cdots$. A word $w_1\cdots w_n$ in $\mathcal E$ is \emph{admissible} if $(w_r,w_{r+1})\in\mathcal R$ for $r<n$, and $[w_1\cdots w_n]$ is the \emph{cylinder} of the $\omega\in\Sigma$ with $\omega_r=w_r$ for $r\le n$. The shift space $\Sigma$ is \emph{finitely irreducible} if there is a finite set $\Lambda$ of admissible words such that, for all $a,b\in\mathcal E$, the word $a\lambda b$ is admissible for some $\lambda\in\Lambda$. If $\mathcal R=\mathcal E\times\mathcal E$, then $\Sigma=\mathcal E^{\N}$ is the \emph{full shift} on $\mathcal E$, which is finitely irreducible, with $\Lambda=\{e\}$ for any $e\in\mathcal E$.

\emph{The shift of the runs.}
Let $\mathcal A=\{(j,k):\ j\in\{1,2,3\},\ k\ge1\}$ be the set of runs, with $\iota(j,k)=j$, $A_{(j,k)}=A_j^k$ and $\psi_{(j,k)}=\psi_j^k$, and let $\Omega$ be the shift space of $\mathcal A$ and of $\mathcal R=\{(a,b):\ \iota(a)\ne\iota(b)\}$:
\[
  \Omega=\bigl\{\omega=\omega_1\omega_2\cdots\in\mathcal A^{\N}:\ \iota(\omega_{r+1})\ne\iota(\omega_r)\text{ for all }r\ge1\bigr\},\qquad
  \operatorname{sh}(\omega_1\omega_2\cdots)=\omega_2\omega_3\cdots .
\]
It is finitely irreducible: for all $a,b\in\mathcal A$ there is $e\in\{(1,1),(2,1),(3,1)\}$ such that $aeb$ is admissible. We identify an admissible word $u=\omega_1\cdots\omega_n$ with the multiindex whose runs, read from the left, are $\omega_1,\dots,\omega_n$, so that $A_u=A_{\omega_1}\cdots A_{\omega_n}$ and $\psi_u=\psi_{\omega_1}\circ\cdots\circ\psi_{\omega_n}$. This is a bijection between the admissible words of length $n$ and the multiindices with $n$ runs; $\omega_1$ is the leftmost run, the last to act, and the last index of $u$ is $\iota(\omega_n)$.

\emph{Functions on shift spaces.}
We use the following classes of functions.

\begin{definition}[{\cite[\S\S2.1--2.2]{MU03}, \cite[\S\S2.1--2.2]{KK16}}]\label{def:shiftfns}
Let $\Sigma$ be a shift space, and let $u:\Sigma\to\R$ be continuous, with Birkhoff sums $u_n=\sum_{r=0}^{n-1}u\circ\operatorname{sh}^r$, $u_0=0$.
\begin{enumerate}
\item[(i)] For $n\ge1$, $\mathrm{var}_nu=\sup\{\abs{u(\omega)-u(\omega')}:\ \omega,\omega'\in\Sigma$ have the same first $n$ terms$\}$. For $0<\theta<1$, $\cF_\theta(\Sigma)$ is the space of the continuous functions $u$ with $\sup_{n\ge1}\theta^{-n}\mathrm{var}_nu<\infty$, and $\cF^b_\theta(\Sigma)$ is the space of the bounded functions in $\cF_\theta(\Sigma)$. The function $u$ is \emph{H\"older continuous} if $u\in\cF_\theta(\Sigma)$ for some $\theta\in(0,1)$. Nothing is required of pairs with $\omega_1\ne\omega'_1$, so a H\"older continuous function may be unbounded.
\item[(ii)] $u$ is \emph{summable} if $\sum_{e\in\mathcal E}\exp\sup_{[e]}u<\infty$. The \emph{pressure} of $u$ is
\[
  P(u)=\lim_{n\to\infty}\frac1n\log Z_n(u),\qquad Z_n(u)=\sum_{\abs w=n}\exp\sup_{[w]}u_n\ \in(0,\infty],
\]
where the sum runs over the admissible words of length $n$. The limit exists in $[-\infty,\infty]$, because $Z_{m+n}(u)\le Z_m(u)\,Z_n(u)$~\cite[Lemmas~2.1.1 and~2.1.2, and (2.1)]{MU03}.
\item[(iii)] A Borel probability measure $m$ on $\Sigma$ is a \emph{Gibbs state} of $u$ if there are constants $C\ge1$ and $P_m\in\R$ such that
\[
  C^{-1}\ \le\ \frac{m([\omega_1\cdots\omega_n])}{\exp\bigl(u_n(\omega)-nP_m\bigr)}\ \le\ C\qquad(\omega\in\Sigma,\ n\ge1);
\]
it is an \emph{invariant Gibbs state} if, in addition, it is shift-invariant.
\end{enumerate}
\end{definition}

\begin{lemma}\label{lem:coding}
Let $\omega\in\Omega$ and $n\ge1$. The cylinder $\psi_{\omega_1\cdots\omega_n}(\cT)$ is contained in $C_{\iota(\omega_1)}$ and has diameter at most $\frac{\sqrt2}4\theta_*^{n-1}$. Hence $\bigcap_n\psi_{\omega_1\cdots\omega_n}(\cT)$ is a single point $x_\omega\in C_{\iota(\omega_1)}$, $x_\omega=\psi_{\omega_1}(x_{\operatorname{sh}\omega})$, and $\abs{x_\omega-x_{\omega'}}\le\frac{\sqrt2}4\theta_*^{n-1}$ for every $\omega'\in[\omega_1\cdots\omega_n]$.
\end{lemma}

\begin{proof}
Since $\psi_j^k(\cT)\subset C_j$ (Lemma~\ref{lem:trap}), the cylinder lies in $\psi_{\omega_1}(\cT)\subset C_{\iota(\omega_1)}$. The triangle $C_i$ is equilateral, with vertices $\frac12(\one-e_i)$, $\frac14(\one+e_j)$, $\frac14(\one+e_l)$, where $j$ and $l$ are the other two indices, and side $\frac{\sqrt2}4$. For $n\ge2$ the cylinder $\psi_{\omega_2\cdots\omega_n}(\cT)$ lies in $C_{\iota(\omega_2)}$, and $\iota(\omega_2)\ne\iota(\omega_1)$; so, by Lemma~\ref{lem:contr}, the diameter of $\psi_{\omega_1\cdots\omega_n}(\cT)=\psi_{\omega_1}\bigl(\psi_{\omega_2\cdots\omega_n}(\cT)\bigr)$ is at most $\theta_*$ times that of $\psi_{\omega_2\cdots\omega_n}(\cT)$, and, by induction, at most $\theta_*^{n-1}\diam C_{\iota(\omega_n)}$. The cylinders are compact and decrease with $n$, because $\psi_{\omega_{n+1}}(\cT)\subset\cT$, so their intersection is a single point $x_\omega$. Moreover $\psi_{\omega_1}(x_{\operatorname{sh}\omega})$ lies in $\psi_{\omega_1}\bigl(\psi_{\omega_2\cdots\omega_n}(\cT)\bigr)=\psi_{\omega_1\cdots\omega_n}(\cT)$ for every $n$, hence equals $x_\omega$; and if $\omega'\in[\omega_1\cdots\omega_n]$, then $x_\omega$ and $x_{\omega'}$ lie in the same cylinder $\psi_{\omega_1\cdots\omega_n}(\cT)$.
\end{proof}

Thus $x_\omega$ is the limit of the soul points $\pr{L_{\omega_1\cdots\omega_n}}=\psi_{\omega_1\cdots\omega_n}(x_0)$, and $\omega\mapsto x_\omega$ is continuous. Define $\xi:\Omega\to\R$ by
\[
  \xi(\omega)=\log\norm{A_{\omega_1}x_{\operatorname{sh}\omega}}_1=\log\bigl(1+2k\,(x_{\operatorname{sh}\omega})_j\bigr)\qquad\bigl(\omega_1=(j,k)\bigr),
\]
let $\xi_n=\sum_{r=0}^{n-1}\xi\circ\operatorname{sh}^r$ be its Birkhoff sums, with $\xi_0=0$, and let $\mathrm{var}_n\xi$ be as in Definition~\ref{def:shiftfns}(i).

\begin{lemma}\label{lem:xi}
Let $\omega\in\Omega$.
\begin{enumerate}
\item[(i)] If $\omega_1=(j,k)$, then $\log(1+\frac k2)\le\xi(\omega)\le\log(1+k)$; in particular $\xi\ge\log\frac32$.
\item[(ii)] $\xi_n(\omega)=\log\norm{A_{\omega_1\cdots\omega_n}x_{\operatorname{sh}^n\omega}}_1$ for every $n\ge0$.
\item[(iii)] $\mathrm{var}_1\xi\le\log2$, and $\mathrm{var}_n\xi\le\sqrt2\,\theta_*^{n-2}$ for $n\ge2$. In particular $\xi$ is continuous and $\sup_{n\ge1}\theta_*^{-n}\mathrm{var}_n\xi<\infty$, that is, $\xi\in\cF_{\theta_*}(\Omega)$ (Definition~\ref{def:shiftfns}(i)).
\end{enumerate}
\end{lemma}

\begin{proof}
(i) The point $x_{\operatorname{sh}\omega}$ lies in $C_{\iota(\omega_2)}$, and $\iota(\omega_2)\ne j$, so $\frac12\le2(x_{\operatorname{sh}\omega})_j\le1$.

(ii) For $n=0$ both sides vanish, since $\norm x_1=1$ on $\cT$. For $n\ge1$ let $y=A_{\omega_2\cdots\omega_n}x_{\operatorname{sh}^n\omega}$. By Lemma~\ref{lem:coding}, applied to $\operatorname{sh}\omega,\dots,\operatorname{sh}^{n-1}\omega$, we have $\pr y=\psi_{\omega_2\cdots\omega_n}(x_{\operatorname{sh}^n\omega})=x_{\operatorname{sh}\omega}$; so, by homogeneity, $\norm{A_{\omega_1\cdots\omega_n}x_{\operatorname{sh}^n\omega}}_1=\norm{A_{\omega_1}y}_1=\norm y_1\,e^{\xi(\omega)}$, and induction on $n$, applied to $\operatorname{sh}\omega$, gives $\log\norm y_1=\xi_{n-1}(\operatorname{sh}\omega)$.

(iii) If $\omega'\in[\omega_1]$, then by (i) $\abs{\xi(\omega)-\xi(\omega')}\le\log\frac{1+k}{1+k/2}<\log2$. If $n\ge2$ and $\omega'\in[\omega_1\cdots\omega_n]$, then $x_{\operatorname{sh}\omega}$ and $x_{\operatorname{sh}\omega'}$ lie in $C_{\iota(\omega_2)}$, at distance at most $\frac{\sqrt2}4\theta_*^{n-2}$ (Lemma~\ref{lem:coding}), and $x\mapsto\log\norm{A_j^kx}_1$ is $4$-Lipschitz on $C_{\iota(\omega_2)}$ by Lemma~\ref{lem:contr}.
\end{proof}

Summable functions on $\Omega$ and their pressure are as in Definition~\ref{def:shiftfns}(ii).

\begin{lemma}\label{lem:pressure}
For $t>1$ the function $-t\xi$ is summable and $P(-t\xi)=p(t)$: the pressure of $-t\xi$ is the growth rate per run of the zeta function of $\GN$ for the norm $\norm\cdot_1$ (Proposition~\ref{prop:runpressure} and Remark~\ref{rem:GNfast}). For $t\le1$, $-t\xi$ is not summable and $P(-t\xi)=\infty$. Hence $\sN$ is the only real number $t$ with $P(-t\xi)=0$.
\end{lemma}

\begin{proof}
By Lemma~\ref{lem:xi}(i), $\sum_e\exp\sup_{[e]}(-t\xi)$ lies between $3\sum_{k\ge1}(1+k)^{-t}$ and $3\sum_{k\ge1}(1+\frac k2)^{-t}$, which are finite if and only if $t>1$. Let $n\ge1$, let $w$ be an admissible word of length $n$, and let $\tau\in[w]$. By Lemma~\ref{lem:xi}(ii), $e^{\xi_n(\tau)}=\norm{A_wx_{\operatorname{sh}^n\tau}}_1$, where $x_{\operatorname{sh}^n\tau}\in C_{\iota(\tau_{n+1})}$ and $\iota(\tau_{n+1})$ differs from the last index $\iota(\tau_n)$ of $w$. By Lemma~\ref{lem:junction}, $\frac14\norm{A_w}_1\le e^{\xi_n(\tau)}\le\frac12\norm{A_w}_1$, so, for $t>0$,
\[
  2^t\,\zeta^{(n)}(t)\ \le\ \sum_{\abs w=n}\exp\sup_{[w]}(-t\xi_n)\ \le\ 4^t\,\zeta^{(n)}(t).
\]
For $t>1=t_0$, Proposition~\ref{prop:runpressure}(i) gives $P(-t\xi)=\lim_n\frac1n\log\zeta^{(n)}(t)=p(t)$. For $t\le1$ the middle sum is infinite for every $n\ge1$: for $t\le0$ it has infinitely many terms, all at least $1$; for $0<t\le1$, let $J$ be a multiindex with $n-1$ runs and $j$ an index different from the first index of $J$ (any $j$ if $J=\emptyset$); the multiindices $j^kJ$, $k\ge1$, have $n$ runs and, by submultiplicativity, $\norm{A_{j^kJ}}_1\le(3+2k)\norm{A_J}_1$, so $\zeta^{(n)}(t)=\infty$. The last statement follows from Proposition~\ref{prop:runpressure} and Remark~\ref{rem:GNfast}: on $(1,\infty)$, $p$ vanishes only at $\sN$.
\end{proof}

\emph{The Gibbs state of the runs.} The measure that we need is given by the following theorem.

\begin{lettertheorem}[Ruelle--Perron--Frobenius~{\cite[Thm.~2.4]{KK16}}; see~{\cite[Lemma~2.4.1, Thms.~2.4.3 and~2.4.6, Cor.~2.7.5]{MU03}}]\label{thm:RPF}
Let $\Sigma$ be a finitely irreducible shift space, let $0<\theta<1$, and let $u\in\cF_\theta(\Sigma)$ be summable. Then $u$ has an ergodic shift-invariant Gibbs state $\mu_u$, with $P_{\mu_u}=P(u)$.
\end{lettertheorem}

By Lemmas~\ref{lem:xi} and~\ref{lem:pressure}, the potential $-\sN\xi$ lies in $\cF_{\theta_*}(\Omega)$, is summable, and has $P(-\sN\xi)=0$. So Theorem~\ref{thm:RPF} applies to $\Sigma=\Omega$, which is finitely irreducible, $\theta=\theta_*$ and $u=-\sN\xi$, and it gives an ergodic shift-invariant Gibbs state $\mu_*=\mu_{-\sN\xi}$ of $-\sN\xi$, with $P_{\mu_*}=0$: for some $C\ge1$,
\[
  C^{-1}\ \le\ \mu_*\bigl([\omega_1\cdots\omega_n]\bigr)\,e^{\sN\xi_n(\omega)}\ \le\ C\qquad(\omega\in\Omega,\ n\ge1).
\]
It is the only invariant Gibbs state of $-\sN\xi$ (Theorem~\ref{thm:MU}(b) below). The Ruelle operator of $-\sN\xi$,
\[
  (\widehat\cL g)(\omega)=\sum_{e\in\mathcal A,\ \iota(e)\ne\iota(\omega_1)}e^{-\sN\xi(e\omega)}\,g(e\omega),
\]
where the sum runs over the preimages $e\omega$ of $\omega$ under $\operatorname{sh}$, is the operator $\cL_{\sN}$ of Section~\ref{sec:fast}, read on the symbolic model: since $\xi(e\omega)=\log\norm{A_ex_\omega}_1$ and $x_{e\omega}=\psi_e(x_\omega)$, if $g(\omega)=f(x_\omega,\iota(\omega_1))$ for a function $f$ on $\mathcal C$, then $(\widehat\cL g)(\omega)=(\cL_{\sN}f)(x_\omega,\iota(\omega_1))$.

\begin{corollary}\label{lem:intxi}
$\log\frac32\le\int\xi\,d\mu_*\le3C\sum_{k\ge1}\bigl(1+\frac k2\bigr)^{-\sN}\log(1+k)<\infty$.
\end{corollary}

\begin{proof}
By the Gibbs property and Lemma~\ref{lem:xi}(i), $\mu_*([(j,k)])\le C\,(1+\frac k2)^{-\sN}$, and $\log\frac32\le\xi\le\log(1+k)$ on $[(j,k)]$; the series converges because $\sN>1$.
\end{proof}

\subsection*{The Lyapunov exponents of the runs}
In this subsection and in the next one we use several results of ergodic theory, that we include for the reader's convenience.

\emph{Entropy of shift-invariant measures.}
Let $\Sigma$ be a shift space, as above. For a shift-invariant Borel probability measure $\nu$ on $\Sigma$ we write $h_\nu=h_\nu(\operatorname{sh})$ for its entropy~\cite[App.~1]{MU03}. For $N\ge1$, let $\alpha_N$ be the partition of $\Sigma$ into the cylinders of the admissible words of length $N$ (for $\Sigma=\Omega$, these words are the multiindices with $N$ runs), and let $H_\nu(\alpha_N)=-\sum_{\abs w=N}\nu([w])\log\nu([w])$ be its entropy. The partition $\alpha_1$ is generating, and $\alpha_N$ is the partition into the sets $[\omega_1]\cap\operatorname{sh}^{-1}[\omega_2]\cap\dots\cap\operatorname{sh}^{-(N-1)}[\omega_N]$. So, by the Kolmogorov--Sinai theorem~\cite[Thm.~A1.0.9]{MU03},
\begin{equation}\label{eq:KS}
  h_\nu=\lim_{N\to\infty}\frac1NH_\nu(\alpha_N)\qquad\text{whenever }H_\nu(\alpha_1)<\infty .
\end{equation}

\emph{Equilibrium states.}
Lemma~\ref{lem:maxent} below identifies $\mu_*$ through a uniqueness theorem of Mauldin and Urba\'nski, and Proposition~\ref{prop:AHS} applies the same theorem to a first-return map. Besides H\"older continuity, the pressure and Gibbs states (Definition~\ref{def:shiftfns}), it uses the following notion.

\begin{definition}[{\cite[\S2.2]{MU03}}]\label{def:equilibrium}
Let $\Sigma$ be a shift space, and let $u:\Sigma\to\R$ be continuous. A shift-invariant Borel probability measure $\nu$ on $\Sigma$ is an \emph{equilibrium state} of $u$ if $\int(-u)\,d\nu<\infty$ and $h_\nu+\int u\,d\nu=P(u)$.
\end{definition}

\begin{lettertheorem}[Mauldin--Urba\'nski~\cite{MU03}]\label{thm:MU}
Let $\Sigma$ be a finitely irreducible shift space, and let $u:\Sigma\to\R$ be H\"older continuous.
\begin{enumerate}
\item[(a)] \cite[Thm.~2.1.7]{MU03} If $u$ is bounded above, then $h_\nu+\int u\,d\nu\le P(u)$ for every shift-invariant Borel probability measure $\nu$ on $\Sigma$ with $\int(-u)\,d\nu<\infty$. If, moreover, $P(u)<\infty$, then $H_\nu(\alpha_q)<\infty$ for some $q\ge1$.
\item[(b)] \cite[Prop.~2.2.2(a) and Thm.~2.2.4]{MU03} If $m$ is a Gibbs state of $u$, then $P_m=P(u)$. If $u$ has a Gibbs state, then it has exactly one invariant Gibbs state, which we denote by $\mu_u$, and $\operatorname{sh}$ is ergodic with respect to $\mu_u$.
\item[(c)] \cite[Thm.~2.2.9]{MU03} If $u$ has a Gibbs state and $\int(-u)\,d\mu_u<\infty$, then $\mu_u$ is the unique equilibrium state of $u$.
\end{enumerate}
\end{lettertheorem}

\begin{lemma}\label{lem:maxent}
Let $\nu$ be a shift-invariant probability measure on $\Omega$ with $\int\xi\,d\nu<\infty$. Then $H_\nu(\alpha_1)<\infty$ and
\[
  h_\nu\ \le\ \sN\int\xi\,d\nu ,
\]
with equality if and only if $\nu=\mu_*$. In particular $h_{\mu_*}=\sN\int\xi\,d\mu_*$.
\end{lemma}

\begin{proof}
We apply Theorem~\ref{thm:MU} to $\Sigma=\Omega$, which is finitely irreducible, and to $u=-\sN\xi$. This $u$ is H\"older continuous (Lemma~\ref{lem:xi}(iii)) and bounded above, $u\le-\sN\log\frac32<0$ (Lemma~\ref{lem:xi}(i)); $\mu_*$ is an invariant Gibbs state of $u$ with $P_{\mu_*}=0$, so that $\mu_u=\mu_*$ and $P(u)=0$ by Theorem~\ref{thm:MU}(b), as in Lemma~\ref{lem:pressure}; and $\int(-u)\,d\mu_*=\sN\int\xi\,d\mu_*<\infty$ (Corollary~\ref{lem:intxi}).
Now let $\nu$ be a shift-invariant probability measure on $\Omega$ with $\int\xi\,d\nu<\infty$, that is, with $\int(-u)\,d\nu<\infty$. By Theorem~\ref{thm:MU}(a), $h_\nu-\sN\int\xi\,d\nu\le P(u)=0$, and, since $P(u)<\infty$, $H_\nu(\alpha_q)<\infty$ for some $q\ge1$; since $\alpha_q$ refines $\alpha_1$, $H_\nu(\alpha_1)\le H_\nu(\alpha_q)<\infty$. Finally, $h_\nu=\sN\int\xi\,d\nu$ means that $\nu$ is an equilibrium state of $u$, and, by Theorem~\ref{thm:MU}(c), $\mu_*$ is an equilibrium state of $u$, and the only one.
\end{proof}

For $\omega\in\Omega$ we set $A^{(N)}(\omega)=A_{\omega_1\cdots\omega_N}$, so that $A^{(N+M)}(\omega)=A^{(N)}(\omega)\,A^{(M)}(\operatorname{sh}^N\omega)$.
The limits $\lambda_k$ in Proposition~\ref{prop:lyap} below are given by Kingman's subadditive ergodic theorem, which we use in the following form.

\begin{lettertheorem}[Kingman~{\cite[Thms.~1, 3 and~5]{Kin68}}]\label{thm:kingman}
Let $T$ be a measurable map of a probability space $(\mathcal X,\mathbb P)$ that preserves $\mathbb P$ and is ergodic, and let $S_N:\mathcal X\to\R$, $N\ge1$, be integrable functions such that $S_{N+M}\le S_N+S_M\circ T^N$ for all $N,M\ge1$ and $\frac1N\int S_N\,d\mathbb P$ is bounded below. Then, as $N\to\infty$, $\frac1NS_N$ converges, almost everywhere and in mean, to the constant $\lim_N\frac1N\int S_N\,d\mathbb P$.
\end{lettertheorem}

In the terms of~\cite[\S1]{Kin68}, the family $S_{st}=S_{t-s}\circ T^s$, $0\le s<t$, is a stationary subadditive process with a time constant.

\begin{proposition}\label{prop:lyap}
For $k=1,2,3$ and $\mu_*$-almost every $\omega$, the limit
$$
\lambda_k=\lim_{N\to\infty}\frac1N\log\sigma_k\bigl(A^{(N)}(\omega)\bigr)
$$
exists and does not depend on $\omega$. Moreover $\lambda_1=\int\xi\,d\mu_*\ge\log\frac32$, $h_{\mu_*}=\sN\lambda_1$, and
\[
  -\frac{\lambda_3}{\lambda_1}=1+\frac{\lambda_2}{\lambda_1}\ \le\ \frac{1+\sA-\sN}{2-\sA}\ <\ 0.7766 .
\]
\end{proposition}

\begin{proof}
    If $\omega_1=(j,k)$, the matrices $A_{\omega_1}^{\pm1}=\id\pm k(\one-e_j)e_j^*$ have Frobenius norm $\sqrt{3+2k^2}\le4(1+\frac k2)\le4e^{\xi(\omega)}$ (Lemma~\ref{lem:xi}(i)); so $\log\norm{A_{\omega_1}^{\pm1}}\le\log4+\xi(\omega)$, which is $\mu_*$-integrable (Corollary~\ref{lem:intxi}). 

    For $i=1,2$ and $N\ge1$, let $S^i_N=\log(\sigma_1\cdots\sigma_i)\bigl(A^{(N)}\bigr)$. Since $\det A_I=1$, we have $\sigma_1(A_I)=\norm{A_I}$ and $(\sigma_1\sigma_2)(A_I)=\sigma_3(A_I)^{-1}=\norm{A_I^{-1}}$, and both are at least $1$; so $S^i_N\ge0$. Since $A^{(N+M)}(\omega)=A^{(N)}(\omega)\,A^{(M)}(\operatorname{sh}^N\omega)$, the inverse of a product is the product of the inverses in the opposite order, and the operator norm is submultiplicative, we have $S^i_{N+M}\le S^i_N+S^i_M\circ\operatorname{sh}^N$. In particular $0\le S^i_N\le\sum_{r=0}^{N-1}S^i_1\circ\operatorname{sh}^r$, where $S^1_1=\log\norm{A_{\omega_1}}$ and $S^2_1=\log\norm{A_{\omega_1}^{-1}}$ are at most $\log4+\xi$; since $\mu_*$ is shift-invariant, every $S^i_N$ is $\mu_*$-integrable. Finally, $\operatorname{sh}$ preserves $\mu_*$ and is ergodic (Theorem~\ref{thm:RPF}).
    So, by Theorem~\ref{thm:kingman}, $\frac1NS^i_N$ converges to a constant $\lambda^i$, $\mu_*$-almost everywhere and in mean. Since $\log\sigma_1=S^1_N$, $\log\sigma_2=S^2_N-S^1_N$ and, because $\sigma_1\sigma_2\sigma_3=1$, $\log\sigma_3=-S^2_N$ at $A^{(N)}$, the limits $\lambda_k$ exist, $\mu_*$-almost everywhere and in mean, and they are the constants
\[
  \lambda_1=\lambda^1,\qquad\lambda_2=\lambda^2-\lambda^1,\qquad\lambda_3=-\lambda^2;
\]
    in particular $\lambda_1+\lambda_2+\lambda_3=0$.

    The difference $\xi_N(\omega)-\log\sigma_1(A^{(N)}(\omega))$ is bounded: $\xi_N(\omega)=\log\norm{A^{(N)}(\omega)x_{\operatorname{sh}^N\omega}}_1$ lies in $[\log\norm{A^{(N)}}_1-\log4,\ \log\norm{A^{(N)}}_1-\log2]$ by Lemma~\ref{lem:xi}(ii) and Lemma~\ref{lem:junction}, and $\log\sigma_1(A^{(N)})$ lies in $[\log\norm{A^{(N)}}_1-\log3,\ \log\norm{A^{(N)}}_1]$ by Lemma~\ref{lem:sv}. 
    By Birkhoff's ergodic theorem~\cite{Bir31}, $\frac1N\xi_N\to\int\xi\,d\mu_*$ $\mu_*$-almost everywhere. So $\lambda_1=\int\xi\,d\mu_*$, which is at least $\log\frac32$ by Lemma~\ref{lem:xi}(i), and Lemma~\ref{lem:maxent} gives $h_{\mu_*}=\sN\int\xi\,d\mu_*=\sN\lambda_1$.

    Let $s\in(\sA,2)$. 
    Then $1<s<2$, so $\varphi^s=\sigma_2^{2-s}\sigma_1^{1-2s}$, and $Z_s=\sum_I\varphi^s(A_I)<\infty$. By the concavity of the logarithm, for $q_w\ge0$ with $\sum_wq_w=1$ and real $a_w$,
\begin{equation}\label{eq:gibbsineq}
  \sum_wq_w(a_w-\log q_w)\ \le\ \log\sum_we^{a_w}.
\end{equation}
    Let $N\ge1$. With $q_w=\mu_*([w])$ and $a_w=\log\varphi^s(A_w)$, where $w$ runs over the admissible words of length $N$, \eqref{eq:gibbsineq} gives
\[
  H_{\mu_*}(\alpha_N)+\int\log\varphi^s\bigl(A^{(N)}\bigr)\,d\mu_*\ \le\ \log Z_s ,
\]
    where the integral is finite, since $\abs{\log\varphi^s(A_I)}\le4\log\norm{A_I}_1+\log3$ (Lemma~\ref{lem:sv}) and $\log\norm{A^{(N)}}_1\le\xi_N+\log4$. By Lemma~\ref{lem:maxent}, $H_{\mu_*}(\alpha_1)<\infty$, so, by~\eqref{eq:KS}, $\frac1NH_{\mu_*}(\alpha_N)\to h_{\mu_*}=\sN\lambda_1$. 
    Since $\log\varphi^s(A^{(N)})=(2-s)\log\sigma_2(A^{(N)})+(1-2s)\log\sigma_1(A^{(N)})=(2-s)S^2_N-(1+s)S^1_N$, the convergence in mean of $\frac1NS^1_N$ and $\frac1NS^2_N$ gives $\frac1N\int\log\varphi^s(A^{(N)})\,d\mu_*\to(2-s)(\lambda_1+\lambda_2)-(1+s)\lambda_1=(2-s)\lambda_2+(1-2s)\lambda_1$. 
    So, dividing by $N$ and letting $N\to\infty$,
\[
  \sN\lambda_1+(2-s)\lambda_2+(1-2s)\lambda_1\ \le\ 0 ,
\]
    that is, $1+\lambda_2/\lambda_1\le(1+s-\sN)/(2-s)$, since $\lambda_1>0$. 
    Now let $s\downarrow\sA$. 
    The limit $(1+\sA-\sN)/(2-\sA)$ increases with $\sA$, since $\sN<3$, and decreases with $\sN$; with $\sA<1.6854$ (Theorem~\ref{thm: dimH}) and $\sN\ge2.4411$ (Theorem~\ref{thm: delta}) it is less than $0.2443/0.3146<0.7766$.
\end{proof}

\subsection*{The renormalization flow of~\cite{AHSinv}}
Write $y\in\Delta$ for the vector of the three lengths, normalized by $y_1+y_2+y_3=1$. A step of the induction of~\cite[\S3.1]{AHSinv} applies to $y$ when $y_j>\sum_{i\ne j}y_i$ for some $j$, and it replaces $y$ by the normalization of $(A_j^*)^{-1}y$; in~\cite{AHSinv} the lengths are kept in decreasing order, so the matrix of a step is $A_1^*$, up to permutation matrices. The accelerated induction groups the consecutive steps with the same winner \cite[Lemma~8]{AHSinv}, so its steps are the runs $(j,k)\in\mathcal A$, and the blocks of its cocycle are $B(n)=A_1^n$ \cite[\S3.6]{AHSinv}. The lengths after a path $\gamma$ of the induction are, up to normalization, $(B_\gamma^T)^{-1}y$. The renormalization flow is built from a complete path $\gamma_*$ and the subsimplex $\Delta_{\gamma_*}$ of the lengths whose induction begins with $\gamma_*$ \cite[Def.~12]{AHSinv}: its base is the first return map of the accelerated induction to $\Delta_{\gamma_*}$, and its roof function is the first return time
\[
  r(y)=-\log\norm{(B_\gamma^T)^{-1}y}_1,
\]
where $\gamma$ is the path followed by $y$ up to its first return. The returns are counted along the accelerated induction, because the cocycle of~\cite[Thm.~27]{AHSinv} is a product of whole blocks $B(n_i)$ and permutation matrices. By the proof of \cite[Lemma~13]{AHSinv}, $r\ge\log3$, while the normalizing factor of a single run of length $1$ is at most $2$ (Lemma~\ref{lem:xi}(i) and Lemma~\ref{lem:dict}(ii) below): $r$ is a return time, not the time of a single step.

Let $J=\one\one^*-\id$, and let $\Pi$ be the projective map $[y]\mapsto[Jy]$. In the $\ell^1$ charts it is the similarity $y\mapsto\frac12(\one-y)$ from $\Delta$ onto $\cT$, and, since $JA_j^*=A_jJ$, it conjugates $[y]\mapsto[A_j^*y]$ to $\psi_j$ (Remark~\ref{rem:conj}). Let $\Omega_0\subset\Omega$ be the shift-invariant set of the sequences in which every index occurs infinitely often. By Birkhoff's ergodic theorem, applied to the indicator functions of the cylinders $[(l,1)]$, $l=1,2,3$, which have positive $\mu_*$-measure by the Gibbs property, $\mu_*$-almost every sequence contains each $(l,1)$ infinitely often; so $\mu_*(\Omega_0)=1$.

\begin{lemma}\label{lem:dict}
\begin{enumerate}
\item[(i)] The map $\omega\mapsto\Pi^{-1}(x_\omega)$ is a bijection from $\Omega_0$ onto the set of the positive lengths at which the accelerated induction can be iterated forever, and the runs of the induction at $\Pi^{-1}(x_\omega)$ are $\omega_1,\omega_2,\dots$. So it conjugates $\operatorname{sh}$ to the accelerated induction.
\item[(ii)] Let $\omega\in\Omega_0$ and $y=\Pi^{-1}(x_\omega)$. For $n\ge0$, $(A_{\omega_1}^*\cdots A_{\omega_n}^*)^{-1}y=e^{-\xi_n(\omega)}\,\Pi^{-1}(x_{\operatorname{sh}^n\omega})$: the logarithm of the normalizing factor of the first $n$ steps of the accelerated induction at $y$ is $\xi_n(\omega)$.
\item[(iii)] With $\omega$ and $y$ as in~(ii) and $n\ge1$, the singular values of the matrix of the cocycle of~\cite{AHSinv} along the first $n$ steps at $y$ lie between $\frac12$ and $2$ times those of $A^{(n)}(\omega)$.
\end{enumerate}
\end{lemma}

\begin{proof}
(i) With $x=\Pi(y)=\frac12(\one-y)$, the lengths are positive if and only if $x$ lies in the interior of $\cT$, and a step won by $j$ applies at $y$ if and only if $x_j<\frac14$; it then maps $y$ to $\Pi^{-1}\bigl(\psi_j^{-1}(x)\bigr)$. At most one coordinate of a point of $\cT$ is less than $\frac14$, because the coordinates add up to $1$ and are at most $\frac12$; so at most one index can win. Since $A_i\frac12(\one-e_l)=\frac12(\one+e_m)$ for $\{i,l,m\}=\{1,2,3\}$, the map $\psi_i$ sends the side $\{x_j=\frac12\}$ of $\cT$ onto the side of $C_i=\{x\in\cT:\ x_i\le\frac14\}$ on which $x_j=\frac12$ if $j\ne i$, and $x_i=\frac14$ if $j=i$. Two facts follow, for $\tau\in\Omega$ and $z\in\cT$:
\begin{enumerate}
\item[(a)] $(x_\tau)_j=\frac12$ if and only if the index $j$ does not occur in $\tau$. Indeed, $\psi_i(z)_j=\frac12$ if and only if $z_j=\frac12$, for $i\ne j$, and $C_j$ does not meet $\{x_j=\frac12\}$. So, if $j$ first occurs in $\tau_m$, then $(x_{\operatorname{sh}^{m-1}\tau})_j\le\frac14$, hence $(x_\tau)_j\ne\frac12$. If $j$ does not occur, every cylinder $\psi_{\tau_1\cdots\tau_n}(\cT)$ meets $\{x_j=\frac12\}$, and so does their intersection.
\item[(b)] $\psi_j^m(z)_j<\frac14$ for every $m\ge1$, unless $m=1$ and $z_j=\frac12$. For $m\ge2$ this holds because $\psi_j^m(\cT)\subset\psi_j(C_j)$ does not meet $\{x_j=\frac14\}$.
\end{enumerate}
Let $\omega\in\Omega_0$, $\omega_1=(j,k)$ and $x'=x_{\operatorname{sh}\omega}\in C_i$, $i=\iota(\omega_2)$. By (a), $x'_j<\frac12$, so by (b) the points $\psi_j^m(x')$, $1\le m\le k$, have $j$-th coordinate less than $\frac14$. On the other hand $x'_j\ge\frac14$, because $C_i\cap C_j$ is a single point, with $j$-th coordinate $\frac14$. So exactly $k$ consecutive steps at $\Pi^{-1}(x_\omega)$ are won by $j$, and they lead to $\Pi^{-1}(x_{\operatorname{sh}\omega})$, whose lengths are positive by~(a). By induction, the runs are $\omega_1,\omega_2,\dots$.

Conversely, suppose the induction can be iterated forever at positive lengths $y$. Every index wins infinitely often: otherwise the unnormalized length of a losing index would be subtracted, from some step on, from the sum of the other two unnormalized lengths at every step. So the runs form a sequence $\omega\in\Omega_0$, and $\Pi(y)\in\psi_{\omega_1\cdots\omega_n}(\cT)$ for every $n$; hence $\Pi(y)=x_\omega$ by Lemma~\ref{lem:coding}.

(ii) Let $n=1$, $\omega_1=(j,k)$, $M=A_j^k$ and $y'=\Pi^{-1}(x_{\operatorname{sh}\omega})$. By~(i), $(M^*)^{-1}y=c\,y'$ with $c=\norm{(M^*)^{-1}y}_1$. On $\Delta$, $Jy=\one-y$, so $x_\omega=\frac12Jy$ and $x_{\operatorname{sh}\omega}=\frac12Jy'$. Since $J(M^*)^{-1}=M^{-1}J$, we get $Mx_{\operatorname{sh}\omega}=c^{-1}x_\omega$, so $\xi(\omega)=\log\norm{Mx_{\operatorname{sh}\omega}}_1=-\log c$. For $n\ge2$, apply this to $\omega$ and then the case $n-1$ to $\operatorname{sh}\omega$, using $\xi_n(\omega)=\xi(\omega)+\xi_{n-1}(\operatorname{sh}\omega)$.

(iii) Undoing the permutations by $\pi A_1\pi^{-1}=A_{\pi(1)}$ (Section~\ref{sec:souls}), the block of the step $(j,k)$ is $A_j^k$. By~\cite[Def.~12]{AHSinv} and~(i), the transpose of the matrix of the cocycle along the first $n$ steps maps the lengths after these steps back to a multiple of $y$. So, up to permutation matrices on both sides, it is $A_{\omega_1}^*\cdots A_{\omega_n}^*=J^{-1}A^{(n)}(\omega)J$. The matrix $J$ is symmetric, with eigenvalues $2,-1,-1$, so $\norm J=2$ and $\norm{J^{-1}}=1$, and $\sigma_k(J^{-1}AJ)$ lies between $\sigma_k(A)/2$ and $2\sigma_k(A)$.
\end{proof}

\emph{First-return maps.}
Let $T$ be a measurable map of a probability space $(\mathcal X,\mathbb P)$ that preserves $\mathbb P$ and is ergodic, and let $A\subset\mathcal X$ be a measurable set with $\mathbb P(A)>0$. By Birkhoff's ergodic theorem, applied to the indicator function of $A$, almost every point of $A$ returns to $A$ infinitely often. So the \emph{first-return time} $\tau_A(x)=\min\{n\ge1:\ T^nx\in A\}$ is finite for almost every $x\in A$, and the \emph{first-return map} $T_A(x)=T^{\tau_A(x)}(x)$ is defined almost everywhere on $A$.

\begin{lettertheorem}[Kac~\cite{Kac47}; see~\cite{Aar97}]\label{thm:kac}
Let $T$, $\mathbb P$ and $A$ be as above. Then $T_A$ preserves the restriction of $\mathbb P$ to $A$, that is, $\mathbb P(T_A^{-1}B)=\mathbb P(B)$ for every measurable set $B\subset A$~\cite[Prop.~1.5.3]{Aar97}. Moreover~\cite[1.5.5 and its proof]{Aar97},
\[
  \mathbb P(B)=\int_A\sum_{k=0}^{\tau_A-1}\mathbf 1_B\circ T^k\,d\mathbb P\qquad\text{for every measurable set }B\subset\mathcal X;
\]
in particular, $\int_A\tau_A\,d\mathbb P=1$.
\end{lettertheorem}

By linearity and monotone convergence, the identity for the sets $B$ gives
\begin{equation}\label{eq:kac}
  \int_A\sum_{k=0}^{\tau_A-1}f\circ T^k\,d\mathbb P=\int_{\mathcal X}f\,d\mathbb P\qquad\text{for every measurable }f:\mathcal X\to[0,\infty] .
\end{equation}

\emph{The flow as a suspension over a first-return map.}
Let $(j_1,k_1),\dots,(j_m,k_m)$ be the runs of the complete path $\gamma_*$. For $\omega\in\Omega_0$, the induction at $\Pi^{-1}(x_\omega)$ begins with $\gamma_*$ if and only if $\omega_r=(j_r,k_r)$ for $r<m$ and $\omega_m=(j_m,k)$ with $k\ge k_m$ (Lemma~\ref{lem:dict}(i); the last run of $\gamma_*$ may continue in $\omega$). Since $\gamma_*$ is complete, every index is the winner of some step of $\gamma_*$~\cite[Def.~8]{AHSinv}, so every index occurs among $j_1,\dots,j_m$, and $m\ge3$. Let $g=(j_1,k_1)\cdots(j_{m-1},k_{m-1})$, let $\mathcal G$ be the set of the words $g\,(j_m,k)$ with $k\ge k_m$, and let $Y_0=\Omega_0\cap\bigcup_{z\in\mathcal G}[z]$. Then $\mu_*(Y_0)>0$, because $\mu_*$ is a Gibbs state, and, by Lemma~\ref{lem:dict}(i), $\omega\mapsto\Pi^{-1}(x_\omega)$ maps $Y_0$ onto the points of $\Delta_{\gamma_*}$ with positive lengths at which the accelerated induction can be iterated forever. We apply the above to $T=\operatorname{sh}$, to $\mathbb P=\mu_*$, with respect to which $\operatorname{sh}$ is ergodic (Theorem~\ref{thm:RPF}), and to $A=Y_0$, and we let $Y$ be the set of the points of $Y_0$ that return to $Y_0$ infinitely often, so that $\mu_*(Y_0\setminus Y)=0$. Thus a point $\omega\in\Omega_0$ lies in $Y$ if and only if its first $m$ terms form a word of $\mathcal G$ and $\operatorname{sh}^n\omega\in Y_0$ for infinitely many $n$; in particular,
\begin{equation}\label{eq:returnY}
  \operatorname{sh}^n\omega\in Y\quad\Longleftrightarrow\quad\omega_{n+1}\cdots\omega_{n+m}\in\mathcal G\qquad(\omega\in Y,\ n\ge1).
\end{equation}
Let $\tau=\tau_Y$ be the first-return time to $Y$, $\operatorname{sh}_Y=\operatorname{sh}^\tau$ the first-return map, and
\[
  \rho=\sum_{j=0}^{\tau-1}\xi\circ\operatorname{sh}^j
\]
the induced roof. Let $\mu_*^Y=\mu_*|_Y/\mu_*(Y)$. By Theorem~\ref{thm:kac}, $\mu_*^Y$ is $\operatorname{sh}_Y$-invariant, and, by~\eqref{eq:kac} with $f=\xi$ and with $f=1$, and since $\int\xi\,d\mu_*<\infty$ (Corollary~\ref{lem:intxi}),
\begin{equation}\label{eq:kacY}
  \int_Y\rho\,d\mu_*^Y=\frac1{\mu_*(Y)}\int_\Omega\xi\,d\mu_*<\infty,\qquad\int_Y\tau\,d\mu_*^Y=\frac1{\mu_*(Y)} .
\end{equation}
By Lemma~\ref{lem:dict}(i), the base map of the flow of~\cite{AHSinv} corresponds to $\operatorname{sh}_Y$; by Lemma~\ref{lem:dict}(ii), since permutation matrices do not change $\ell^1$ norms, its roof function $r$ corresponds to $\rho$; and by Lemma~\ref{lem:dict}(iii), its cocycle along $q$ returns at $\omega\in Y$ has the singular values of $A^{(\tau_q(\omega))}(\omega)$, up to factors between $\frac12$ and $2$, where $\tau_q=\sum_{i<q}\tau\circ\operatorname{sh}_Y^i$. So the flow of~\cite{AHSinv} is the suspension of $(Y,\operatorname{sh}_Y)$ with roof $\rho$, that is, the suspension flow of $(\Omega_0,\operatorname{sh})$ with roof $\xi$, restricted to the orbits that meet $Y$.

\emph{The induced shift.}
Let $\mathcal W$ be the set of the first-return words $\omega_1\cdots\omega_{\tau(\omega)}$, $\omega\in Y$. For $\omega\in Y$ and $i\ge1$, let $u_i(\omega)\in\mathcal W$ be the first-return word of $\operatorname{sh}_Y^{i-1}\omega$; since $\omega$ returns to $Y$ infinitely often, $\omega$ is the concatenation $u_1(\omega)u_2(\omega)\cdots$. For $u_1,\dots,u_n\in\mathcal W$, let $[u_1\cdots u_n]_Y$ be the set of the $\omega\in Y$ with $u_i(\omega)=u_i$ for $i\le n$.

\begin{lemma}\label{lem:inducedshift}
\begin{enumerate}
\item[(i)] Let $u_1,\dots,u_n\in\mathcal W$, and let $L$ be the length of the word $w=u_1\cdots u_n$ in $\mathcal A$. Then $[u_1\cdots u_n]_Y$ is the set of the $\omega\in\Omega_0\cap[w]$ with $\operatorname{sh}^L\omega\in Y$.
\item[(ii)] For every sequence $(u_i)_{i\ge1}$ in $\mathcal W$, the concatenation $u_1u_2\cdots$ lies in $Y$, and its first-return words are $u_1,u_2,\dots$.
\end{enumerate}
So $\omega\mapsto(u_i(\omega))_{i\ge1}$ is a bijection from $Y$ onto $\mathcal W^{\N}$, which conjugates $\operatorname{sh}_Y$ to the shift of $\mathcal W^{\N}$ and maps the sets $[u_1\cdots u_n]_Y$ onto the cylinders.
\end{lemma}

\begin{proof}
Every point of $Y$ begins with $g$, a word of length $m-1$. We first show:
\begin{enumerate}
\item[(a)] if $y\in Y$ has first-return word $u$, of length $t$, and $x\in\Omega_0\cap[u]$ satisfies $\operatorname{sh}^tx\in Y$, then $x\in Y$, and its first-return word is $u$.
\end{enumerate}
Indeed, $x$ and $y$ both begin with $ug$, so they have the same first $t+m-1$ terms, and, for $0\le s<t$, the first $m$ terms of $\operatorname{sh}^sx$ are those of $\operatorname{sh}^sy$. For $s=0$ they form a word of $\mathcal G$, and $x$ returns to $Y_0$ infinitely often because $\operatorname{sh}^tx\in Y$; so $x\in Y$. For $1\le s<t$ they do not form a word of $\mathcal G$, by~\eqref{eq:returnY} for $y$; so $\operatorname{sh}^sx\notin Y$, by~\eqref{eq:returnY} for $x$. Hence $\tau(x)=t$.

(i) If $\omega\in[u_1\cdots u_n]_Y$, then $\omega\in[w]$ and $\operatorname{sh}^L\omega=\operatorname{sh}_Y^n\omega\in Y$. Conversely, let $x\in\Omega_0\cap[w]$ with $\operatorname{sh}^Lx\in Y$. For $n=1$, (a) gives $x\in[u_1]_Y$. For $n\ge2$, by induction on $n$, $\operatorname{sh}^{t}x\in[u_2\cdots u_n]_Y\subset Y$, where $t$ is the length of $u_1$, and (a) gives $x\in Y$, with first-return word $u_1$; so $x\in[u_1\cdots u_n]_Y$.

(ii) For each $i$ let $y_i\in Y$ have first-return word $u_i$. In $y_i$, the last term of $u_i$ is followed by $(j_1,k_1)$, the first term of $g$, which is also the first term of $u_{i+1}$; so $x=u_1u_2\cdots$ is admissible. We claim that, for $i\ge1$, $r\ge0$ and $z\in\mathcal G$, the first $m$ terms of the word $u_i\cdots u_{i+r}z$ form a word of $\mathcal G$. For $r=0$: $y_i$ begins with $u_iz_0$, where $z_0\in\mathcal G$ is formed by the first $m$ terms of $\operatorname{sh}^{t_i}y_i\in Y$, $t_i$ the length of $u_i$. If $t_i\ge m$, the first $m$ terms of $u_iz$ and of $y_i$ coincide. If $t_i<m$, the first $m$ terms of $u_iz$ are $u_i$ followed by the first $m-t_i\le m-1$ terms of $z$, which are those of $g$, hence those of $z_0$; so again they are the first $m$ terms of $y_i$. In both cases they form a word of $\mathcal G$, because $y_i\in Y$. For $r\ge1$, by induction on $r$ the first $m$ terms of $v=u_{i+1}\cdots u_{i+r}z$ form a word $z'\in\mathcal G$, and the first $m$ terms of $u_iv$ are those of $u_iz'$, which form a word of $\mathcal G$ by the case $r=0$. Choosing $r$ so large that $u_i\cdots u_{i+r}$ has length at least $m$, we get that the first $m$ terms of every tail $u_iu_{i+1}\cdots$ of $x$ form a word of $\mathcal G$. Such a word contains every index, because $j_1,\dots,j_m$ are all the indices; so $x\in\Omega_0$, and all the tails $\operatorname{sh}^{L_n}x$ lie in $Y_0$, where $L_n$ is the length of $u_1\cdots u_n$. Hence $\operatorname{sh}^{L_n}x\in Y$ for every $n$, and, by~(i), $x\in[u_1\cdots u_n]_Y$ for every $n$.

The last statement follows from~(i) and~(ii), since every $\omega\in Y$ is the concatenation of its first-return words.
\end{proof}

So the induced system is the full shift $\Sigma_Y=\mathcal W^{\N}$ on the countable alphabet $\mathcal W$, with shift $\operatorname{sh}_Y$ and cylinders $[u_1\cdots u_n]_Y$; it is finitely irreducible. We write $P_Y$ for the pressure of Definition~\ref{def:shiftfns}(ii) on $\Sigma_Y$, $\rho_n=\sum_{i<n}\rho\circ\operatorname{sh}_Y^i$ for the Birkhoff sums of $\rho$ with respect to $\operatorname{sh}_Y$, and, for $n\ge1$,
\[
  \mathrm{var}^Y_n\rho=\sup\bigl\{\abs{\rho(\omega)-\rho(\omega')}:\ \omega,\omega'\in Y\text{ have the same first $n$ first-return words}\bigr\}.
\]
The roof $\rho$ has the following two properties.
\begin{enumerate}
\item[(R1)] $\rho\ge\tau\log\frac32\ge c_0$, where $c_0=\log\frac32$, by Lemma~\ref{lem:xi}(i). So $\rho_n\ge nc_0$.
\item[(R2)] $\mathrm{var}^Y_n\rho\le C_\xi\,\theta_*^{n-1}$ for $n\ge1$, where $C_\xi=\log2+\sqrt2/(\theta_*(1-\theta_*))$. Indeed, let $\omega,\omega'\in Y$ have the same first $n$ first-return words $u_1,\dots,u_n$, and let $t$ and $L$ be the lengths of $u_1$ and of $u_1\cdots u_n$ as words in $\mathcal A$. Then $\tau(\omega)=\tau(\omega')=t$, $L\ge t+n-1$, and, for $0\le j<t$, $\operatorname{sh}^j\omega$ and $\operatorname{sh}^j\omega'$ have the same first $L-j\ge n$ terms. So
\[
  \abs{\rho(\omega)-\rho(\omega')}\ \le\ \sum_{j=0}^{t-1}\mathrm{var}_{L-j}\,\xi\ \le\ \sum_{i\ge n}\mathrm{var}_i\,\xi\ \le\ C_\xi\,\theta_*^{n-1},
\]
where the last inequality follows from Lemma~\ref{lem:xi}(iii) by summing a geometric series.
\end{enumerate}
\emph{The Gibbs property of $\mu_*^Y$.}
Let $C$ be the constant of the Gibbs property of $\mu_*$, and let $K=C^3e^{2\sN C_\xi}$. Then, for every admissible word $wv$ in $\mathcal A$,
\begin{equation}\label{eq:quasiB}
  K^{-1}\mu_*([w])\,\mu_*([v])\ \le\ \mu_*([wv])\ \le\ K\mu_*([w])\,\mu_*([v]).
\end{equation}
Indeed, let $L$ be the length of $w$, and let $x\in[wv]$, $y\in[w]$ and $z\in[v]$. Then $\xi_{L+\abs v}(x)=\xi_L(x)+\xi_{\abs v}(\operatorname{sh}^Lx)$. Since $x$ and $y$ have the same first $L$ terms, $\abs{\xi_L(x)-\xi_L(y)}\le\sum_{r=0}^{L-1}\mathrm{var}_{L-r}\,\xi\le C_\xi$, as in~(R2); in the same way, $\abs{\xi_{\abs v}(\operatorname{sh}^Lx)-\xi_{\abs v}(z)}\le C_\xi$. The Gibbs property bounds $\mu_*([wv])$, $\mu_*([w])$ and $\mu_*([v])$ by $e^{-\sN\xi_{L+\abs v}(x)}$, $e^{-\sN\xi_L(y)}$ and $e^{-\sN\xi_{\abs v}(z)}$, up to the factor $C^{\pm1}$, and~\eqref{eq:quasiB} follows.

Now let $u_1,\dots,u_n\in\mathcal W$, let $w=u_1\cdots u_n$, of length $L$, and let $\omega\in[u_1\cdots u_n]_Y$. By Lemma~\ref{lem:inducedshift}(i), $[u_1\cdots u_n]_Y$ differs from the union of the pairwise disjoint cylinders $[wz]$, $z\in\mathcal G$, by a set of measure zero, because $\mu_*(\Omega_0)=1$, $\mu_*(Y_0\setminus Y)=0$ and $\mu_*$ is shift-invariant. The words $wz$ are admissible, because in $\omega$ the last term of $w$ is followed by $(j_1,k_1)$. So, by~\eqref{eq:quasiB},
\[
  K^{-1}\mu_*([w])\,\mu_*(Y)\ \le\ \mu_*\bigl([u_1\cdots u_n]_Y\bigr)=\sum_{z\in\mathcal G}\mu_*([wz])\ \le\ K\mu_*([w])\,\mu_*(Y),
\]
since $\sum_{z\in\mathcal G}\mu_*([z])=\mu_*(Y_0)=\mu_*(Y)$. Dividing by $\mu_*(Y)$, $K^{-1}\mu_*([w])\le\mu_*^Y([u_1\cdots u_n]_Y)\le K\mu_*([w])$. Moreover, the first $n$ returns of $\omega$ to $Y$ occur at the lengths of $u_1$, $u_1u_2$, \dots, $w$, so $\rho_n(\omega)=\xi_L(\omega)$, and $C^{-1}\le\mu_*([w])\,e^{\sN\xi_L(\omega)}\le C$ because $\omega\in[w]$. Hence, with $C_Y=CK$,
\begin{equation}\label{eq:gibbsY}
  C_Y^{-1}\ \le\ \mu_*^Y\bigl([u_1\cdots u_n]_Y\bigr)\,e^{\sN\rho_n(\omega)}\ \le\ C_Y\qquad(\omega\in[u_1\cdots u_n]_Y,\ n\ge1).
\end{equation}
That is, $\mu_*^Y$ is a Gibbs state of $-\sN\rho$ on $\Sigma_Y$, with $P_{\mu_*^Y}=0$.

\emph{The pressure of~\cite{AHSinv}.}
In~\cite{AHSinv}, the pressure of $-\kappa r$ is the Gurevich--Sarig pressure~\cite[Def.~15]{AHSinv}, \cite{Sar99}, and, by~\cite[Lemma~23]{AHSinv}, it is also the growth rate of $L^n_{-\kappa r}1$, where $L_{-\kappa r}$ is the Ruelle operator of~\cite[\S5.1]{AHSinv}. Read on $\Sigma_Y$ through Lemmas~\ref{lem:dict} and~\ref{lem:inducedshift}, these quantities are, for $\varphi=-\kappa\rho$, $a\in\mathcal W$ and $\omega\in\Sigma_Y$,
\[
  Z_n(\varphi,a)=\sum_{\substack{\omega'\in[a]_Y\\ \operatorname{sh}_Y^n\omega'=\omega'}}e^{\varphi_n(\omega')},\qquad
  L^n_\varphi1(\omega)=\sum_{\operatorname{sh}_Y^n\omega'=\omega}e^{\varphi_n(\omega')},
\]
where $\varphi_n=-\kappa\rho_n$, and the Gurevich--Sarig pressure of $\varphi$ is $\lim_n\frac1n\log Z_n(\varphi,a)$. The next lemma identifies these two growth rates with $P_Y(\varphi)$, without using~\cite{Sar99}.

\begin{lemma}\label{lem:gurevich}
Let $\kappa>0$ and $\varphi=-\kappa\rho$. For every $a\in\mathcal W$ and every $\omega\in\Sigma_Y$, both $\frac1n\log Z_n(\varphi,a)$ and $\frac1n\log L^n_\varphi1(\omega)$ converge to $P_Y(\varphi)$ as $n\to\infty$. In particular, the pressure of $-\kappa r$ in~\cite{AHSinv} is $P_Y(-\kappa\rho)$, also when it is infinite.
\end{lemma}

\begin{proof}
Let $V=\kappa C_\xi/(1-\theta_*)$, and let $Z_n=Z_n(\varphi)$ be as in Definition~\ref{def:shiftfns}(ii). If $\omega'$ and $\omega''$ lie in the same cylinder $[w]_Y$ with $\abs w=n$, then, for $0\le i<n$, $\operatorname{sh}_Y^i\omega'$ and $\operatorname{sh}_Y^i\omega''$ have the same first $n-i$ first-return words; so, by~(R2),
\[
  \abs{\varphi_n(\omega')-\varphi_n(\omega'')}\ \le\ \kappa\sum_{i=0}^{n-1}\mathrm{var}^Y_{n-i}\rho\ \le\ V .
\]
Since $\Sigma_Y$ is a full shift (Lemma~\ref{lem:inducedshift}), the preimages of $\omega$ under $\operatorname{sh}_Y^n$ are the points $w\omega$, one in each cylinder $[w]_Y$ with $\abs w=n$; so $e^{-V}Z_n\le L^n_\varphi1(\omega)\le Z_n$. The fixed points of $\operatorname{sh}_Y^n$ in $[a]_Y$ are the periodic points $avav\cdots$, one in each cylinder $[av]_Y$ with $\abs v=n-1$; so
\[
  e^{-V}\sum_{\abs v=n-1}\exp\sup_{[av]_Y}\varphi_n\ \le\ Z_n(\varphi,a)\ \le\ Z_n .
\]
Since $\varphi_n(a\omega')=\varphi(a\omega')+\varphi_{n-1}(\omega')$, we have $\sup_{[av]_Y}\varphi_n\ge\inf_{[a]_Y}\varphi+\sup_{[v]_Y}\varphi_{n-1}$, and $\inf_{[a]_Y}\varphi=-\kappa\sup_{[a]_Y}\rho$ is finite by~(R2). Hence $e^{-V+\inf_{[a]_Y}\varphi}\,Z_{n-1}\le Z_n(\varphi,a)\le Z_n$ for $n\ge2$. Since $\frac1n\log Z_n\to P_Y(\varphi)$ (Definition~\ref{def:shiftfns}(ii)), both claims follow.
\end{proof}

\begin{proposition}[$\kappa_0=\sN$]\label{prop:AHS}
\begin{enumerate}
\item[(i)] $P_Y(-\kappa\rho)\ge(\sN-\kappa)c_0>0$ for $\kappa<\sN$, $P_Y(-\sN\rho)=0$, and $P_Y(-\kappa\rho)\le-(\kappa-\sN)c_0<0$ for $\kappa>\sN$. Hence $\kappa_0=\sN$, and it is the only zero of the pressure of $-\kappa r$.
\item[(ii)] The potential $-\sN\rho$ is summable, $\sum_{u\in\mathcal W}\exp\sup_{[u]_Y}(-\sN\rho)<\infty$, so $\norm{L_{-\sN r}1}_\infty<\infty$ for the Ruelle operator of~\cite[\S5.1]{AHSinv}. The invariant Gibbs measure $\mu_0$ of $-\sN r$ exists, is unique, and is the image of $\mu_*^Y$ under $\omega\mapsto\Pi^{-1}(x_\omega)$, and $\int r\,d\mu_0<\infty$. Moreover, for every $\operatorname{sh}_Y$-invariant probability measure $\nu$ on $Y$ with $\int_Y\rho\,d\nu<\infty$,
\begin{equation}\label{eq:induced-entropy-bound}
  h_\nu(\operatorname{sh}_Y)\ \le\ \sN\int_Y\rho\,d\nu ,
\end{equation}
with equality if and only if $\nu=\mu_*^Y$. So $\mu_0$ is the unique equilibrium state of $-\sN r$ for the base map of the flow of~\cite{AHSinv}, $h_{\mu_0}=\sN\int r\,d\mu_0$, and $\sN$ is the maximum of $h_\nu/\int r\,d\nu$ over the invariant probability measures $\nu$ of the base map with $\int r\,d\nu<\infty$, attained only at $\nu=\mu_0$.
\item[(iii)] Let $B$ be the matrix of the cocycle of~\cite{AHSinv} along the path to the first return. Then $\log\norm B$ and $\log\norm{(B^{-1})^T}$ are $\mu_0$-integrable. For $\mu_0$-almost every point, the Lyapunov exponents of the cocycle, normalized by the number of blocks as in~\cite[Thm.~27]{AHSinv}, are those of Proposition~\ref{prop:lyap}; normalized by the number of returns, they are multiplied by $\int_Y\tau\,d\mu_*^Y=1/\mu_*(Y)$. In either case, the ratio $-\lambda_3/\lambda_1$ is that of Proposition~\ref{prop:lyap}.
\end{enumerate}
\end{proposition}

\begin{proof}
Theorem~\ref{thm:MU} applies to $\Sigma=\Sigma_Y$ and $u=-\sN\rho$. Indeed, $\Sigma_Y$ is a full shift (Lemma~\ref{lem:inducedshift}), hence finitely irreducible; $u\in\cF_{\theta_*}(\Sigma_Y)$ by~(R2), since $\mathrm{var}_nu$, taken on $\Sigma_Y$, is $\sN\,\mathrm{var}^Y_n\rho$, and $u\le-\sN c_0$ is bounded above, by~(R1); $\mu_*^Y$ is $\operatorname{sh}_Y$-invariant (Theorem~\ref{thm:kac}) and is a Gibbs state of $u$ with $P_{\mu_*^Y}=0$, by~\eqref{eq:gibbsY}, so that $P_Y(-\sN\rho)=0$ and $\mu_u=\mu_*^Y$ by Theorem~\ref{thm:MU}(b); and $\int(-u)\,d\mu_*^Y=\sN\int_Y\rho\,d\mu_*^Y<\infty$, by~\eqref{eq:kacY}.

(i) Let $\kappa<\sN$. By~(R1), $-\kappa\rho_n=-\sN\rho_n+(\sN-\kappa)\rho_n\ge-\sN\rho_n+n(\sN-\kappa)c_0$ on $\Sigma_Y$. So $Z_n(-\kappa\rho)\ge e^{n(\sN-\kappa)c_0}Z_n(-\sN\rho)$ for every $n$, and $P_Y(-\kappa\rho)\ge P_Y(-\sN\rho)+(\sN-\kappa)c_0=(\sN-\kappa)c_0$. For $\kappa>\sN$, in the same way, $-\kappa\rho_n\le-\sN\rho_n-n(\kappa-\sN)c_0$, and $P_Y(-\kappa\rho)\le-(\kappa-\sN)c_0$. By Lemma~\ref{lem:gurevich}, $P_Y(-\kappa\rho)$ is the pressure of $-\kappa r$ in~\cite{AHSinv}; so this pressure vanishes only at $\kappa=\sN$, and $\kappa_0=\sN$.

(ii) By~\eqref{eq:gibbsY} with $n=1$, $\exp\sup_{[u]_Y}(-\sN\rho)\le C_Y\,\mu_*^Y([u]_Y)$ for every $u\in\mathcal W$, and these measures add up to $1$; so $-\sN\rho$ is summable. Since $L_{-\sN r}1$ reads on $\Sigma_Y$ as $\omega\mapsto\sum_{u\in\mathcal W}e^{-\sN\rho(u\omega)}$, also $\norm{L_{-\sN r}1}_\infty\le C_Y$. The Gibbs measures of~\cite{AHSinv} are the Gibbs states of Definition~\ref{def:shiftfns}(iii), read on $\Sigma_Y$, and $\kappa_0=\sN$ by~(i); by Theorem~\ref{thm:MU}(b), $\mu_*^Y$ is the only invariant Gibbs state of $-\sN\rho$. So $\mu_0$ exists, is unique, and is the image of $\mu_*^Y$, and $\int r\,d\mu_0=\int_Y\rho\,d\mu_*^Y<\infty$, by~\eqref{eq:kacY}. Now let $\nu$ be as in~(ii). By Theorem~\ref{thm:MU}(a), $h_\nu(\operatorname{sh}_Y)-\sN\int_Y\rho\,d\nu\le P_Y(-\sN\rho)=0$, which is~\eqref{eq:induced-entropy-bound}. Equality in~\eqref{eq:induced-entropy-bound} means that $\nu$ is an equilibrium state of $-\sN\rho$, and, by Theorem~\ref{thm:MU}(c), $\mu_*^Y$ is one, and the only one. Since $\int_Y\rho\,d\nu\ge c_0>0$ by~(R1), \eqref{eq:induced-entropy-bound} also says that $h_\nu(\operatorname{sh}_Y)/\int_Y\rho\,d\nu\le\sN$, with equality if and only if $\nu=\mu_*^Y$. Since $\omega\mapsto\Pi^{-1}(x_\omega)$ conjugates $\operatorname{sh}_Y$ to the base map of the flow and carries $\rho$ to $r$ (Lemma~\ref{lem:dict}), the same holds for the base map, $r$ and $\mu_0$.

(iii) At $\omega\in Y$, $\log\norm B\le\log2+\log\sigma_1(A^{(\tau)}(\omega))\le\log8+\rho(\omega)$ by Lemma~\ref{lem:dict}(iii) and the proof of Proposition~\ref{prop:lyap}, and $\norm{(B^{-1})^T}=\sigma_1\sigma_2(B)\le\norm B^2$, since $\det B=\pm1$. Both are integrable, because $\int_Y\rho\,d\mu_*^Y<\infty$ by~\eqref{eq:kacY}. Since $\mu_*^Y\ll\mu_*$, the limits of Proposition~\ref{prop:lyap} hold at $\mu_*^Y$-almost every $\omega$, and, by Lemma~\ref{lem:dict}(iii), the cocycle along $N$ blocks at $\Pi^{-1}(x_\omega)$ has the singular values of $A^{(N)}(\omega)$, up to factors between $\frac12$ and $2$. Along $q$ returns, $N=\tau_q(\omega)$. By Birkhoff's ergodic theorem, applied to the indicator function of $Y$, for $\mu_*$-almost every $\omega$ the number of the $n<N$ with $\operatorname{sh}^n\omega\in Y$, divided by $N$, tends to $\mu_*(Y)$ as $N\to\infty$. For $\omega\in Y$ and $N=\tau_q(\omega)$, these $n$ are $0=\tau_0(\omega)<\tau_1(\omega)<\dots<\tau_{q-1}(\omega)$; so $q/\tau_q(\omega)\to\mu_*(Y)$, that is, $\tau_q/q\to1/\mu_*(Y)=\int_Y\tau\,d\mu_*^Y$ by~\eqref{eq:kacY}, $\mu_*^Y$-almost everywhere.
\end{proof}

\begin{remark}[the measure of maximal entropy of the flow]\label{rem:maxentflow}
Theorem~25 of~\cite{AHSinv} states that the lift of $\mu_0$ is the unique measure of maximal entropy of the flow. Proposition~\ref{prop:AHS}(ii) is the corresponding statement for the base map $T$ of the flow: by~\eqref{eq:induced-entropy-bound}, $h_\nu(T)/\int r\,d\nu\le\sN$ for its invariant probability measures $\nu$ with $\int r\,d\nu<\infty$, with equality only for $\nu=\mu_0$. It gives the statement for the flow, by the argument of~\cite[proof of Thm.~25]{AHSinv}, once two facts on suspension flows, which we do not prove, are granted: every ergodic invariant probability measure of the flow is the lift of an ergodic invariant probability measure $\nu$ of $T$ with $\int r\,d\nu<\infty$; and, by Abramov's formula~\cite{Abr59b}, the entropy of the flow with respect to this lift is $h_\nu(T)/\int r\,d\nu$. Abramov's formula is classical for flows built over invertible maps, while $T$ is not invertible. Theorem~\ref{cor:diffusion} does not use this remark, and neither do the statements on the maximal entropy per unit time in the abstract, in the introduction and in Proposition~\ref{prop:AHS}(ii), which concern the base map $T$.
\end{remark}

We are finally ready to prove our last main result.

\thmdiffusion*

\begin{proof}
By Proposition~\ref{prop:AHS}, $\kappa_0=\sN$, the Gibbs measure $\mu_0$ of $-\sN r$ exists, $\norm{L_{-\sN r}1}_\infty<\infty$ and the cocycle is log-integrable for $\mu_0$; so \cite[Lemma~26, Thm.~27 and Thm.~30]{AHSinv} apply (see the beginning of this section), and $\delta=-\lambda_3/\lambda_1$ for the Lyapunov exponents of their cocycle for $\mu_0$. By Proposition~\ref{prop:AHS}(iii), this is the ratio of Proposition~\ref{prop:lyap}, which gives the bound.
\end{proof}

Our argument does not use the strict inequalities $\lambda_2<0$ and $\lambda_2>\lambda_3$
of~\cite{AHSinv}. In fact, for every ergodic invariant measure on the shift of the runs for which the cocycle is log-integrable, the non-strict bounds $\frac12\le1+\lambda_2/\lambda_1\le1$ hold. Indeed, $\lambda_1\ge\log\frac32>0$ as in the proof of Proposition~\ref{prop:lyap}, and, since $\sigma_1\sigma_2(A_I)=\norm{A_I^{-1}}$, we have $1+\lambda_2/\lambda_1=\lim_N\log\norm{(A^{(N)})^{-1}}/\log\norm{A^{(N)}}$. So $\norm{A_I^{-1}}_1\le\norm{A_I}_1$ (Lemma~\ref{lem:geom}(ii)) gives $1+\lambda_2/\lambda_1\le1$, and $\sigma_2(A_I)\ge\sigma_1(A_I)^{-1/2}$ (Lemma~\ref{lem:sv}) gives $1+\lambda_2/\lambda_1\ge\frac12$. The bounds $\frac12<\delta<1$ of~\cite{AHSinv} are their strict forms for $\mu_0$, which use the Pisot property and the simplicity of the Lyapunov spectrum.

\begin{remark}[numerics]\label{rem:diffnum}
The bound of Theorem~\ref{cor:diffusion} is probably not sharp. By Proposition~\ref{prop:AHS}, $\delta$ is the ratio $-\lambda_3/\lambda_1$ for the Gibbs state $\mu_*$ of the shift of the runs, which gives to the cylinder of a multiindex $I$ a mass comparable to $\norm{A_I}_1^{-\sN}$ (by its Gibbs property and Lemma~\ref{lem:junction}), nearly constant on a dyadic shell of $\norm{A_I}_1$. As a heuristic stand-in for $\mu_*$ we use the uniform distribution on the multiindices with
$\norm{A_I}_1$ in a dyadic shell; we do not prove that the resulting estimates converge to $\delta$. Over the $3.0\times10^{10}$ multiindices with
$\norm{A_I}_1\le65536$, consider the increments between consecutive shells of the averages
of $\log\norm{A_I^{-1}}$ and of $\log\norm{A_I}$ (Frobenius norms). Their ratio decreases
from $0.750$ at $2^{12}$ to $0.744$ at $2^{16}$, so numerically $\delta\approx0.74$. For
the central values of $\sA$ and $\sN$, the bound is $(1+\sA-\sN)/(2-\sA)\approx0.77$.
\end{remark}

\section*{Acknowledgments}
The author is thankful to Ivan Dynnikov for pointing out to him the work of Meester where the gasket first appeared in literature.
This material is based upon work supported in part by the National Science Foundation under Grant Number DMS-2308225.

\section*{Tools use}
The programs of the certificates, and drafts of parts of the proofs and of the text, were written with the assistance of Claude (Anthropic) and ChatGPT (OpenAI), two AI models. 
The author is responsible for all the content of this article.

\bibliographystyle{amsplain}
\bibliography{\jobname}

\providecommand{\bysame}{\leavevmode\hbox to3em{\hrulefill}\thinspace}
\providecommand{\MR}{\relax\ifhmode\unskip\space\fi MR }
\providecommand{\MRhref}[2]{%
  \href{http://www.ams.org/mathscinet-getitem?mr=#1}{#2}
}
\providecommand{\href}[2]{#2}
\begin{thebibliography}{10}

\bibitem{Aar97}
Jon Aaronson, \emph{An introduction to infinite ergodic theory}, Mathematical
  Surveys and Monographs, vol.~50, American Mathematical Society, Providence,
  RI, 1997.

\bibitem{Abr59b}
Leonid~M. Abramov, \emph{On the entropy of a flow}, Dokl. Akad. Nauk SSSR
  \textbf{128} (1959), 873--875, English transl.: Amer. Math. Soc. Transl. (2)
  \textbf{49} (1966), 167--170.

\bibitem{AR91}
Pierre Arnoux and G{\'e}rard Rauzy, \emph{Repr{\'e}sentation
  g{\'e}om{\'e}trique de suites de complexit{\'e} $2n+1$}, Bull. Soc. Math.
  France \textbf{119} (1991), 199--215.

\bibitem{AS13}
Pierre Arnoux and {\v{S}}t{\v{e}}p{\'a}n Starosta, \emph{The {R}auzy gasket},
  Further Developments in Fractals and Related Fields, Birkh{\"a}user, 2013,
  pp.~1--23.

\bibitem{AHSinv}
Artur Avila, Pascal Hubert, and Alexandra Skripchenko, \emph{Diffusion for
  chaotic plane sections of 3-periodic surfaces}, Invent. Math. \textbf{206}
  (2016), no.~1, 109--146.

\bibitem{AHS16}
\bysame, \emph{On the {H}ausdorff dimension of the {R}auzy gasket}, Bull. Soc.
  Math. France \textbf{144} (2016), 539--568.

\bibitem{Bar98}
Arthur Baragar, \emph{The exponent for the {M}arkoff--{H}urwitz equations},
  Pacific J. Math. \textbf{182} (1998), no.~1, 1--21.

\bibitem{BV12}
Michael~F. Barnsley and Andrew Vince, \emph{Real projective iterated function
  systems}, J. Geom. Anal. \textbf{22} (2012), no.~4, 1137--1172,
  arXiv:1003.3473.

\bibitem{BV13}
\bysame, \emph{The {C}onley attractors of an iterated function system}, Bull.
  Aust. Math. Soc. \textbf{88} (2013), 267--279.

\bibitem{BJJP13}
Val{\'e}rie Berth{\'e}, Damien Jamet, Timo Jolivet, and Xavier Proven{\c{c}}al,
  \emph{Critical connectedness of thin arithmetical discrete planes}, Discrete
  Geometry for Computer Imagery (DGCI 2013), Lecture Notes in Comput. Sci.,
  vol. 7749, Springer, 2013, Extended version: arXiv:1312.7820, pp.~107--118.

\bibitem{Bir31}
George~D. Birkhoff, \emph{Proof of the ergodic theorem}, Proc. Nat. Acad. Sci.
  U.S.A. \textbf{17} (1931), 656--660.

\bibitem{Bow79}
Rufus Bowen, \emph{Hausdorff dimension of quasi-circles}, Inst. Hautes
  {\'E}tudes Sci. Publ. Math. \textbf{50} (1979), 11--25.

\bibitem{Col42}
Lothar Collatz, \emph{Einschlie{\ss}ungssatz f{\"u}r die charakteristischen
  {Z}ahlen von {M}atrizen}, Math. Z. \textbf{48} (1942), 221--226.

\bibitem{DL03}
Roberto De~Leo, \emph{Characterization of the set of ``ergodic directions'' in
  {N}ovikov's problem of quasi-electron orbits in normal metals}, Russian Math.
  Surveys \textbf{58} (2003), no.~5, 1042--1043, arXiv:math/0207234.

\bibitem{DL05}
\bysame, \emph{Proof of a {D}ynnikov conjecture on the {N}ovikov problem of
  plane sections of periodic surfaces}, Russian Math. Surveys \textbf{60}
  (2005), no.~3, 566--567.

\bibitem{DL08}
\bysame, \emph{On a generalized {S}ierpi{\'n}ski fractal in $\mathbb{RP}^n$},
  2008, arXiv:0804.1154.

\bibitem{DL12}
\bysame, \emph{Exponential growth of norms in semigroups of linear
  automorphisms and {H}ausdorff dimension of self-projective {IFS}}, 2012,
  arXiv:1204.0250v1.

\bibitem{DL15EM}
\bysame, \emph{A conjecture on the {H}ausdorff dimension of attractors of real
  self-projective {I}terated {F}unction {S}ystems}, Exp. Math. \textbf{24}
  (2015), no.~3, 270--288.

\bibitem{DL15JGA}
\bysame, \emph{On the exponential growth of norms in semigroups of linear
  endomorphisms and the {H}ausdorff dimension of attractors of projective
  {I}terated {F}unction {S}ystems}, J. Geom. Anal. \textbf{25} (2015), no.~3,
  1798--1827.

\bibitem{DL26parity}
\bysame, \emph{A parity constraint on the topological numbers of
  self-complementary triply-periodic surfaces}, 2026, arXiv:2610.18416.

\bibitem{DLdata}
\bysame, \emph{Programs, logs and test functions for the computer-assisted
  bounds on the {L}evitt--{Y}occoz gasket}, 2026, Version 1.1, Zenodo,
  \href{https://doi.org/10.5281/zenodo.24140876}{doi:10.5281/zenodo.24140876}.

\bibitem{DLstreams}
\bysame, \emph{Recurrent and gradient dynamics of iterated function systems},
  2026, arXiv:2604.17588.

\bibitem{DD07}
Roberto De~Leo and Ivan~A. Dynnikov, \emph{An example of a fractal set of plane
  directions having chaotic intersections with a fixed 3-periodic surface},
  Russian Math. Surveys \textbf{62} (2007), no.~5, 990--992.

\bibitem{DD09}
\bysame, \emph{Geometry of plane sections of the infinite regular skew
  polyhedron $\{4,6\,|\,4\}$}, Geom. Dedicata \textbf{138} (2009), 51--67.

\bibitem{DLY26DCDS}
Roberto De~Leo and James~A. Yorke, \emph{Streams, graphs and global attractors
  of dynamical systems on locally compact spaces}, Discrete and Continuous
  Dynamical Systems \textbf{47} (2026), 308--340.

\bibitem{DeLeoYorke2026JDE}
\bysame, \emph{Streams, graphs and global attractors of dynamical systems with
  applications to parabolic {PDE}s}, Journal of Differential Equations
  \textbf{467} (2026), 114309.

\bibitem{DJT09}
{\'E}ric Domenjoud, Damien Jamet, and Jean-Luc Toutant, \emph{On the connecting
  thickness of arithmetical discrete planes}, Discrete Geometry for Computer
  Imagery (DGCI 2009), 2009, pp.~362--372.

\bibitem{DHS23}
Ivan Dynnikov, Pascal Hubert, and Alexandra Skripchenko, \emph{Dynamical
  systems around the {R}auzy gasket and their ergodic properties}, Int. Math.
  Res. Not. IMRN (2023), no.~8, 6461--6503.

\bibitem{DS15}
Ivan Dynnikov and Alexandra Skripchenko, \emph{Symmetric band complexes of thin
  type and chaotic sections which are not quite chaotic}, Trans. Moscow Math.
  Soc. \textbf{76} (2015), no.~2, 251--269.

\bibitem{Dyn99}
Ivan~A. Dynnikov, \emph{The geometry of stability regions in {N}ovikov's
  problem on the semiclassical motion of an electron}, Russian Math. Surveys
  \textbf{54} (1999), no.~1, 21--59.

\bibitem{Dyn08}
\bysame, \emph{Interval identification systems and plane sections of 3-periodic
  surfaces}, Proc. Steklov Inst. Math. \textbf{263} (2008), 65--77.

\bibitem{Ede62}
Michael Edelstein, \emph{On fixed and periodic points under contractive
  mappings}, J. London Math. Soc. \textbf{37} (1962), 74--79.

\bibitem{Fal88}
Kenneth~J. Falconer, \emph{The {H}ausdorff dimension of self-affine fractals},
  Math. Proc. Cambridge Philos. Soc. \textbf{103} (1988), 339--350.

\bibitem{FalBook}
\bysame, \emph{Fractal geometry: Mathematical foundations and applications},
  3rd ed., Wiley, Chichester, 2014.

\bibitem{FN18}
Richard~S. Falk and Roger~D. Nussbaum, \emph{{$C^m$} eigenfunctions of
  {P}erron--{F}robenius operators and a new approach to numerical computation
  of {H}ausdorff dimension: applications in $\mathbb{R}^1$}, J. Fractal Geom.
  \textbf{5} (2018), no.~3, 279--337.

\bibitem{Fed69}
Herbert Federer, \emph{Geometric measure theory}, Grundlehren Math. Wiss., vol.
  153, Springer, New York, 1969.

\bibitem{Fek23}
Michael Fekete, \emph{{\"U}ber die {V}erteilung der {W}urzeln bei gewissen
  algebraischen {G}leichungen mit ganzzahligen {K}oeffizienten}, Math. Z.
  \textbf{17} (1923), 228--249.

\bibitem{FKN11}
Robbert Fokkink, Cor Kraaikamp, and Hitoshi Nakada, \emph{On {S}chweiger's
  problems on fully subtractive algorithms}, Israel J. Math. \textbf{186}
  (2011), 285--296.

\bibitem{Fou20}
Charles Fougeron, \emph{Dynamical properties of simplicial systems and
  continued fraction algorithms}, 2020, Preprint, arXiv:2001.01367.

\bibitem{GMR19}
Alex Gamburd, Michael Magee, and Ryan Ronan, \emph{An asymptotic formula for
  integer points on {M}arkoff--{H}urwitz varieties}, Ann. of Math. (2)
  \textbf{190} (2019), no.~3, 751--809.

\bibitem{GW08}
Andreas Griewank and Andrea Walther, \emph{Evaluating derivatives: Principles
  and techniques of algorithmic differentiation}, 2nd ed., SIAM, Philadelphia,
  PA, 2008.

\bibitem{GRM20}
Rodolfo Guti{\'e}rrez-Romo and Carlos Matheus, \emph{Lower bounds on the
  dimension of the {R}auzy gasket}, Bull. Soc. Math. France \textbf{148}
  (2020), 321--327.

\bibitem{HW08}
Godfrey~H. Hardy and Edward~M. Wright, \emph{An introduction to the theory of
  numbers}, 6th ed., Oxford University Press, Oxford, 2008.

\bibitem{HJ85}
R.A. Horn and C.R. Johnson, \emph{Matrix analysis}, 2nd ed., Cambridge
  University Press, Cambridge, 2013.

\bibitem{HN17}
Yi~Huang and Paul Norbury, \emph{Simple geodesics and {M}arkoff quads}, Geom.
  Dedicata \textbf{186} (2017), 113--148.

\bibitem{HPR19}
Pascal Hubert and Olga Paris-Romaskevich, \emph{Triangle tiling billiards and
  the exceptional family of their escaping trajectories: circumcenters and
  {R}auzy gasket}, Exp. Math. (2019), Published online 21 May 2019;
  arXiv:1804.00181.

\bibitem{Hut81}
John~E. Hutchinson, \emph{Fractals and self similarity}, Indiana Univ. Math. J.
  \textbf{30} (1981), no.~5, 713--747.

\bibitem{IEEE754}
{IEEE}, \emph{{IEEE} standard for floating-point arithmetic}, 2019, IEEE Std
  754-2019.

\bibitem{JT09}
Damien Jamet and Jean-Luc Toutant, \emph{Minimal arithmetic thickness
  connecting discrete planes}, Discrete Appl. Math. \textbf{157} (2009), no.~3,
  500--509.

\bibitem{JLPX23}
Yuxiang Jiao, Jialun Li, Wenyu Pan, and Disheng Xu, \emph{On the dimension of
  limit sets on {$\mathbb P(\mathbb R^3)$} via stationary measures: variational
  principles and applications}, Int. Math. Res. Not. IMRN (2024), no.~19,
  13015--13045.

\bibitem{Arb}
Fredrik Johansson, \emph{Arb: efficient arbitrary-precision midpoint-radius
  interval arithmetic}, IEEE Trans. Comput. \textbf{66} (2017), no.~8,
  1281--1292.

\bibitem{Jur26}
Natalia Jurga, \emph{Hausdorff dimension of the {R}auzy gasket}, J. Eur. Math.
  Soc. (2026), Online first; arXiv:2312.04999, whose numbering we use.

\bibitem{Kac47}
Mark Kac, \emph{On the notion of recurrence in discrete stochastic processes},
  Bull. Amer. Math. Soc. \textbf{53} (1947), 1002--1010.

\bibitem{KK16}
Marc Kesseb{\"o}hmer and Sabrina Kombrink, \emph{A complex
  {R}uelle--{P}erron--{F}robenius theorem for infinite {M}arkov shifts with
  applications to renewal theory}, Discrete Contin. Dyn. Syst. Ser. S
  \textbf{10} (2017), no.~2, 335--352, arXiv:1604.08252; we use the numbering
  of the arXiv version.

\bibitem{Kie02}
Bernd Kieninger, \emph{Iterated function systems on compact {H}ausdorff
  spaces}, Shaker Verlag, Aachen, 2002.

\bibitem{Kin68}
John F.~C. Kingman, \emph{The ergodic theory of subadditive stochastic
  processes}, J. Roy. Statist. Soc. Ser. B \textbf{30} (1968), 499--510.

\bibitem{KO11}
Alex Kontorovich and Hee Oh, \emph{Apollonian circle packings and closed
  horospheres on hyperbolic 3-manifolds}, J. Amer. Math. Soc. \textbf{24}
  (2011), no.~3, 603--648.

\bibitem{KM95}
Cor Kraaikamp and Ronald Meester, \emph{Ergodic properties of a dynamical
  system arising from percolation theory}, Ergodic Theory Dynam. Systems
  \textbf{15} (1995), no.~4, 653--661.

\bibitem{Lal89}
Steven~P. Lalley, \emph{Renewal theorems in symbolic dynamics, with
  applications to geodesic flows, non-{E}uclidean tessellations and their
  fractal limits}, Acta Math. \textbf{163} (1989), 1--55.

\bibitem{Lev93}
Gilbert Levitt, \emph{La dynamique des pseudogroupes de rotations}, Invent.
  Math. \textbf{113} (1993), 633--670.

\bibitem{LPX23}
Jialun Li, Wenyu Pan, and Disheng Xu, \emph{On the dimension of limit sets on
  $\mathbb{P}(\mathbb{R}^3)$ via stationary measures: the theory and
  applications}, 2023, arXiv:2311.10265.

\bibitem{Mag20}
Michael Magee, \emph{Counting one-sided simple closed geodesics on {F}uchsian
  thrice punctured projective planes}, Int. Math. Res. Not. IMRN (2020),
  no.~13, 3886--3901.

\bibitem{MN04}
Andrei~Ya. Maltsev and Sergei~P. Novikov, \emph{Dynamical systems, topology,
  and conductivity in normal metals}, J. Stat. Phys. \textbf{115} (2004),
  no.~1--2, 31--46, arXiv:cond-mat/0304471.

\bibitem{MU00}
R.~Daniel Mauldin and Mariusz Urba{\'n}ski, \emph{Parabolic iterated function
  systems}, Ergodic Theory Dynam. Systems \textbf{20} (2000), no.~5,
  1423--1447.

\bibitem{MU03}
\bysame, \emph{Graph directed {M}arkov systems: Geometry and dynamics of limit
  sets}, Cambridge Tracts in Math., vol. 148, Cambridge Univ. Press, 2003.

\bibitem{MW06}
Richard~P. McGehee and Tamas Wiandt, \emph{Conley decomposition for closed
  relations}, J. Difference Equ. Appl. \textbf{12} (2006), no.~1, 1--47.

\bibitem{Mee89}
Ronald W.~J. Meester, \emph{An algorithm for calculating critical probabilities
  and percolation functions in percolation models defined by rotations},
  Ergodic Theory Dynam. Systems \textbf{9} (1989), 495--509.

\bibitem{MN89}
Ronald W.~J. Meester and Tomasz Nowicki, \emph{Infinite clusters and critical
  values in two-dimensional circle percolation}, Israel J. Math. \textbf{68}
  (1989), no.~1, 63--81.

\bibitem{MKC09}
Ramon~E. Moore, R.~Baker Kearfott, and Michael~J. Cloud, \emph{Introduction to
  interval analysis}, SIAM, Philadelphia, PA, 2009.

\bibitem{PS23}
Mark Pollicott and Benedict Sewell, \emph{An upper bound on the dimension of
  the {R}auzy gasket}, Bull. Soc. Math. France \textbf{151} (2023), no.~4,
  595--611.

\bibitem{PS24}
\bysame, \emph{An elementary proof that the {R}auzy gasket is fractal}, Ergodic
  Theory Dynam. Systems \textbf{44} (2024), no.~7, 1913--1922.

\bibitem{PS98}
George P{\'o}lya and G{\'a}bor Szeg{\H o}, \emph{Problems and theorems in
  analysis {I}}, Classics in Mathematics, Springer, Berlin, 1998.

\bibitem{Sar99}
Omri~M. Sarig, \emph{Thermodynamic formalism for countable {M}arkov shifts},
  Ergodic Theory Dynam. Systems \textbf{19} (1999), no.~6, 1565--1593.

\bibitem{Sch14}
Rolf Schneider, \emph{Convex bodies: the {B}runn--{M}inkowski theory}, 2nd ed.,
  Encyclopedia of Mathematics and its Applications, vol. 151, Cambridge
  University Press, Cambridge, 2014.

\bibitem{Sch00}
Fritz Schweiger, \emph{Multidimensional continued fractions}, Oxford University
  Press, Oxford, 2000.

\bibitem{Sew23}
Benedict Sewell, \emph{A formula for the upper box-counting dimension of
  self-projective sets}, 2023, arXiv:2306.03047.

\bibitem{Skr12}
Alexandra Skripchenko, \emph{Symmetric interval identification systems of order
  three}, Discrete Contin. Dyn. Syst. \textbf{32} (2012), no.~2, 643--656.

\bibitem{FLINT}
{The FLINT team}, \emph{{FLINT}: {F}ast {L}ibrary for {N}umber {T}heory}, 2026,
  Version 3.6.0, \url{https://flintlib.org}.

\bibitem{Tuc11}
Warwick Tucker, \emph{Validated numerics: A short introduction to rigorous
  computations}, Princeton University Press, Princeton, NJ, 2011.

\bibitem{Vin13}
Andrew Vince, \emph{M\"obius iterated function systems}, Trans. Amer. Math.
  Soc. \textbf{365} (2013), no.~1, 491--509.

\bibitem{VW24}
Polina Vytnova and Caroline Wormell, \emph{Hausdorff dimension of the
  {A}pollonian gasket}, Invent. Math. \textbf{239} (2025), no.~3, 909--946,
  arXiv:2406.04922.

\bibitem{Wie50}
Helmut Wielandt, \emph{Unzerlegbare, nicht negative {M}atrizen}, Math. Z.
  \textbf{52} (1950), 642--648.

\end{thebibliography}

\end{document}